\documentclass[reqno]{amsart}
\usepackage[active]{srcltx}
\usepackage{amssymb}
\usepackage[active]{srcltx}
\usepackage{amsmath}
\usepackage{times}
\usepackage[T1]{fontenc}
\usepackage{textcomp}
\usepackage{txfonts}
\usepackage{tikz}
\usetikzlibrary{arrows,decorations.markings}
\usetikzlibrary{shapes}
\usepackage{mathrsfs}
\usepackage{multicol}
\usepackage{hyperref}
\usepackage{amsthm}
\usepackage{cases}
\usepackage{leftidx}

\usepackage{titlesec}

\titleformat{\section}{\Large\bfseries}{\thesection.}{0.5em}{}
\titleformat{\subsection}{\large\itshape}{\thesubsection.}{0.5em}{}

\usepackage{pifont}

\usepackage{bm}

\DeclareMathAlphabet{\mathpzc}{OT1}{pzc}{m}{it}

{\catcode`\|=\active
  \gdef\set#1{\mathinner{\lbrace\,{\mathcode`\|"8000%
                                   \let|\midvert #1}\,\rbrace}}
}
\def\midvert{\egroup\mid\bgroup}

\usepackage{aliascnt}
\def\NewTheorem#1{%
  \newaliascnt{#1}{equation}
  \newtheorem{#1}[#1]{#1}
  \aliascntresetthe{#1}
  \expandafter\def\csname #1autorefname\endcsname{#1}
}

\newcounter{main}
\theoremstyle{plain}

\numberwithin{equation}{section}
\NewTheorem{Assumption}
\NewTheorem{Proposition}
\NewTheorem{Theorem}
\NewTheorem{Corollary}
\NewTheorem{Conjecture}
\NewTheorem{Lemma}
\NewTheorem{Definition}
\theoremstyle{definition}
\NewTheorem{Remark}
\NewTheorem{Remarks}
\NewTheorem{Example}

\def\@thm#1#2#3{%
  \ifhmode\unskip\unskip\par\fi
  \normalfont
  \trivlist
  \let\thmheadnl\relax
  \let\thm@swap\@gobble
  \let\thm@indent\noindent
  \thm@headfont{\SMStheoremfont}
  \thm@notefont{\fontseries\mddefault\upshape}%
  \thm@headpunct{.}
  \thm@headsep 5\p@ plus\p@ minus\p@\relax
  \thm@space@setup
  #1
  \@topsep \thm@preskip               
  \@topsepadd \thm@postskip           
  \def\@tempa{#2}\ifx\@empty\@tempa
    \def\@tempa{\@oparg{\@begintheorem{#3}{}}[]}%
  \else
    \refstepcounter{#2}%
    \def\@tempa{\@oparg{\@begintheorem{#3}{\csname the#2\endcsname}}[]}%
  \fi
  \@tempa
}

\renewcommand{\P}{\mathscr{P}^\Lambda}
\newcommand{\Z}{\mathbb{Z}}
\newcommand{\N}{\mathbb{N}}
\newcommand{\R}{\mathscr{R}^\Lambda_n}
\newcommand{\Q}{\mathbb{Q}}

\newcommand{\B}{\mathscr B_n^{\inftyweight}}
\newcommand{\Bi}{\mathscr B_\alpha^{\inftyweight}}

\newcommand{\affineP}{\mathscr P^{\kappa}_n}
\newcommand{\inftyweight}{\Lambda^\infty}

\newcommand{\F}{\mathbb{F}}

\newcommand{\Rn}{\mathscr R_n}
\newcommand{\Ri}{\mathscr R_{\alpha}^\Lambda}

\newcommand{\Rni}{\mathscr R_{\alpha}}

\def\ten{10}
\def\eleven{11}
\def\twelve{12}
\def\thirteen{13}
\def\fourteen{14}
\def\fifteen{15}
\def\sixteen{16}
\def\seventeen{17}
\def\eighteen{18}
\def\nineteen{19}
\def\twenty{20}
\def\twentyone{21}
\def\twentytwo{22}
\def\twentythree{23}
\def\twentyfour{24}

\def\hh{\leavevmode\newline\hspace*{6.5mm}}
\def\bi{\mathbf{i}}
\def\bj{\mathbf{j}}

\def\t{\mathsf{t}}
\def\s{\mathsf{s}}
\def\u{\mathsf{u}}
\def\v{\mathsf{v}}
\def\T{\hat\t}
\def\S{\hat\s}
\def\U{\hat\u}
\def\V{\hat\v}
\def\x{\mathsf{x}}
\def\y{\mathsf{y}}
\def\w{\mathsf{w}}
\def\l{\ell}

\newcommand\BZ{R^\Lambda_n}
\newcommand\Blam[1][\lambda]{R^{>#1}_n}
\newcommand\Bgelam[1][\lambda]{R^{\ge #1}_n}

\usepackage{tikz}  

\colorlet{darkgreen}{green!50!black}
\tikzset{dots/.style={ultra thick,loosely dotted},
         belt/.style={draw,blue,thick,fill=blue!50},
         greendot/.style={fill,circle,color=darkgreen,inner sep=1.5pt,outer sep=0}
}
\newenvironment{braid}{
  \begin{tikzpicture}[baseline=6mm,blue,line width=1pt, xscale = 0.35, yscale=0.4,
                      draw/.append style={rounded corners},
                      every node/.append style={font=\fontsize{5}{5}\selectfont}]%
  }{\end{tikzpicture}
}

\def\Grid(#1,#2){
  \draw[very thin,gray,step=2mm] (0,0)grid(#1,#2);
  \draw[very thin,darkgreen,step=10mm] (0,0)grid(#1,#2);
}

\newcommand\Tableau[2][\relax]{
  \begin{tikzpicture}[scale=0.5,draw/.append style={thick,black},baseline=-3mm]
    \ifx\relax#1\relax%
    \else 
      \foreach\box in {#1} {
        \ifx\box\relax\else
          \filldraw[blue!10]\box+(-.5,-.5) rectangle ++(.5,.5);
        \fi
      }
    \fi;
    \newcount\row\newcount\col
    \row=0
    \foreach \Row in {#2} {
       \col=1
       \foreach \k in \Row {
          \draw (\the\col,\the\row) +(-.5,-.5) rectangle ++(.5,.5);
          \ifnum\k<0
            \draw(\the\col,\the\row) node[fill=gray!20]{-\k};
          \else \draw (\the\col,\the\row) node{\k};
          \fi
          \global\advance\col by 1
       }
       \global\advance\row by -1
    }
  \end{tikzpicture}
}
\def\Tritab(#1|#2|#3){\Bigg(\hspace*{1mm}\Tableau{#1}\hspace*{1mm}\Bigg|\hspace*{1mm}%
    \Tableau{#2}\hspace*{1mm}\Bigg|\hspace*{1mm}\Tableau{#3}\hspace*{1mm}\Bigg)
}

\newcommand\YoungDiagram[2][\relax]{
  \begin{tikzpicture}[scale=0.5,draw/.append style={thick,black},baseline=-1mm]
    \ifx\relax#1\relax%
    \else 
    \foreach\box in {#1} {
      \filldraw[blue!10]\box rectangle ++(1,1);
    }
    \fi
    \newcount\row
    \row=0
    \foreach \col in {#2} {
       \draw(1,\the\row)grid ++(\col,1);
       \global\advance\row by -1
    }
  \end{tikzpicture}
}

\def\Tridiag(#1|#2|#3){\Bigg(\hspace*{1mm}\YoungDiagram{#1}\hspace*{1mm}\Bigg|\hspace*{1mm}%
    \YoungDiagram{#2}\hspace*{1mm}\Bigg|\hspace*{1mm}\YoungDiagram{#3}\hspace*{1mm}\Bigg)
}

\DeclareMathOperator{\res}{res} 
 \DeclareMathOperator{\node}{node}

\DeclareMathOperator{\Std}{Std}
\DeclareMathOperator{\Shape}{Shape}

\newcommand{\map}[2]{\,{:}\,#1\!\longrightarrow\!#2}

\usepackage[enableskew,vcentermath]{youngtab}
\def\tab(#1){\mbox{\tiny$\young(#1)$}\,}
\def\ydiag(#1){\mbox{\tiny$\yng(#1)$}\,}

\newpage

\begin{document}    
\bibliographystyle{andrew}

\title{Integral Basis Theorem of cyclotomic Khovanov-Lauda-Rouquier Algebras of type A}


\author{Ge Li}

\address{School of Mathematics and Statistics\\
		University of Sydney\\
		Sydney, NSW 2006}
		
\email{geli@maths.usyd.edu.au}

\begin{abstract}
In this paper we prove that the cyclotomic Khovanov-Lauda-Rouquier algebras in type A, $\R$, are $\Z$-free. We then extend the graded cellular basis of $\R$ constructed by Hu and Mathas to $\Rn$ and use this basis to give a classification of all irreducible $\Rn$-modules.

\end{abstract}

\keywords{Cyclotomic Hecke algebras, Khovanov-Lauda-Rouquier algebras, Cellular basis}
%

\maketitle          
\pagenumbering{arabic}


\section{Introduction}

Khovanov and Lauda~\cite{KhovLaud:diagI,KhovLaud:diagII} and Rouquier~\cite{Rouq:2KM} have introduced a remarkable new family of algebras $\mathscr R_n$, the \textbf{quiver Hecke algebras}, for each oriented quiver. They showed that these algebras categorify the positive part of the enveloping algebras of the corresponding quantum groups. The algebras $\mathscr R_n$ are naturally $\Z$-graded. Varagnolo and Vasserot~\cite{VV:KhovLauda} proved that, under this categorification, the canonical basis of the positive part of the quantum group corresponds to the image of the projective indecomposable modules in the Grothendieck rings of the quiver Hecke algebras when the Cartan matrix is symmetric.

The algebra $\Rn$ is infinite dimensional and for every highest weight vector in the corresponding Kac-Moody algebra there is an associated finite dimensional 'cyclotomic quotient' $\R$ of $\Rn$. The cyclotomic quiver algebras $\R$ were originally defined by Khovanov and Lauda~\cite{KhovLaud:diagI, KhovLaud:diagII} and Rouquier~\cite{Rouq:2KM} who conjectured that these algebras should categorify the irreducible representations of the corresponding quantum group. Lauda and Vazirani~\cite{LV:Crystals} proved that, up to shift, the simple $\Rn$-modules are indexed by the vertices of the corresponding crystal graph, and Kang and Kashiwara~\cite{KK:HighestWeight} proved the full conjecture by showing that the images of the projective irreducible modules in the Grothendieck ring $Rep(\R)$ correspond to the canonical basis of the corresponding highest weight module. Prior to this work, Brundan and Stropple~\cite{BrundanStroppel:KhovanovIII} proved this conjecture in the special the case when $\Lambda$ is a dominat weight of level $2$ and $\Gamma$ is the linear quiver and Brundan and Kleshchev~\cite{BK:GradedDecomp} established the conjecture for all $\Lambda$ when $\Gamma$ is a quiver of type $A$.

Let $\Gamma$ be the quiver of type $A_e$, for $e\in\{0,2,3,4,\dots\}$. Brundan and Kleshchev \cite{BK:GradedKL} proved that every degenerate and non-degenerate cyclotomic Hecke algebra $H^\Lambda_n$ of type $G(r,1,n)$ over a field is isomorphic to a cyclotomic quiver Hecke algebra $\R$ of type $A$. They did this by constructing an explicit
isomorphisms between these two algebras.

The algebras $\R$ are defined by generators and relations and so these algebras are defined over any integral domain. Hu and Mathas~\cite{HuMathas:GradedCellular} defined a homogeneous basis $\{\psi_{\s\t}\}$ of the cyclotomic quiver algebras $\R$  (see \autoref{basis with field} below), and they showed that $\R$ is $\Z$-free whenever $e=0$ or $e$ is invertible in the ground ring. They asked whether the algebra $\R$ is always $\Z$-free. Kleshchev-Mathas-Ram~\cite{KMR:UniversalSpecht} defined $\Z$-free Specht modules for the cyclotomic KLR algebras of type A (and the affine KLR algebras of type A), but that the existence of these modules does not imply that the cyclotomic KLR algebras are torsion free. The main result of this paper shows that this is always the case. More precisely, we prove the following.

\begin{Theorem} \label{main theorem intro}
Let $\R(\Z)$ be a cyclotomic Khovanov-Lauda-Rouquier algebra of type $A$ over $\Z$, where $\Lambda$ is a dominant weight of height~$\ell$. Then $\R(\Z)$ is a graded cellular algebra, with respect to the dominance order, with homogeneous cellular basis
$
\{ \psi_{\s\t}\ |\ \lambda\in\mathscr P_n^\Lambda\text{ and $\s,\t\in\Std(\lambda)$}\}.
$
In particular, $\R(\Z)$ is $\Z$-free of rank $\ell^nn!$.
\end{Theorem}

If $\mathscr O$ is any integral domain then $\R(\mathscr O)\cong\R(\Z)\otimes_\Z\mathscr O$, so it follows that $\R(\mathscr O)$ is free over~$\mathscr O$.

The proof of our main theorem is long and technical, requiring a delicate multistage induction. Fortunately, by \cite[Theorem~5.14]{HuMathas:GradedCellular} we may assume that $e\ne2$. Even though our arguments should apply in this case, being able to assume that $e\ne2$ dramatically simplifies our arguments because the quiver of type $A_e$ is simply laced when $e\ne2$.

The starting point for our arguments is the observation that the definition of Hu and Mathas' the homogeneous elements $\psi_{\s\t}$ makes sense over any ring. Consequently, the linearly independents elements $\{\psi_{\s\t}\}$ span a $\Z$-free submodule $R^\Lambda_n$ of $\R$. To prove our Main Theorem it is therefore enough to show that $R^\Lambda_n$ is closed multiplication by the generators of~$\R$ and that the identity element of $\R$ belongs to~$R^\Lambda_n$.

The algebra $\R$ is generated by elements $y_r$, $\psi_s$ and $e(\bi)$, where $1\le r\le n$, $1\le s<n$ and $\bi\in(Z/e\Z)^n$. To prove our Main Theorem we need to treat these three classes of generators separately. The cellular basis element $\psi_{\s\t}$ is indexed by two standard $\lambda$-tableaux where $\lambda$ is a multipartition of~$n$; the definitions of these terms are recalled in Section~1. We argue by simultaneous induction on~$n$, and on the lexicographic orderings on the set of multipartitions, to show that multiplication by the KLR generators always sends $\psi_{\s\t}$ to a $\Z$-linear combination of terms $\psi_{\u\v}$ which are larger in the lexicographic order. Multiplication by~$y_r$ is the hardest case, partly because once this case is understood it can be used to describe the action of $\psi_r$ and $e(\bi)$ on the $\psi$-basis of~$\R$.

As a consequence of \autoref{main theorem intro}, we obtain a graded cellular basis of $\R(\Z)$. We then extend this basis to obtain a graded cellular basis of $\Rn$, and hence $\Rn$ is a graded cellular algebra. Hence we can use similar arguments to Graham-Lehrer~\cite{GL} to give a complete set of non-isomorphic graded irreducible $\Rn$-modules. Koenig and Xi~\cite{KonigXi:AffineCA} introduced the notion of affine cellular algebras and they have shown that the affine Hecke algebra of type A is an affine cellular algebra. The work of Koenig and Xi predates the KLR grading, but nonetheless it gives a classification of the irreducible representation of affine Hecke algebras. Our approach is very different to that of Koenig and Xi in that, first, we incorporate the grading and, secondly, we obtain a new labeling of the irreducible representations that is compatible the labeling of the irreducible modules of the cyclotomic quotients.

In more detail, this paper is organized as follows. In Section 1 we summarise the background material from the representation theory of the cyclotomic Khovanov-Lauda-Rouquier algebras that we need, including the theory of (graded) cellular algebras and the combinatorics of multipartitions and tableaux. In Section 2 considers the special case where $\lambda$ is a multicomposition that has at most two rows. Once this case is understood we show for an arbitrary multipartition $\lambda$ that $\psi_{\t^\lambda\t^\lambda}y_r$ is a $\Z$-linear combination of higher terms, where $\t^\lambda$ is the `initial' $\lambda$-tableau. Section 3 begins by proving, again by induction, that $\psi_{\s\t}y_r$ is a linear combination of bigger terms in~$R^\Lambda_n$. By considering the Garnir tableau of two-rowed multipartition we then show that $\psi_{\s\t}\psi_r$ can be written in the required form. This result is then extended to multipartitions of arbitrary shape. Finally, we deduce that $e(\bi)\in R^\Lambda_n$, for all $\bi\in (\Z/e\Z)^n$, which completes the proof of our main result. In Section 4 we define a sequence of weights $\inftyweight$ and using it to extend the graded cellular basis of $\Ri$ to $\Rni$ and hence define a graded cellular basis for $\Rn$. By adapting the arguments of Graham-Lehrer~\cite{GL} we give a complete set of non-isomorphic graded simple $\Rn$-modules.

Finally, we remark that the calculations in Sections~2 and~3 gives an algorithm inductively for computing $\psi_{\s\t}y_r$ and $\psi_{\s\t}\psi_r$. We extend the techniques developed in this paper to give a KLR grading of the Brauer algebras~\cite{Li:GradedBrauer}.

\textbf{Acknowledgements.} I would like to thank Professor Andrew Mathas for his support and comments on this paper. I also want to thank Matthew Gibson, who gave me advice on writing paper in English. Further thanks go to Jun Hu, Jonathan Brundan and Alexander Kleshchev for many helpful discussions with me about cyclotomic Khovanov-Lauda-Rouquier algebras.

\section{Khovanov-Lauda-Rouquier Algebras}

In this section we are going to introduce the necessary background for our work. First we will define our principal object of study --- the (cyclotomic) Khovanov-Lauda-Rouquier algebras $\R$. Then we give a brief introduction to (graded) cellular algebras and symmetric groups. Finally after explaining tableaux combinatorics we describe a graded cellular basis for the cyclotomic KLR algebra, found by Hu and Mathas~\cite{HuMathas:GradedCellular}.

\subsection{The cyclotomic Khovanov-Lauda-Rouquier algebras}

Fix an integer $e\in \{0,2,3,4\ldots\}$ and $I = \Z/e\Z$. Let $\Gamma_e$ be the oriented quiver with vertex set $I$ and directed edges $i\rightarrow i+1$, for $i\in I$. Thus, $\Gamma_e$ is the quiver of type $A_{\infty}$ if $e = 0$ and if $e \geq 2$ then it is a cyclic quiver of type $A_e^{(1)}$:
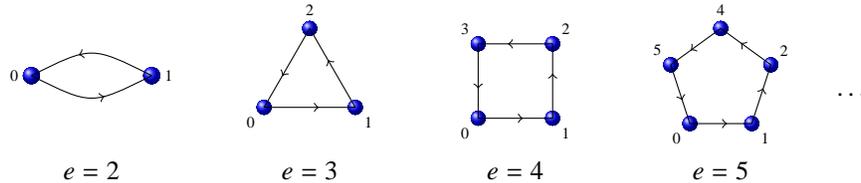
\begin{center}
\begin{tabular}{*5c}
  \begin{tikzpicture}[scale=0.8,decoration={curveto, markings,
            mark=at position 0.6 with {\arrow{>}}
    }]
    \useasboundingbox (-1.7,-0.7) rectangle (1.7,0.7);
    \foreach \x in {0, 180} {
      \shade[ball color=blue] (\x:1cm) circle(4pt);
    }
    \tikzstyle{every node}=[font=\tiny]
    \draw[postaction={decorate}] (0:1cm) .. controls (90:5mm) .. (180:1cm)
           node[left,xshift=-0.5mm]{$0$};
    \draw[postaction={decorate}] (180:1cm) .. controls (270:5mm) .. (0:1cm)
           node[right,xshift=0.5mm]{$1$};
  \end{tikzpicture}
& 
  \begin{tikzpicture}[scale=0.7,decoration={ markings,
            mark=at position 0.6 with {\arrow{>}}
    }]
    \useasboundingbox (-1.7,-0.7) rectangle (1.7,1.4);
    \foreach \x in {90,210,330} {
      \shade[ball color=blue] (\x:1cm) circle(4pt);
    }
    \tikzstyle{every node}=[font=\tiny]
    \draw[postaction={decorate}]( 90:1cm)--(210:1cm) node[below left]{$0$};
    \draw[postaction={decorate}](210:1cm)--(330:1cm) node[below right]{$1$};
    \draw[postaction={decorate}](330:1cm)--( 90:1cm)
           node[above,yshift=.5mm]{$2$};
  \end{tikzpicture}
& 
  \begin{tikzpicture}[scale=0.7,decoration={ markings,
            mark=at position 0.6 with {\arrow{>}}
    }]
    \useasboundingbox (-1.7,-0.7) rectangle (1.7,0.7);
    \foreach \x in {45, 135, 225, 315} {
      \shade[ball color=blue] (\x:1cm) circle(4pt);
    }
    \tikzstyle{every node}=[font=\tiny]
    \draw[postaction={decorate}] (135:1cm) -- (225:1cm) node[below left]{$0$};
    \draw[postaction={decorate}] (225:1cm) -- (315:1cm) node[below right]{$1$};
    \draw[postaction={decorate}] (315:1cm) -- (45:1cm) node[above right]{$2$};
    \draw[postaction={decorate}] (45:1cm) -- (135:1cm) node[above left]{$3$};
  \end{tikzpicture}
& 
  \begin{tikzpicture}[scale=0.7,decoration={ markings,
            mark=at position 0.6 with {\arrow{>}}
    }]
    \useasboundingbox (-1.7,-0.7) rectangle (1.7,0.7);
    \foreach \x in {18,90,162,234,306} {
      \shade[ball color=blue] (\x:1cm) circle(4pt);
    }
    \tikzstyle{every node}=[font=\tiny]
    \draw[postaction={decorate}] (162:1cm) -- (234:1cm)node[below left]{$0$};
    \draw[postaction={decorate}] (234:1cm) -- (306:1cm) node[below right]{$1$};
    \draw[postaction={decorate}] (306:1cm) -- (18:1cm) node[above right]{$2$};
    \draw[postaction={decorate}] (18:1cm) -- (90:1cm)
           node[above,yshift=.5mm]{$4$};
    \draw[postaction={decorate}] (90:1cm) -- (162:1cm) node[above left]{$5$};
  \end{tikzpicture}
&\raisebox{3mm}{$\dots$}
\\[4mm]
  $e=2$&$e=3$&$e=4$&$e=5$&
\end{tabular}
\end{center}

Let $(a_{i,j})_{i,j\in I}$ be the symmetric Cartan matrix associated with $\Gamma_e$, so that
\begin{equation*} \label{notation: cartan}
a_{i,j} = \begin{cases}
2, & \text{ if $i = j$,}\\
0, & \text{ if $i\neq j \pm 1$,}\\
-1, & \text{ if $e\neq 2$ and $i = j \pm 1$,}\\
-2, & \text{ if $e = 2$ and $i = j + 1$.}
\end{cases}
\end{equation*}

To the quiver $\Gamma_e$ attach the standard Lie theoretic data of a Cartan matrix $(a_{ij})_{i,j\in I}$, fundamental weights $\{\Lambda_i| i\in I\}$, positive weights $P_+ = \sum_{i\in I}\N \Lambda_i$, positive roots $Q_+ = \bigoplus_{i\in I}\N\alpha_i$ and let $(\cdot,\cdot)$ be the bilinear form determined by
$$
(\alpha_i,\alpha_j)=a_{ij}\qquad\text{and}\qquad
          (\Lambda_i,\alpha_j)=\delta_{ij},\qquad\text{for }i,j\in I.
$$
Fix a \textbf{weight} \label{notation: weight} $\Lambda = \sum_{i\in I} a_i \Lambda_i \in P_+$. Then $\Lambda$ is a weight of
\textbf{level} $l(\Lambda) = \l = \sum_{i\in I} a_i$. A \textbf{multicharge} \label{notation: multicharge} for $\Lambda$ is a sequence
$\kappa_\Lambda = (\kappa_1,\ldots,\kappa_\l) \in I^\l$ such that
$$
(\Lambda,\alpha_i) = a_i = \#\set{ 1\leq s\leq \l|\kappa_s \equiv i\pmod{e}}
$$
for any $i\in I$.

The following algebras were introduced by Khovanov and Lauda and Rouquier who defined KLR algebras for arbitrary oriented quivers.

\begin{Definition}[Khovanov and Lauda~\cite{KhovLaud:diagI,KhovLaud:diagII} and Rouquier~\cite{Rouq:2KM}] \label{notation: KLR}
Suppose $\mathscr O$ is an integral ring and $n$ is a positive integer. The \textbf{Khovanov-Lauda--Rouquier algebra}, $\mathscr{R}_n(\mathscr O)$ of type $\Gamma_e$ is the unital associative $\mathscr O$-algebra with generators
  $$\{\hat \psi_1,\dots,\hat \psi_{n-1}\} \cup
         \{ \hat y_1,\dots,\hat y_n \} \cup \set{\hat e(\bi)|\bi\in I^n}$$
  and relations
\begin{align}
\hat e(\bi) \hat e(\bj) &= \delta_{\bi\bj}\hat e(\bi),
&{\textstyle\sum_{\bi \in I^n}} \hat e(\bi)&= 1,\label{alg:trivial1}\\
\hat y_r \hat e(\bi) &= \hat e(\bi) y_r,
&\hat \psi_r \hat e(\bi)&= \hat e(s_r{\cdot}\bi) \hat \psi_r,
&\hat y_r \hat y_s &= \hat y_s \hat y_r,\label{alg:trivial2}
\end{align}
\vskip-21pt
\begin{align}
\hat \psi_r \hat y_s  &= \hat y_s \hat \psi_r,&\text{if }s \neq r,r+1,\label{alg:trivial3}\\
\hat \psi_r \hat \psi_s &= \hat \psi_s \hat \psi_r,&\text{if }|r-s|>1,\label{alg:trivial4}
\end{align}
\vskip-21pt
\begin{align}
  \hat \psi_r \hat y_{r+1} \hat e(\bi) &= \begin{cases}
      (\hat y_r\hat \psi_r+1)\hat e(\bi),\hspace*{18mm} &\text{if $i_r=i_{r+1}$},\\
    \hat y_r\hat \psi_r \hat e(\bi),&\text{if $i_r\neq i_{r+1}$}
  \end{cases} \label{alg:psi-y com}\\
  \hat y_{r+1} \hat \psi_r \hat e(\bi) &= \begin{cases}
      (\hat \psi_r \hat y_r+1) \hat e(\bi),\hspace*{18mm} &\text{if $i_r=i_{r+1}$},\\
    \hat \psi_r \hat y_r \hat e(\bi), &\text{if $i_r\neq i_{r+1}$}
  \end{cases}\label{alg:y-psi com}\\
  \hat \psi_r^2 \hat e(\bi) &= \begin{cases}
       0,&\text{if $i_r = i_{r+1}$},\\
      \hat e(\bi),&\text{if $i_r \ne i_{r+1}\pm1$},\\
      (\hat y_{r+1}-\hat y_r)\hat e(\bi),&\text{if  $e\ne2$ and $i_{r+1}=i_r+1$},\\
       (\hat y_r - \hat y_{r+1})\hat e(\bi),&\text{if $e\ne2$ and  $i_{r+1}=i_r-1$},\\
      (\hat y_{r+1} - \hat y_{r})(\hat y_{r}-\hat y_{r+1}) \hat e(\bi),&\text{if $e=2$ and
$i_{r+1}=i_r+1$}
\end{cases}\label{alg:psipsi}\\
\hat \psi_{r}\hat \psi_{r+1} \hat \psi_{r} \hat e(\bi) &= \begin{cases}
    (\hat \psi_{r+1} \hat \psi_{r} \hat \psi_{r+1} +1)\hat e(\bi),\hspace*{7mm}
       &\text{if $e\ne2$ and $i_{r+2}=i_r=i_{r+1}-1$},\\
  (\hat \psi_{r+1} \hat \psi_{r} \hat \psi_{r+1} -1)\hat e(\bi),
       &\text{if $e\ne2$ and $i_{r+2}=i_r=i_{r+1}+1$},\\
  \big(\hat \psi_{r+1} \hat \psi_{r} \hat \psi_{r+1} +\hat y_r\\
  \qquad -2\hat y_{r+1}+\hat y_{r+2}\big)\hat e(\bi),
    &\text{if $e=2$ and $i_{r+2}=i_r=i_{r+1}+1$},\\
  \hat \psi_{r+1} \hat \psi_{r} \hat \psi_{r+1} \hat e(\bi),&\text{otherwise.}
\end{cases}\label{alg:braid}
\end{align}
for $\bi,\bj\in I^n$ and all admissible $r$ and $s$. Moreover, $\mathscr{R}_n(\mathscr O)$ is naturally
$\Z$-graded with degree function determined by
$$\deg \hat e(\bi)=0,\qquad \deg \hat y_r=2\qquad\text{and}\qquad \deg
  \hat \psi_s \hat e(\bi)=-a_{i_s,i_{s+1}},$$
for $1\le r\le n$, $1\le s<n$ and $\bi\in I^n$.
\end{Definition}

Notice that the relations depend on the quiver $\Gamma_e$. By \cite[Theorem~5.14]{HuMathas:GradedCellular}, if $\mathscr O$ is a commutative integral domain and suppose either $e = 0$, $e$ is non-zero prime, or that $e{\cdot}1_{\mathscr O}$ is invertible in $\mathscr O$, $\mathscr R_n^\Lambda(\mathscr O)$ is an $\mathscr O$-free algebra.

 Following Khovanov and Lauda~\cite{KhovLaud:diagI}, we will frequently use
diagrammatic analogues of the relations of $\mathscr{R}_n(\mathscr O)$ in order
to simplify our calculations. To do this we associate to each generator
of $\mathscr{R}_n(\mathscr O)$ an $I$-labelled decorated planar diagram on~$2n$ points in
the following way:
$$
e(\bi)=
\begin{braid}\tikzset{baseline=7mm}
  \draw (0,4)node[above]{$i_1$}--(0,0);
  \draw (1,4)node[above]{$i_2$}--(1,0);
  \draw[dots] (1.2,4)--(3.8,4);
  \draw[dots] (1.2,0)--(3.8,0);
  \draw (4,4)node[above]{$i_{n}$}--(4,0);
\end{braid},
\quad \psi_re(\bi)=
\begin{braid}\tikzset{baseline=7mm}
  \draw (0,4)node[above]{$i_1$}--(0,0);
  \draw[dots] (0.2,4)--(1.8,4);
  \draw[dots] (0.2,0)--(1.8,0);
  \draw (2,4)node[above]{$i_{r-1}$}--(2,0);
  \draw (3,4)node[above]{$i_r$}--(4,0);
  \draw (4,4)node[above]{$i_{r+1}$}--(3,0);
  \draw (5,4)--(5,0);
  \draw[dots] (5.2,4)--(6.8,4);
  \draw[dots] (5.2,0)--(6.8,0);
  \draw (7,4)node[above]{$i_{n}$}--(7,0);
\end{braid},
\quad\text{and}\quad
y_se(\bi) =
\begin{braid}\tikzset{baseline=7mm}
  \draw (0,4)node[above]{$i_1$}--(0,0);
  \draw[dots] (0.2,4)--(1.8,4);
  \draw[dots] (0.2,0)--(1.8,0);
  \draw (2,4)node[above]{$i_{s-1}$}--(2,0);
  \draw (3,4)node[above]{$i_s$}--(3,0);
  \node[greendot] at (3,2){};
  \draw (4,4)node[above]{$i_s$}--(4,0);
  \draw[dots] (4.2,4)--(5.8,4);
  \draw[dots] (4.2,0)--(5.8,0);
  \draw (6,4)node[above]{$i_{n}$}--(6,0);
\end{braid},
$$
for $\bi\in I^n$, $1\le r<n$ and $1\le s\le n$. The $r$-th string of the diagram is the string labelled with $i_r$.

Diagrams are considered up to isotopy, and multiplication of diagrams is given
by concatenation, subject to the relations \eqref{alg:trivial1}--\eqref{alg:braid}. In
more detail, if $D_1$ and $D_2$ are two diagrams then the diagrammatic analogue of the relation
$e(\bi)e(\bj)=\delta_{\bi\bj}e(\bi)$ is
$$D_1 \cdot D_2 =
 \begin{tikzpicture}[baseline=9mm,blue,line width=1pt,scale=0.4,
                      draw/.append style={rounded corners},
                      every node/.append style={font=\fontsize{5}{5}\selectfont}]%
  \draw[thick](-0.4,6)rectangle(4.4,3.5)node[midway]{$D_1$};
  \draw (0,3.5)node[below right=-1mm]{$i_1$}--(0,3);
  \draw (1,3.5)node[below right=-1mm]{$i_2$}--(1,3);
  \draw (4,3.5)node[below right=-1mm]{$i_n$}--(4,3);
  \draw[dotted](0,3)--(0,2);
  \draw[dotted](1,3)--(1,2);
  \draw[dotted](4,3)--(4,2);
  \draw (0,2)--(0,1.5)node[above right=-1mm]{$j_1$};
  \draw (1,2)--(1,1.5)node[above right=-1mm]{$j_2$};
  \draw (4,2)--(4,1.5)node[above right=-1mm]{$j_n$};
  \draw[thick](-0.4,1.5)rectangle(4.4,-1)node[midway]{$D_2$};
 \end{tikzpicture}
 =\delta_{\bi\bj}
 \begin{tikzpicture}[baseline=11mm,blue,line width=1pt, scale=0.4,
                      draw/.append style={rounded corners},
                      every node/.append style={font=\fontsize{5}{5}\selectfont}]%
  \draw[thick](-0.4,6)rectangle(4.4,3.5)node[midway]{$D_1$};
  \draw (0,3.5)--(0,3)node[right=-1mm,pos=0.9]{$i_1$};
  \draw (1,3.5)--(1,3)node[right=-1mm,pos=0.9]{$i_2$};
  \draw (4,3.5)--(4,3)node[right=-1mm,pos=0.9]{$i_n$};
  \draw (0,3)--(0,2.5);
  \draw (1,3)--(1,2.5);
  \draw (4,3)--(4,2.5);
  \draw[thick](-0.4,2.5)rectangle(4.4,0)node[midway]{$D_2$};
 \end{tikzpicture}
$$
That is, $D_1\cdot D_2=0$ unless the labels of the strings on the bottom of $D_1$
match the corresponding labels on the top of the strings in~$D_2$ in which case
we just concatenate the two diagrams.

Multiplication by $y_r$ simply adds a decorative dot to the $r$-th string, reading left to right, so relations \eqref{alg:trivial2}--\eqref{alg:trivial4}
become self when written in terms of diagrams. Ignoring the extraneous strings on the left and right, and setting $i=i_r$ and $j=i_{r+1}$, the diagrammatic analogue of relations~\eqref{alg:psi-y com} and \eqref{alg:y-psi com} is
\begin{equation}\label{dia:y-psi com}
\begin{braid}
 \draw (0,4)node[above]{$i$}--(4,0);
 \draw (4,4)node[above]{$i$}--(0,0);
  \node[greendot] at (3,3){};
 \end{braid}
 \ \  -\ \
 \begin{braid}
 \draw (0,4)node[above]{$i$}--(4,0);
 \draw (4,4)node[above]{$i$}--(0,0);
  \node[greendot] at (1,1){};
 \end{braid}
 =\delta_{ij}
 \begin{braid}
 \draw (0,4)node[above]{$i$}--(0,0);
 \draw (2,4)node[above]{$i$}--(2,0);
 \end{braid}
=
 \begin{braid}
 \draw (0,4)node[above]{$i$}--(4,0);
 \draw (4,4)node[above]{$i$}--(0,0);
  \node[greendot] at (3,1){};
 \end{braid}
 \ \ -\ \
 \begin{braid}
 \draw (0,4)node[above]{$i$}--(4,0);
 \draw (4,4)node[above]{$i$}--(0,0);
  \node[greendot] at (1,3){};
 \end{braid}.
 \end{equation}
 Similarly, if $e\ne2$ then relation~\eqref{alg:psipsi} becomes
\begin{align}
\begin{braid}
 \draw (0,4)node[above]{$i$}--(2,2)--(0,0);
 \draw (2,4)node[above]{$j$}--(0,2)--(2,0);
 \end{braid} &\ \  =\ \
 \begin{cases}
\ \ \ \ \ \ \ 0, & \text{if $i = j$,}\\
\ \ \ \begin{braid}
 \draw (0,4)node[above]{$i$}--(0,0);
 \draw (2,4)node[above]{$j$}--(2,0);
\end{braid}, & \text{if $i\neq j\pm 1$,}\\
\ \ \ \pm \begin{braid}
 \draw (0,4)node[above]{$i$}--(0,0);
 \draw (2,4)node[above]{$j$}--(2,0);
 \node[greendot] at (2,2){};
\end{braid}
\ \ \mp\ \
\begin{braid}
 \draw (0,4)node[above]{$i$}--(0,0);
 \draw (2,4)node[above]{$j$}--(2,0);
 \node[greendot] at (0,2){};
\end{braid}, & \text{if $j = i\pm1$.}
\end{cases}\label{dia:psipsi}
\end{align}
and if $e\ne2$ then the diagrammatic analogue of relation~\eqref{alg:braid} is
\begin{align}
\begin{braid}
 \draw (0,4)node[above]{$i$}--(4,0);
 \draw (2,4)node[above]{$j$}--(4,2)--(2,0);
 \draw (4,4)node[above]{$k$}--(0,0);
\end{braid}
 \ \ -\ \
\begin{braid}
 \draw (0,4)node[above]{$i$}--(4,0);
 \draw (2,4)node[above]{$j$}--(0,2)--(2,0);
 \draw (4,4)node[above]{$k$}--(0,0);
\end{braid}
&\ \ =\ \
\delta_{i,k}(\delta_{i,j+1} - \delta_{i,j-1})\begin{braid}
 \draw (0,4)node[above]{$i$}--(0,0);
 \draw (2,4)node[above]{$j$}--(2,0);
 \draw (4,4)node[above]{$k$}--(4,0);
\end{braid}.
&\label{dia:braid}
\end{align}

Using the relations in $\mathscr{R}_n(\mathscr O)$ it is easy to verify the following identity which we record for future use:
\begin{equation}
\hat e(\bi)\hat y_r^k \hat y_{r+1}^k \hat \psi_r = \hat e(\bi)\hat \psi_r \hat y_r^k \hat y_{r+1}^k \label{alg:ypsi}
\end{equation}
for any $\bi$. Clearly it is enough to prove this relation when $k=1$ when,
diagrammatically, this identity takes the form
\begin{equation}
\begin{braid}
 \draw (0,4)node[above]{$i$}--(4,0);
 \draw (4,4)node[above]{$j$}--(0,0);
 \node[greendot] at (1,1){};
 \node[greendot] at (3,1){};
\end{braid}
=
\begin{braid}
 \draw (0,4)node[above]{$i$}--(4,0);
 \draw (4,4)node[above]{$j$}--(0,0);
 \node[greendot] at (1,3){};
 \node[greendot] at (3,3){};
\end{braid}\label{dia:ypsi}
\end{equation}
locally on the $r$ and $r+1$-th strings and where we set $i=i_r$ and $j=i_{r+1}$.

Three more easy, and very useful, consequences of the relations are the
following:

\begin{align}
\begin{braid}
 \draw (0,4)node[above]{$i$}--(0,0);
 \draw (4,4)node[above]{$i$}--(4,0);
\end{braid}
&=
-\begin{braid}
 \draw (0,4)node[above]{$i$}--(1,3)--(1,1)--(0,0);
 \draw (1,4)node[above]{$i$}--(0,3)--(0,1)--(1,0);
 \node[greendot] at (0,2.5){};
 \node[greendot] at (0,1.5){};
\end{braid}\label{dia:ii}
-
\begin{braid}
 \draw (0,4)node[above]{$i$}--(4,0);
 \draw (4,4)node[above]{$i$}--(0,0);
 \node[greendot] at (1,1){};
\end{braid}
-
\begin{braid}
 \draw (0,4)node[above]{$i$}--(4,0);
 \draw (4,4)node[above]{$i$}--(0,0);
 \node[greendot] at (1,3){};
\end{braid}
\\
\begin{braid}
 \draw (0,4)node[above]{$i$}--(0,0);
 \draw (4,4)node[above]{$i$}--(4,0);
 \node[greendot] at (4,2){};
\end{braid}
&=
-\begin{braid}
 \draw (0,4)node[above]{$i$}--(1,3)--(1,1)--(0,0);
 \draw (1,4)node[above]{$i$}--(0,3)--(0,1)--(1,0);
 \node[greendot] at (0,2.5){};
 \node[greendot] at (0,1.5){};
 \node[greendot] at (0.8,0.2){};
\end{braid}
-
\begin{braid}
 \draw (0,4)node[above]{$i$}--(4,0);
 \draw (4,4)node[above]{$i$}--(0,0);
 \node[greendot] at (1,1){};
 \node[greendot] at (3,1){};
\end{braid}
-
\begin{braid}
 \draw (0,4)node[above]{$i$}--(4,0);
 \draw (4,4)node[above]{$i$}--(0,0);
 \node[greendot] at (0.5,3.5){};
 \node[greendot] at (1.5,2.5){};
\end{braid}
-
\begin{braid}
 \draw (0,4)node[above]{$i$}--(0,0);
 \draw (4,4)node[above]{$i$}--(4,0);
 \node[greendot] at (0,2){};
\end{braid}
\label{dia:ii2}
\\
\label{dia:ii3}
\begin{braid}
 \draw (0,4)node[above]{$i$}--(4,0);
 \draw (4,4)node[above]{$i$}--(0,0);
\end{braid}
&\overset{(\ref{dia:y-psi com})}=
\begin{braid}
 \draw (0,4)node[above]{$i$}--(4,2)--(0,0);
 \draw (4,4)node[above]{$i$}--(0,2)--(4,0);
 \node[greendot] at (4,4){};
\end{braid}
-
\begin{braid}
 \draw (0,4)node[above]{$i$}--(4,2)--(0,0);
 \draw (4,4)node[above]{$i$}--(0,2)--(4,0);
 \node[greendot] at (0,2){};
\end{braid}
\overset{(\ref{dia:psipsi})}=
-\begin{braid}
 \draw (0,4)node[above]{$i$}--(4,2)--(0,0);
 \draw (4,4)node[above]{$i$}--(0,2)--(4,0);
 \node[greendot] at (0,2){};
\end{braid}\\
&\overset{(\ref{dia:y-psi com})}=
\begin{braid}
 \draw (0,4)node[above]{$i$}--(4,2)--(0,0);
 \draw (4,4)node[above]{$i$}--(0,2)--(4,0);
 \node[greendot] at (4,2){};
\end{braid}
-
\begin{braid}
 \draw (0,4)node[above]{$i$}--(4,2)--(0,0);
 \draw (4,4)node[above]{$i$}--(0,2)--(4,0);
 \node[greendot] at (0,4){};
\end{braid}
\overset{(\ref{dia:psipsi})}=
\begin{braid}
 \draw (0,4)node[above]{$i$}--(4,2)--(0,0);
 \draw (4,4)node[above]{$i$}--(0,2)--(4,0);
 \node[greendot] at (4,2){};
\end{braid} \notag
\end{align}

Note that (\ref{dia:ii2}) follows by multiplicating (\ref{dia:ii}) by $y_{r+1}$ and expanding.

In the rest of the paper we will play around with these diagrammatic notations a lot. In order to make the reader easy to follow our calculation we will use dotted strands to represent moving strands and arrows to represent moving dots. If we are going to move a dot then we will also write the strand which the dot is on dotted so the reader can see the arrow clearly. For example, we will write
$$
\begin{braid}
 \draw (0,4)[densely dotted] node[above]{$1$}--(2,0);
 \draw (1,4)[densely dotted] node[above]{$2$}--(2,2)--(1,0);
 \draw (2,4)[densely dotted] node[above]{$1$}--(0,0);
 \draw (3,4)node[above]{$3$}--(3,0);
 \draw (4,4)node[above]{$0$}--(4,0);
 \draw (5,4)node[above]{$1$}--(5,0);
\end{braid}
\overset{(\ref{dia:braid})}=
\begin{braid}
 \draw (0,4)node[above]{$1$}--(2,0);
 \draw (1,4)node[above]{$2$}--(0,2)--(1,0);
 \draw (2,4)node[above]{$1$}--(0,0);
 \draw (3,4)node[above]{$3$}--(3,0);
 \draw (4,4)node[above]{$0$}--(4,0);
 \draw (5,4)node[above]{$1$}--(5,0);
\end{braid}
-
\begin{braid}
 \draw (0,4)node[above]{$1$}--(0,0);
 \draw (1,4)node[above]{$2$}--(1,0);
 \draw (2,4)node[above]{$1$}--(2,0);
 \draw (3,4)node[above]{$3$}--(3,0);
 \draw (4,4)node[above]{$0$}--(4,0);
 \draw (5,4)node[above]{$1$}--(5,0);
\end{braid}
$$
to signify the application of relation~(\ref{dia:braid}) and
$$
\begin{braid}
 \draw (0,4)node[above]{$1$}--(0,0);
 \draw (1,4)node[above]{$2$}--(1,0);
 \draw (2,4)node[above]{$1$}--(2,0);
 \draw (3,4)[densely dotted] node[above]{$3$}--(4,0);
 \draw (4,4)[densely dotted] node[above]{$3$}--(3,0);
 \draw (5,4)node[above]{$0$}--(5,0);
 \node[greendot] at (3.25,1){};
 \node[greendot] at (3.75,1){};
 \draw[->] (3.25,1) -- (3.75,3);
 \draw[->] (3.75,1) -- (3.25,3);
\end{braid}
\overset{(\ref{dia:ypsi})}=
\begin{braid}
 \draw (0,4)node[above]{$1$}--(0,0);
 \draw (1,4)node[above]{$2$}--(1,0);
 \draw (2,4)node[above]{$1$}--(2,0);
 \draw (3,4) node[above]{$3$}--(4,0);
 \draw (4,4) node[above]{$3$}--(3,0);
 \draw (5,4)node[above]{$0$}--(5,0);
 \node[greendot] at (3.25,3){};
 \node[greendot] at (3.75,3){};
\end{braid}
$$
to signify the application of relation~(\ref{dia:ypsi}).

We can define a linear map $*\map{\Rn}{\Rn}$ by swapping the diagrams of $\Rn$ up-side-down. For example,
$$
\left(\begin{braid}
 \draw (0,4)node[above]{$0$}--(1,0);
 \draw (1,4)node[above]{$1$}--(3,0);
 \draw (2,4)node[above]{$3$}--(0,0);
 \draw (3,4) node[above]{$2$}--(4,2)--(2,0);
 \draw (4,4) node[above]{$2$}--(3,2)--(4,0);
 \node[greendot] at (0.8,0.8){};
 \node[greendot] at (2,2){};
 \node[greendot] at (1.8,3.6){};
 \node[greendot] at (3,1){};
\end{braid}\right)^* =
\begin{braid}
 \draw (0,4)node[above]{$3$}--(2,0);
 \draw (1,4)node[above]{$0$}--(0,0);
 \draw (2,4)node[above]{$2$}--(4,2)--(3,0);
 \draw (3,4) node[above]{$1$}--(1,0);
 \draw (4,4) node[above]{$2$}--(3,2)--(4,0);
 \node[greendot] at (0.8,3.2){};
 \node[greendot] at (2,2){};
 \node[greendot] at (1.8,0.4){};
 \node[greendot] at (3,3){};
\end{braid}.
$$

It is obvious that $*$ is an anti-isomorphism and it preserves the generators of $\Rn$.\\

Fix a weight $\Lambda = \sum_{i\in I} a_i \Lambda_i$ with $a_i \in \N$. Let $N_n^\Lambda(\mathscr O)$ be the two-sided ideal of $\mathscr{R}_n$ generated by the elements with form $e(\bi)y_1^{(\Lambda,\alpha_{i_1})}$. We can now define the main object of study in this paper, the cyclotomic Khovanov-Lauda-Rouquier algebras, which were introduced by Khovanov and Lauda~\cite[Section 3.4]{KhovLaud:diagI}.

\begin{Definition} \label{notation: cyclotomic KLR}
The \textbf{cyclotomic Khovanov-Lauda-Rouquier algebras} of weight $\Lambda$ and type $\Gamma_e$ is the algebra $\R(\mathscr O) = \mathscr{R}_n(\mathscr O)/ N_n^\Lambda(\mathscr O)$.
\end{Definition}

Therefore, if we write $e(\bi) = \hat e(\bi) + N_n^\Lambda(\mathscr O)$, $y_r = \hat y_r + N_n^\Lambda(\mathscr O)$ and $\psi_s = \hat \psi_s + N_n^\Lambda(\mathscr O)$, the algebra $\R(\mathscr O)$ is the unital $\mathscr O$-algebra generated by
$$
\{\psi_1,\dots,\psi_{n-1}\} \cup \{ y_1,\dots,y_n \} \cup \set{e(\bi)|\bi\in I^n}
$$
subject to the relations \eqref{alg:trivial1}--\eqref{alg:braid} of $\mathscr{R}_n(\mathscr O)$ together with the additional relation
\begin{equation}\label{alg:additional}
e(\bi)y_1^{(\Lambda,\alpha_{i_1})} = 0, \hspace*{10mm}\text{for each $\bi\in I^n$}.
\end{equation}

\subsection{The (graded) cellular algebras and the symmetric groups}

Following Graham and Lehrer \cite{GL}, we now introduce the graded cellular algebras. Reader may also refer to Hu-Mathas~\cite{HuMathas:GradedCellular}. Let $\mathscr O$ be a commutative ring with $1$ and let $A$ be a unital $\mathscr O$-algebra.

\begin{Definition} \label{definition: graded cellular basis}
A \textbf{graded cell datum} for $A$ is a triple $(\Lambda,T,C,\deg)$ where $\Lambda = (\Lambda,>)$ is a poset, either finite or infinite, and $T(\lambda)$ is a finite set for each $\lambda \in \Lambda$, $\deg$ is a function from $\coprod_\lambda T(\lambda)$ to $\Z$, and
$$
C\map{\prod_{\lambda\in\Lambda} T(\lambda)\times T(\lambda)}{A}
$$
is an injective map which sends $(s,t)$ to $a_{st}^\lambda$ such that:

(a)  $\set{ a_{\s\t}^\lambda|\lambda\in\Lambda,s,t\in T(\lambda)}$ is an $\mathscr O$-free basis of $A$;

(b) for any $r\in A$ and $\t\in T(\lambda)$, there exists scalars $c_{\t}^v(r)$ such that, for any $\s\in T(\lambda)$,
$$
a_{\s\t}^\lambda{\cdot}r \equiv \sum_{\v\in T(\lambda)} c_{\t}^{\v}(r) a_{\s\v}^\lambda \mod{A^{>\lambda}}
$$
where $A^{>\lambda}$ is the $\mathscr O$-submodule of $A$ spanned by $\set{ a_{\x\y}^\mu|\mu>\lambda, \x,\y\in T(\mu)}$;

(c) the $\mathscr O$-linear map $*\map{A}{A}$ which sends $a_{\s\t}^\lambda$ to $a_{\t\s}^\lambda$, for all $\lambda\in \Lambda$ and $\s,\t\in T(\lambda)$, is an anti-isomorphism of~$A$.

(d) each basis element $a_{\s\t}^\lambda$ is homogeneous of degree $\deg a_{\s\t}^\lambda = \deg(\s) + \deg(\t)$, for $\lambda\in \Lambda$ and all $\s,\t\in T(\lambda)$.
\end{Definition}

If a graded cell datum exists for $A$ then $A$ is a \textbf{graded cellular algebra}. Similarly, by forgetting the grading we can define a \textbf{cell datum} and hence a \textbf{cellular algebra}.

Suppose $A$ is a graded cellular algebra with graded cell datum $(\Lambda,T,C,\deg)$. For any $\lambda \in \Lambda$, define $A^{\geq \lambda}$ to be the $\mathscr O$-submodule of $A$ spanned by
$$
\set{c_{\s\t}^\mu|\mu \geq \lambda, \s,\t\in T(\mu)}.
$$

Then $A^{>\lambda}$ is an ideal of $A^{\geq\lambda}$ and hence $A^{\geq\lambda}/A^{>\lambda}$ is a $A$-module. For any $\s\in T(\lambda)$ we define $C_\s^\lambda$ to be the $A$-submodule of $A^{\geq\lambda}/A^{>\lambda}$ with basis $\set{a^\lambda_{\s\t} + A^{>\lambda}|\t\in T(\lambda)}$. By the cellularity of $A$ we have $C_\s^\lambda \cong C_\t^\lambda$ for any $\s,\t\in T(\lambda)$.

\begin{Definition} \label{cell module}
Suppose $\lambda \in \P_n$. Define the \textbf{cell module} of $A$ to be $C^\lambda = C_\s^\lambda$ for any $\s\in T(\lambda)$, which has basis $\set{a^\lambda_\t|\t\in T(\lambda)}$ and for any $r \in A$,
$$
a^\lambda_\t{\cdot}r= \sum_{\u\in T(\lambda)} c_\u^r a_\u^\lambda
$$
where $c_\u^r$ are determined by
$$
a^\lambda_{\s\t}{\cdot}r = \sum_{\u\in T(\lambda)} c_\u^r a^\lambda_{\s\u} + A^{>\lambda}.
$$
\end{Definition}

We can define a bilinear map $\langle\cdot,\cdot\rangle\map{C^\lambda\times C^\lambda}{\Z}$ such that
$$
\langle a^\lambda_\s,a^\lambda_\t\rangle a^\lambda_{\u\v} = a^\lambda_{\u\s}a^\lambda_{\t\v} + A^{>\lambda}
$$
and let $\text{rad }C^\lambda = \set{\s\in C^\lambda|\langle \s,\t\rangle = 0\text{ for all }\t\in C^\lambda}$. The $\text{rad }C^\lambda$ is a graded $A$-submodule of $C^\lambda$.

\begin{Definition}\label{D lambda}
Suppose $\lambda \in \P_n$. Let $D^\lambda = C^\lambda/\text{rad }C^\lambda$ as a graded $A$-module.
\end{Definition}

Exactly as in the ungraded case~\cite[Theorem 3.4]{GL} or~\cite[Theorem 2.10]{HuMathas:GradedCellular}, we obtain the following:\\

\begin{Theorem}~\label{all graded simple in CA}
The set $\set{D^\lambda\langle k\rangle|\lambda\in \Lambda, D^\lambda\neq 0, k\in\Z}$ is a complete set of pairwise non-isomorphic graded simple $A$-modules.
\end{Theorem}

We give an example of graded cellular algebras here, which is called the cyclotomic Hecke algebras.

Let $\F_p$ be a fixed field of characteristic $p \geq 0$ with $q\in \F_p^\times$. Let $e$ be the smallest positive integer such that $1 + q + \ldots + q^{e-1} = 0$ and setting $e = 0$ if no such integer exists. Then define $I = \Z/e\Z$ if $e > 0$ and $I = \Z$ if $e = 0$.

For $n \geq 0$, assume that $q = 1$. Let $H_n$ be the \textbf{degenerate affine Hecke algebra}, working over $\F_p$. So $H_n$ has generators
$$
\{x_1,\ldots,x_n\}\cup\{s_1,\ldots,s_{n-1}\}
$$
subject to the following relations
$$
\begin{array}{cc}
x_r x_s = x_s x_r; & \\
s_r x_{r+1} = x_r s_r + 1, & s_r s_x = x_s s_r \hspace*{5mm}\text{if $s \neq r,r+1$}\\
s_r^2 = 1; \\
s_r s_{r+1}s_r = s_{r+1}s_r s_{r+1}, & s_r s_t = s_t s_r \hspace*{5mm}\text{if $|r - t| > 1$}
\end{array}
$$

Now we assume that $q\neq 1$ and $H_n$ be the \textbf{non-degenerate affine Hecke algebra} over $\F_p$. So $H_n$ has generators
$$
\{ X_1^{\pm 1},\ldots,X_n^{\pm 1}\}\cup\{T_1,\ldots,T_{n-1}\}
$$
subject to the following relations
$$
\begin{array}{cc}
X_r^{\pm 1} X_s^{\pm 1} = X_s^{\pm 1} X_r^{\pm1}, & X_r X_r^{-1} = 1;\\
T_r X_r T_r = q X_{r+1}, & T_r X_s = X_s T_r\hspace*{5mm}\text{if $s\neq r,r+1$};\\
T_r^2 = (q-1)T_r + q; &\\
T_r T_{r+1}T_r = T_{r+1}T_rT_{r+1},& T_r T_s = T_s T_r\hspace*{5mm}\text{if $|r-s| > 1$}.
\end{array}
$$

Then for any $\Lambda \in P_+$, we define
\begin{equation} \label{equation: quotient hecke}
H_n^\Lambda = \begin{cases}
H_n/ \langle \prod_{i\in I}(X_1 - q^i)^{(\Lambda,\alpha_i)}\rangle, & \text{if $q\neq 1$,}\\
H_n/ \langle \prod_{i\in I}(x_1 - i)^{(\Lambda,\alpha_i)}\rangle, & \text{if $q = 1$.}
\end{cases}
\end{equation}
and we call $H_n^\Lambda$ the \textbf{degenerate cyclotomic Hecke algebra} if $q = 1$ and \textbf{non-degenerate cyclotomic Hecke algebra} if $q\neq 1$.

By the definitions, degenerate and non-degenerate cyclotomic Hecke algebras are similar with some minor difference. In order to minimize their difference we define
\begin{equation} \label{equation: q_i}
q_i = \begin{cases}
i, & \text{if $q = 1$,}\\
q^i, & \text{if $q\neq 1$.}
\end{cases}
\end{equation}
and use $x_r$ instead of $X_r$ when we don't have to distinguish which case we are working with. Hence we can re-write (\ref{equation: quotient hecke}) as
\begin{equation} \label{equation: cyclotomic Hecke algebras}
H_n^\Lambda = H_n/\langle \prod_{i\in I}(x_1 - q_i)^{(\Lambda,\alpha_i)}\rangle.
\end{equation}

Murphy~\cite{Murphy:basis} gave a set of cellular basis for $H_n^\Lambda$ which shows that $H_n^\Lambda$ is a cellular algebra. Brundan and Kleshchev~\cite{BK:GradedKL} proved the remarkable result that every $H_n^\Lambda$ over $\F_p$ is isomorphic to $\R(\F_p)$ introduced in \autoref{notation: KLR}, where in both algebras $\Lambda$ and $e$ are the same. Therefore when $H_n^\Lambda$ is over a field it is a graded cellular algebra.

Let $\mathfrak{S}_n$ be the symmetric group on $\{1,2,\ldots, n\}$. Then $\mathfrak{S}_n$ is a Coxeter group and $\{s_1, \ldots, s_{n-1}\}$ is its standard set of Coxeter generators, where $s_i = (i,i+1)$ for $i = 1, 2, \ldots, n-1$. Suppose $w$ is an element of $\mathfrak{S}_n$ and $w = s_{i_1} s_{i_2}\ldots s_{i_m}$. If $m$ is minimal we say that $w$ has \textbf{length} $m$ \label{notation: length} and write $l(w) = m$. In this case we say $s_{i_1}s_{i_2}\ldots s_{i_m}$ is a \textbf{reduced expression} of $w$. In general an element of $\mathfrak{S}_n$ has more than one reduced expressions. For example, we have $w = s_1 s_2 s_1 = s_2 s_1 s_2$. Nonetheless, all the reduced expression of an element have the same length.

In this paper we let $\mathfrak{S}_n$ act on $\{1,2,\ldots,n\}$ from right. For example, $(i) s_i s_{i+1} = (i + 1) s_{i+1} = i+2$. The following result is well-known. See, for example, \cite[Corollary~1.4]{M:ULect}.

\begin{Proposition} \label{change length of w}
Suppose that $w\in \mathfrak{S}_n$. For $i = 1,2,\ldots,n-1$,
$$
l(w s_i) = \begin{cases}
  l(w) + 1, & \text{ if $(i)w < (i+1)w$,}\\
  l(w) - 1, & \text{ if $(i)w > (i+1)w$.}
\end{cases}
$$
\end{Proposition}

We recall the definition of the \textbf{Bruhat order} $\leq$ on $\mathfrak{S}_n$. For $u,w\in \mathfrak{S}_n$ define $u \leq w$ if $u = s_{r_{a_1}}s_{r_{a_2}} \ldots s_{r_{a_b}}$ for some $1 \leq a_1 < a_2 < \ldots < a_b \leq m$, where $w = s_{r_1}s_{r_2}\ldots s_{r_m}$ is a reduced expression for $w$.

%
%
%

\subsection{Tableaux combinatorics}

In this subsection we recall the combinatorics of (multi)partitions and (multi)tableaux that we will need in this paper.

  Let $n$ be a positive integer. A \textbf{composition} of~$n$ is an ordered sequence of nonnegative integers $\lambda = (\lambda_1,\lambda_2,\ldots)$ and $|\lambda|=\sum_{i = 1}^\infty \lambda_i = n$. We say $\lambda$ is a \textbf{partition} of $n$ if $\lambda = (\lambda_1,\lambda_2,\ldots)$ is a composition and $\lambda_1 \geq \lambda_2 \geq \lambda_3 \geq \ldots $. We can then identify $\lambda$ with a sequence $(\lambda_1,\ldots,\lambda_k)$ whenever $\lambda_i = 0$ for $i > k$.

As we now recall, there is a natural partial ordering on the set of compositions of $n$. Suppose $\lambda = (\lambda_1,\lambda_2,\dots,)$ and $\mu = (\mu_1,\mu_2,\dots)$ are compositions of~$n$. Then
$\lambda$ \textbf{dominates} $\mu$, and we write $\lambda \unrhd \mu$, if
$$
\sum_{i = 1}^k \lambda_i \geq \sum_{i = 1}^k \mu_i
$$
for any $k$. We write $\lambda \rhd \mu$ if $\lambda \unrhd \mu$ and $\lambda \neq \mu$. The dominance ordering can be extended to a total ordering~$\ge$, called the \textbf{lexicographic ordering}. We write $\lambda > \mu$ if there exist some $k$, such that $\lambda_i = \mu_i$ for all $i < k$ and $\lambda_k > \mu_k$. Define $\lambda \geq \mu$ if $\lambda > \mu$ or $\lambda = \mu$. Then $\lambda \unrhd \mu$ implies $\lambda \geq \mu$.

A \textbf{multicomposition} of $n$ of \textbf{level} $\l$ is an ordered sequence $\lambda = (\lambda^{(1)},\ldots,\lambda^{(\l)})$ of compositions such that $\sum_{i = 1}^\l|\lambda^{(i)}| = n$. Similarly, a \textbf{multipartition} of level~$\l$ is multicomposition $\lambda = (\lambda^{(1)},\ldots,\lambda^{(\l)})$ of~$n$ such that each $\lambda^{(i)}$ is a partition. We will identify multicompositions and multipartitions of level~$1$ with compositions and partitions in the obvious way.

Let $\mathscr{C}_n^\Lambda$ \label{notation: multicomposition}be the set of all multicomposition of $n$ and $\mathscr{P}_n^\Lambda$ \label{notation: multipartition}be the set of all multipartitions of $n$. We can extend the dominance ordering to $\mathscr C_n^\Lambda$ by defining $\lambda \unrhd \mu$ if
$$
\sum_{i = 1}^{k-1}|\lambda^{(i)}| + \sum_{j = 1}^s \lambda^{(k)}_j \geq\sum_{i = 1}^{k-1}|\mu^{(i)}| + \sum_{j = 1}^s \mu^{(k)}_j
$$
for any $1 \leq k \leq \l$ and all $s \geq 1$. Again, we write $\lambda \rhd \mu$ if $\lambda \unrhd \mu$ and $\lambda \neq \mu$. Similarly, we extend the lexicographic ordering $\lambda > \mu$ and $\lambda \geq \mu$ to $\mathscr{C}^\Lambda_n$ in the obvious way way.

The \textbf{Young diagram} of a multicomposition $\lambda$ of level~$\l$ is the set
$$
[\lambda] = \{ (r,c,l)\ |\ 1 \leq c \leq \lambda_r^{(l)}, r \geq 0\text{ and }1\leq l \leq \l\}
$$
which we think of as an ordered $\l$-tuple of the diagrams of the partitions $\lambda^{(1)},\ldots,\lambda^{(\l)}$. The triple $(r,c,l)\in[\lambda]$ is \textbf{node} of~$\lambda$ in row~$r$, column~$c$ and component~$l$. A \textbf{$\lambda$-tableau} is any bijection $\t\map{[\lambda]}{\{1,2,\ldots,n\}}$. We identify a $\lambda$-tableau $\t$ with a labeling of the diagram of $\lambda$. That is, we label the node $(r,c,l)\in[\lambda]$ with the integer $\t(r,c,l)$. For example,
$$\Bigg(\hspace*{1mm}\tab(1234,567,8)\hspace*{1mm}\Bigg|\hspace*{1mm}%
    \tab(9\ten,\eleven\twelve,\thirteen)\hspace*{1mm}\Bigg|\hspace*{1mm}\tab(\fourteen\fifteen\sixteen)\hspace*{1mm}\Bigg)
$$
is a $(4,3,1|2^2,1|3)$-tableaux.
If $\t$ is a $\lambda$-tableau then the \textbf{shape} of $\t$ is the multicomposition $\lambda$ and we write $\Shape(\t) = \lambda$. A $\lambda$-tableau $\t$ is \textbf{standard} if
$\lambda=\Shape(\t)$ is a multipartition and if, in each component, the entries increase along each row and down each column. More precisely, if $(r,c,l)\in[\lambda]$ then
$\t(r,c,l)<\t(r+1,c,l)$ whenever $(r+1,c,l)\in[\lambda]$ and
$\t(r,c,l)<\t(r,c+1,l)$ whenever $(r,c+1,l)\in[\lambda]$. Let $\Std(\lambda)$ be the set of all standard $\lambda$-tableaux and $\Std(>\lambda)$ be the set of all standard $\mu$-tableaux with $\mu > \lambda$. We can define $\Std(\geq\lambda)$ similarly. Note that if $\t$ is standard then so is $\t|_k$ for $1\le k\le n$.

If $\t\in\Std(\lambda)$ and $1\le k\le n$ define $\t|_k$ to be the subtableau of $\t$ obtained by removing all the nodes containing an entry greater than $k$. We define an analogue of the dominance ordering for standard tableaux by defining $\t \unrhd \s$ if $\Shape(\t|_k) \unrhd \Shape(\s|_k)$, for $1\le k\le n$. As with the dominance ordering,  if $\t\unrhd\s$ then we write $\s\unlhd\t$ and if $\s\ne\t$ then write $\t\rhd\s$ and $\s\lhd\t$. We also define $(\s,\t) \rhd (\u,\v)$ if $\s\unrhd \u$, $\t \unrhd \v$ and $(\s,\t) \neq (\u,\v)$.

For any multicomposition $\lambda$, define $\t^\lambda$ \label{notation: tlambda}to be the unique $\lambda$-tableau such that $\t^\lambda \unrhd \t$ for all standard $\lambda$-tableau $\t$. For example, if $\lambda=(4,3,1|2^2,1|3)$ then $\t^\lambda$ is the tableau displayed above.

The symmetric group acts on the set of all $\lambda$-tableaux. Let $\t$ be a $\lambda$-tableau, then $\t{\cdot}s_r$ is the tableau obtained by exchanging the entries $r$ and $r+1$ in $\t$, i.e. $(r)\t^{-1} = (r+1)(\t{\cdot}s_r)^{-1}$, $(r+1)\t^{-1} = (r)(\t{\cdot}s_r)^{-1}$, and $(k)\t^{-1} = (k)(\t{\cdot}s_r)^{-1}$ for $k\neq r,r+1$. Then for each $\lambda$-tableau $\t$ let $d(\t)$ \label{notation: dt}be the permutation in $\mathfrak{S}_n$ such that $\t^\lambda{\cdot}d(\t) = \t$.

Recall the Bruhat order $\leq$ on $\mathfrak{S}_n$ from subsection~1.1. The following result, which goes back to work of Ehresmann and James, is part of the folklore for these algebras. The proof for level $1$ can be found from \cite[Lemma~3.7]{M:ULect}. The higher level cases follow easily.

\begin{Lemma}
Suppose $\lambda\in\mathscr P_n^\Lambda$ and $\s$ and $\t$ are standard $\lambda$-tableaux. Then $\s \unrhd \t$ if and only if $d(\s) \leq d(\t)$.
\end{Lemma}

Suppose $\lambda$ is a multicomposition and $\gamma = (r,c,l) \in [\lambda]$ and recall from subsection~1.1 that $\kappa_\Lambda = (\kappa_1,\kappa_2,\ldots,\kappa_\l)$ is a fixed multicharge of~$\Lambda$. The \textbf{residue} of $\gamma$ associate to $\kappa_\Lambda$ is
$$
\res(\gamma) \equiv r - c + \kappa_l \pmod{e}.
$$

If $\t$ is a standard $\lambda$-tableau and the \textbf{residue sequence} of $\t$ is $\res(\t) = \bi_t = (i_1, i_2, i_3, \ldots, i_n)$, where $i_k = \res(\gamma_k)$ and $\gamma_k$ is the unique node in $[\lambda]$ such that $\t(\gamma_k) = k$. In particular, we write $\bi_{\t^\lambda} = \bi_\lambda$
and $\res_\t(k)=\res(\gamma_k)$. \label{notation: bi lambda}

Recall that for each standard tableau $\t$, we can define a permutation $d(\t) \in \mathfrak S_n$ such that $\t = \t^\lambda{\cdot}d(\t)$. For each permutation we may have more than one reduced expression. Here we fix a choice of the reduced expression of $d(\t)$.

For any standard $\lambda$-tableau $\t$ and $1 \leq i \leq n+1$, define $\lambda^{(i)} = \Shape(\t|_{i-1})$ and $\t^{(i)}$ to be a standard $\lambda$-tableau where $\t^{(i)}|_{i-1} = \t^{\lambda^{(i)}}$, and ${\t^{(i)}}^{-1}(k) = \t^{-1}(k)$ for any $i \leq k \leq n$. In particular, $\t^{(1)} = \t$ and $\t^{(n+1)} = \t^\lambda$. Therefore we have a series of standard $\lambda$-tableau
$$
\t^\lambda = \t^{(n+1)}, \t^{(n)}, \t^{(n-1)},\ldots, \t^{(2)},\t^{(1)} = \t.
$$

Define $w_i$ to be the unique permutation in $\mathfrak{S}_n$ such that $\t^{(i+1)} w_i = \t^{(i)}$. For each $w_i \neq 1$, we can write $w_i = s_{a_i}s_{a_i+1}s_{a_i+2}\ldots s_{i-2}s_{i-1}$ for some $a_i \leq i-1$. Notice that
\begin{equation} \label{remark: standard expression}
(k)(\t^{(i+1)})^{-1} =
\begin{cases}
(i)(\t^{(i)})^{-1}, & \text{ if $k = a_i$,}\\
(k-1)(\t^{(i)})^{-1}, & \text{ if $a_i < k \leq i$,}\\
(k)(\t^{(i)})^{-1}, & \text{ otherwise,}
\end{cases}
\end{equation}
and $l(w_i)$ is always greater than or equal to the length of the row containing $i$ in $\t^{(i+1)}$. Also for each $i$, if $\Shape(\t^{(i)}|_{i-1}) = \lambda$, then $\t^{(i)}|_{i-1} = \t^\lambda$.

\begin{Example}
Suppose $\t^{(11)} = \tab(1234\twelve,567\eleven,89\ten\thirteen,\fourteen\fifteen)$ and $\t^{(10)} = \tab(123\ten\twelve,456\eleven,789\thirteen,\fourteen\fifteen)$. Therefore we have $w_{10} = s_4 s_5 s_6 s_7 s_8 s_9$ such that $\t^{(11)}{\cdot}w_{10} = \t^{(10)}$.

Notice that in this case, $i = 10$ and $a_{10} = 4$. So
$$
(k)(\t^{(11)})^{-1} = \begin{cases}
(10)(\t^{(10)})^{-1}, &\text{ if $k = a_i = 4$,}\\
(k-1)(\t^{(10)})^{-1}, &\text{ if $4 = a_i < k \leq i = 10$,}\\
(k)(\t^{(10)})^{-1}, &\text{ otherwise.}
\end{cases}
$$

Furthermore, $\t^{(11)}|_{10} = \tab(1234,567,89\ten) = \t^{(4,3,3)}$ and $\t^{(10)}|_9 = \tab(123,456,789) = \t^{(3,3,3)}$.
\end{Example}

One can see that $d(\t) = w_n w_{n-1} \ldots w_2 w_1$ and define it to be the \textbf{standard expression} of $d(\t)$. By ~\autoref{change length of w} one can verify that standard expression is a reduced expression of $d(\t)$. In the rest of this paper, we fix $d(\t)$ to be its standard expression.

\begin{Lemma}\label{standard: short}
Suppose $\t$ is a standard $\lambda$-tableau and $d(\t) = s_{r_1}\ldots s_{r_m}$ is the standard expression. For any $1 \leq k \leq m$, define $\s = \t^\lambda {\cdot}s_{r_1}s_{r_2}\ldots s_{r_k}$. Then $\s$ is a standard $\lambda$-tableau.
\end{Lemma}

\proof The proof is trivial by the definition of the standard expression. \endproof

\begin{Example}
Suppose $\t = \tab(12467,35)$. Then we have $d(\t) = s_5 s_6{\cdot}s_4 s_5{\cdot}s_3$. Then
\begin{eqnarray*}
\t^\lambda{\cdot}s_5 & = & \tab(12346,57),\\
\t^\lambda{\cdot}s_5s_6 & = & \tab(12347,56),\\
\t^\lambda{\cdot}s_5s_6s_4 & = & \tab(12357,46),\\
\t^\lambda{\cdot}s_5s_6s_4s_5 & = & \tab(12367,45),
\end{eqnarray*}
and the above tableaux are all standard.
\end{Example}

\subsection{Graded cellular basis of KLR algebras over a field}\label{S:PsiBasis}

Suppose $\mathscr O$ is a field. Hu and Mathas \cite[Theorem 5.8]{HuMathas:GradedCellular} constructed a graded cellular basis of $\R(\mathscr O)$. Here we give an equivalent definition of their basis. For any multicomposition $\lambda$, recall $\t^\lambda$ to be the unique standard $\lambda$-tableau such that $\t^\lambda \unrhd \t$ for all standard $\lambda$-tableau $\t$, and $\bi_\lambda$ is the residue sequence of $\t^\lambda$. We define $\hat e_\lambda = \hat e(\bi_\lambda)$. \label{notation: e lambda}

Suppose $\lambda$ is a multicomposition. A node $(r,c,l)$ is an \textbf{addable node} of $\lambda$ if $(r,c,l)\not\in [\lambda]$ and $[\lambda]\cup \{(r,c,l)\}$ is the Young diagram of a multipartition. Similarly, a node $(r,c,l)$ is a \textbf{removable node} of $\lambda$ if $(r,c,l)\in [\lambda]$ and $[\lambda]\backslash \{(r,c,l)\}$ is the Young diagram of a multipartition. Given two nodes $\alpha = (r,c,l)$ and $\beta = (s,t,m)$ then $\alpha$ is \textbf{below} $\beta$ if either $l > m$, or $l = m$ and $r > s$.

Suppose that $\s\in\Std(\lambda)$. Let $\mathscr{A}_\s(k)$ \label{notation: addable1}be the set of addable nodes of the multicomposition $\Shape(\s|_k)$ which are below $\s^{-1}(k)$ and let
\begin{equation*}\label{notation: addable2}
\mathscr{A}_\s^\Lambda(k) = \{ \alpha \in \mathscr{A}_\s(k)\ |\ \res(\alpha) = \res_\t(k)\}.
\end{equation*}

Similarly as in \cite[Definition~4.12]{HuMathas:GradedCellular}, define
\begin{equation*} \label{notation: y lambda}
\hat y_\lambda = \prod_{k = 1}^n
\hat y_k^{|\mathscr{A}_{\t^\lambda}^\Lambda(k)|}\in\Rn(\mathscr O).
\end{equation*}

For example, if $\lambda = (3,1|4^2,2|5,1)$, $e = 4$ and $\Lambda = 3\Lambda_0$ then
$$
\t^\lambda = \Bigg(\hspace*{1mm}\tab(123,4)\hspace*{1mm}\Bigg|\hspace*{1mm}%
    \tab(5678,9\ten\eleven\twelve,\thirteen\fourteen)\hspace*{1mm}\Bigg|\hspace*{1mm}\tab(\fifteen\sixteen\seventeen\eighteen\nineteen,\twenty)\hspace*{1mm}\Bigg)
$$
and $\hat y_\lambda = \hat y_1^2 \hat y_5\hat y_8\hat y_{10}\hat y_{12}\hat y_{18}$. Therefore,
$$
\hat e_\lambda \hat y_\lambda = \hat e(0123012330122012303)\hat y_1^2 \hat y_5\hat y_8\hat y_{10}\hat y_{12}\hat y_{18}.
$$

We define a particular kind of element in $\Rn(\mathscr O)$. Suppose $w \in \mathfrak S_n$ has length $\l$ and $s_{i_1} s_{i_2} \ldots s_{i_\l}$ is a reduced expression for $w$ in $\mathfrak{S}_n$. Recall that $\Rn(\mathscr O)$ has a unique anti-isomorphism $*$ which fixes all of the KLR generators. Define
$$
\hat\psi_{w} = \hat\psi_{i_1} \hat\psi_{i_2} \ldots \hat\psi_{i_\l}\in \Rn(\mathscr O) \hspace*{5mm} \text{and} \hspace*{5mm}\hat\psi_{w}^* = \hat\psi_{i_\l}\hat\psi_{i_{\l-1}}\ldots \hat\psi_{i_2}\hat\psi_{i_1}\in \Rn(\mathscr O).
$$

Notice that $\hat \psi_{w}$ and $\hat \psi^*_{w}$ depend on the choice of the reduced expression of $w$, even though in $\mathfrak{S}_n$ all reduced expressions of $w$ are the same. For example, $s_1 s_2 s_1$ and $s_2 s_1 s_2$ are equal to the same element of $\mathfrak{S}_n$, but in general $\hat \psi_1\hat \psi_2\hat \psi_1 \neq \hat \psi_2\hat \psi_1\hat \psi_2$ in $\Rn(\mathscr O)$. Define $l(\hat\psi_{w}) = l(\hat\psi_{w}^*) = l(w)$ for any standard tableau $\t$. Similarly we can define
$$
\psi_{w} = \psi_{i_1} \psi_{i_2} \ldots \psi_{i_\l}\in \R(\mathscr O) \hspace*{5mm} \text{and} \hspace*{5mm}\psi_{w}^* = \psi_{i_\l}\psi_{i_{\l-1}}\ldots \psi_{i_2}\psi_{i_1}\in \R(\mathscr O)
$$
and $\psi_w$ and $\psi_w^*$ depends on the choice of reduced expressions of $w$ as well.

Suppose $l(d(\t)) = \l$ and $d(\t) = s_{i_1} s_{i_2} \ldots s_{i_\l}$ is the standard expression of $d(\t)$ where $\t^\lambda{\cdot}d(\t) = \t$. Define $\hat \psi_{d(\t)} = \hat  \psi_{i_1}\hat \psi_{i_2}\ldots \hat \psi_{i_\l}$ and $\hat \psi^*_{d(\t)} = \hat \psi_{i_\l}\hat \psi_{i_{\l-1}}\ldots \hat \psi_{i_2}\hat \psi_{i_1}$.

\begin{Definition} \label{notation: psi basis}
Suppose $\Lambda\in P_+$, $\lambda \in \mathscr{P}_n^\Lambda$ and $\s,\t$ are two standard $\lambda$-tableaux. We define
$$
\hat\psi^{\mathscr O}_{\s\t} = \hat\psi_{d(\s)}^* \hat e_\lambda \hat y_\lambda \hat \psi_{d(\t)} \in \Rn(\mathscr O),
$$
and hence
$$
\psi^{\mathscr O}_{\s\t} = \hat\psi^{\mathscr O}_{\s\t} + N_n^\Lambda \in \R(\mathscr O).
$$
\end{Definition}

\begin{Remark}
Notice that \protect{Hu and Mathas~\cite[Definition~5.1]{HuMathas:GradedCellular}} defined $\psi_{\s\t}^{\mathscr O}$ differently. Actually if we define $e_\lambda$, $y_\lambda$ and $\psi_w$ in $\R(\mathscr O)$ as analogues of $\hat e_\lambda$, $\hat y_\lambda$ and $\hat\psi_w$, and define $\psi_{\s\t}^{\mathscr O} = \psi_{d(\s)}^* e_\lambda y_\lambda \psi_{d(\t)}$ for $\s,\t\in\Std(\lambda)$, it is equivalent to \autoref{notation: psi basis}. We define $\psi_{\s\t}^{\mathscr O}$ as in \autoref{notation: psi basis} because we need to work in $\Rn(\mathscr O)$ later.
\end{Remark}

\begin{Remark}
By construction, then this $\psi_{\s\t}^{\mathscr O}$ is well defined as an element of $\R(\mathscr O)$ for any ring~$\mathscr O$. Many of the calculations in this paper depend heavily on the choice of $\mathscr O$ so we write $\psi_{\s\t}^{\mathscr O}$ to emphasize that $\psi_{\s\t}^{\mathscr O}$ is an element of~$\R(\mathscr O)$. Most of the time, however, we will work in $\R(\Z)$ so for convenience we set $\psi_{\s\t}=\psi_{\s\t}^\Z$.
\end{Remark}

\begin{Lemma} [\protect{Hu and Mathas~\cite[Lemma~5.2]{HuMathas:GradedCellular}~\cite[Corollary 3.11,3.12]{HuMathas:GradedInduction}}] \label{residue sequence has to be the same}
Suppose $\mathscr O$ is a field and $\s$ and $\t$ are standard $\lambda$-tableaux and $1\leq r\leq n$,
\begin{eqnarray*}
\psi_{\s\t} \psi_r & = & \begin{cases}\sum_{(\u,\v)\rhd(\s,\t)}c_{\u\v} \psi_{\u\v}, & \text{if $\t{\cdot}s_r$ is not standard}\\
&\hspace*{5mm}\text{or $d(\t){\cdot}s_r$ is not reduced,}\\
\psi_{\s\v} + \sum_{(\u,\v)\rhd (\s,\t)}c_{\u\v} \psi_{\u\v}, & \text{if $\v = \t{\cdot}s_r$ standard and $d(\t){\cdot}s_r = d(\v)$.}\end{cases}
\end{eqnarray*}
for $c_{\u\v}\in\mathscr O$, and $c_{\u\v} \neq 0$ only if $\res(\s) = \res(\u)$ and $\res(\t{\cdot}s_r) = \res(\v)$. Similarly, we have
$$
\psi_{\s\t}^{\mathscr O}y_r = \sum_{(\u,\v)\rhd (\s,\t)} c_{\u\v}\psi_{\u\v}^{\mathscr O}
$$
for $c_{\u\v}\in\mathscr O$, and $c_{\u\v}\neq 0$ only if $\res(\s) = \res(\u)$ and $\res(\t) = \res(\v)$.

\end{Lemma}

\begin{Theorem}[\protect{Hu and Mathas~\cite[Theorem 5.14]{HuMathas:GradedCellular}}] \label{basis with field}
Suppose $\mathscr O$ is an integral domain and that either $e=0$, $e$ is a prime
or $e$ is a non-zero non-prime integer such that $e\cdot 1_{\mathscr O}$ is
invertible in $\mathscr O$. Then
$$
\{ \psi_{\s\t}^{\mathscr O} \ |\ \s,\t \in \text{Std}(\lambda) \text{ for $\lambda\in\mathscr{P}_n^\Lambda$} \}
$$
is a graded cellular basis of $\R(\mathscr O)$. In
particular, $\R(\mathscr O)$ is free as an $\mathscr
O$-module of rank $\l^n n!$.
\end{Theorem}

The main purpose of this paper is to prove that $\R(\Z)$ is free of rank~$\l^nn!$. To do this we will show that $\{\psi_{\s\t}^\Z\ |\ \s,\t \in \text{Std}(\lambda) \text{ for $\lambda\in\mathscr{P}_n^\Lambda$} \}$ is a homogeneous basis of $\R(\Z)$.\\

We define some notation for future use.

\begin{Definition} \label{Definition: B}
Suppose $\lambda$ is a multipartition of $\mathscr{P}_n^\Lambda$. Define:
\begin{align*}
  \BZ&=\langle\psi_{\s\t}\mid \s,\t\in\text{Std}(\mu)\text{ for }\mu\in\P_n\rangle_\Z ,\\
\Bgelam&=\langle\psi_{\s\t}\mid \s,\t\in\text{Std}(\mu)\text{ and }\mu\ge\lambda
                \text{ for }\mu\in\P_n\rangle_\Z ,\\
  \Blam & = \langle\psi_{\s\t}\mid \s,\t\in\text{Std}(\mu)\text{ and }\mu>\lambda
        \text{ for }\mu\in\P_n\rangle_\Z .
\end{align*}
where $\Blam \subseteq \Bgelam \subseteq \BZ \subseteq \R(\Z)$
\end{Definition}

This subsection closes with an important Proposition:


Consider the quiver Hecke algebra  $\mathscr{R}_n(\Q)$ defined over the rational field $\Q$. We have $\mathscr{R}_n(\Q) \cong \mathscr{R}_n(\Z)\otimes\Q$ and we can define a linear map $f\map{\mathscr{R}_n(\Z)}{\mathscr{R}_n(\Q)}$ by sending $x\in \mathscr{R}_n(\Z)$ to $x\otimes 1$.

\begin{Lemma}

The linear map $f\map{\mathscr{R}_n(\Z)}{\mathscr{R}_n(\Q)}$ is an injection.

\end{Lemma}

\begin{proof}
In \cite{KhovLaud:diagI}\cite{KhovLaud:diagII}, Khovanov and Lauda constructed a basis of $\Rn(\mathscr O)$
\begin{equation} \label{basis of KLR}
\set{ \hat e(\bi)\hat y_1^{\l_1} \hat y_2^{\l_2} \ldots \hat y_n^{\l_n} \hat \psi_w |\bi\in I^n, w\in\mathfrak{S}_n, \l_1,\l_2,\ldots,\l_n \geq 0 }
\end{equation}
for any ring $\mathscr O$. Hence that $\mathscr{R}_n(\Z)$ is free over $\Z$. The Lemma follows immediately.
\end{proof}

From the definitions, it is evident that $f(N^\Lambda_n(\Z)) \subseteq N^\Lambda_n(\Q)$. Hence, $f$ induces a homomorphism,
$$f\map{\R(\Z)}\R(\Q);
      x+N^\Lambda_n(\Z) \mapsto f(x)+N^\Lambda_n(\Q),$$
which by abuse of notation we also denote by~$f$. In particular, observe that
$f(\psi_{\s\t}^\Z) = \psi_{\s\t}^\Q$. The main Theorem of this paper is equivalently to prove that $f\map{\R(\Z)}\R(\Q)$ is an injection.

We then introduce an important special case where we already know that~$f$ is injective.

\begin{Proposition} \label{basis to hecke}
The homomorphism $f\map{\R(\Z)}\R(\Q)$ restricts to an injective map from $\BZ$ to $\R(\Q)$.
\end{Proposition}

\begin{proof}As we have already noted above, $f(\psi_{\s\t}^\Z) =\psi_{\s\t}^\Q$ for all $\s,\t\in\Std(\lambda)$ and $\lambda\in\mathscr{P}^\Lambda_n$. Hence, \autoref{basis with field} implies the result.
\end{proof}

\begin{Corollary} \label{linearly independent}
The elements $\{ \psi_{\s\t}^\Z \ |\ \s,\t \in \text{Std}(\lambda) \text{ for $\lambda\in\mathscr{P}_n^\Lambda$} \}$ are a linearly independent subset of $\R(\Z)$.
\end{Corollary}

\begin{Remark}
  Proposition~\ref{basis to hecke} is quite crucial. In this paper we prove that $\psi_{\s\t}^\Z{\cdot} \psi_r \in \BZ$ whenever $d(\t){\cdot}s_r$ is not reduced or $\t{\cdot} s_r$ is not standard in $\R(\Z)$. We can only have
$$
\psi_{\s\t}^\Z {\cdot}\psi_r = \sum_{\u,\v} c^\Z_{\u\v} \psi_{\u\v}^\Z.
$$

In $\R(\Q)$, however, by Lemma~\ref{residue sequence has to be the same}, under these conditions we have
$$
\psi_{\s\t}^\Q{\cdot}\psi_r = \sum_{(\u,\v)\rhd (\s,\t)}c^\Q_{\u\v} \psi_{\u\v}^\Q
$$
for some $c^\Q_{\u\v}\in\Q$, where $(\u,\v)\rhd (\s,\t)$ if $\u \unrhd \s$, $\v \unrhd \t$ and $(\u,\v) \neq (\s,\t)$. Therefore, $c^\Q_{\u\v} = c^\Z_{\u\v}$ by Proposition~\ref{basis to hecke} and we see that $c^\Z_{\u\v}\ne0$ only if $(\u,\v)\rhd (\s,\t)$. In such case we have much more information about $\u$ and $\v$ with $c^\Z_{\u\v}\neq 0$. Similar remarks apply to the products $\psi_{\s\t}^\Z{\cdot}y_r$.
\end{Remark}

\section{Integral Basis Theorem I}




In the next two sections we will prove that $\R$ is $\Z$-free by showing that $\set{\psi_{\s\t} | \s,\t \in \Std(\lambda), \lambda \in \P_n}$ spans $\R$ over $\Z$. Then by ~\autoref{linearly independent} it is a basis of $\R$ over $\Z$.

In the rest of this paper we write $\mathscr{R}_n(\Z)$ as $\mathscr{R}_n$ and $\mathscr{R}_n^\Lambda(\Z)$ as $\mathscr{R}_n^\Lambda$. Fix a weight $\Lambda$, a multicharge $\kappa_\Lambda = (\kappa_1,\ldots,\kappa_l)$ corresponding to $\Lambda$ and an integer $e > 2$. In this and the next section we mainly work with the algebra $\mathscr{R}^\Lambda_n$.

\subsection{The base step of the induction}

In this subsection we set up the notations and inductive machinery that we use in the next two sections to prove our main theorem. We then consider the base case of our induction which is when $\lambda = (n|\emptyset|\ldots|\emptyset)$. Finally we develop some technical Lemmas which will be useful later.

\begin{Definition} \label{notation: overrightarrow lambda}
  Suppose that $\lambda$ is a multipartition of~$n$. Let
 $\lambda^+$ be the multicomposition of $n+1$ obtained by adding a node at the end of the last non-empty row of $\lambda$, and $\lambda_- = \lambda|_{n-1}$ be the multipartition of $n-1$ obtained by removing the last node from $\lambda$.
\end{Definition}

For example, if $\lambda = (4,3|3,3)$ then $\lambda^+ = (4,3|3,4)$ and $\lambda_- = (4,3|3,2)$. Notice that in general, $\lambda^+$ will be a multicomposition rather than a multipartition.

For $k \in I$ and $\lambda \in \P_n$, define $\mathscr A_{\t^\lambda}^k = \set{\alpha \in \mathscr A_{\t^\lambda}(n)|\res(\alpha) = k}$. Recall $\bi_\lambda=\res(\t^\lambda)$ and $e_\lambda=e(\bi_\lambda)$ from
\autoref{S:PsiBasis}.

\begin{Definition} \label{notation: b k lambda}
Suppose that $\lambda\in \P_n$ and $k\in I$. Define the integer $b_k^\lambda$
by
$$
b_k^\lambda = \begin{cases}
        |\mathscr{A}_{\t^\lambda}^k| + 1,& \text{if $\lambda^+$ is a
                            multipartition and $i_n+1 = k$},\\
                            |\mathscr{A}_{\t^\lambda}^k|,& \text{otherwise.}
       \end{cases}
$$
\end{Definition}

If $\bi = (i_1, i_2, \ldots, i_n) \in I^n$ and $k\in I$ then define $\bi\vee k = (i_1, i_2, \ldots, i_n, k) \in I^{n+1}$. \label{notation: concatenate residues}

\begin{Lemma} \label{I-problem: downstair cases}
  Suppose that $\lambda\in \P_n$ and $k\in I$. Then for each
  integer $b$ with $0\le b<b^\lambda_k$, there exists a multipartition
  $\nu=\nu(b)$ such that $e_\nu y_\nu = e(\bi_\lambda \vee k)y_\lambda y_{n+1}^b$.
\end{Lemma}

\proof The definitions of $\lambda$ and $b_k^\lambda$ ensure that there are $b_k^\lambda$ addable nodes of residue $k$ below $(\t^\lambda)^{-1}(n)$. Suppose those nodes are $(r_1,c_1,l_1)$, $(r_2,c_2,l_2), \ldots, (r_{b_k^\lambda},c_{b_k^\lambda},l_{b_k^\lambda})$, where $l_1 \geq l_2 \geq l_3 \geq \ldots \geq l_{b_k^\lambda}$, and if $l_i = l_{i+1}$ then $r_i \geq r_{i+1}$. In another word, $(r_i, c_i,l_i)$ is a node below $(r_{i+1},c_{i+1},l_{i+1})$.

For any $b$ with $0 \leq b < b_k^\lambda$, we define $\nu$ to be the multipartition obtained by adding the node $(r_{b+1},c_{b+1},l_{b+1})$ on to $\lambda$. Then $y_\nu = y_\lambda y_{n+1}^b$ and $e_\nu = e(\bi_\lambda\vee k) = e(\bi\vee k)$. This completes the proof. \endproof

\begin{Example} \label{I-problem: example: add below}
Suppose that $\lambda = (4,3|2,1|0|0)$ with $e = 4$ and $\kappa_\Lambda = (0,0,2,1)$. Then $e(\bi_\lambda) = e(0123301013)$ and $y_\lambda = y_1 y_2 y_3 y_4 y_6 y_7 y_9$. Then $b_0^\lambda = 1$, $b_1^\lambda = 1$, $b_2^\lambda = 2$ and $b_3^\lambda = 0$ and the proof of
\autoref{I-problem: downstair cases} shows that:
\begin{eqnarray*}
e(01233010130)y_1 y_2 y_3 y_4 y_6 y_7 y_9 & = & e_{\mu_1} y_{\mu_1},\\
e(01233010131)y_1 y_2 y_3 y_4 y_6 y_7 y_9 & = & e_{\mu_2} y_{\mu_2},\\
e(01233010132)y_1 y_2 y_3 y_4 y_6 y_7 y_9 & = & e_{\mu_3} y_{\mu_3},\\
e(01233010132)y_1 y_2 y_3 y_4 y_6 y_7 y_9 y_{11} & = & e_{\mu_4} y_{\mu_4},
\end{eqnarray*}
where $\mu_1 = (4,3|2,2|\emptyset|\emptyset)$, $\mu_2 = (4,3|2,1|\emptyset|1)$, $\mu_3 = (4,3|2,1|1|\emptyset)$ and $\mu_4 = (4,3|2,1,1|\emptyset|\emptyset)$.
\end{Example}

\begin{Definition} \label{Definition: P}
Let $\P = \cup_{n \geq 0} \P_n$. Define three sets $\P_I$, $\P_y$ and $\P_\psi$ of multipartitions by:
\begin{align*}
\P_I&= \{ \lambda\in\P \mid |\lambda| = n \text{ and } e(\bi_{\lambda_-}\vee k)y_{\lambda_-}y_n^{b_k^{\lambda_-}} \in \Blam \text{ for all }k\in I\},\\
\P_y&=\{\lambda\in\P \mid |\lambda| = n  \text{ and } \psi_{\s\t}y_r\in\BZ
		 \text{ whenever $\s,\t\in\Std(\lambda)$ and }1\le r\le n\},\\
\P_\psi &=\{\lambda\in\P \mid |\lambda| = n \text{ and }\psi_{\s\t}\psi_r\in\BZ
                \text{ whenever $\s,\t\in\Std(\lambda)$ and }1\le r<n\}.
\end{align*}
\end{Definition}

\begin{Remark} \label{two side by involution}
Notice that if for some $\s,\t\in \Std(\lambda)$ and $1 \leq r \leq n$ we have $\psi_{\s\t}y_r \in \BZ$, then $y_r \psi_{\s\t} \in \BZ$ as well. Similar property holds for $\psi_{\s\t}\psi_r$. Therefore we can write
\begin{align*}
\P_y&=\{\lambda\in\P \mid |\lambda| = n  \text{ and } y_r\psi_{\s\t}\in\BZ
		 \text{ whenever $\s,\t\in\Std(\lambda)$ and }1\le r\le n\},\\
\P_\psi &=\{\lambda\in\P \mid |\lambda| = n \text{ and }\psi_r\psi_{\s\t}\in\BZ
                \text{ whenever $\s,\t\in\Std(\lambda)$ and }1\le r<n\}
\end{align*}
as well.
\end{Remark}

By \autoref{basis to hecke} if one of $e(\bi_\nu)y_{\lambda_-}y_n^n$, $\psi_{\s\t}y_r$ or $\psi_{\s\t}\psi_r$ belongs to $\BZ$ then it can be written in a unique way as an (integral) linear combination of the $\psi$-basis elements. In particular, these linear combinations must satisfy the restrictions imposed by \autoref{residue sequence has to be the same}.

We define a total ordering on $\P$ extended by lexicographic ordering. Suppose $\lambda$ and $\mu$ are two multipartitions, not necessarily of the same integer.  Define $\mu \prec \lambda$ if $|\mu| < |\lambda|$, or $|\mu| = |\lambda|$ and $l(\mu) < l(\lambda)$, or $|\mu| = |\lambda|$, $l(\mu) = l(\lambda)$ and $\lambda < \mu$.

\begin{Definition} \label{Definition: S}
Define
$
\mathscr S_n^\Lambda = \{\lambda \in \P_n\mid
\mu \in \mathscr P^{\Lambda'}_I\cap \mathscr P^{\Lambda'}_y \cap \mathscr P^{\Lambda'}_\psi\text{ whenever $\mu\in\mathscr P_m^{\Lambda'}$ and }\mu\prec\lambda\}
$

\end{Definition}

We note that $\P_n \subseteq \P_I \cap \P_y \cap \P_\psi$ implies $R_n^\Lambda = \R$. Recall that $R_n^\Lambda$ is the $\Z$-vector space spanned by $\{\psi_{\s\t}\}$. Hence $\R$ is a $\Z$-span of $\{\psi_{\s\t}\}$. The main goal of this and the next sections is to prove that $\P_n \subseteq \P_I \cap \P_y \cap \P_\psi$ and we will use induction to prove the argument.

Now we can state the main result of this section.

\begin{Theorem} \label{I-problem: final}

Suppose $\lambda \in \mathscr S_n^\Lambda$. Then we have $\lambda\in \P_I$.

\end{Theorem}

As we mentioned before we are going to apply induction on $\lambda$ to prove the main Theorem. \autoref{onerow: I}, \autoref{onerow: y} and \autoref{onerow: psi} give the base case of the induction. Recall that $e \neq 2$.

\begin{Lemma}\label{onerow: I}
Suppose that $n\geq 1$ and $\lambda = (n|0|\ldots|0)\in \mathscr P_n^\Lambda$. Then
$
e(\bi_{\lambda_-}\vee k)y_{\lambda_-} y_n^{b_k^{\lambda_-}} \in \Blam
$
for any $k\in I$.
\end{Lemma}
\proof As $\lambda$ is the maximal element of $\P_n$, $\Blam  = \{0\}$. Therefore the Lemma is equivalent to the claim that $e(\bi_{\lambda_-}\vee k)y_{\lambda_-} y_{n}^{b_k^{\lambda_-}} = 0$. We prove this by induction on~$n$.

If $n =1$ then it is easy to see that $b_k^{\lambda_-} = (\Lambda,\alpha_{i_1})$. Therefore, $e(\bi_{\lambda_-} \vee k) y_{\lambda_-} y_n^{b_k^{\lambda_-}} = e(k) y_1^{(\Lambda,\alpha_{i_1})} = 0$ by~\eqref{alg:additional}.

Suppose now that the Lemma holds for any $n' <n$. Notice that for any $m \geq 1$, set $\gamma = (m|0|\ldots|0)$ then $|\mathscr{A}_{\t^\gamma}^k|$ is independent to the value of $m$. For the rest of the proof we set $a_k = |\mathscr A_{\t^{\lambda_-}}^k|$.

In order to simplify the notations, for the rest of the proof we will omit $i_1i_2\ldots i_{n-3}$ and simply write $e(\bi) = e(i_{n-2},i_{n-1},i_n)$. We will also suppress $y_\nu$, where $\nu = \lambda|_{n-3}$.

We consider four cases separately, depending on the value of $k$.

\medskip\textbf{Case ~\ref{onerow: I}a: $k=i_{n-1}$.} Then
$e(\textbf{i}_{\lambda_-}\vee k) y_{\lambda_-} y_{n}^{b_k^{\lambda_-}} = e(\textbf{i}|_{n-3},k-1,k,k)y_{\lambda|_{n-3}} y_{n-2}^{a_{k-1}} y_{n-1}^{a_k}y_{n}^{a_k}$. In this case we have
\begin{align*}
e(k-1,k,k)y_{n-2}^{a_{k-1}}y_{n-1}^{a_k}y_n^{a_k} &\overset{(\ref{alg:psi-y com})}= -e(k-1,k,k)y_{n-2}^{a_{k-1}}y_{n-1}^{a_k+1}y_n^{a_k}\psi_{n-1} + e(k-1,k,k)y_{n-2}^{a_{k-1}}y_{n-1}^{a_k}y_n^{a_k}\psi_{n-1} y_n\\
&\overset{(\ref{alg:ypsi})}= \psi_{n-1}e(k-1,k,k)y_{n-2}^{a_{k-1}}y_{n-1}^{a_k}y_n^{a_k}y_n,
\end{align*}
where $e(k-1,k,k)y_{n-2}^{a_{k-1}}y_{n-1}^{a_k+1}y_n^{a_k} = 0$ by induction. Therefore,
\begin{align*}
e(k-1,k,k)y_{n-2}^{a_{k-1}}y_{n-1}^{a_k}y_n^{a_k} & = \psi_{n-1}e(k-1,k,k)y_{n-2}^{a_{k-1}}y_{n-1}^{a_k}y_n^{a_k}y_n\\
& = \psi_{n-1}^2e(k-1,k,k)y_{n-2}^{a_{k-1}}y_{n-1}^{a_k}y_n^{a_k}y_n^2 = 0
\end{align*}
by relation (\ref{alg:psipsi}).

\medskip\medskip\textbf{Case ~\ref{onerow: I}b: $k=i_{n-1}+1$}. Now,
$
e(\textbf{i}_{\lambda_-}\vee k) y_{\lambda_-} y_{n}^{b_k^{\lambda_-}} = e(\textbf{i}|_{n-3},k-2,k-1,k)y_{\lambda|_{n-3}} y_{n-2}^{a_{k-2}} y_{n-1}^{a_{k-1}} y_{n}^{a_k+1}
$.
Therefore,
\begin{eqnarray*}
&&e(k-2,k-1,k)y_{n-2}^{a_{k-2}}y_{n-1}^{a_{k-1}}y_n^{a_k+1}\\
 & \overset{(\ref{alg:psipsi})} = & e(k-2,k-1,k)y_{n-2}^{a_{k-2}}y_{n-1}^{a_{k-1}+1}y_n^{a_k} + e(k-2,k-1,k)y_{n-2}^{a_{k-2}}y_{n-1}^{a_{k-1}}y_n^{a_k}\psi^2_{n-1}\\
 &\overset{\substack{(\ref{alg:psi-y com})\\(\ref{alg:y-psi com})}}= & e(k-2,k-1,k)y_{n-2}^{a_{k-2}}y_{n-1}^{a_{k-1}+1}y_n^{a_k} + \psi_{n-1} e(k-2,k,k-1)y_{n-2}^{a_{k-2}}y_{n-1}^{a_k}y_n^{a_{k-1}}\psi_{n-1} = 0,
\end{eqnarray*}
where the last equality follows by induction.

\medskip\medskip\textbf{Case ~\ref{onerow: I}c: $k=i_{n-1}-1$.} If $n = 2$ then
$e(\textbf{i}_{\lambda_-}\vee k) y_{\lambda_-} y_{n}^{b_k^{\lambda_-}}= e(k,k-1)y_1^{a_k}y_2^{a_{k-1}}$. Then $a_{k-1} \geq 1$. Therefore,
$$
e(k,k-1)y_1^{a_k}y_2^{a_{k-1}} \overset{(\ref{alg:psipsi})}= e(k,k-1)y_1^{a_k+1}y_2^{a_{k-1}-1} - \psi_1e(k-1,k)y_1^{a_{k-1}-1}y_2^{a_k}\psi_1 = 0,
$$
using relation (\ref{alg:additional}) and induction. Hence, the lemma follows in this case when $n=2$.

If $n > 2$ then $e(\textbf{i}_{\lambda_-}\vee k) y_{\lambda_-} y_{n}^{b_k^{\lambda_-}} = e(\textbf{i}|_{n-3},k,k+1,k)y_{\lambda|_{n-3}}y_{n-2}^{a_k} y_{n-1}^{a_{k+1}} y_{n}^{a_k}$. Hence,
\begin{align*}
e(k,k+1,k)y_{n-2}^{a_k}
y_{n-1}^{a_{k+1}} y_{n}^{a_k}
&\overset{(\ref{alg:braid})}= \psi_{n-2}\psi_{n-1}\psi_{n-2}e(k,k+1,k)y_{n-2}^{a_k} y_{n-1}^{a_{k+1}} y_{n}^{a_k} \\
&\qquad- \psi_{n-1}\psi_{n-2}\psi_{n-1}e(k,k+1,k)y_{n-2}^{a_k} y_{n-1}^{a_{k+1}} y_{n}^{a_k}\\
&=  \psi_{n-2}\psi_{n-1} e(k+1,k,k)y_{n-2}^{a_{k+1}}y_{n-1}^{a_k} y_{n}^{a_k}\psi_{n-2}\\
&\qquad- \psi_{n-1}\psi_{n-2}e(k,k,k+1)y_{n-2}^{a_k} y_{n-1}^{a_k} y_{n}^{a_{k+1}}\psi_{n-1}\\
&= 0,
\end{align*}
where the last equality follows by induction.

\medskip\medskip\textbf{Case ~\ref{onerow: I}d: $|k-i_{n-1}|>1$.} Because $i_{n-2} = i_{n-1} - 1$, we have $i_{n-2} \neq k$. Therefore we have
\begin{align*}
e(\textbf{i}_{\lambda_-}\vee k) y_{\lambda_-} y_{n}^{b_k^{\lambda_-}}
  &= e(\textbf{i}|_{n-3},i_{n-2},i_{n-1},k)y_{\lambda|_{n-3}} y_{n-2}^{a_{i_{n-2}}} y_{n-1}^{a_{i_{n-1}}} y_{n}^{a_k}\\
  &\overset{(\ref{alg:psipsi})}= \psi_{n-1} e(\textbf{i}|_{n-3},i_{n-2},k,i_{n-1})y_{{\lambda}|_{n-3}} y_{n-2}^{a_{i_{n-2}}} y_{n-1}^{a_k} y_{n}^{a_{i_{n-1}}} \psi_{n-1} = 0
\end{align*}
by induction. This completes the proof. \endproof

\autoref{onerow: I} has two immediate Corollaries:

\begin{Corollary}\label{onerow: y}
Suppose $n\geq 2$ and $\lambda = (n|0|\ldots|0)$. Then
$
e_\lambda y_\lambda y_r \in \Blam
$
for any $1\leq r\leq n$.
\end{Corollary}

\begin{Corollary}\label{onerow: psi}
Suppose that $n\geq 2$ and $\lambda = (n|0|\ldots|0)$. Then
$
e_\lambda y_\lambda \psi_r \in \Blam
$,
for any $1\leq r\leq n-1$.
\end{Corollary}
\proof Write $y_\lambda = y_1^{l_1}y_2^{l_2}\ldots y_n^{l_n}$ and $\bi_\lambda = (i_1i_2\ldots i_n)$,
$$
e(i_1\ldots i_n)y_1^{l_1}\ldots y_n^{l_n}\psi_r = \psi_r e(i_1\ldots i_{r-1}i_{r+1}i_r\ldots i_n)y_1^{l_1}\ldots y_{r-1}^{l_{r-1}}y_r^{l_{r+1}}y_{r+1}^{l_r}\ldots y_n^{l_n} = 0,
$$
by \autoref{onerow: I}.\endproof

The results in the rest of the subsection will be used frequently in the later proofs.

Recall that for any multipartition $\lambda$, $\Blam $ is the subspace of $\R$ spanned by all of the elements $\psi_{\s\t}$, where $\Shape(\s) = \Shape(\t) > \lambda$.

\begin{Lemma} \label{B_lambda is an ideal}
  Suppose $\lambda \in \mathscr{S}_n^\Lambda$. Then $\Blam $ is a two-sided ideal of $\mathscr{R}_n^\Lambda$. More precisely, $\Blam[\mu]$ is a two-sided ideal of~$\R$ whenever $\mu\prec\lambda$.
\end{Lemma}

\proof The Lemma follows directly from the definition of the set $\mathscr S_n^\Lambda$,
$\P_y$, $\P_\psi$ and \autoref{two side by involution}.
\endproof

In order to simplify the notation, for each $i \in I$ define
$\theta_i\map{\mathscr{R}_n^\Lambda}{\mathscr{R}_{n+1}^\Lambda}$ to be the unique $\Z$-linear map which sends  \label{notation: theta} $e(\bi)$ to $e(\bi\vee i)$, $y_r$ to $y_r$ and $\psi_r$ to $\psi_r$. It is easy to see that $\theta_i$ respects the relations in~$\R$, so $\theta_i$ is a
$\Z$-algebra homomorphism.

\begin{Lemma}\label{send:help}
Suppose $\lambda\in \mathscr S_n^\Lambda$ and $\u,\v\in\Std(\mu)$, where $\mu \in \P_m$ with $m<n$ such that $\mu > \lambda|_m$. Let $\sigma=\lambda|_{m+1} \in \P_{m+1}$. Then
$\theta_i(\psi_{\u\v}) \in \Blam[\sigma]$, for any $i\in I$.
\end{Lemma}

\proof Write $\mu = (\mu^{(1)},\ldots,\mu^{(\l)})$ and $\mu^{(\l)} = (\mu^{(\l)}_1,\ldots,\mu^{(\l)}_k)$ and define $\gamma = (\mu^{(1)},\ldots,\mu^{(\l-1)},\gamma^{(\l)})$ where

$$
\gamma^{(\l)} =
\begin{cases}
(\mu^{(\l)}_1,\ldots,\mu^{(\l)}_{k-1},\mu^{(\l)}_k+1),& \text{if $\mu^{(\l)}_{k-1} > \mu^{(\l)}_k$,}\\
(\mu^{(\l)}_1,\ldots,\mu^{(\l)}_{k-1},\mu^{(\l)}_k,1),& \text{if $\mu^{(\l)}_{k-1} = \mu^{(\l)}_k$.}
\end{cases}
$$

Then $\gamma$ is a multipartition of $m+1$ and $\gamma|_m = \mu$. Since $m < n$,
if $m = n-1$, then $\gamma|_{n-1} = \mu > \lambda_-$, so that $\gamma >
\lambda$. On the other hand, if $m < n-1$ then $|\gamma| = m+1 < n = |\lambda|$.
So we always have $\gamma \prec \lambda$.  Therefore, $\gamma \in \P_I \cap \P_y \cap \P_\psi$
because $\lambda\in \mathscr S_n^\Lambda$.

As $\gamma|_m = \mu$, we have
$
\theta_i(\psi_{\u\v}) = \theta_i(\psi_{d(\u)}^* e_\mu y_\mu \psi_{d(\v)})
   = \psi_{d(\u)}^* e(\bi_{\gamma|_m}\vee i) y_\mu \psi_{d(\v)}.
$
First suppose that $b_i^\mu = 0$. Then using the definition of $\P_I$, we have
$
e(\bi_{\gamma|_m} \vee i) y_\mu \in \Blam[\gamma] \subseteq \Blam[\sigma ]
$.
Hence, by \autoref{B_lambda is an ideal}, we have
$
\theta_i(\psi_{\u\v}) = \psi_{d(\u)}^* e(\bi_{\gamma|_m}\vee i) y_\mu \psi_{d(\v)} \in \Blam[\sigma ]
$.

Now supose that $b_i^\mu > 0$. By \autoref{I-problem: downstair cases} there exists a multipartition
$\nu$ with $\nu|_m = \mu$ such that
$
e(\bi_{\gamma|_m} \vee i) y_\mu = e_\nu y_\nu
$.
Further, as $\nu|_m = \mu$, there exist two standard $\nu$-tableaux $\s$ and $\t$ such that $\s|_m = \u$ and $\t|_m = \v$. That is, $d(\s) = d(\u)$ and $d(\t) = d(\v)$. Therefore,
$$
\theta_i(\psi_{\u\v}) = \psi_{d(\u)}^* e(\bi_{\gamma|_m}\vee i) y_\mu \psi_{d(\v)} = \psi_{d(\s)}^* e_\nu y_\nu \psi_{d(\t)} = \psi_{\s\t} \in \Bgelam[\nu] \subseteq \Blam[\sigma]
$$
because $\nu  > \sigma $. This completes the proof.
\endproof

If $\bi=(i_1,\dots,i_n)\in I^n$ and $1\le m\le n$ let $\bi_m=(i_1\dots i_m)$.

\begin{Lemma}\label{send}
Suppose $\lambda\in\mathscr S_n^\Lambda$, $m \leq n$ and $\sigma = \lambda|_m$. For any $\bi = (i_1,i_2,\ldots,i_{n-m})$ we have $\mathscr R_n^\Lambda \theta_\bi(\Blam[\sigma])\mathscr R_n^\Lambda \subseteq \Blam$.
\end{Lemma}

\proof Suppose $r\in \Blam[\sigma]$, we have that
$$
r = \sum_{\u,\v \in \Std( > \sigma)}c_{\u\v}\psi_{\u\v}
$$
for some $c_{\u\v} \in \Z$. For any $i\in I$,
$$
\theta_i(r) = \sum_{\u,\v \in \Std( > \sigma)}c_{\u\v} \theta_i(\psi_{\u\v}).
$$

By \autoref{send:help}, $\theta_i(\psi_{\u\v})\in \Blam[\lambda|_{m+1}]$. Hence $\theta_i(r) \in \Blam[\lambda|_{m+1}]$. By induction we have $\theta_\bi(r) \in \Blam$.
By \autoref{B_lambda is an ideal}, $\Blam$ is an ideal. Therefore
$
\mathscr R_n^\Lambda \theta_\bi(\Blam[\sigma])\mathscr R_n^\Lambda \subseteq \Blam
$
which completes the proof.  \endproof

\subsection{The action of $y_r$ on two-rowed partitions}

Recall that the main result of this section is to prove that if $\lambda\in \mathscr S_n$, then
$$
e(\bi_{\lambda_-}\vee k)y_{\lambda_-} y_n^{b_k^{\lambda_-}} \in R_n^{>\lambda}.
$$

In the inductive process we consider different types of multipartitions $\lambda$ and a residue $k\in I$. We will consider the more difficult case first, namely when $\lambda = (\lambda^{(1)},\ldots,\lambda^{(\l)})$ and $\lambda^{(\l)} = (\lambda^{(\l)}_1,\ldots,\lambda^{(\l)}_l,1) \neq \emptyset$ with $l \geq 2$, $\lambda^{(\l)}_{l-1} = \lambda^{(\l)}_l = m$ and $k \equiv \kappa_\l - l + m + 1\pmod{e}$. In this subsection we assume that $\l = 1$ and $l = 2$. We will extend the result to the general case in the next subsection. Notice that in this case $\lambda = (m,m,1)$ for some integer $m$ and $k \equiv \kappa_1 - 1 + m\pmod{e}$. Then $e(\bi_{\lambda_-}\vee k)y_{\lambda_-} y_n^{b_k^{\lambda_-}} = e_\gamma y_\gamma$ where $\gamma = (m,m+1)$. It is very hard to prove that $e_\gamma y_\gamma \in R_n^{>\lambda}$ directly, so we are going to work with $\gamma$ which is in a more general form.

In this subsection we fix $\Lambda = \Lambda_j$ for some $j\in I$, $\gamma = (\gamma_1,\gamma_2)$ and $\lambda = (\gamma_1,\gamma_2-1,1)$ with $\gamma_2 > 1$ and $\gamma_2 - \gamma_1 \equiv 1\pmod{e}$. We will prove that if $\gamma_1 + 1 = \gamma_2$ and $\lambda\in \mathscr S_n^\Lambda$ then $e_\gamma y_\gamma \in \Blam[\gamma]$.

Without loss of generality we can assume that $\Lambda = \Lambda_0$. Define $i \equiv \gamma_2-2 \pmod{e}$, which is the residue of $(2,\gamma_2,1)$. Because $\gamma_2 \equiv \gamma_1 + 1\pmod{e}$, it is also the residue of the node $(1,\gamma_1,1)$. In diagrammatic notation, we have
$$
e_\gamma y_\gamma =
\begin{braid}
 \draw(2,-0.5)node{$\underbrace{\hspace*{16mm}}_{\gamma_1}$};
 \draw(7,-0.5)node{$\underbrace{\hspace*{16mm}}_{\gamma_2}$};
 \draw (0,4)node[above]{$0$}--(0,0);
 \draw (1,4)node[above]{$1$}--(1,0);
 \draw[dots] (1.2,4)--(2.8,4);
 \draw[dots] (1.2,0)--(2.8,0);
 \draw (3,4)node[above]{$i-1$}--(3,0);
 \draw (4,4)node[above]{$i$}--(4,0);
 \draw (5,4)node[above]{$e-1$}--(5,0);
 \draw (6,4)node[above]{$0$}--(6,0);
 \draw[dots] (6.2,4)--(7.8,4);
 \draw[dots] (6.2,0)--(7.8,0);
 \draw (8,4)node[above]{$i-1$}--(8,0);
 \draw (9,4)node[above]{$i$}--(9,0);
 \node[greendot] at (0,2){};
 \node[above right,scale = 0.7] at (0,2){$l_1$};
 \node[greendot] at (1,2){};
 \node[above right,scale = 0.7] at (1,2){$l_2$};
 \node[greendot] at (3,2){};
 \node[above right,scale = 0.7] at (3,2){$l_{m-1}$};
 \node[greendot] at (4,2){};
 \node[above right,scale = 0.7] at (4,2){$l_m$};
 \node[greendot] at (5,2){};
 \node[above right,scale = 0.7] at (5,2){$l_{m+1}$};
 \node[greendot] at (6,2){};
 \node[above right,scale = 0.7] at (6,2){$l_{m+2}$};
 \node[greendot] at (8,2){};
 \node[above right,scale = 0.7] at (8,2){$l_{2m}$};
 \node[greendot] at (9,2){};
 \node[above right,scale = 0.7] at (9,2){$l_{2m+1}$};
\end{braid}
$$
where $\bi_\gamma = (i_1, i_2, \ldots, i_n)$ and $l_k = |\mathscr{A}_{\t^\gamma|_k}^{i_k}|$ is the multiplicity of the green dot on the $k$-th string. For the rest of this subsection, for clarity we will omit extraneous dots when they do not play an important role in the argument.

Next we introduce an important equivalent relation $=_\gamma$. For $\gamma\in\mathscr S_n^\Lambda$, and $r_1,r_2\in \mathscr{R}_n^\Lambda$, we write $r_1 =_\gamma r_2$ if $r_1 \pm r_2 \in \Blam[\gamma]$. It is clearly an equivalent relation. Moreover, by \autoref{B_lambda is an ideal}, for any $r \in \R$ we have $r_1{\cdot}r =_\gamma r_2{\cdot}r$ if $r_1 =_\gamma r_2$. This will be helpful for us to simplify the notations and calculations.

Recall that $\gamma_2 > 1$. We can write $\gamma_2 = k{\cdot}e + t$ for some nonnegative integer $k$ and $2 \leq t \leq e + 1$. We will first prove

\begin{equation} \label{I-problem: non-addable: help}
e_\gamma y_\gamma
=_\gamma
\begin{cases}
\begin{braid}
 \draw(2,-0.5)node{$\underbrace{\hspace*{16mm}}_{\gamma_1}$};
 \draw(7,-0.5)node{$\underbrace{\hspace*{16mm}}_{\gamma_2}$};
 \draw (0,4)node[above]{$0$}--(0,0);
 \draw (1,4)node[above]{$1$}--(1,0);
 \draw[dots] (1.2,4)--(2.8,4);
 \draw[dots] (1.2,0)--(2.8,0);
 \draw (3,4)node[above]{$i-1$}--(3,0);
 \draw (4,4)node[above]{$i$}--(8,2)--(9,0);
 \draw (5,4)node[above]{$e-1$}--(4,2)--(5,0);
 \draw (6,4)node[above]{$0$}--(5,2)--(6,0);
 \draw[dots] (6.2,4)--(7.8,4);
 \draw[dots] (6.2,0)--(7.8,0);
 \draw (8,4)node[above]{$i-1$}--(7,2)--(8,0);
 \draw (9,4)node[above]{$i$}--(8,2)--(4,0);
\end{braid}, & \text{if $i \neq e-1$},\\
\begin{braid}
 \draw(2,-0.5)node{$\underbrace{\hspace*{16mm}}_{\gamma_1}$};
 \draw(7,-0.5)node{$\underbrace{\hspace*{16mm}}_{\gamma_2}$};
 \draw (0,4)node[above]{$0$}--(0,0);
 \draw (1,4)node[above]{$1$}--(1,0);
 \draw[dots] (1.2,4)--(2.8,4);
 \draw[dots] (1.2,0)--(2.8,0);
 \draw (3,4)node[above]{$e-2$}--(3,0);
 \draw (4,4)node[above]{$e-1$}--(4,0);
 \draw (5,4)node[above]{$e-1$}--(8,2)--(9,0);
 \draw (6,4)node[above]{$0$}--(5,2)--(6,0);
 \draw[dots] (6.2,4)--(7.8,4);
 \draw[dots] (6.2,0)--(7.8,0);
 \draw (8,4)node[above]{$e-2$}--(7,2)--(8,0);
 \draw (9,4)node[above]{$e-1$}--(8,2)--(5,0);
\end{braid}, & \text{if $i = e-1$}.
\end{cases}
\end{equation}
by induction on $k$, which can imply $e_\gamma y_\gamma \in R_n^{>\lambda}$ easily.

In order to clarify the meaning of the diagrams in (\ref{I-problem: non-addable: help}), let us give two examples below. In these
examples for convenience we fix $e = 4$.

\begin{Example}
Suppose $\gamma = (8,5)$, then $\gamma = \tab(\ \ \ \ \ \ \ \ ,\ \ \ \ \ )$ and $i = 3$. Then we are trying to prove that
$$
e_\gamma y_\gamma  =
\begin{braid}
 \draw (0,4)node[above]{$0$}--(0,0);
 \draw (1,4)node[above]{$1$}--(1,0);
 \draw (2,4)node[above]{$2$}--(2,0);
 \draw (3,4)node[above]{$3$}--(3,0);
 \draw (4,4)node[above]{$0$}--(4,0);
 \draw (5,4)node[above]{$1$}--(5,0);
 \draw (6,4)node[above]{$2$}--(6,0);
 \draw (7,4)node[above]{$3$}--(7,0);
 \draw (8,4)node[above]{$3$}--(8,0);
 \draw (9,4)node[above]{$0$}--(9,0);
 \draw (10,4)node[above]{$1$}--(10,0);
 \draw (11,4)node[above]{$2$}--(11,0);
 \draw (12,4)node[above]{$3$}--(12,0);
 \node[greendot] at (3,2){};
 \node[greendot] at (7,2){};
 \node[greendot] at (11,2){};
\end{braid}
 =_\gamma
\begin{braid}
 \draw (0,4)node[above]{$0$}--(0,0);
 \draw (1,4)node[above]{$1$}--(1,0);
 \draw (2,4)node[above]{$2$}--(2,0);
 \draw (3,4)node[above]{$3$}--(3,0);
 \draw (4,4)node[above]{$0$}--(4,0);
 \draw (5,4)node[above]{$1$}--(5,0);
 \draw (6,4)node[above]{$2$}--(6,0);
 \draw (7,4)node[above]{$3$}--(7,0);
 \draw (8,4)node[above]{$3$}--(11,2) -- (12,0);
 \draw (9,4)node[above]{$0$}--(8,2)--(9,0);
 \draw (10,4)node[above]{$1$}--(9,2)--(10,0);
 \draw (11,4)node[above]{$2$}--(10,2)--(11,0);
 \draw (12,4)node[above]{$3$}--(11,2)--(8,0);
 \node[greendot] at (3,2){};
 \node[greendot] at (10,2){};
\end{braid}.
$$
\end{Example}

\begin{Example}
Suppose $\gamma = (9,10)$, then $\gamma = \tab(\ \ \ \ \ \ \ \ \ ,\ \ \ \ \ \ \ \ \ \ )$ and $i = 0$. We are trying to prove that
\begin{eqnarray*}
e_\gamma y_\gamma & = &
\begin{braid}
 \draw (0,4)node[above]{$0$}--(0,0);
 \draw (1,4)node[above]{$1$}--(1,0);
 \draw (2,4)node[above]{$2$}--(2,0);
 \draw (3,4)node[above]{$3$}--(3,0);
 \draw (4,4)node[above]{$0$}--(4,0);
 \draw (5,4)node[above]{$1$}--(5,0);
 \draw (6,4)node[above]{$2$}--(6,0);
 \draw (7,4)node[above]{$3$}--(7,0);
 \draw (8,4)node[above]{$0$}--(8,0);
 \draw (9,4)node[above]{$3$}--(9,0);
 \draw (10,4)node[above]{$0$}--(10,0);
 \draw (11,4)node[above]{$1$}--(11,0);
 \draw (12,4)node[above]{$2$}--(12,0);
 \draw (13,4)node[above]{$3$}--(13,0);
 \draw (14,4)node[above]{$0$}--(14,0);
 \draw (15,4)node[above]{$1$}--(15,0);
 \draw (16,4)node[above]{$2$}--(16,0);
 \draw (17,4)node[above]{$3$}--(17,0);
 \draw (18,4)node[above]{$0$}--(18,0);
 \node[greendot] at (3,2){};
 \node[greendot] at (7,2){};
 \node[greendot] at (12,2){};
 \node[greendot] at (16,2){};
\end{braid}\\
& =_\gamma &
\begin{braid}
 \draw (0,4)node[above]{$0$}--(0,0);
 \draw (1,4)node[above]{$1$}--(1,0);
 \draw (2,4)node[above]{$2$}--(2,0);
 \draw (3,4)node[above]{$3$}--(3,0);
 \draw (4,4)node[above]{$0$}--(4,0);
 \draw (5,4)node[above]{$1$}--(5,0);
 \draw (6,4)node[above]{$2$}--(6,0);
 \draw (7,4)node[above]{$3$}--(7,0);
 \draw (8,4)node[above]{$0$}--(17,2) -- (18,0);
 \draw (9,4)node[above]{$3$}--(8,2)--(9,0);
 \draw (10,4)node[above]{$0$}--(9,2)--(10,0);
 \draw (11,4)node[above]{$1$}--(10,2)--(11,0);
 \draw (12,4)node[above]{$2$}--(11,2)--(12,0);
 \draw (13,4)node[above]{$3$}--(12,2)--(13,0);
 \draw (14,4)node[above]{$0$}--(13,2)--(14,0);
 \draw (15,4)node[above]{$1$}--(14,2)--(15,0);
 \draw (16,4)node[above]{$2$}--(15,2)--(16,0);
 \draw (17,4)node[above]{$3$}--(16,2)--(17,0);
 \draw (18,4)node[above]{$0$}--(17,2)--(8,0);
 \node[greendot] at (3,2){};
 \node[greendot] at (7,2){};
 \node[greendot] at (11,2){};
 \node[greendot] at (15,2){};
\end{braid}.
\end{eqnarray*}
\end{Example}

The next Proposition is the base case of the induction. When $k = 0$, we have $2\leq \gamma_2 \leq e + 1$.\\

\begin{Proposition}~\label{I-problem: induction1}
Suppose $\gamma = (\gamma_1,\gamma_2) \in \mathscr C_n^\Lambda$ with $\gamma_2 > 1$ and $\gamma_2 - \gamma_1 \equiv 1\pmod{e}$ and $\lambda = (\gamma_1,\gamma_2-1,1) \in \mathscr S_n^\Lambda$. Define $i$ to be the residue of the node at position $(1,\gamma_1,1)$ or $(2,\gamma_2,1)$. When $2 \leq \gamma_2 \leq e + 1$, (\ref{I-problem: non-addable: help}) holds.
\end{Proposition}

Before proving \autoref{I-problem: induction1} we first give a useful lemma.

\begin{Lemma}~\label{I-problem: short}
For any $i \in I$, we have
$$
\begin{braid}
 \draw (0,4)node[above]{$i$}--(4,2)--(5,0);
 \draw (1,4)node[above]{$i+1$}--(0,2)--(1,0);
 \draw (2,4)node[above]{$i+2$}--(1,2)--(2,0);
 \draw[dots] (2.2,4)--(3.8,4);
 \draw[dots] (2.2,0)--(3.8,0);
 \draw (4,4)node[above]{$i-1$}--(3,2)--(4,0);
 \draw (5,4)node[above]{$i$}--(4,2)--(0,0);
\end{braid}
=
\begin{braid}
 \draw (0,4)node[above]{$i$}--(0,0);
 \draw (1,4)node[above]{$i+1$}--(1,0);
 \draw (2,4)node[above]{$i+2$}--(2,0);
 \draw[dots] (2.2,4)--(3.8,4);
 \draw[dots] (2.2,0)--(3.8,0);
 \draw (4,4)node[above]{$i-1$}--(4,0);
 \draw (5,4)node[above]{$i$}--(5,0);
 \node[greendot] at (5,2){};
\end{braid}
-
\begin{braid}
 \draw (0,4)node[above]{$i$}--(0,0);
 \draw (1,4)node[above]{$i+1$}--(1,0);
 \draw (2,4)node[above]{$i+2$}--(2,0);
 \draw[dots] (2.2,4)--(3.8,4);
 \draw[dots] (2.2,0)--(3.8,0);
 \draw (4,4)node[above]{$i-1$}--(4,0);
 \draw (5,4)node[above]{$i$}--(5,0);
 \node[greendot] at (4,2){};
\end{braid}
-
\begin{braid}
 \draw (0,4)node[above]{$i$}--(0,0);
 \draw (1,4)node[above]{$i+1$}--(1,0);
 \draw (2,4)node[above]{$i+2$}--(2,0);
 \draw[dots] (2.2,4)--(3.8,4);
 \draw[dots] (2.2,0)--(3.8,0);
 \draw (4,4)node[above]{$i-1$}--(4,0);
 \draw (5,4)node[above]{$i$}--(5,0);
 \node[greendot] at (1,2){};
\end{braid}
+
\begin{braid}
 \draw (0,4)node[above]{$i$}--(0,0);
 \draw (1,4)node[above]{$i+1$}--(1,0);
 \draw (2,4)node[above]{$i+2$}--(2,0);
 \draw[dots] (2.2,4)--(3.8,4);
 \draw[dots] (2.2,0)--(3.8,0);
 \draw (4,4)node[above]{$i-1$}--(4,0);
 \draw (5,4)node[above]{$i$}--(5,0);
 \node[greendot] at (0,2){};
\end{braid}
+
\begin{braid}
 \draw (0,4)node[above]{$i$}--(1,2)--(5,0);
 \draw (1,4)node[above]{$i+1$}--(2,2)--(1,0);
 \draw (2,4)node[above]{$i+2$}--(3,2)--(2,0);
 \draw[dots] (2.2,4)--(3.8,4);
 \draw[dots] (2.2,0)--(3.8,0);
 \draw (4,4)node[above]{$i-1$}--(5,2)--(4,0);
 \draw (5,4)node[above]{$i$}--(1,2)--(0,0);
\end{braid}.
$$
\end{Lemma}

\proof The Lemma follows by directly applying braid relations on the left hand side of the equation.\endproof

Now we are ready to prove \autoref{I-problem: induction1}.

\proof We prove the Proposition by considering four different cases depending
upon the value of~$i$. Notice that in this Proposition, we have $\gamma_1 \geq \gamma_2-1$ because $2 \leq \gamma_2 \leq e+1$ and $\gamma_2 - \gamma_1 \equiv 1\pmod{e}$.

\medskip\textbf{Case ~\ref{I-problem: induction1}a:} $i = 0$, i.e. $\gamma_2 = 2$.
\begin{eqnarray*}
e_\gamma y_\gamma & = &
\begin{braid}
 \draw(2,-0.5)node{$\underbrace{\hspace*{16mm}}_{\gamma_1}$};
 \draw (0,4)node[above]{$0$}--(0,0);
 \draw (1,4)node[above]{$1$}--(1,0);
 \draw[dots] (1.2,4)--(2.8,4);
 \draw[dots] (1.2,0)--(2.8,0);
 \draw (3,4)node[above]{$e-1$}--(3,0);
 \draw (4,4)[densely dotted] node[above]{$0$}--(4,0);
 \draw (5,4)[densely dotted] node[above]{$e-1$}--(5,0);
 \draw (6,4)[densely dotted] node[above]{$0$}--(6,0);
\end{braid}
\overset{(\ref{dia:braid})}=
\begin{braid}
 \draw (0,4)node[above]{$0$}--(0,0);
 \draw (1,4)node[above]{$1$}--(1,0);
 \draw[dots] (1.2,4)--(2.8,4);
 \draw[dots] (1.2,0)--(2.8,0);
 \draw (3,4)node[above]{$e-1$}--(3,0);
 \draw (4,4)[densely dotted] node[above]{$0$}--(6,0);
 \draw (5,4)node[above]{$e-1$}--(6,2)--(5,0);
 \draw (6,4)[densely dotted] node[above]{$0$}--(4,0);
\end{braid}
-
\begin{braid}
 \draw (0,4)node[above]{$0$}--(0,0);
 \draw (1,4)node[above]{$1$}--(1,0);
 \draw[dots] (1.2,4)--(2.8,4);
 \draw[dots] (1.2,0)--(2.8,0);
 \draw (3,4)node[above]{$e-1$}--(3,0);
 \draw (4,4)node[above]{$0$}--(6,0);
 \draw (5,4)node[above]{$e-1$}--(4,2)--(5,0);
 \draw (6,4)node[above]{$0$}--(4,0);
 \node[greendot] at (3,2){};
\end{braid}\\
& \overset{(\ref{dia:ii3})} = &
-\begin{braid}
 \draw (0,4)node[above]{$0$}--(0,0);
 \draw (1,4)node[above]{$1$}--(1,0);
 \draw[dots] (1.2,4)--(2.8,4);
 \draw[dots] (1.2,0)--(2.8,0);
 \draw (3,4)node[above]{$e-1$}--(3,0);
 \draw (4,4)node[above]{$0$}--(5,2)--(4,0);
 \draw (5,4)node[above]{$e-1$}--(6,2)--(5,0);
 \draw (6,4)node[above]{$0$}--(4,2)--(6,0);
 \node[greendot] at (3,2){};
 \node[greendot] at (4,2){};
\end{braid}
-
\begin{braid}
 \draw (0,4)node[above]{$0$}--(0,0);
 \draw (1,4)node[above]{$1$}--(1,0);
 \draw[dots] (1.2,4)--(2.8,4);
 \draw[dots] (1.2,0)--(2.8,0);
 \draw (3,4)node[above]{$e-1$}--(3,0);
 \draw (4,4)node[above]{$0$}--(6,0);
 \draw (5,4)node[above]{$e-1$}--(4,2)--(5,0);
 \draw (6,4)node[above]{$0$}--(4,0);
 \node[greendot] at (3,2){};
\end{braid}
 = 
-\begin{braid}
 \draw (0,4)node[above]{$0$}--(0,0);
 \draw (1,4)node[above]{$1$}--(1,0);
 \draw[dots] (1.2,4)--(2.8,4);
 \draw[dots] (1.2,0)--(2.8,0);
 \draw (3,4)node[above]{$e-1$}--(3,0);
 \draw (4,4)node[above]{$0$}--(5,3)--(5,1)--(4,0);
 \draw (5,4)node[above]{$e-1$}--(6,3)--(6,1)--(5,0);
 \draw (6,4)node[above]{$0$}--(4,3)--(4,1)--(6,0);
 \node[greendot] at (3,2){};
 \node[greendot] at (4,2){};
 \draw[loosely dashed,red] (-0.5,3) -- (4.5,3)--(4.5,1)--(-0.5,1)--(-0.5,3);
\end{braid}
-
\begin{braid}
 \draw (0,4)node[above]{$0$}--(0,0);
 \draw (1,4)node[above]{$1$}--(1,0);
 \draw[dots] (1.2,4)--(2.8,4);
 \draw[dots] (1.2,0)--(2.8,0);
 \draw (3,4)node[above]{$e-1$}--(3,0);
 \draw (4,4)node[above]{$0$}--(6,0);
 \draw (5,4)node[above]{$e-1$}--(4,2)--(5,0);
 \draw (6,4)node[above]{$0$}--(4,0);
 \node[greendot] at (3,2){};
\end{braid}.
\end{eqnarray*}

Because
$$
\begin{braid}
 \draw (0,4)node[above]{$0$}--(0,0);
 \draw (1,4)node[above]{$1$}--(1,0);
 \draw[dots] (1.2,4)--(2.8,4);
 \draw[dots] (1.2,0)--(2.8,0);
 \draw (3,4)node[above]{$e-1$}--(3,0);
 \draw (4,4)node[above]{$0$}--(5,3)--(5,1)--(4,0);
 \draw (5,4)node[above]{$e-1$}--(6,3)--(6,1)--(5,0);
 \draw (6,4)node[above]{$0$}--(4,3)--(4,1)--(6,0);
 \node[greendot] at (3,2){};
 \node[greendot] at (4,2){};
 \draw[loosely dashed,red] (-0.5,3) -- (4.5,3)--(4.5,1)--(-0.5,1)--(-0.5,3);
\end{braid} =
\begin{braid}
 \draw (0,4)node[above]{$0$}--(0,0);
 \draw (1,4)node[above]{$1$}--(1,0);
 \draw[dots] (1.2,4)--(2.8,4);
 \draw[dots] (1.2,0)--(2.8,0);
 \draw (3,4)node[above]{$e-1$}--(3,0);
 \draw (4,4)node[above]{$0$}--(5,0);
 \draw (5,4)node[above]{$e-1$}--(6,0);
 \draw (6,4)node[above]{$0$}--(4,0);
\end{braid}
{\cdot}
\begin{braid}
 \draw (0,4)node[above]{$0$}--(0,0);
 \draw (1,4)node[above]{$1$}--(1,0);
 \draw[dots] (1.2,4)--(2.8,4);
 \draw[dots] (1.2,0)--(2.8,0);
 \draw (3,4)node[above]{$e-1$}--(3,0);
 \draw (4,4)node[above]{$0$}--(4,0);
 \draw (5,4)node[above]{$0$}--(5,0);
 \draw (6,4)node[above]{$e-1$}--(6,0);
 \node[greendot] at (3,2){};
 \node[greendot] at (4,2){};
 \draw[loosely dashed,red] (-0.5,3) -- (4.5,3)--(4.5,1)--(-0.5,1)--(-0.5,3);
\end{braid}
{\cdot}
\begin{braid}
 \draw (0,4)node[above]{$0$}--(0,0);
 \draw (1,4)node[above]{$1$}--(1,0);
 \draw[dots] (1.2,4)--(2.8,4);
 \draw[dots] (1.2,0)--(2.8,0);
 \draw (3,4)node[above]{$e-1$}--(3,0);
 \draw (4,4)node[above]{$0$}--(6,0);
 \draw (5,4)node[above]{$0$}--(4,0);
 \draw (6,4)node[above]{$e-1$}--(5,0);
\end{braid}
$$
and if we define $\nu = (\gamma_1,1) = \lambda|_{\gamma_1+1}$, then as $\lambda \in \mathscr S_n^\Lambda$ and $|\nu| = \gamma_1+1 < n = \gamma_1 + \gamma_2 = |\lambda|$, $\nu \in \P_I$. Moreover as $b_0^{\nu_-} = 1$,
$$
\begin{braid}
 \draw (0,4)node[above]{$0$}--(0,0);
 \draw (1,4)node[above]{$1$}--(1,0);
 \draw[dots] (1.2,4)--(2.8,4);
 \draw[dots] (1.2,0)--(2.8,0);
 \draw (3,4)node[above]{$e-1$}--(3,0);
 \draw (4,4)node[above]{$0$}--(4,0);
 \node[greendot] at (3,2){};
 \node[greendot] at (4,2){};
\end{braid} = e(\bi_{\nu_-}\vee 0) y_{\nu_-} y_{|\nu|}^1 \in \Blam[\nu].
$$

Then by \autoref{send},
$$
\begin{braid}
 \draw (0,4)node[above]{$0$}--(0,0);
 \draw (1,4)node[above]{$1$}--(1,0);
 \draw[dots] (1.2,4)--(2.8,4);
 \draw[dots] (1.2,0)--(2.8,0);
 \draw (3,4)node[above]{$e-1$}--(3,0);
 \draw (4,4)node[above]{$0$}--(5,3)--(5,1)--(4,0);
 \draw (5,4)node[above]{$e-1$}--(6,3)--(6,1)--(5,0);
 \draw (6,4)node[above]{$0$}--(4,3)--(4,1)--(6,0);
 \node[greendot] at (3,2){};
 \node[greendot] at (4,2){};
 \draw[loosely dashed,red] (-0.5,3) -- (4.5,3)--(4.5,1)--(-0.5,1)--(-0.5,3);
\end{braid}\in \Blam[\gamma].
$$

Therefore,
$$
e_\gamma y_\gamma =_\gamma \begin{braid}
 \draw (0,4)node[above]{$0$}--(0,0);
 \draw (1,4)node[above]{$1$}--(1,0);
 \draw[dots] (1.2,4)--(2.8,4);
 \draw[dots] (1.2,0)--(2.8,0);
 \draw (3,4)node[above]{$e-1$}--(3,0);
 \draw (4,4)node[above]{$0$}--(6,0);
 \draw (5,4)node[above]{$e-1$}--(4,2)--(5,0);
 \draw (6,4)node[above]{$0$}--(4,0);
 \node[greendot] at (3,2){};
\end{braid},
$$
which gives the proposition in this case.

\medskip\textbf{Case ~\ref{I-problem: induction1}b:} $1\leq i \leq e-3$, i.e. $3\leq \gamma_2 \leq e-1$.
\begin{eqnarray*}
\begin{braid}
 \draw(2,-0.5)node{$\underbrace{\hspace*{16mm}}_{\gamma_1}$};
 \draw (0,4)node[above]{$0$}--(0,0);
 \draw (1,4)node[above]{$1$}--(1,0);
 \draw[dots] (1.2,4)--(2.8,4);
 \draw[dots] (1.2,0)--(2.8,0);
 \draw (3,4)node[above]{$i-1$}--(3,0);
 \draw (4,4)[densely dotted] node[above]{$i$}--(9,2)--(10,0);
 \draw (5,4)node[above]{$e-1$}--(4,2)--(5,0);
 \draw (6,4)node[above]{$0$}--(5,2)--(6,0);
 \draw[dots] (6.2,4)--(7.8,4);
 \draw[dots] (6.2,0)--(7.8,0);
 \draw (8,4)node[above]{$i-2$}--(7,2)--(8,0);
 \draw (9,4)[densely dotted] node[above]{$i-1$}--(8,2)--(9,0);
 \draw (10,4)[densely dotted] node[above]{$i$}--(9,2)--(4,0);
\end{braid}
& \overset{(\ref{dia:braid})}= &
\begin{braid}
 \draw (0,4)node[above]{$0$}--(0,0);
 \draw (1,4)node[above]{$1$}--(1,0);
 \draw[dots] (1.2,4)--(2.8,4);
 \draw[dots] (1.2,0)--(2.8,0);
 \draw (3,4)node[above]{$i-1$}--(3,0);
 \draw (4,4)[densely dotted] node[above]{$i$}--(9,2)--(10,0);
 \draw (5,4)[densely dotted] node[above]{$e-1$}--(4,2)--(5,0);
 \draw (6,4)[densely dotted] node[above]{$0$}--(5,2)--(6,0);
 \draw[dots] (6.2,4)--(7.8,4);
 \draw[dots] (6.2,0)--(7.8,0);
 \draw (8,4)[densely dotted] node[above]{$i-2$}--(7,2)--(8,0);
 \draw (9,4)node[above]{$i-1$}--(10,2)--(9,0);
 \draw (10,4)[densely dotted] node[above]{$i$}--(9,2)--(4,0);
\end{braid}
-
\begin{braid}
 \draw (0,4)node[above]{$0$}--(0,0);
 \draw (1,4)node[above]{$1$}--(1,0);
 \draw[dots] (1.2,4)--(2.8,4);
 \draw[dots] (1.2,0)--(2.8,0);
 \draw (3,4)node[above]{$i-1$}--(3,0);
 \draw (4,4)node[above]{$i$}--(8,2)--(4,0);
 \draw (5,4)node[above]{$e-1$}--(4,2)--(5,0);
 \draw (6,4)node[above]{$0$}--(5,2)--(6,0);
 \draw[dots] (6.2,4)--(7.8,4);
 \draw[dots] (6.2,0)--(7.8,0);
 \draw (8,4)node[above]{$i-2$}--(7,2)--(8,0);
 \draw (9,4)node[above]{$i-1$}--(9,0);
 \draw (10,4)node[above]{$i$}--(10,0);
\end{braid}\\
& \overset{(\ref{dia:braid})} = &
\begin{braid}
 \draw (0,4)node[above]{$0$}--(0,0);
 \draw (1,4)node[above]{$1$}--(1,0);
 \draw[dots] (1.2,4)--(2.8,4);
 \draw[dots] (1.2,0)--(2.8,0);
 \draw (3,4)node[above]{$i-1$}--(3,0);
 \draw (4,4)node[above]{$i$}--(5,2)--(10,0);
 \draw (5,4)node[above]{$e-1$}--(6,2)--(5,0);
 \draw (6,4)node[above]{$0$}--(7,2)--(6,0);
 \draw[dots] (6.2,4)--(7.8,4);
 \draw[dots] (6.2,0)--(7.8,0);
 \draw (8,4)node[above]{$i-2$}--(9,2)--(8,0);
 \draw (9,4)node[above]{$i-1$}--(10,2)--(9,0);
 \draw (10,4)node[above]{$i$}--(5,2)--(4,0);
\end{braid}
-
\begin{braid}
 \draw (0,4)node[above]{$0$}--(0,0);
 \draw (1,4)node[above]{$1$}--(1,0);
 \draw[dots] (1.2,4)--(2.8,4);
 \draw[dots] (1.2,0)--(2.8,0);
 \draw (3,4)node[above]{$i-1$}--(3,0);
 \draw (4,4)[densely dotted] node[above]{$i$}--(8,2)--(4,0);
 \draw (5,4)[densely dotted] node[above]{$e-1$}--(4,2)--(5,0);
 \draw (6,4)[densely dotted] node[above]{$0$}--(5,2)--(6,0);
 \draw[dots] (6.2,4)--(7.8,4);
 \draw[dots] (6.2,0)--(7.8,0);
 \draw (8,4)[densely dotted] node[above]{$i-2$}--(7,2)--(8,0);
 \draw (9,4)node[above]{$i-1$}--(9,0);
 \draw (10,4)node[above]{$i$}--(10,0);
\end{braid}\\
& \overset{(\ref{dia:psipsi})}= &
\begin{braid}
 \draw (0,4)node[above]{$0$}--(0,0);
 \draw (1,4)node[above]{$1$}--(1,0);
 \draw[dots] (1.2,4)--(2.8,4);
 \draw[dots] (1.2,0)--(2.8,0);
 \draw (3,4)node[above]{$i-1$}--(3,0);
 \draw (4,4)node[above]{$i$}--(5,2)--(10,0);
 \draw (5,4)node[above]{$e-1$}--(6,2)--(5,0);
 \draw (6,4)node[above]{$0$}--(7,2)--(6,0);
 \draw[dots] (6.2,4)--(7.8,4);
 \draw[dots] (6.2,0)--(7.8,0);
 \draw (8,4)node[above]{$i-2$}--(9,2)--(8,0);
 \draw (9,4)node[above]{$i-1$}--(10,2)--(9,0);
 \draw (10,4)node[above]{$i$}--(5,2)--(4,0);
\end{braid}
-
\begin{braid}
 \draw (0,4)node[above]{$0$}--(0,0);
 \draw (1,4)node[above]{$1$}--(1,0);
 \draw[dots] (1.2,4)--(2.8,4);
 \draw[dots] (1.2,0)--(2.8,0);
 \draw (3,4)node[above]{$i-1$}--(3,0);
 \draw (4,4)node[above]{$i$}--(4,0);
 \draw (5,4)node[above]{$e-1$}--(5,0);
 \draw (6,4)node[above]{$0$}--(6,0);
 \draw[dots] (6.2,4)--(7.8,4);
 \draw[dots] (6.2,0)--(7.8,0);
 \draw (8,4)node[above]{$i-2$}--(8,0);
 \draw (9,4)node[above]{$i-1$}--(9,0);
 \draw (10,4)node[above]{$i$}--(10,0);
\end{braid}\\
& = &
\begin{braid}
 \draw (0,4)node[above]{$0$}--(0,0);
 \draw (1,4)node[above]{$1$}--(1,0);
 \draw[dots] (1.2,4)--(2.8,4);
 \draw[dots] (1.2,0)--(2.8,0);
 \draw (3,4)node[above]{$i-1$}--(3,0);
 \draw (4,4)node[above]{$i$}--(4,1)--(4.5,0.5)--(10,0);
 \draw (5,4)node[above]{$e-1$}--(6,3)--(6,1)--(5,0);
 \draw (6,4)node[above]{$0$}--(7,3)--(7,1)--(6,0);
 \draw[dots] (6.2,4)--(7.8,4);
 \draw[dots] (6.2,0)--(7.8,0);
 \draw (8,4)node[above]{$i-2$}--(9,3)--(9,1)--(8,0);
 \draw (9,4)node[above]{$i-1$}--(10,3)--(10,1)--(9,0);
 \draw (10,4)node[above]{$i$}--(5,3)--(5,1)--(4.5,0.5)--(4,0);
 \draw[loosely dashed,red] (-0.5,3) -- (5.5,3)--(5.5,1)--(-0.5,1)--(-0.5,3);
\end{braid}
-
\begin{braid}
 \draw(2,-0.5)node{$\underbrace{\hspace*{16mm}}_{\gamma_1}$};
 \draw (0,4)node[above]{$0$}--(0,0);
 \draw (1,4)node[above]{$1$}--(1,0);
 \draw[dots] (1.2,4)--(2.8,4);
 \draw[dots] (1.2,0)--(2.8,0);
 \draw (3,4)node[above]{$i-1$}--(3,0);
 \draw (4,4)node[above]{$i$}--(4,0);
 \draw (5,4)node[above]{$e-1$}--(5,0);
 \draw (6,4)node[above]{$0$}--(6,0);
 \draw[dots] (6.2,4)--(7.8,4);
 \draw[dots] (6.2,0)--(7.8,0);
 \draw (8,4)node[above]{$i-2$}--(8,0);
 \draw (9,4)node[above]{$i-1$}--(9,0);
 \draw (10,4)node[above]{$i$}--(10,0);
\end{braid}.
\end{eqnarray*}

For the same reason as in Case ~\ref{I-problem: induction1}a,
$$
\begin{braid}
 \draw (0,4)node[above]{$0$}--(0,0);
 \draw (1,4)node[above]{$1$}--(1,0);
 \draw[dots] (1.2,4)--(2.8,4);
 \draw[dots] (1.2,0)--(2.8,0);
 \draw (3,4)node[above]{$i-1$}--(3,0);
 \draw (4,4)node[above]{$i$}--(5,2)--(10,0);
 \draw (5,4)node[above]{$e-1$}--(6,2)--(5,0);
 \draw (6,4)node[above]{$0$}--(7,2)--(6,0);
 \draw[dots] (6.2,4)--(7.8,4);
 \draw[dots] (6.2,0)--(7.8,0);
 \draw (8,4)node[above]{$i-2$}--(9,2)--(8,0);
 \draw (9,4)node[above]{$i-1$}--(10,2)--(9,0);
 \draw (10,4)node[above]{$i$}--(5,2)--(4,0);
\end{braid}
\in \Blam[\gamma],
$$
which implies the proposition in this case.

\medskip\textbf{Case ~\ref{I-problem: induction1}c:} $i = e-2$, i.e. $\gamma_2 = e$. By \autoref{I-problem: short},
\begin{eqnarray*}
\begin{braid}
 \draw(2,-0.5)node{$\underbrace{\hspace*{16mm}}_{\gamma_1}$};
 \draw (0,4)node[above]{$0$}--(0,0);
 \draw (1,4)node[above]{$1$}--(1,0);
 \draw[dots] (1.2,4)--(2.8,4);
 \draw[dots] (1.2,0)--(2.8,0);
 \draw (3,4)node[above]{$e-3$}--(3,0);
 \draw (4,4)node[above]{$e-2$}--(8,2)--(9,0);
 \draw (5,4)node[above]{$e-1$}--(4,2)--(5,0);
 \draw (6,4)node[above]{$0$}--(5,2)--(6,0);
 \draw[dots] (6.2,4)--(7.8,4);
 \draw[dots] (6.2,0)--(7.8,0);
 \draw (8,4)node[above]{$e-3$}--(7,2)--(8,0);
 \draw (9,4)node[above]{$e-2$}--(8,2)--(4,0);
\end{braid}
& = &
\begin{braid}
 \draw(2,-0.5)node{$\underbrace{\hspace*{16mm}}_{\gamma_1}$};
 \draw (0,4)node[above]{$0$}--(0,0);
 \draw (1,4)node[above]{$1$}--(1,0);
 \draw[dots] (1.2,4)--(2.8,4);
 \draw[dots] (1.2,0)--(2.8,0);
 \draw (3,4)node[above]{$e-3$}--(3,0);
 \draw (4,4)node[above]{$e-2$}--(4,0);
 \draw (5,4)node[above]{$e-1$}--(5,0);
 \draw (6,4)node[above]{$0$}--(6,0);
 \draw[dots] (6.2,4)--(7.8,4);
 \draw[dots] (6.2,0)--(7.8,0);
 \draw (8,4)node[above]{$e-3$}--(8,0);
 \draw (9,4)node[above]{$e-2$}--(9,0);
 \node[greendot] at (9,2){};
\end{braid}
-
\begin{braid}
 \draw (0,4)node[above]{$0$}--(0,0);
 \draw (1,4)node[above]{$1$}--(1,0);
 \draw[dots] (1.2,4)--(2.8,4);
 \draw[dots] (1.2,0)--(2.8,0);
 \draw (3,4)node[above]{$e-3$}--(3,0);
 \draw (4,4)node[above]{$e-2$}--(4,0);
 \draw (5,4)node[above]{$e-1$}--(5,0);
 \draw (6,4)node[above]{$0$}--(6,0);
 \draw[dots] (6.2,4)--(7.8,4);
 \draw[dots] (6.2,0)--(7.8,0);
 \draw (8,4)node[above]{$e-3$}--(8,0);
 \draw (9,4)node[above]{$e-2$}--(9,0);
 \node[greendot] at (8,2){};
\end{braid}\\
& - &
\begin{braid}
 \draw (0,4)node[above]{$0$}--(0,0);
 \draw (1,4)node[above]{$1$}--(1,0);
 \draw[dots] (1.2,4)--(2.8,4);
 \draw[dots] (1.2,0)--(2.8,0);
 \draw (3,4)node[above]{$e-3$}--(3,0);
 \draw (4,4)node[above]{$e-2$}--(4,0);
 \draw (5,4)node[above]{$e-1$}--(5,0);
 \draw (6,4)node[above]{$0$}--(6,0);
 \draw[dots] (6.2,4)--(7.8,4);
 \draw[dots] (6.2,0)--(7.8,0);
 \draw (8,4)node[above]{$e-3$}--(8,0);
 \draw (9,4)node[above]{$e-2$}--(9,0);
 \node[greendot] at (5,2){};
\end{braid}
+
\begin{braid}
 \draw (0,4)node[above]{$0$}--(0,0);
 \draw (1,4)node[above]{$1$}--(1,0);
 \draw[dots] (1.2,4)--(2.8,4);
 \draw[dots] (1.2,0)--(2.8,0);
 \draw (3,4)node[above]{$e-3$}--(3,0);
 \draw (4,4)node[above]{$e-2$}--(4,0);
 \draw (5,4)node[above]{$e-1$}--(5,0);
 \draw (6,4)node[above]{$0$}--(6,0);
 \draw[dots] (6.2,4)--(7.8,4);
 \draw[dots] (6.2,0)--(7.8,0);
 \draw (8,4)node[above]{$e-3$}--(8,0);
 \draw (9,4)node[above]{$e-2$}--(9,0);
 \node[greendot] at (4,2){};
\end{braid}
+
\begin{braid}
 \draw (0,4)node[above]{$0$}--(0,0);
 \draw (1,4)node[above]{$1$}--(1,0);
 \draw[dots] (1.2,4)--(2.8,4);
 \draw[dots] (1.2,0)--(2.8,0);
 \draw (3,4)node[above]{$e-3$}--(3,0);
 \draw (4,4)node[above]{$e-2$}--(5,2)--(9,0);
 \draw (5,4)node[above]{$e-1$}--(6,2)--(5,0);
 \draw (6,4)node[above]{$0$}--(7,2)--(6,0);
 \draw[dots] (6.2,4)--(7.8,4);
 \draw[dots] (6.2,0)--(7.8,0);
 \draw (8,4)node[above]{$e-3$}--(9,2)--(8,0);
 \draw (9,4)node[above]{$e-2$}--(5,2)--(4,0);
\end{braid}.
\end{eqnarray*}

Set $\nu = (\gamma_1,\gamma_2-1) = \gamma|_{n-1}$. As $\gamma_1 \geq \gamma_2-1$, we have $\nu\in \P_{n-1}$. As $\lambda\in\mathscr S_n^\Lambda$ and $|\nu| < |\lambda|$, we have $\nu \in \P_I$. It is not hard to see that $b_{e-3}^{\nu_-} = 1$. Hence
$$
\begin{braid}
 \draw(2,-0.5)node{$\underbrace{\hspace*{16mm}}_{\gamma_1}$};
 \draw(7,-0.5)node{$\underbrace{\hspace*{16mm}}_{\gamma_2-1}$};
 \draw (0,4)node[above]{$0$}--(0,0);
 \draw (1,4)node[above]{$1$}--(1,0);
 \draw[dots] (1.2,4)--(2.8,4);
 \draw[dots] (1.2,0)--(2.8,0);
 \draw (3,4)node[above]{$e-3$}--(3,0);
 \draw (4,4)node[above]{$e-2$}--(4,0);
 \draw (5,4)node[above]{$e-1$}--(5,0);
 \draw (6,4)node[above]{$0$}--(6,0);
 \draw[dots] (6.2,4)--(7.8,4);
 \draw[dots] (6.2,0)--(7.8,0);
 \draw (8,4)node[above]{$e-4$}--(8,0);
 \draw (9,4)node[above]{$e-3$}--(9,0);
 \node[greendot] at (9,2){};
\end{braid} = e(\bi_{\nu_-}\vee e-3)y_{\nu_-} y_{n-1}^1 \in \Blam[\nu].
$$

Then by \autoref{send},
$$
\begin{braid}
 \draw (0,4)node[above]{$0$}--(0,0);
 \draw (1,4)node[above]{$1$}--(1,0);
 \draw[dots] (1.2,4)--(2.8,4);
 \draw[dots] (1.2,0)--(2.8,0);
 \draw (3,4)node[above]{$e-3$}--(3,0);
 \draw (4,4)node[above]{$e-2$}--(4,0);
 \draw (5,4)node[above]{$e-1$}--(5,0);
 \draw (6,4)node[above]{$0$}--(6,0);
 \draw[dots] (6.2,4)--(7.8,4);
 \draw[dots] (6.2,0)--(7.8,0);
 \draw (8,4)node[above]{$e-3$}--(8,0);
 \draw (9,4)node[above]{$e-2$}--(9,0);
 \node[greendot] at (8,2){};
\end{braid} \in \Blam[\gamma].
$$

Similarly, we have
$$
\begin{braid}
 \draw (0,4)node[above]{$0$}--(0,0);
 \draw (1,4)node[above]{$1$}--(1,0);
 \draw[dots] (1.2,4)--(2.8,4);
 \draw[dots] (1.2,0)--(2.8,0);
 \draw (3,4)node[above]{$e-3$}--(3,0);
 \draw (4,4)node[above]{$e-2$}--(4,0);
 \draw (5,4)node[above]{$e-1$}--(5,0);
 \draw (6,4)node[above]{$0$}--(6,0);
 \draw[dots] (6.2,4)--(7.8,4);
 \draw[dots] (6.2,0)--(7.8,0);
 \draw (8,4)node[above]{$e-3$}--(8,0);
 \draw (9,4)node[above]{$e-2$}--(9,0);
 \node[greendot] at (5,2){};
\end{braid}
,
\begin{braid}
 \draw (0,4)node[above]{$0$}--(0,0);
 \draw (1,4)node[above]{$1$}--(1,0);
 \draw[dots] (1.2,4)--(2.8,4);
 \draw[dots] (1.2,0)--(2.8,0);
 \draw (3,4)node[above]{$e-3$}--(3,0);
 \draw (4,4)node[above]{$e-2$}--(4,0);
 \draw (5,4)node[above]{$e-1$}--(5,0);
 \draw (6,4)node[above]{$0$}--(6,0);
 \draw[dots] (6.2,4)--(7.8,4);
 \draw[dots] (6.2,0)--(7.8,0);
 \draw (8,4)node[above]{$e-3$}--(8,0);
 \draw (9,4)node[above]{$e-2$}--(9,0);
 \node[greendot] at (4,2){};
\end{braid} \in \Blam[\gamma],
$$
and for the similar argument as in Case ~\ref{I-problem: induction1}a, we have
$$
\begin{braid}
 \draw (0,4)node[above]{$0$}--(0,0);
 \draw (1,4)node[above]{$1$}--(1,0);
 \draw[dots] (1.2,4)--(2.8,4);
 \draw[dots] (1.2,0)--(2.8,0);
 \draw (3,4)node[above]{$e-3$}--(3,0);
 \draw (4,4)node[above]{$e-2$}--(5,2)--(9,0);
 \draw (5,4)node[above]{$e-1$}--(6,2)--(5,0);
 \draw (6,4)node[above]{$0$}--(7,2)--(6,0);
 \draw[dots] (6.2,4)--(7.8,4);
 \draw[dots] (6.2,0)--(7.8,0);
 \draw (8,4)node[above]{$e-3$}--(9,2)--(8,0);
 \draw (9,4)node[above]{$e-2$}--(5,2)--(4,0);
\end{braid} \in \Blam[\gamma].
$$

Therefore,
$$
\begin{braid}
 \draw (0,4)node[above]{$0$}--(0,0);
 \draw (1,4)node[above]{$1$}--(1,0);
 \draw[dots] (1.2,4)--(2.8,4);
 \draw[dots] (1.2,0)--(2.8,0);
 \draw (3,4)node[above]{$e-3$}--(3,0);
 \draw (4,4)node[above]{$e-2$}--(8,2)--(9,0);
 \draw (5,4)node[above]{$e-1$}--(4,2)--(5,0);
 \draw (6,4)node[above]{$0$}--(5,2)--(6,0);
 \draw[dots] (6.2,4)--(7.8,4);
 \draw[dots] (6.2,0)--(7.8,0);
 \draw (8,4)node[above]{$e-3$}--(7,2)--(8,0);
 \draw (9,4)node[above]{$e-2$}--(8,2)--(4,0);
\end{braid} =_\gamma
\begin{braid}
 \draw (0,4)node[above]{$0$}--(0,0);
 \draw (1,4)node[above]{$1$}--(1,0);
 \draw[dots] (1.2,4)--(2.8,4);
 \draw[dots] (1.2,0)--(2.8,0);
 \draw (3,4)node[above]{$e-3$}--(3,0);
 \draw (4,4)node[above]{$e-2$}--(4,0);
 \draw (5,4)node[above]{$e-1$}--(5,0);
 \draw (6,4)node[above]{$0$}--(6,0);
 \draw[dots] (6.2,4)--(7.8,4);
 \draw[dots] (6.2,0)--(7.8,0);
 \draw (8,4)node[above]{$e-3$}--(8,0);
 \draw (9,4)node[above]{$e-2$}--(9,0);
 \node[greendot] at (9,2){};
\end{braid} = e_\gamma y_\gamma,
$$
which follows the Proposition.

\medskip\textbf{Case ~\ref{I-problem: induction1}d:} $i = e-1$, i.e. $\gamma_2 = e + 1$. By \autoref{I-problem: short},
\begin{eqnarray}
&&
\begin{braid}
 \draw(2,-0.5)node{$\underbrace{\hspace*{16mm}}_{\gamma_1}$};
 \draw (0,4)node[above]{$0$}--(0,0);
 \draw (1,4)node[above]{$1$}--(1,0);
 \draw[dots] (1.2,4)--(2.8,4);
 \draw[dots] (1.2,0)--(2.8,0);
 \draw (3,4)node[above]{$e-2$}--(3,0);
 \draw (4,4)node[above]{$e-1$}--(4,0);
 \draw (5,4)node[above]{$e-1$}--(10,2)--(11,0);
 \draw (6,4)node[above]{$0$}--(5,2)--(6,0);
 \draw (7,4)node[above]{$1$}--(6,2)--(7,0);
 \draw[dots] (7.2,4)--(8.8,4);
 \draw[dots] (7.2,0)--(8.8,0);
 \draw (9,4)node[above]{$e-3$}--(8,2)--(9,0);
 \draw (10,4)node[above]{$e-2$}--(9,2)--(10,0);
 \draw (11,4)node[above]{$e-1$}--(10,2)--(5,0);
 \node[greendot] at (9,2){};
\end{braid}
 =
\begin{braid}
 \draw(2,-0.5)node{$\underbrace{\hspace*{16mm}}_{\gamma_1}$};
 \draw (0,4)node[above]{$0$}--(0,0);
 \draw (1,4)node[above]{$1$}--(1,0);
 \draw[dots] (1.2,4)--(2.8,4);
 \draw[dots] (1.2,0)--(2.8,0);
 \draw (3,4)node[above]{$e-2$}--(3,0);
 \draw (4,4)node[above]{$e-1$}--(4,0);
 \draw (5,4)node[above]{$e-1$}--(5,0);
 \draw (6,4)node[above]{$0$}--(6,0);
 \draw (7,4)node[above]{$1$}--(7,0);
 \draw[dots] (7.2,4)--(8.8,4);
 \draw[dots] (7.2,0)--(8.8,0);
 \draw (9,4)node[above]{$e-3$}--(9,0);
 \draw (10,4)node[above]{$e-2$}--(10,0);
 \draw (11,4)node[above]{$e-1$}--(11,0);
 \node[greendot] at (10,2){};
 \node[greendot] at (11,2){};
\end{braid}\notag\\
&&-
\begin{braid}
 \draw (0,4)node[above]{$0$}--(0,0);
 \draw (1,4)node[above]{$1$}--(1,0);
 \draw[dots] (1.2,4)--(2.8,4);
 \draw[dots] (1.2,0)--(2.8,0);
 \draw (3,4)node[above]{$e-2$}--(3,0);
 \draw (4,4)node[above]{$e-1$}--(4,0);
 \draw (5,4)node[above]{$e-1$}--(5,0);
 \draw (6,4)node[above]{$0$}--(6,0);
 \draw (7,4)node[above]{$1$}--(7,0);
 \draw[dots] (7.2,4)--(8.8,4);
 \draw[dots] (7.2,0)--(8.8,0);
 \draw (9,4)node[above]{$e-3$}--(9,0);
 \draw (10,4)node[above]{$e-2$}--(10,0);
 \draw (11,4)node[above]{$e-1$}--(11,0);
 \node[greendot] at (10,3){};
 \node[greendot] at (10,1){};
\end{braid}
-
\begin{braid}
 \draw (0,4)node[above]{$0$}--(0,0);
 \draw (1,4)node[above]{$1$}--(1,0);
 \draw[dots] (1.2,4)--(2.8,4);
 \draw[dots] (1.2,0)--(2.8,0);
 \draw (3,4)node[above]{$e-2$}--(3,0);
 \draw (4,4)node[above]{$e-1$}--(4,0);
 \draw (5,4)node[above]{$e-1$}--(5,0);
 \draw (6,4)node[above]{$0$}--(6,0);
 \draw (7,4)node[above]{$1$}--(7,0);
 \draw[dots] (7.2,4)--(8.8,4);
 \draw[dots] (7.2,0)--(8.8,0);
 \draw (9,4)node[above]{$e-3$}--(9,0);
 \draw (10,4)node[above]{$e-2$}--(10,0);
 \draw (11,4)node[above]{$e-1$}--(11,0);
 \node[greendot] at (10,2){};
 \node[greendot] at (6,2){};
\end{braid} \notag \\
&& +
\begin{braid}
 \draw (0,4)node[above]{$0$}--(0,0);
 \draw (1,4)node[above]{$1$}--(1,0);
 \draw[dots] (1.2,4)--(2.8,4);
 \draw[dots] (1.2,0)--(2.8,0);
 \draw (3,4)node[above]{$e-2$}--(3,0);
 \draw (4,4)node[above]{$e-1$}--(4,0);
 \draw (5,4)node[above]{$e-1$}--(5,0);
 \draw (6,4)node[above]{$0$}--(6,0);
 \draw (7,4)node[above]{$1$}--(7,0);
 \draw[dots] (7.2,4)--(8.8,4);
 \draw[dots] (7.2,0)--(8.8,0);
 \draw (9,4)node[above]{$e-3$}--(9,0);
 \draw (10,4)node[above]{$e-2$}--(10,0);
 \draw (11,4)node[above]{$e-1$}--(11,0);
 \node[greendot] at (5,2){};
 \node[greendot] at (10,2){};
\end{braid}
+
\begin{braid}
 \draw (0,4)node[above]{$0$}--(0,0);
 \draw (1,4)node[above]{$1$}--(1,0);
 \draw[dots] (1.2,4)--(2.8,4);
 \draw[dots] (1.2,0)--(2.8,0);
 \draw (3,4)node[above]{$e-2$}--(3,0);
 \draw (4,4)node[above]{$e-1$}--(4,0);
 \draw (5,4)node[above]{$e-1$}--(6,2)--(11,0);
 \draw (6,4)node[above]{$0$}--(7,2)--(6,0);
 \draw (7,4)node[above]{$1$}--(8,2)--(7,0);
 \draw[dots] (7.2,4)--(8.8,4);
 \draw[dots] (7.2,0)--(8.8,0);
 \draw (9,4)node[above]{$e-3$}--(10,2)--(9,0);
 \draw (10,4)node[above]{$e-2$}--(11,2)--(10,0);
 \draw (11,4)node[above]{$e-1$}--(6,2)--(5,0);
 \node[greendot] at (11,2){};
\end{braid} \label{I-problem: short: equation1}.
\end{eqnarray}

For the first two terms of (\ref{I-problem: short: equation1}),
\begin{align*}
& \begin{braid}
 \draw (0,4)node[above]{$0$}--(0,0);
 \draw (1,4)node[above]{$1$}--(1,0);
 \draw[dots] (1.2,4)--(2.8,4);
 \draw[dots] (1.2,0)--(2.8,0);
 \draw (3,4)node[above]{$e-2$}--(3,0);
 \draw (4,4)node[above]{$e-1$}--(4,0);
 \draw (5,4)node[above]{$e-1$}--(5,0);
 \draw (6,4)node[above]{$0$}--(6,0);
 \draw (7,4)node[above]{$1$}--(7,0);
 \draw[dots] (7.2,4)--(8.8,4);
 \draw[dots] (7.2,0)--(8.8,0);
 \draw (9,4)node[above]{$e-3$}--(9,0);
 \draw (10,4)node[above]{$e-2$}--(10,0);
 \draw (11,4)node[above]{$e-1$}--(11,0);
 \node[greendot] at (10,2){};
 \node[greendot] at (11,2){};
\end{braid} -
\begin{braid}
 \draw (0,4)node[above]{$0$}--(0,0);
 \draw (1,4)node[above]{$1$}--(1,0);
 \draw[dots] (1.2,4)--(2.8,4);
 \draw[dots] (1.2,0)--(2.8,0);
 \draw (3,4)node[above]{$e-2$}--(3,0);
 \draw (4,4)node[above]{$e-1$}--(4,0);
 \draw (5,4)node[above]{$e-1$}--(5,0);
 \draw (6,4)node[above]{$0$}--(6,0);
 \draw (7,4)node[above]{$1$}--(7,0);
 \draw[dots] (7.2,4)--(8.8,4);
 \draw[dots] (7.2,0)--(8.8,0);
 \draw (9,4)node[above]{$e-3$}--(9,0);
 \draw (10,4)node[above]{$e-2$}--(10,0);
 \draw (11,4)node[above]{$e-1$}--(11,0);
 \node[greendot] at (10,3){};
 \node[greendot] at (10,1){};
\end{braid}
\overset{(\ref{dia:psipsi})}=
\begin{braid}
 \draw (0,4)node[above]{$0$}--(0,0);
 \draw (1,4)node[above]{$1$}--(1,0);
 \draw[dots] (1.2,4)--(2.8,4);
 \draw[dots] (1.2,0)--(2.8,0);
 \draw (3,4)node[above]{$e-2$}--(3,0);
 \draw (4,4)[densely dotted] node[above]{$e-1$}--(4,0);
 \draw (5,4)[densely dotted] node[above]{$e-1$}--(5,0);
 \draw (6,4)node[above]{$0$}--(6,0);
 \draw (7,4)node[above]{$1$}--(7,0);
 \draw[dots] (7.2,4)--(8.8,4);
 \draw[dots] (7.2,0)--(8.8,0);
 \draw (9,4)node[above]{$e-3$}--(9,0);
 \draw (10,4)node[above]{$e-2$}--(11,2)--(10,0);
 \draw (11,4)node[above]{$e-1$}--(10,2)--(11,0);
 \node[greendot] at (11,2){};
\end{braid}\\
\overset{(\ref{dia:ii})}= &
-\begin{braid}
 \draw (0,4)node[above]{$0$}--(0,0);
 \draw (1,4)node[above]{$1$}--(1,0);
 \draw[dots] (1.2,4)--(2.8,4);
 \draw[dots] (1.2,0)--(2.8,0);
 \draw (3,4)node[above]{$e-2$}--(3,0);
 \draw (4,4)node[above]{$e-1$}--(5,3)--(5,1)--(4,0);
 \draw (5,4)node[above]{$e-1$}--(4,3)--(4,1)--(5,0);
 \draw (6,4)node[above]{$0$}--(6,0);
 \draw (7,4)node[above]{$1$}--(7,0);
 \draw[dots] (7.2,4)--(8.8,4);
 \draw[dots] (7.2,0)--(8.8,0);
 \draw (9,4)node[above]{$e-3$}--(9,0);
 \draw (10,4)node[above]{$e-2$}--(11,2)--(10,0);
 \draw (11,4)node[above]{$e-1$}--(10,2)--(11,0);
 \node[greendot] at (4,2.7){};
 \node[greendot] at (4,1.3){};
 \node[greendot] at (11,2){};
 \draw[loosely dashed,red] (-0.5,3) -- (4.5,3)--(4.5,1)--(-0.5,1)--(-0.5,3);
\end{braid} -
\begin{braid}
 \draw (0,4)node[above]{$0$}--(0,0);
 \draw (1,4)node[above]{$1$}--(1,0);
 \draw[dots] (1.2,4)--(2.8,4);
 \draw[dots] (1.2,0)--(2.8,0);
 \draw (3,4)node[above]{$e-2$}--(3,0);
 \draw (4,4)node[above]{$e-1$}--(4,1)--(5,0);
 \draw (5,4)node[above]{$e-1$}--(5,1)--(4,0);
 \draw (6,4)node[above]{$0$}--(6,0);
 \draw (7,4)node[above]{$1$}--(7,0);
 \draw[dots] (7.2,4)--(8.8,4);
 \draw[dots] (7.2,0)--(8.8,0);
 \draw (9,4)node[above]{$e-3$}--(9,0);
 \draw (10,4)node[above]{$e-2$}--(11,3)--(11,1)--(10,0);
 \draw (11,4)node[above]{$e-1$}--(10,3)--(10,1)--(11,0);
 \node[greendot] at (4,2){};
 \node[greendot] at (11,2){};
 \draw[loosely dashed,red] (-0.5,3) -- (10.5,3)--(10.5,1)--(-0.5,1)--(-0.5,3);
\end{braid} -
\begin{braid}
 \draw (0,4)node[above]{$0$}--(0,0);
 \draw (1,4)node[above]{$1$}--(1,0);
 \draw[dots] (1.2,4)--(2.8,4);
 \draw[dots] (1.2,0)--(2.8,0);
 \draw (3,4)node[above]{$e-2$}--(3,0);
 \draw (4,4)node[above]{$e-1$}--(5,3)--(5,0);
 \draw (5,4)node[above]{$e-1$}--(4,3)--(4,0);
 \draw (6,4)node[above]{$0$}--(6,0);
 \draw (7,4)node[above]{$1$}--(7,0);
 \draw[dots] (7.2,4)--(8.8,4);
 \draw[dots] (7.2,0)--(8.8,0);
 \draw (9,4)node[above]{$e-3$}--(9,0);
 \draw (10,4)node[above]{$e-2$}--(11,3)--(11,1)--(10,0);
 \draw (11,4)node[above]{$e-1$}--(10,3)--(10,1)--(11,0);
 \node[greendot] at (4,2){};
 \node[greendot] at (11,2){};
 \draw[loosely dashed,red] (-0.5,3) -- (10.5,3)--(10.5,1)--(-0.5,1)--(-0.5,3);
\end{braid}.
\end{align*}

Using the similar argument as in Case ~\ref{I-problem: induction1}a,
$$
\begin{braid}
 \draw (0,4)node[above]{$0$}--(0,0);
 \draw (1,4)node[above]{$1$}--(1,0);
 \draw[dots] (1.2,4)--(2.8,4);
 \draw[dots] (1.2,0)--(2.8,0);
 \draw (3,4)node[above]{$e-2$}--(3,0);
 \draw (4,4)node[above]{$e-1$}--(5,3)--(5,1)--(4,0);
 \draw (5,4)node[above]{$e-1$}--(4,3)--(4,1)--(5,0);
 \draw (6,4)node[above]{$0$}--(6,0);
 \draw (7,4)node[above]{$1$}--(7,0);
 \draw[dots] (7.2,4)--(8.8,4);
 \draw[dots] (7.2,0)--(8.8,0);
 \draw (9,4)node[above]{$e-3$}--(9,0);
 \draw (10,4)node[above]{$e-2$}--(10,0);
 \draw (11,4)node[above]{$e-1$}--(11,0);
 \node[greendot] at (4,2.7){};
 \node[greendot] at (4,1.3){};
 \node[greendot] at (10,2){};
 \node[greendot] at (11,2){};
 \draw[loosely dashed,red] (-0.5,3) -- (4.5,3)--(4.5,1)--(-0.5,1)--(-0.5,3);
\end{braid} \in \Blam[\gamma],
$$
and using the similar argument in Case~\ref{I-problem: induction1}c,
$$
-\begin{braid}
 \draw (0,4)node[above]{$0$}--(0,0);
 \draw (1,4)node[above]{$1$}--(1,0);
 \draw[dots] (1.2,4)--(2.8,4);
 \draw[dots] (1.2,0)--(2.8,0);
 \draw (3,4)node[above]{$e-2$}--(3,0);
 \draw (4,4)node[above]{$e-1$}--(4,1)--(5,0);
 \draw (5,4)node[above]{$e-1$}--(5,1)--(4,0);
 \draw (6,4)node[above]{$0$}--(6,0);
 \draw (7,4)node[above]{$1$}--(7,0);
 \draw[dots] (7.2,4)--(8.8,4);
 \draw[dots] (7.2,0)--(8.8,0);
 \draw (9,4)node[above]{$e-3$}--(9,0);
 \draw (10,4)node[above]{$e-2$}--(11,3)--(11,1)--(10,0);
 \draw (11,4)node[above]{$e-1$}--(10,3)--(10,1)--(11,0);
 \node[greendot] at (4,2){};
 \node[greendot] at (11,2){};
 \draw[loosely dashed,red] (-0.5,3) -- (10.5,3)--(10.5,1)--(-0.5,1)--(-0.5,3);
\end{braid} -
\begin{braid}
 \draw (0,4)node[above]{$0$}--(0,0);
 \draw (1,4)node[above]{$1$}--(1,0);
 \draw[dots] (1.2,4)--(2.8,4);
 \draw[dots] (1.2,0)--(2.8,0);
 \draw (3,4)node[above]{$e-2$}--(3,0);
 \draw (4,4)node[above]{$e-1$}--(5,3)--(5,0);
 \draw (5,4)node[above]{$e-1$}--(4,3)--(4,0);
 \draw (6,4)node[above]{$0$}--(6,0);
 \draw (7,4)node[above]{$1$}--(7,0);
 \draw[dots] (7.2,4)--(8.8,4);
 \draw[dots] (7.2,0)--(8.8,0);
 \draw (9,4)node[above]{$e-3$}--(9,0);
 \draw (10,4)node[above]{$e-2$}--(11,3)--(11,1)--(10,0);
 \draw (11,4)node[above]{$e-1$}--(10,3)--(10,1)--(11,0);
 \node[greendot] at (4,2){};
 \node[greendot] at (11,2){};
 \draw[loosely dashed,red] (-0.5,3) -- (10.5,3)--(10.5,1)--(-0.5,1)--(-0.5,3);
\end{braid}\in \Blam[\gamma].
$$

Hence the first two terms of (\ref{I-problem: short: equation1}) is in $\Blam[\gamma]$. For the third term of (\ref{I-problem: short: equation1}),
\begin{align*}
&\begin{braid}
 \draw (0,4)node[above]{$0$}--(0,0);
 \draw (1,4)node[above]{$1$}--(1,0);
 \draw[dots] (1.2,4)--(2.8,4);
 \draw[dots] (1.2,0)--(2.8,0);
 \draw (3,4)node[above]{$e-2$}--(3,0);
 \draw (4,4)[densely dotted] node[above]{$e-1$}--(4,0);
 \draw (5,4)[densely dotted] node[above]{$e-1$}--(5,0);
 \draw (6,4)node[above]{$0$}--(6,0);
 \draw (7,4)node[above]{$1$}--(7,0);
 \draw[dots] (7.2,4)--(8.8,4);
 \draw[dots] (7.2,0)--(8.8,0);
 \draw (9,4)node[above]{$e-3$}--(9,0);
 \draw (10,4)node[above]{$e-2$}--(10,0);
 \draw (11,4)node[above]{$e-1$}--(11,0);
 \node[greendot] at (10,2){};
 \node[greendot] at (6,2){};
\end{braid}\\
\overset{(\ref{dia:ii})}= &
-\begin{braid}
 \draw (0,4)node[above]{$0$}--(0,0);
 \draw (1,4)node[above]{$1$}--(1,0);
 \draw[dots] (1.2,4)--(2.8,4);
 \draw[dots] (1.2,0)--(2.8,0);
 \draw (3,4)node[above]{$e-2$}--(3,0);
 \draw (4,4)node[above]{$e-1$}--(5,3)--(5,1)--(4,0);
 \draw (5,4)node[above]{$e-1$}--(4,3)--(4,1)--(5,0);
 \draw (6,4)node[above]{$0$}--(6,0);
 \draw (7,4)node[above]{$1$}--(7,0);
 \draw[dots] (7.2,4)--(8.8,4);
 \draw[dots] (7.2,0)--(8.8,0);
 \draw (9,4)node[above]{$e-3$}--(9,0);
 \draw (10,4)node[above]{$e-2$}--(10,0);
 \draw (11,4)node[above]{$e-1$}--(11,0);
 \node[greendot] at (4,2.7){};
 \node[greendot] at (4,1.3){};
 \node[greendot] at (10,2){};
 \node[greendot] at (6,2){};
 \draw[loosely dashed,red] (-0.5,3) -- (4.5,3)--(4.5,1)--(-0.5,1)--(-0.5,3);
\end{braid} -
\begin{braid}
 \draw (0,4)node[above]{$0$}--(0,0);
 \draw (1,4)node[above]{$1$}--(1,0);
 \draw[dots] (1.2,4)--(2.8,4);
 \draw[dots] (1.2,0)--(2.8,0);
 \draw (3,4)node[above]{$e-2$}--(3,0);
 \draw (4,4)node[above]{$e-1$}--(4,1)--(5,0);
 \draw (5,4)node[above]{$e-1$}--(5,1)--(4,0);
 \draw (6,4)node[above]{$0$}--(6,0);
 \draw (7,4)node[above]{$1$}--(7,0);
 \draw[dots] (7.2,4)--(8.8,4);
 \draw[dots] (7.2,0)--(8.8,0);
 \draw (9,4)node[above]{$e-3$}--(9,0);
 \draw (10,4)node[above]{$e-2$}--(10,0);
 \draw (11,4)node[above]{$e-1$}--(11,0);
 \node[greendot] at (4,2){};
 \node[greendot] at (10,2){};
 \node[greendot] at (6,2){};
 \draw[loosely dashed,red] (-0.5,3) -- (6.5,3)--(6.5,1)--(-0.5,1)--(-0.5,3);
\end{braid} -
\begin{braid}
 \draw (0,4)node[above]{$0$}--(0,0);
 \draw (1,4)node[above]{$1$}--(1,0);
 \draw[dots] (1.2,4)--(2.8,4);
 \draw[dots] (1.2,0)--(2.8,0);
 \draw (3,4)node[above]{$e-2$}--(3,0);
 \draw (4,4)node[above]{$e-1$}--(5,3)--(5,0);
 \draw (5,4)node[above]{$e-1$}--(4,3)--(4,0);
 \draw (6,4)node[above]{$0$}--(6,0);
 \draw (7,4)node[above]{$1$}--(7,0);
 \draw[dots] (7.2,4)--(8.8,4);
 \draw[dots] (7.2,0)--(8.8,0);
 \draw (9,4)node[above]{$e-3$}--(9,0);
 \draw (10,4)node[above]{$e-2$}--(10,0);
 \draw (11,4)node[above]{$e-1$}--(11,0);
 \node[greendot] at (10,2){};
 \node[greendot] at (6,2){};
 \node[greendot] at (4,2){};
 \draw[loosely dashed,red] (-0.5,3) -- (6.5,3)--(6.5,1)--(-0.5,1)--(-0.5,3);
\end{braid} \in \Blam[\gamma].
\end{align*}

For the fourth term of (\ref{I-problem: short: equation1}), by (\ref{dia:ii2}),
\begin{align*}
\begin{braid}
 \draw (0,4)node[above]{$0$}--(0,0);
 \draw (1,4)node[above]{$1$}--(1,0);
 \draw[dots] (1.2,4)--(2.8,4);
 \draw[dots] (1.2,0)--(2.8,0);
 \draw (3,4)node[above]{$e-2$}--(3,0);
 \draw (4,4)[densely dotted] node[above]{$e-1$}--(4,0);
 \draw (5,4)[densely dotted] node[above]{$e-1$}--(5,0);
 \draw (6,4)node[above]{$0$}--(6,0);
 \draw (7,4)node[above]{$1$}--(7,0);
 \draw[dots] (7.2,4)--(8.8,4);
 \draw[dots] (7.2,0)--(8.8,0);
 \draw (9,4)node[above]{$e-3$}--(9,0);
 \draw (10,4)node[above]{$e-2$}--(10,0);
 \draw (11,4)node[above]{$e-1$}--(11,0);
 \node[greendot] at (5,2){};
 \node[greendot] at (10,2){};
\end{braid}
& =
-\begin{braid}
 \draw (0,4)node[above]{$0$}--(0,0);
 \draw (1,4)node[above]{$1$}--(1,0);
 \draw[dots] (1.2,4)--(2.8,4);
 \draw[dots] (1.2,0)--(2.8,0);
 \draw (3,4)node[above]{$e-2$}--(3,0);
 \draw (4,4)node[above]{$e-1$}--(5,3)--(5,1)--(4,0);
 \draw (5,4)node[above]{$e-1$}--(4,3)--(4,1)--(5,0);
 \draw (6,4)node[above]{$0$}--(6,0);
 \draw (7,4)node[above]{$1$}--(7,0);
 \draw[dots] (7.2,4)--(8.8,4);
 \draw[dots] (7.2,0)--(8.8,0);
 \draw (9,4)node[above]{$e-3$}--(9,0);
 \draw (10,4)node[above]{$e-2$}--(10,0);
 \draw (11,4)node[above]{$e-1$}--(11,0);
 \node[greendot] at (4,2.7){};
 \node[greendot] at (4,1.3){};
 \node[greendot] at (5,0){};
 \node[greendot] at (10,2){};
 \draw[loosely dashed,red] (-0.5,3) -- (4.5,3)--(4.5,1)--(-0.5,1)--(-0.5,3);
\end{braid}
-
\begin{braid}
 \draw (0,4)node[above]{$0$}--(0,0);
 \draw (1,4)node[above]{$1$}--(1,0);
 \draw[dots] (1.2,4)--(2.8,4);
 \draw[dots] (1.2,0)--(2.8,0);
 \draw (3,4)node[above]{$e-2$}--(3,0);
 \draw (4,4)node[above]{$e-1$}--(5,3)--(5,0);
 \draw (5,4)node[above]{$e-1$}--(4,3)--(4,0);
 \draw (6,4)node[above]{$0$}--(6,0);
 \draw (7,4)node[above]{$1$}--(7,0);
 \draw[dots] (7.2,4)--(8.8,4);
 \draw[dots] (7.2,0)--(8.8,0);
 \draw (9,4)node[above]{$e-3$}--(9,0);
 \draw (10,4)node[above]{$e-2$}--(10,0);
 \draw (11,4)node[above]{$e-1$}--(11,0);
 \node[greendot] at (4,2){};
 \node[greendot] at (5,2){};
 \node[greendot] at (10,2){};
 \draw[loosely dashed,red] (-0.5,3) -- (5.5,3)--(5.5,1)--(-0.5,1)--(-0.5,3);
\end{braid}\\
& -
\begin{braid}
 \draw (0,4)node[above]{$0$}--(0,0);
 \draw (1,4)node[above]{$1$}--(1,0);
 \draw[dots] (1.2,4)--(2.8,4);
 \draw[dots] (1.2,0)--(2.8,0);
 \draw (3,4)node[above]{$e-2$}--(3,0);
 \draw (4,4)node[above]{$e-1$}--(4,1)--(5,0);
 \draw (5,4)node[above]{$e-1$}--(5,1)--(4,0);
 \draw (6,4)node[above]{$0$}--(6,0);
 \draw (7,4)node[above]{$1$}--(7,0);
 \draw[dots] (7.2,4)--(8.8,4);
 \draw[dots] (7.2,0)--(8.8,0);
 \draw (9,4)node[above]{$e-3$}--(9,0);
 \draw (10,4)node[above]{$e-2$}--(10,0);
 \draw (11,4)node[above]{$e-1$}--(11,0);
 \node[greendot] at (4,1.3){};
 \node[greendot] at (4,2.7){};
 \node[greendot] at (10,2){};
 \draw[loosely dashed,red] (-0.5,3) -- (4.5,3)--(4.5,1)--(-0.5,1)--(-0.5,3);
\end{braid}
-
\begin{braid}
 \draw (0,4)node[above]{$0$}--(0,0);
 \draw (1,4)node[above]{$1$}--(1,0);
 \draw[dots] (1.2,4)--(2.8,4);
 \draw[dots] (1.2,0)--(2.8,0);
 \draw (3,4)node[above]{$e-2$}--(3,0);
 \draw (4,4)node[above]{$e-1$}--(4,0);
 \draw (5,4)node[above]{$e-1$}--(5,0);
 \draw (6,4)node[above]{$0$}--(6,0);
 \draw (7,4)node[above]{$1$}--(7,0);
 \draw[dots] (7.2,4)--(8.8,4);
 \draw[dots] (7.2,0)--(8.8,0);
 \draw (9,4)node[above]{$e-3$}--(9,0);
 \draw (10,4)node[above]{$e-2$}--(10,0);
 \draw (11,4)node[above]{$e-1$}--(11,0);
 \node[greendot] at (4,2){};
 \node[greendot] at (10,2){};
\end{braid}\\
&=_\gamma
\begin{braid}
 \draw (0,4)node[above]{$0$}--(0,0);
 \draw (1,4)node[above]{$1$}--(1,0);
 \draw[dots] (1.2,4)--(2.8,4);
 \draw[dots] (1.2,0)--(2.8,0);
 \draw (3,4)node[above]{$e-2$}--(3,0);
 \draw (4,4)node[above]{$e-1$}--(4,0);
 \draw (5,4)node[above]{$e-1$}--(5,0);
 \draw (6,4)node[above]{$0$}--(6,0);
 \draw (7,4)node[above]{$1$}--(7,0);
 \draw[dots] (7.2,4)--(8.8,4);
 \draw[dots] (7.2,0)--(8.8,0);
 \draw (9,4)node[above]{$e-3$}--(9,0);
 \draw (10,4)node[above]{$e-2$}--(10,0);
 \draw (11,4)node[above]{$e-1$}--(11,0);
 \node[greendot] at (4,2){};
 \node[greendot] at (10,2){};
\end{braid} = e_\gamma y_\gamma,
\end{align*}
and for the last term of (\ref{I-problem: short: equation1}),
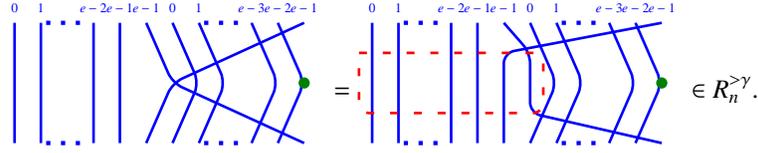
$$
\begin{braid}
 \draw (0,4)node[above]{$0$}--(0,0);
 \draw (1,4)node[above]{$1$}--(1,0);
 \draw[dots] (1.2,4)--(2.8,4);
 \draw[dots] (1.2,0)--(2.8,0);
 \draw (3,4)node[above]{$e-2$}--(3,0);
 \draw (4,4)node[above]{$e-1$}--(4,0);
 \draw (5,4)node[above]{$e-1$}--(6,2)--(11,0);
 \draw (6,4)node[above]{$0$}--(7,2)--(6,0);
 \draw (7,4)node[above]{$1$}--(8,2)--(7,0);
 \draw[dots] (7.2,4)--(8.8,4);
 \draw[dots] (7.2,0)--(8.8,0);
 \draw (9,4)node[above]{$e-3$}--(10,2)--(9,0);
 \draw (10,4)node[above]{$e-2$}--(11,2)--(10,0);
 \draw (11,4)node[above]{$e-1$}--(6,2)--(5,0);
 \node[greendot] at (11,2){};
\end{braid}
=
\begin{braid}
 \draw (0,4)node[above]{$0$}--(0,0);
 \draw (1,4)node[above]{$1$}--(1,0);
 \draw[dots] (1.2,4)--(2.8,4);
 \draw[dots] (1.2,0)--(2.8,0);
 \draw (3,4)node[above]{$e-2$}--(3,0);
 \draw (4,4)node[above]{$e-1$}--(4,0);
 \draw (5,4)node[above]{$e-1$}--(6,3)--(6,1)--(11,0);
 \draw (6,4)node[above]{$0$}--(7,2)--(6,0);
 \draw (7,4)node[above]{$1$}--(8,2)--(7,0);
 \draw[dots] (7.2,4)--(8.8,4);
 \draw[dots] (7.2,0)--(8.8,0);
 \draw (9,4)node[above]{$e-3$}--(10,2)--(9,0);
 \draw (10,4)node[above]{$e-2$}--(11,2)--(10,0);
 \draw (11,4)node[above]{$e-1$}--(5,3)--(5,0);
 \node[greendot] at (11,2){};
 \draw[loosely dashed,red] (-0.5,3) -- (6.5,3)--(6.5,1)--(-0.5,1)--(-0.5,3);
\end{braid} \in \Blam[\gamma].
$$

Combine the result above, we have
$$
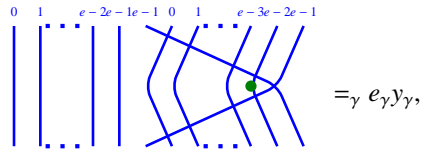
\begin{braid}
 \draw (0,4)node[above]{$0$}--(0,0);
 \draw (1,4)node[above]{$1$}--(1,0);
 \draw[dots] (1.2,4)--(2.8,4);
 \draw[dots] (1.2,0)--(2.8,0);
 \draw (3,4)node[above]{$e-2$}--(3,0);
 \draw (4,4)node[above]{$e-1$}--(4,0);
 \draw (5,4)node[above]{$e-1$}--(10,2)--(11,0);
 \draw (6,4)node[above]{$0$}--(5,2)--(6,0);
 \draw (7,4)node[above]{$1$}--(6,2)--(7,0);
 \draw[dots] (7.2,4)--(8.8,4);
 \draw[dots] (7.2,0)--(8.8,0);
 \draw (9,4)node[above]{$e-3$}--(8,2)--(9,0);
 \draw (10,4)node[above]{$e-2$}--(9,2)--(10,0);
 \draw (11,4)node[above]{$e-1$}--(10,2)--(5,0);
 \node[greendot] at (9,2){};
\end{braid}
=_\gamma e_\gamma y_\gamma,
$$
which completes the proof.

\endproof

\begin{Remark} \label{Remark: application}

The technique of applying \autoref{send} in proving \autoref{I-problem: induction1} will be used many times in the rest of the paper. Although the process is straightforward, individual details will vary from case to case, thus in order to simplify the process we will omit details in the future.

\end{Remark}

Recall $\gamma_2 = k{\cdot}e + t$ where $k$ is a nonnegative integer and $2 \leq t \leq e + 1$. Now we remove the restriction on $\gamma_2$ by applying the induction on $k$.

\begin{Proposition} \label{I-problem: induction2}
Suppose $\gamma = (\gamma_1,\gamma_2) \in \mathscr C_n^\Lambda$ with $\gamma_2 > 1$ and $\gamma_2 - \gamma_1 \equiv 1\pmod{e}$ and $\lambda = (\gamma_1,\gamma_2-1,1) \in \mathscr S_n^\Lambda$. Define $i$ to be the residue of the node at position $(1,m,1)$. Then (\ref{I-problem: non-addable: help}) holds.
\end{Proposition}

\proof We prove this Proposition by induction. As we claimed before that we can write $\gamma_2 = k{\cdot}e + t$ with $2 \leq t \leq e+1$ and we will apply induction on $k$. \autoref{I-problem: induction1} implies that for $k = 0$ the Proposition holds. Assume that for $k \leq k'$ the Proposition holds. For $k = k'$, we consider two different cases, which are $i = e-2$, $i = e-1$ and $i\neq e-2,e-1$. Recall that $i$ is the residue of the node at $(1,m,1)$ or $(2,m+1,1)$.

\textbf{Case \ref{I-problem: induction2}a:} $i \neq e-2, e-1$.

\begin{eqnarray*}
&&
\begin{braid}
 \draw(2,-0.5)node{$\underbrace{\hspace*{16mm}}_{\gamma_1}$};
 \draw(6.5,-0.5)node{$\underbrace{\hspace*{12mm}}_{\gamma_2-e}$};
 \draw(11.5,-0.5)node{$\underbrace{\hspace*{20mm}}_{e}$};
 \draw (0,4)node[above]{$0$}--(0,0);
 \draw (1,4)node[above]{$1$}--(1,0);
 \draw[dots] (1.2,4)--(2.8,4);
 \draw[dots] (1.2,0)--(2.8,0);
 \draw (3,4)node[above]{$i-1$}--(3,0);
 \draw (4,4)node[above]{$i$}--(13,2)--(14,0);
 \draw (5,4)node[above]{$e-1$}--(4,2)--(5,0);
 \draw[dots] (5.2,4)--(6.8,4);
 \draw[dots] (5.2,0)--(6.8,0);
 \draw (7,4)node[above]{$i-1$}--(6,2)--(7,0);
 \draw (8,4)node[above]{$i$}--(7,2)--(8,0);
 \draw (9,4)node[above]{$i+1$}--(8,2)--(9,0);
 \draw[dots] (9.2,4)--(10.8,4);
 \draw[dots] (9.2,0)--(10.8,0);
 \draw (11,4)[densely dotted] node[above]{$e-2$}--(10,2)--(11,0);
 \draw[dots] (11.2,4)--(12.8,4);
 \draw[dots] (11.2,0)--(12.8,0);
 \draw (13,4)node[above]{$i-1$}--(12,2)--(13,0);
 \draw (14,4)node[above]{$i$}--(13,2)--(4,0);
 \node[greendot] at (10,2){};
 \draw[->] (10,2) -- (10.5,1);
\end{braid}
\overset{(\ref{dia:y-psi com})}=
\begin{braid}
 \draw(2,-0.5)node{$\underbrace{\hspace*{16mm}}_{\gamma_1}$};
 \draw(6.5,-0.5)node{$\underbrace{\hspace*{12mm}}_{\gamma_2-e}$};
 \draw(11.5,-0.5)node{$\underbrace{\hspace*{20mm}}_{e}$};
 \draw (0,4)node[above]{$0$}--(0,0);
 \draw (1,4)node[above]{$1$}--(1,0);
 \draw[dots] (1.2,4)--(2.8,4);
 \draw[dots] (1.2,0)--(2.8,0);
 \draw (3,4)node[above]{$i-1$}--(3,0);
 \draw (4,4)node[above]{$i$}--(13,2)--(14,0);
 \draw (5,4)node[above]{$e-1$}--(4,2)--(5,0);
 \draw[dots] (5.2,4)--(6.8,4);
 \draw[dots] (5.2,0)--(6.8,0);
 \draw (7,4)node[above]{$i-1$}--(6,2)--(7,0);
 \draw (8,4)node[above]{$i$}--(7,2)--(8,0);
 \draw (9,4)node[above]{$i+1$}--(8,2)--(9,0);
 \draw[dots] (9.2,4)--(10.8,4);
 \draw[dots] (9.2,0)--(10.8,0);
 \draw (11,4)node[above]{$e-2$}--(10,2)--(11,0);
 \draw[dots] (11.2,4)--(12.8,4);
 \draw[dots] (11.2,0)--(12.8,0);
 \draw (13,4)node[above]{$i-1$}--(12,2)--(13,0);
 \draw (14,4)node[above]{$i$}--(13,2)--(4,0);
 \node[greendot] at (11,0){};
\end{braid}\\
& = &
\begin{braid}
 \draw(2,-0.5)node{$\underbrace{\hspace*{16mm}}_{\gamma_1}$};
 \draw(6.5,-0.5)node{$\underbrace{\hspace*{12mm}}_{\gamma_2-e}$};
 \draw(11.5,-0.5)node{$\underbrace{\hspace*{20mm}}_{e}$};
 \draw (0,4)node[above]{$0$}--(0,0);
 \draw (1,4)node[above]{$1$}--(1,0);
 \draw[dots] (1.2,4)--(2.8,4);
 \draw[dots] (1.2,0)--(2.8,0);
 \draw (3,4)node[above]{$i-1$}--(3,0);
 \draw (4,4)node[above]{$i$}--(8,3)--(13,2)--(14,1)--(14,0);
 \draw (5,4)node[above]{$e-1$}--(4,3)--(4,1)--(5,0);
 \draw[dots] (5.2,4)--(6.8,4);
 \draw[dots] (5.2,0)--(6.8,0);
 \draw (7,4)node[above]{$i-1$}--(6,3)--(6,1)--(7,0);
 \draw (8,4)node[above]{$i$}--(7,3)--(7,1)--(8,0);
 \draw (9,4)node[above]{$i+1$}--(9,3)--(8,2)--(9,1)--(9,0);
 \draw[dots] (9.2,4)--(10.8,4);
 \draw[dots] (9.2,0)--(10.8,0);
 \draw (11,4)node[above]{$e-2$}--(11,3)--(10,2)--(11,1)--(11,0);
 \draw[dots] (11.2,4)--(12.8,4);
 \draw[dots] (11.2,0)--(12.8,0);
 \draw (13,4)node[above]{$i-1$}--(13,3)--(12,2)--(13,1)--(13,0);
 \draw (14,4)node[above]{$i$}--(14,3)--(13,2)--(8,1)--(4,0);
 \node[greendot] at (11,0){};
 \draw[loosely dashed,red] (7.5,3) -- (14.5,3)--(14.5,1)--(7.5,1)--(7.5,3);
\end{braid}\hspace*{11mm} \text{ by \autoref{I-problem: short}}\\
& = &
\begin{braid}
 \draw(2,-0.5)node{$\underbrace{\hspace*{16mm}}_{\gamma_1}$};
 \draw(6.5,-0.5)node{$\underbrace{\hspace*{12mm}}_{\gamma_2-e}$};
 \draw(11.5,-0.5)node{$\underbrace{\hspace*{20mm}}_{e}$};
 \draw (0,4)node[above]{$0$}--(0,0);
 \draw (1,4)node[above]{$1$}--(1,0);
 \draw[dots] (1.2,4)--(2.8,4);
 \draw[dots] (1.2,0)--(2.8,0);
 \draw (3,4)node[above]{$i-1$}--(3,0);
 \draw (4,4)[densely dotted] node[above]{$i$}--(8,2)--(4,0);
 \draw (5,4)node[above]{$e-1$}--(4,2)--(5,0);
 \draw[dots] (5.2,4)--(6.8,4);
 \draw[dots] (5.2,0)--(6.8,0);
 \draw (7,4)node[above]{$i-1$}--(6,2)--(7,0);
 \draw (8,4)[densely dotted] node[above]{$i$}--(7,2)--(8,0);
 \draw (9,4)node[above]{$i+1$}--(9,0);
 \draw[dots] (9.2,4)--(10.8,4);
 \draw[dots] (9.2,0)--(10.8,0);
 \draw (11,4)node[above]{$e-2$}--(11,0);
 \draw[dots] (11.2,4)--(12.8,4);
 \draw[dots] (11.2,0)--(12.8,0);
 \draw (13,4)node[above]{$i-1$}--(13,0);
 \draw (14,4)node[above]{$i$}--(14,0);
 \node[greendot] at (11,0){};
 \node[greendot] at (14,2){};
\end{braid}
-
\begin{braid}
 \draw(2,-0.5)node{$\underbrace{\hspace*{16mm}}_{\gamma_1}$};
 \draw(6.5,-0.5)node{$\underbrace{\hspace*{12mm}}_{\gamma_2-e}$};
 \draw(11.5,-0.5)node{$\underbrace{\hspace*{20mm}}_{e}$};
 \draw (0,4)node[above]{$0$}--(0,0);
 \draw (1,4)node[above]{$1$}--(1,0);
 \draw[dots] (1.2,4)--(2.8,4);
 \draw[dots] (1.2,0)--(2.8,0);
 \draw (3,4)node[above]{$i-1$}--(3,0);
 \draw (4,4)[densely dotted] node[above]{$i$}--(8,2)--(4,0);
 \draw (5,4)node[above]{$e-1$}--(4,2)--(5,0);
 \draw[dots] (5.2,4)--(6.8,4);
 \draw[dots] (5.2,0)--(6.8,0);
 \draw (7,4)node[above]{$i-1$}--(6,2)--(7,0);
 \draw (8,4)[densely dotted] node[above]{$i$}--(7,2)--(8,0);
 \draw (9,4)node[above]{$i+1$}--(9,0);
 \draw[dots] (9.2,4)--(10.8,4);
 \draw[dots] (9.2,0)--(10.8,0);
 \draw (11,4)node[above]{$e-2$}--(11,0);
 \draw[dots] (11.2,4)--(12.8,4);
 \draw[dots] (11.2,0)--(12.8,0);
 \draw (13,4)node[above]{$i-1$}--(13,0);
 \draw (14,4)node[above]{$i$}--(14,0);
 \node[greendot] at (11,0){};
 \node[greendot] at (13,2){};
\end{braid}\\
&&
-
\begin{braid}
 \draw(2,-0.5)node{$\underbrace{\hspace*{16mm}}_{\gamma_1}$};
 \draw(6.5,-0.5)node{$\underbrace{\hspace*{12mm}}_{\gamma_2-e}$};
 \draw(11.5,-0.5)node{$\underbrace{\hspace*{20mm}}_{e}$};
 \draw (0,4)node[above]{$0$}--(0,0);
 \draw (1,4)node[above]{$1$}--(1,0);
 \draw[dots] (1.2,4)--(2.8,4);
 \draw[dots] (1.2,0)--(2.8,0);
 \draw (3,4)node[above]{$i-1$}--(3,0);
 \draw (4,4)[densely dotted] node[above]{$i$}--(8,2)--(4,0);
 \draw (5,4)node[above]{$e-1$}--(4,2)--(5,0);
 \draw[dots] (5.2,4)--(6.8,4);
 \draw[dots] (5.2,0)--(6.8,0);
 \draw (7,4)node[above]{$i-1$}--(6,2)--(7,0);
 \draw (8,4)[densely dotted] node[above]{$i$}--(7,2)--(8,0);
 \draw (9,4)node[above]{$i+1$}--(9,0);
 \draw[dots] (9.2,4)--(10.8,4);
 \draw[dots] (9.2,0)--(10.8,0);
 \draw (11,4)node[above]{$e-2$}--(11,0);
 \draw[dots] (11.2,4)--(12.8,4);
 \draw[dots] (11.2,0)--(12.8,0);
 \draw (13,4)node[above]{$i-1$}--(13,0);
 \draw (14,4)node[above]{$i$}--(14,0);
 \node[greendot] at (11,0){};
 \node[greendot] at (9,2){};
\end{braid}
+
\begin{braid}
 \draw(2,-0.5)node{$\underbrace{\hspace*{16mm}}_{\gamma_1}$};
 \draw(6.5,-0.5)node{$\underbrace{\hspace*{12mm}}_{\gamma_2-e}$};
 \draw(11.5,-0.5)node{$\underbrace{\hspace*{20mm}}_{e}$};
 \draw (0,4)node[above]{$0$}--(0,0);
 \draw (1,4)node[above]{$1$}--(1,0);
 \draw[dots] (1.2,4)--(2.8,4);
 \draw[dots] (1.2,0)--(2.8,0);
 \draw (3,4)node[above]{$i-1$}--(3,0);
 \draw (4,4)[densely dotted] node[above]{$i$}--(8,2)--(4,0);
 \draw (5,4)node[above]{$e-1$}--(4,2)--(5,0);
 \draw[dots] (5.2,4)--(6.8,4);
 \draw[dots] (5.2,0)--(6.8,0);
 \draw (7,4)node[above]{$i-1$}--(6,2)--(7,0);
 \draw (8,4)[densely dotted] node[above]{$i$}--(7,2)--(8,0);
 \draw (9,4)node[above]{$i+1$}--(9,0);
 \draw[dots] (9.2,4)--(10.8,4);
 \draw[dots] (9.2,0)--(10.8,0);
 \draw (11,4)node[above]{$e-2$}--(11,0);
 \draw[dots] (11.2,4)--(12.8,4);
 \draw[dots] (11.2,0)--(12.8,0);
 \draw (13,4)node[above]{$i-1$}--(13,0);
 \draw (14,4)node[above]{$i$}--(14,0);
 \node[greendot] at (11,0){};
 \node[greendot] at (8,2){};
 \draw[->] (8,2) -- (7,2.5);
\end{braid}\\
&&
+
\begin{braid}
 \draw(2,-0.5)node{$\underbrace{\hspace*{16mm}}_{\gamma_1}$};
 \draw(6.5,-0.5)node{$\underbrace{\hspace*{12mm}}_{\gamma_2-e}$};
 \draw(11.5,-0.5)node{$\underbrace{\hspace*{20mm}}_{e}$};
 \draw (0,4)node[above]{$0$}--(0,0);
 \draw (1,4)node[above]{$1$}--(1,0);
 \draw[dots] (1.2,4)--(2.8,4);
 \draw[dots] (1.2,0)--(2.8,0);
 \draw (3,4)node[above]{$i-1$}--(3,0);
 \draw (4,4)node[above]{$i$}--(9,2)--(14,0);
 \draw (5,4)node[above]{$e-1$}--(4,2)--(5,0);
 \draw[dots] (5.2,4)--(6.8,4);
 \draw[dots] (5.2,0)--(6.8,0);
 \draw (7,4)node[above]{$i-1$}--(6,2)--(7,0);
 \draw (8,4)node[above]{$i$}--(7,2)--(8,0);
 \draw (9,4)node[above]{$i+1$}--(10,2)--(9,0);
 \draw[dots] (9.2,4)--(10.8,4);
 \draw[dots] (9.2,0)--(10.8,0);
 \draw (11,4)node[above]{$e-2$}--(12,2)--(11,0);
 \draw[dots] (11.2,4)--(12.8,4);
 \draw[dots] (11.2,0)--(12.8,0);
 \draw (13,4)node[above]{$i-1$}--(14,2)--(13,0);
 \draw (14,4)node[above]{$i$}--(9,2)--(4,0);
 \node[greendot] at (11,0){};
\end{braid}\\
& \overset{\substack{(\ref{dia:psipsi})\\(\ref{dia:ii3})}}= &
\begin{braid}
 \draw(2,-0.5)node{$\underbrace{\hspace*{16mm}}_{\gamma_1}$};
 \draw(6.5,-0.5)node{$\underbrace{\hspace*{12mm}}_{\gamma_2-e}$};
 \draw(11.5,-0.5)node{$\underbrace{\hspace*{20mm}}_{e}$};
 \draw (0,4)node[above]{$0$}--(0,0);
 \draw (1,4)node[above]{$1$}--(1,0);
 \draw[dots] (1.2,4)--(2.8,4);
 \draw[dots] (1.2,0)--(2.8,0);
 \draw (3,4)node[above]{$i-1$}--(3,0);
 \draw (4,4)node[above]{$i$}--(7,2)--(8,0);
 \draw (5,4)node[above]{$e-1$}--(4,2)--(5,0);
 \draw[dots] (5.2,4)--(6.8,4);
 \draw[dots] (5.2,0)--(6.8,0);
 \draw (7,4)node[above]{$i-1$}--(6,2)--(7,0);
 \draw (8,4)node[above]{$i$}--(7,2)--(4,0);
 \draw (9,4)node[above]{$i+1$}--(9,0);
 \draw[dots] (9.2,4)--(10.8,4);
 \draw[dots] (9.2,0)--(10.8,0);
 \draw (11,4)node[above]{$e-2$}--(11,0);
 \draw[dots] (11.2,4)--(12.8,4);
 \draw[dots] (11.2,0)--(12.8,0);
 \draw (13,4)node[above]{$i-1$}--(13,0);
 \draw (14,4)node[above]{$i$}--(14,0);
 \node[greendot] at (11,0){};
\end{braid}
+
\begin{braid}
 \draw(2,-0.5)node{$\underbrace{\hspace*{16mm}}_{\gamma_1}$};
 \draw(6.5,-0.5)node{$\underbrace{\hspace*{12mm}}_{\gamma_2-e}$};
 \draw(11.5,-0.5)node{$\underbrace{\hspace*{20mm}}_{e}$};
 \draw (0,4)node[above]{$0$}--(0,0);
 \draw (1,4)node[above]{$1$}--(1,0);
 \draw[dots] (1.2,4)--(2.8,4);
 \draw[dots] (1.2,0)--(2.8,0);
 \draw (3,4)node[above]{$i-1$}--(3,0);
 \draw (4,4)node[above]{$i$}--(9,2)--(14,0);
 \draw (5,4)node[above]{$e-1$}--(4,2)--(5,0);
 \draw[dots] (5.2,4)--(6.8,4);
 \draw[dots] (5.2,0)--(6.8,0);
 \draw (7,4)node[above]{$i-1$}--(6,2)--(7,0);
 \draw (8,4)node[above]{$i$}--(7,2)--(8,0);
 \draw (9,4)node[above]{$i+1$}--(10,2)--(9,0);
 \draw[dots] (9.2,4)--(10.8,4);
 \draw[dots] (9.2,0)--(10.8,0);
 \draw (11,4)[densely dotted] node[above]{$e-2$}--(12,2)--(11,0);
 \draw[dots] (11.2,4)--(12.8,4);
 \draw[dots] (11.2,0)--(12.8,0);
 \draw (13,4)node[above]{$i-1$}--(14,2)--(13,0);
 \draw (14,4)node[above]{$i$}--(9,2)--(4,0);
 \node[greendot] at (11,0){};
 \draw[->] (11,0) -- (12,2);
\end{braid}\\
& \overset{(\ref{dia:y-psi com})}= &
\begin{braid}
 \draw(2,-0.5)node{$\underbrace{\hspace*{16mm}}_{\gamma_1}$};
 \draw(6.5,-0.5)node{$\underbrace{\hspace*{12mm}}_{\gamma_2-e}$};
 \draw(11.5,-0.5)node{$\underbrace{\hspace*{20mm}}_{e}$};
 \draw (0,4)node[above]{$0$}--(0,0);
 \draw (1,4)node[above]{$1$}--(1,0);
 \draw[dots] (1.2,4)--(2.8,4);
 \draw[dots] (1.2,0)--(2.8,0);
 \draw (3,4)node[above]{$i-1$}--(3,0);
 \draw (4,4)node[above]{$i$}--(7,2)--(8,0);
 \draw (5,4)node[above]{$e-1$}--(4,2)--(5,0);
 \draw[dots] (5.2,4)--(6.8,4);
 \draw[dots] (5.2,0)--(6.8,0);
 \draw (7,4)node[above]{$i-1$}--(6,2)--(7,0);
 \draw (8,4)node[above]{$i$}--(7,2)--(4,0);
 \draw (9,4)node[above]{$i+1$}--(9,0);
 \draw[dots] (9.2,4)--(10.8,4);
 \draw[dots] (9.2,0)--(10.8,0);
 \draw (11,4)node[above]{$e-2$}--(11,0);
 \draw[dots] (11.2,4)--(12.8,4);
 \draw[dots] (11.2,0)--(12.8,0);
 \draw (13,4)node[above]{$i-1$}--(13,0);
 \draw (14,4)node[above]{$i$}--(14,0);
 \node[greendot] at (11,2){};
\end{braid}
+
\begin{braid}
 \draw(2,-0.5)node{$\underbrace{\hspace*{16mm}}_{\gamma_1}$};
 \draw(6.5,-0.5)node{$\underbrace{\hspace*{12mm}}_{\gamma_2-e}$};
 \draw(11.5,-0.5)node{$\underbrace{\hspace*{20mm}}_{e}$};
 \draw (0,4)node[above]{$0$}--(0,0);
 \draw (1,4)node[above]{$1$}--(1,0);
 \draw[dots] (1.2,4)--(2.8,4);
 \draw[dots] (1.2,0)--(2.8,0);
 \draw (3,4)node[above]{$i-1$}--(3,0);
 \draw (4,4)node[above]{$i$}--(4,3)--(8,2)--(9,1)--(14,0);
 \draw (5,4)node[above]{$e-1$}--(5,3)--(4,2)--(5,1)--(5,0);
 \draw[dots] (5.2,4)--(6.8,4);
 \draw[dots] (5.2,0)--(6.8,0);
 \draw (7,4)node[above]{$i-1$}--(7,3)--(6,2)--(7,1)--(7,0);
 \draw (8,4)node[above]{$i$}--(8,3)--(7,2)--(8,1)--(8,0);
 \draw (9,4)node[above]{$i+1$}--(10,3)--(10,1)--(9,0);
 \draw[dots] (9.2,4)--(10.8,4);
 \draw[dots] (9.2,0)--(10.8,0);
 \draw (11,4)node[above]{$e-2$}--(12,3)--(12,1)--(11,0);
 \draw[dots] (11.2,4)--(12.8,4);
 \draw[dots] (11.2,0)--(12.8,0);
 \draw (13,4)node[above]{$i-1$}--(14,3)--(14,1)--(13,0);
 \draw (14,4)node[above]{$i$}--(9,3)--(8,2)--(4,1)--(4,0);
 \node[greendot] at (12,2){};
 \draw[loosely dashed,red] (-0.5,3) -- (9.5,3)--(9.5,1)--(-0.5,1)--(-0.5,3);
\end{braid}.
\end{eqnarray*}

Then by induction and \autoref{send}, we have
$$
\begin{braid}
 \draw(2,-0.5)node{$\underbrace{\hspace*{16mm}}_{\gamma_1}$};
 \draw(6.5,-0.5)node{$\underbrace{\hspace*{12mm}}_{\gamma_2-e}$};
 \draw(11.5,-0.5)node{$\underbrace{\hspace*{20mm}}_{e}$};
 \draw (0,4)node[above]{$0$}--(0,0);
 \draw (1,4)node[above]{$1$}--(1,0);
 \draw[dots] (1.2,4)--(2.8,4);
 \draw[dots] (1.2,0)--(2.8,0);
 \draw (3,4)node[above]{$i-1$}--(3,0);
 \draw (4,4)node[above]{$i$}--(4,3)--(8,2)--(9,1)--(14,0);
 \draw (5,4)node[above]{$e-1$}--(5,3)--(4,2)--(5,1)--(5,0);
 \draw[dots] (5.2,4)--(6.8,4);
 \draw[dots] (5.2,0)--(6.8,0);
 \draw (7,4)node[above]{$i-1$}--(7,3)--(6,2)--(7,1)--(7,0);
 \draw (8,4)node[above]{$i$}--(8,3)--(7,2)--(8,1)--(8,0);
 \draw (9,4)node[above]{$i+1$}--(10,3)--(10,1)--(9,0);
 \draw[dots] (9.2,4)--(10.8,4);
 \draw[dots] (9.2,0)--(10.8,0);
 \draw (11,4)node[above]{$e-2$}--(12,3)--(12,1)--(11,0);
 \draw[dots] (11.2,4)--(12.8,4);
 \draw[dots] (11.2,0)--(12.8,0);
 \draw (13,4)node[above]{$i-1$}--(14,3)--(14,1)--(13,0);
 \draw (14,4)node[above]{$i$}--(9,3)--(8,2)--(4,1)--(4,0);
 \node[greendot] at (12,2){};
 \draw[loosely dashed,red] (-0.5,3) -- (9.5,3)--(9.5,1)--(-0.5,1)--(-0.5,3);
\end{braid}
=_\gamma
\begin{braid}
 \draw(2,-0.5)node{$\underbrace{\hspace*{16mm}}_{\gamma_1}$};
 \draw(6.5,-0.5)node{$\underbrace{\hspace*{12mm}}_{\gamma_2-e}$};
 \draw(11.5,-0.5)node{$\underbrace{\hspace*{20mm}}_{e}$};
 \draw (0,4)node[above]{$0$}--(0,0);
 \draw (1,4)node[above]{$1$}--(1,0);
 \draw[dots] (1.2,4)--(2.8,4);
 \draw[dots] (1.2,0)--(2.8,0);
 \draw (3,4)node[above]{$i-1$}--(3,0);
 \draw (4,4)node[above]{$i$}--(4,0);
 \draw (5,4)node[above]{$e-1$}--(5,0);
 \draw[dots] (5.2,4)--(6.8,4);
 \draw[dots] (5.2,0)--(6.8,0);
 \draw (7,4)node[above]{$i-1$}--(7,0);
 \draw (8,4)node[above]{$i$}--(8,0);
 \draw (9,4)node[above]{$i+1$}--(10,3)--(10,1)--(9,0);
 \draw[dots] (9.2,4)--(10.8,4);
 \draw[dots] (9.2,0)--(10.8,0);
 \draw (11,4)node[above]{$e-2$}--(12,3)--(12,1)--(11,0);
 \draw[dots] (11.2,4)--(12.8,4);
 \draw[dots] (11.2,0)--(12.8,0);
 \draw (13,4)node[above]{$i-1$}--(14,3)--(14,1)--(13,0);
 \draw (14,4)node[above]{$i$}--(9,3)--(9,1)--(14,0);
 \node[greendot] at (12,2){};
\end{braid} = 0,
$$
which implies that
$$
\begin{braid}
 \draw(2,-0.5)node{$\underbrace{\hspace*{16mm}}_{\gamma_1}$};
 \draw(6.5,-0.5)node{$\underbrace{\hspace*{12mm}}_{\gamma_2-e}$};
 \draw(11.5,-0.5)node{$\underbrace{\hspace*{20mm}}_{e}$};
 \draw (0,4)node[above]{$0$}--(0,0);
 \draw (1,4)node[above]{$1$}--(1,0);
 \draw[dots] (1.2,4)--(2.8,4);
 \draw[dots] (1.2,0)--(2.8,0);
 \draw (3,4)node[above]{$i-1$}--(3,0);
 \draw (4,4)node[above]{$i$}--(9,2)--(14,0);
 \draw (5,4)node[above]{$e-1$}--(4,2)--(5,0);
 \draw[dots] (5.2,4)--(6.8,4);
 \draw[dots] (5.2,0)--(6.8,0);
 \draw (7,4)node[above]{$i-1$}--(6,2)--(7,0);
 \draw (8,4)node[above]{$i$}--(7,2)--(8,0);
 \draw (9,4)node[above]{$i+1$}--(10,2)--(9,0);
 \draw[dots] (9.2,4)--(10.8,4);
 \draw[dots] (9.2,0)--(10.8,0);
 \draw (11,4)node[above]{$e-2$}--(12,2)--(11,0);
 \draw[dots] (11.2,4)--(12.8,4);
 \draw[dots] (11.2,0)--(12.8,0);
 \draw (13,4)node[above]{$i-1$}--(14,2)--(13,0);
 \draw (14,4)node[above]{$i$}--(9,2)--(4,0);
 \node[greendot] at (12,2){};
\end{braid}
\in \Blam[\gamma].
$$

Hence by induction and \autoref{send}
$$
\begin{braid}
 \draw(2,-0.5)node{$\underbrace{\hspace*{16mm}}_{\gamma_1}$};
 \draw(6.5,-0.5)node{$\underbrace{\hspace*{12mm}}_{\gamma_2-e}$};
 \draw(11.5,-0.5)node{$\underbrace{\hspace*{20mm}}_{e}$};
 \draw (0,4)node[above]{$0$}--(0,0);
 \draw (1,4)node[above]{$1$}--(1,0);
 \draw[dots] (1.2,4)--(2.8,4);
 \draw[dots] (1.2,0)--(2.8,0);
 \draw (3,4)node[above]{$i-1$}--(3,0);
 \draw (4,4)node[above]{$i$}--(13,2)--(14,0);
 \draw (5,4)node[above]{$e-1$}--(4,2)--(5,0);
 \draw[dots] (5.2,4)--(6.8,4);
 \draw[dots] (5.2,0)--(6.8,0);
 \draw (7,4)node[above]{$i-1$}--(6,2)--(7,0);
 \draw (8,4)node[above]{$i$}--(7,2)--(8,0);
 \draw (9,4)node[above]{$i+1$}--(8,2)--(9,0);
 \draw[dots] (9.2,4)--(10.8,4);
 \draw[dots] (9.2,0)--(10.8,0);
 \draw (11,4)node[above]{$e-2$}--(10,2)--(11,0);
 \draw[dots] (11.2,4)--(12.8,4);
 \draw[dots] (11.2,0)--(12.8,0);
 \draw (13,4)node[above]{$i-1$}--(12,2)--(13,0);
 \draw (14,4)node[above]{$i$}--(13,2)--(4,0);
 \node[greendot] at (10,2){};
\end{braid}
=_\gamma
\begin{braid}
 \draw(2,-0.5)node{$\underbrace{\hspace*{16mm}}_{\gamma_1}$};
 \draw(6.5,-0.5)node{$\underbrace{\hspace*{12mm}}_{\gamma_2-e}$};
 \draw(11.5,-0.5)node{$\underbrace{\hspace*{20mm}}_{e}$};
 \draw (0,4)node[above]{$0$}--(0,0);
 \draw (1,4)node[above]{$1$}--(1,0);
 \draw[dots] (1.2,4)--(2.8,4);
 \draw[dots] (1.2,0)--(2.8,0);
 \draw (3,4)node[above]{$i-1$}--(3,0);
 \draw (4,4)node[above]{$i$}--(7,2)--(8,0);
 \draw (5,4)node[above]{$e-1$}--(4,2)--(5,0);
 \draw[dots] (5.2,4)--(6.8,4);
 \draw[dots] (5.2,0)--(6.8,0);
 \draw (7,4)node[above]{$i-1$}--(6,2)--(7,0);
 \draw (8,4)node[above]{$i$}--(7,2)--(4,0);
 \draw (9,4)node[above]{$i+1$}--(9,0);
 \draw[dots] (9.2,4)--(10.8,4);
 \draw[dots] (9.2,0)--(10.8,0);
 \draw (11,4)node[above]{$e-2$}--(11,0);
 \draw[dots] (11.2,4)--(12.8,4);
 \draw[dots] (11.2,0)--(12.8,0);
 \draw (13,4)node[above]{$i-1$}--(13,0);
 \draw (14,4)node[above]{$i$}--(14,0);
 \node[greendot] at (11,2){};
\end{braid}
=_\gamma e_\gamma y_\gamma.
$$

\textbf{Case \ref{I-problem: induction2}b:} $i = e-2$.

\begin{eqnarray}
&&\begin{braid}
 \draw(2,-0.5)node{$\underbrace{\hspace*{16mm}}_{\gamma_1}$};
 \draw(6.5,-0.5)node{$\underbrace{\hspace*{12mm}}_{\gamma_2-e}$};
 \draw(10.5,-0.5)node{$\underbrace{\hspace*{12mm}}_{e}$};
 \draw (0,4)node[above]{$0$}--(0,0);
 \draw (1,4)node[above]{$1$}--(1,0);
 \draw[dots] (1.2,4)--(2.8,4);
 \draw[dots] (1.2,0)--(2.8,0);
 \draw (3,4)node[above]{$e-3$}--(3,0);
 \draw (4,4)node[above]{$e-2$}--(11,2)--(12,0);
 \draw (5,4)node[above]{$e-1$}--(4,2)--(5,0);
 \draw[dots] (5.2,4)--(6.8,4);
 \draw[dots] (5.2,0)--(6.8,0);
 \draw (7,4)node[above]{$e-3$}--(6,2)--(7,0);
 \draw (8,4)node[above]{$e-2$}--(7,2)--(8,0);
 \draw (9,4)node[above]{$e-1$}--(8,2)--(9,0);
 \draw[dots] (9.2,4)--(10.8,4);
 \draw[dots] (9.2,0)--(10.8,0);
 \draw (11,4)node[above]{$e-3$}--(10,2)--(11,0);
 \draw (12,4)node[above]{$e-2$}--(11,2)--(4,0);
 \node[greendot] at (7,2){};
\end{braid}
 =
\begin{braid}
 \draw(2,-0.5)node{$\underbrace{\hspace*{16mm}}_{\gamma_1}$};
 \draw(6.5,-0.5)node{$\underbrace{\hspace*{12mm}}_{\gamma_2-e}$};
 \draw(10.5,-0.5)node{$\underbrace{\hspace*{12mm}}_{e}$};
 \draw (0,4)node[above]{$0$}--(0,0);
 \draw (1,4)node[above]{$1$}--(1,0);
 \draw[dots] (1.2,4)--(2.8,4);
 \draw[dots] (1.2,0)--(2.8,0);
 \draw (3,4)node[above]{$e-3$}--(3,0);
 \draw (4,4)node[above]{$e-2$}--(8,3)--(11,2)--(12,1)--(12,0);
 \draw (5,4)node[above]{$e-1$}--(4,3)--(4,1)--(5,0);
 \draw[dots] (5.2,4)--(6.8,4);
 \draw[dots] (5.2,0)--(6.8,0);
 \draw (7,4)node[above]{$e-3$}--(6,3)--(6,1)--(7,0);
 \draw (8,4)node[above]{$e-2$}--(7,3)--(7,1)--(8,0);
 \draw (9,4)node[above]{$e-1$}--(9,3)--(8,2)--(9,1)--(9,0);
 \draw[dots] (9.2,4)--(10.8,4);
 \draw[dots] (9.2,0)--(10.8,0);
 \draw (11,4)node[above]{$e-3$}--(11,3)--(10,2)--(11,1)--(11,0);
 \draw (12,4)node[above]{$e-2$}--(12,3)--(11,2)--(8,1)--(4,0);
 \node[greendot] at (7,2){};
 \draw[loosely dashed,red] (7.5,3) -- (12.5,3)--(12.5,1)--(7.5,1)--(7.5,3);
\end{braid} \notag \hspace*{11mm}\text{by \autoref{I-problem: short}}\\
& = &
\begin{braid}
 \draw(2,-0.5)node{$\underbrace{\hspace*{16mm}}_{\gamma_1}$};
 \draw(6.5,-0.5)node{$\underbrace{\hspace*{12mm}}_{\gamma_2-e}$};
 \draw(10.5,-0.5)node{$\underbrace{\hspace*{12mm}}_{e}$};
 \draw (0,4)node[above]{$0$}--(0,0);
 \draw (1,4)node[above]{$1$}--(1,0);
 \draw[dots] (1.2,4)--(2.8,4);
 \draw[dots] (1.2,0)--(2.8,0);
 \draw (3,4)node[above]{$e-3$}--(3,0);
 \draw (4,4)node[above]{$e-2$}--(8,2)--(4,0);
 \draw (5,4)node[above]{$e-1$}--(4,2)--(5,0);
 \draw[dots] (5.2,4)--(6.8,4);
 \draw[dots] (5.2,0)--(6.8,0);
 \draw (7,4)node[above]{$e-3$}--(6,2)--(7,0);
 \draw (8,4)[densely dotted] node[above]{$e-2$}--(7,2)--(8,0);
 \draw (9,4)node[above]{$e-1$}--(9,0);
 \draw[dots] (9.2,4)--(10.8,4);
 \draw[dots] (9.2,0)--(10.8,0);
 \draw (11,4)node[above]{$e-3$}--(11,0);
 \draw (12,4)node[above]{$e-2$}--(12,0);
 \node[greendot] at (7,2){};
 \node[greendot] at (12,2){};
 \draw[->] (7,2) -- (7.5,3);
\end{braid}
-
\begin{braid}
 \draw(2,-0.5)node{$\underbrace{\hspace*{16mm}}_{\gamma_1}$};
 \draw(6.5,-0.5)node{$\underbrace{\hspace*{12mm}}_{\gamma_2-e}$};
 \draw(10.5,-0.5)node{$\underbrace{\hspace*{12mm}}_{e}$};
 \draw (0,4)node[above]{$0$}--(0,0);
 \draw (1,4)node[above]{$1$}--(1,0);
 \draw[dots] (1.2,4)--(2.8,4);
 \draw[dots] (1.2,0)--(2.8,0);
 \draw (3,4)node[above]{$e-3$}--(3,0);
 \draw (4,4)node[above]{$e-2$}--(8,2)--(4,0);
 \draw (5,4)node[above]{$e-1$}--(4,2)--(5,0);
 \draw[dots] (5.2,4)--(6.8,4);
 \draw[dots] (5.2,0)--(6.8,0);
 \draw (7,4)node[above]{$e-3$}--(6,2)--(7,0);
 \draw (8,4)[densely dotted] node[above]{$e-2$}--(7,2)--(8,0);
 \draw (9,4)node[above]{$e-1$}--(9,0);
 \draw[dots] (9.2,4)--(10.8,4);
 \draw[dots] (9.2,0)--(10.8,0);
 \draw (11,4)node[above]{$e-3$}--(11,0);
 \draw (12,4)node[above]{$e-2$}--(12,0);
 \node[greendot] at (7,2){};
 \node[greendot] at (11,2){};
 \draw[->] (7,2) -- (7.5,3);
\end{braid}\notag \\
&&-
\begin{braid}
 \draw(2,-0.5)node{$\underbrace{\hspace*{16mm}}_{\gamma_1}$};
 \draw(6.5,-0.5)node{$\underbrace{\hspace*{12mm}}_{\gamma_2-e}$};
 \draw(10.5,-0.5)node{$\underbrace{\hspace*{12mm}}_{e}$};
 \draw (0,4)node[above]{$0$}--(0,0);
 \draw (1,4)node[above]{$1$}--(1,0);
 \draw[dots] (1.2,4)--(2.8,4);
 \draw[dots] (1.2,0)--(2.8,0);
 \draw (3,4)node[above]{$e-3$}--(3,0);
 \draw (4,4)node[above]{$e-2$}--(8,2)--(4,0);
 \draw (5,4)node[above]{$e-1$}--(4,2)--(5,0);
 \draw[dots] (5.2,4)--(6.8,4);
 \draw[dots] (5.2,0)--(6.8,0);
 \draw (7,4)node[above]{$e-3$}--(6,2)--(7,0);
 \draw (8,4)[densely dotted] node[above]{$e-2$}--(7,2)--(8,0);
 \draw (9,4)node[above]{$e-1$}--(9,0);
 \draw[dots] (9.2,4)--(10.8,4);
 \draw[dots] (9.2,0)--(10.8,0);
 \draw (11,4)node[above]{$e-3$}--(11,0);
 \draw (12,4)node[above]{$e-2$}--(12,0);
 \node[greendot] at (7,2){};
 \node[greendot] at (9,2){};
 \draw[->] (7,2) -- (7.5,3);
\end{braid}
+
\begin{braid}
 \draw(2,-0.5)node{$\underbrace{\hspace*{16mm}}_{\gamma_1}$};
 \draw(6.5,-0.5)node{$\underbrace{\hspace*{12mm}}_{\gamma_2-e}$};
 \draw(10.5,-0.5)node{$\underbrace{\hspace*{12mm}}_{e}$};
 \draw (0,4)node[above]{$0$}--(0,0);
 \draw (1,4)node[above]{$1$}--(1,0);
 \draw[dots] (1.2,4)--(2.8,4);
 \draw[dots] (1.2,0)--(2.8,0);
 \draw (3,4)node[above]{$e-3$}--(3,0);
 \draw (4,4)[densely dotted] node[above]{$e-2$}--(8,2)--(4,0);
 \draw (5,4)node[above]{$e-1$}--(4,2)--(5,0);
 \draw[dots] (5.2,4)--(6.8,4);
 \draw[dots] (5.2,0)--(6.8,0);
 \draw (7,4)node[above]{$e-3$}--(6,2)--(7,0);
 \draw (8,4)[densely dotted] node[above]{$e-2$}--(7,2)--(8,0);
 \draw (9,4)node[above]{$e-1$}--(9,0);
 \draw[dots] (9.2,4)--(10.8,4);
 \draw[dots] (9.2,0)--(10.8,0);
 \draw (11,4)node[above]{$e-3$}--(11,0);
 \draw (12,4)node[above]{$e-2$}--(12,0);
 \node[greendot] at (7,2){};
 \node[greendot] at (8,2){};
 \draw[->] (7,2) -- (7.5,3);
 \draw[->] (8,2) -- (6.5,2.75);
\end{braid}\notag \\
&& +
\begin{braid}
 \draw(2,-0.5)node{$\underbrace{\hspace*{16mm}}_{\gamma_1}$};
 \draw(6.5,-0.5)node{$\underbrace{\hspace*{12mm}}_{\gamma_2-e}$};
 \draw(10.5,-0.5)node{$\underbrace{\hspace*{12mm}}_{e}$};
 \draw (0,4)node[above]{$0$}--(0,0);
 \draw (1,4)node[above]{$1$}--(1,0);
 \draw[dots] (1.2,4)--(2.8,4);
 \draw[dots] (1.2,0)--(2.8,0);
 \draw (3,4)node[above]{$e-3$}--(3,0);
 \draw (4,4)node[above]{$e-2$}--(9,2)--(12,0);
 \draw (5,4)node[above]{$e-1$}--(4,2)--(5,0);
 \draw[dots] (5.2,4)--(6.8,4);
 \draw[dots] (5.2,0)--(6.8,0);
 \draw (7,4)node[above]{$e-3$}--(6,2)--(7,0);
 \draw (8,4)node[above]{$e-2$}--(7,2)--(8,0);
 \draw (9,4)node[above]{$e-1$}--(10,2)--(9,0);
 \draw[dots] (9.2,4)--(10.8,4);
 \draw[dots] (9.2,0)--(10.8,0);
 \draw (11,4)node[above]{$e-3$}--(12,2)--(11,0);
 \draw (12,4)node[above]{$e-2$}--(9,2)--(4,0);
 \node[greendot] at (7,2){};
\end{braid}\notag \\
& \overset{\substack{(\ref{dia:y-psi com})\\(\ref{dia:psipsi})}}= &
-\begin{braid}
 \draw(2,-0.5)node{$\underbrace{\hspace*{16mm}}_{\gamma_1}$};
 \draw(6.5,-0.5)node{$\underbrace{\hspace*{12mm}}_{\gamma_2-e}$};
 \draw(10.5,-0.5)node{$\underbrace{\hspace*{12mm}}_{e}$};
 \draw (0,4)node[above]{$0$}--(0,0);
 \draw (1,4)node[above]{$1$}--(1,0);
 \draw[dots] (1.2,4)--(2.8,4);
 \draw[dots] (1.2,0)--(2.8,0);
 \draw (3,4)node[above]{$e-3$}--(3,0);
 \draw (4,4)node[above]{$e-2$}--(7,2)--(8,0);
 \draw (5,4)node[above]{$e-1$}--(4,2)--(5,0);
 \draw[dots] (5.2,4)--(6.8,4);
 \draw[dots] (5.2,0)--(6.8,0);
 \draw (7,4)node[above]{$e-3$}--(6,2)--(7,0);
 \draw (8,4)node[above]{$e-2$}--(7,2)--(4,0);
 \draw (9,4)node[above]{$e-1$}--(9,0);
 \draw[dots] (9.2,4)--(10.8,4);
 \draw[dots] (9.2,0)--(10.8,0);
 \draw (11,4)node[above]{$e-3$}--(11,0);
 \draw (12,4)node[above]{$e-2$}--(12,0);
 \node[greendot] at (12,2){};
\end{braid}
+
\begin{braid}
 \draw(2,-0.5)node{$\underbrace{\hspace*{16mm}}_{\gamma_1}$};
 \draw(6.5,-0.5)node{$\underbrace{\hspace*{12mm}}_{\gamma_2-e}$};
 \draw(10.5,-0.5)node{$\underbrace{\hspace*{12mm}}_{e}$};
 \draw (0,4)node[above]{$0$}--(0,0);
 \draw (1,4)node[above]{$1$}--(1,0);
 \draw[dots] (1.2,4)--(2.8,4);
 \draw[dots] (1.2,0)--(2.8,0);
 \draw (3,4)node[above]{$e-3$}--(3,0);
 \draw (4,4)node[above]{$e-2$}--(7,2)--(8,0);
 \draw (5,4)node[above]{$e-1$}--(4,2)--(5,0);
 \draw[dots] (5.2,4)--(6.8,4);
 \draw[dots] (5.2,0)--(6.8,0);
 \draw (7,4)node[above]{$e-3$}--(6,2)--(7,0);
 \draw (8,4)node[above]{$e-2$}--(7,2)--(4,0);
 \draw (9,4)node[above]{$e-1$}--(9,0);
 \draw[dots] (9.2,4)--(10.8,4);
 \draw[dots] (9.2,0)--(10.8,0);
 \draw (11,4)node[above]{$e-3$}--(11,0);
 \draw (12,4)node[above]{$e-2$}--(12,0);
 \node[greendot] at (11,2){};
\end{braid} \label{I-problem: induction2: example2: equation1}\\
&& +
\begin{braid}
 \draw(2,-0.5)node{$\underbrace{\hspace*{16mm}}_{\gamma_1}$};
 \draw(6.5,-0.5)node{$\underbrace{\hspace*{12mm}}_{\gamma_2-e}$};
 \draw(10.5,-0.5)node{$\underbrace{\hspace*{12mm}}_{e}$};
 \draw (0,4)node[above]{$0$}--(0,0);
 \draw (1,4)node[above]{$1$}--(1,0);
 \draw[dots] (1.2,4)--(2.8,4);
 \draw[dots] (1.2,0)--(2.8,0);
 \draw (3,4)node[above]{$e-3$}--(3,0);
 \draw (4,4)node[above]{$e-2$}--(7,2)--(8,0);
 \draw (5,4)node[above]{$e-1$}--(4,2)--(5,0);
 \draw[dots] (5.2,4)--(6.8,4);
 \draw[dots] (5.2,0)--(6.8,0);
 \draw (7,4)node[above]{$e-3$}--(6,2)--(7,0);
 \draw (8,4)node[above]{$e-2$}--(7,2)--(4,0);
 \draw (9,4)node[above]{$e-1$}--(9,0);
 \draw[dots] (9.2,4)--(10.8,4);
 \draw[dots] (9.2,0)--(10.8,0);
 \draw (11,4)node[above]{$e-3$}--(11,0);
 \draw (12,4)node[above]{$e-2$}--(12,0);
 \node[greendot] at (9,2){};
\end{braid}
+
\begin{braid}
 \draw(2,-0.5)node{$\underbrace{\hspace*{16mm}}_{\gamma_1}$};
 \draw(6.5,-0.5)node{$\underbrace{\hspace*{12mm}}_{\gamma_2-e}$};
 \draw(10.5,-0.5)node{$\underbrace{\hspace*{12mm}}_{e}$};
 \draw (0,4)node[above]{$0$}--(0,0);
 \draw (1,4)node[above]{$1$}--(1,0);
 \draw[dots] (1.2,4)--(2.8,4);
 \draw[dots] (1.2,0)--(2.8,0);
 \draw (3,4)node[above]{$e-3$}--(3,0);
 \draw (4,4)node[above]{$e-2$}--(9,2)--(12,0);
 \draw (5,4)node[above]{$e-1$}--(4,2)--(5,0);
 \draw[dots] (5.2,4)--(6.8,4);
 \draw[dots] (5.2,0)--(6.8,0);
 \draw (7,4)node[above]{$e-3$}--(6,2)--(7,0);
 \draw (8,4)node[above]{$e-2$}--(7,2)--(8,0);
 \draw (9,4)node[above]{$e-1$}--(10,2)--(9,0);
 \draw[dots] (9.2,4)--(10.8,4);
 \draw[dots] (9.2,0)--(10.8,0);
 \draw (11,4)node[above]{$e-3$}--(12,2)--(11,0);
 \draw (12,4)node[above]{$e-2$}--(9,2)--(4,0);
 \node[greendot] at (7,2){};
\end{braid}\notag.
\end{eqnarray}

By induction and \autoref{send}, the second and the third terms of (\ref{I-problem: induction2: example2: equation1}) are both in $\Blam[\gamma]$. Now for the last term.
\begin{eqnarray}
&&\begin{braid}
 \draw(2,-0.5)node{$\underbrace{\hspace*{16mm}}_{\gamma_1}$};
 \draw(6.5,-0.5)node{$\underbrace{\hspace*{12mm}}_{\gamma_2-e}$};
 \draw(10.5,-0.5)node{$\underbrace{\hspace*{12mm}}_{e}$};
 \draw (0,4)node[above]{$0$}--(0,0);
 \draw (1,4)node[above]{$1$}--(1,0);
 \draw[dots] (1.2,4)--(2.8,4);
 \draw[dots] (1.2,0)--(2.8,0);
 \draw (3,4)node[above]{$e-3$}--(3,0);
 \draw (4,4)node[above]{$e-2$}--(9,2)--(12,0);
 \draw (5,4)node[above]{$e-1$}--(4,2)--(5,0);
 \draw[dots] (5.2,4)--(6.8,4);
 \draw[dots] (5.2,0)--(6.8,0);
 \draw (7,4)node[above]{$e-3$}--(6,2)--(7,0);
 \draw (8,4)[densely dotted] node[above]{$e-2$}--(7,2)--(8,0);
 \draw (9,4)node[above]{$e-1$}--(10,2)--(9,0);
 \draw[dots] (9.2,4)--(10.8,4);
 \draw[dots] (9.2,0)--(10.8,0);
 \draw (11,4)node[above]{$e-3$}--(12,2)--(11,0);
 \draw (12,4)node[above]{$e-2$}--(9,2)--(4,0);
 \node[greendot] at (7,2){};
 \draw[->] (7,2) -- (7.5,3);
\end{braid}\notag\\
& \overset{(\ref{dia:y-psi com})}= &
\begin{braid}
 \draw(2,-0.5)node{$\underbrace{\hspace*{16mm}}_{\gamma_1}$};
 \draw(6.5,-0.5)node{$\underbrace{\hspace*{12mm}}_{\gamma_2-e}$};
 \draw(10.5,-0.5)node{$\underbrace{\hspace*{12mm}}_{e}$};
 \draw (0,4)node[above]{$0$}--(0,0);
 \draw (1,4)node[above]{$1$}--(1,0);
 \draw[dots] (1.2,4)--(2.8,4);
 \draw[dots] (1.2,0)--(2.8,0);
 \draw (3,4)node[above]{$e-3$}--(3,0);
 \draw (4,4)node[above]{$e-2$}--(9,2)--(12,0);
 \draw (5,4)node[above]{$e-1$}--(4,2)--(5,0);
 \draw[dots] (5.2,4)--(6.8,4);
 \draw[dots] (5.2,0)--(6.8,0);
 \draw (7,4)node[above]{$e-3$}--(6,2)--(7,0);
 \draw (8,4)node[above]{$e-2$}--(7,2)--(8,0);
 \draw (9,4)node[above]{$e-1$}--(10,2)--(9,0);
 \draw[dots] (9.2,4)--(10.8,4);
 \draw[dots] (9.2,0)--(10.8,0);
 \draw (11,4)node[above]{$e-3$}--(12,2)--(11,0);
 \draw (12,4)node[above]{$e-2$}--(9,2)--(4,0);
 \node[greendot] at (8,4){};
\end{braid}
-
\begin{braid}
 \draw(2,-0.5)node{$\underbrace{\hspace*{16mm}}_{\gamma_1}$};
 \draw(6.5,-0.5)node{$\underbrace{\hspace*{12mm}}_{\gamma_2-e}$};
 \draw(10.5,-0.5)node{$\underbrace{\hspace*{12mm}}_{e}$};
 \draw (0,4)node[above]{$0$}--(0,0);
 \draw (1,4)node[above]{$1$}--(1,0);
 \draw[dots] (1.2,4)--(2.8,4);
 \draw[dots] (1.2,0)--(2.8,0);
 \draw (3,4)node[above]{$e-3$}--(3,0);
 \draw (4,4)node[above]{$e-2$}--(7,2)--(8,0);
 \draw (5,4)node[above]{$e-1$}--(4,2)--(5,0);
 \draw[dots] (5.2,4)--(6.8,4);
 \draw[dots] (5.2,0)--(6.8,0);
 \draw (7,4)node[above]{$e-3$}--(6,2)--(7,0);
 \draw (8,4)node[above]{$e-2$}--(9,2)--(12,0);
 \draw (9,4)node[above]{$e-1$}--(10,2)--(9,0);
 \draw[dots] (9.2,4)--(10.8,4);
 \draw[dots] (9.2,0)--(10.8,0);
 \draw (11,4)node[above]{$e-3$}--(12,2)--(11,0);
 \draw (12,4)node[above]{$e-2$}--(9,2)--(4,0);
\end{braid} \label{I-problem: induction2: example2: equation2}.
\end{eqnarray}

Substitute (\ref{I-problem: induction2: example2: equation2}) to (\ref{I-problem: induction2: example2: equation1}), let $n = \gamma_1 + \gamma_2$, we have
\begin{eqnarray}
&&\begin{braid}
 \draw(2,-0.5)node{$\underbrace{\hspace*{16mm}}_{\gamma_1}$};
 \draw(6.5,-0.5)node{$\underbrace{\hspace*{12mm}}_{\gamma_2-e}$};
 \draw(10.5,-0.5)node{$\underbrace{\hspace*{12mm}}_{e}$};
 \draw (0,4)node[above]{$0$}--(0,0);
 \draw (1,4)node[above]{$1$}--(1,0);
 \draw[dots] (1.2,4)--(2.8,4);
 \draw[dots] (1.2,0)--(2.8,0);
 \draw (3,4)node[above]{$e-3$}--(3,0);
 \draw (4,4)node[above]{$e-2$}--(11,2)--(12,0);
 \draw (5,4)node[above]{$e-1$}--(4,2)--(5,0);
 \draw[dots] (5.2,4)--(6.8,4);
 \draw[dots] (5.2,0)--(6.8,0);
 \draw (7,4)node[above]{$e-3$}--(6,2)--(7,0);
 \draw (8,4)node[above]{$e-2$}--(7,2)--(8,0);
 \draw (9,4)node[above]{$e-1$}--(8,2)--(9,0);
 \draw[dots] (9.2,4)--(10.8,4);
 \draw[dots] (9.2,0)--(10.8,0);
 \draw (11,4)node[above]{$e-3$}--(10,2)--(11,0);
 \draw (12,4)node[above]{$e-2$}--(11,2)--(4,0);
 \node[greendot] at (7,2){};
\end{braid}
=
\begin{braid}
 \draw(2,-0.5)node{$\underbrace{\hspace*{16mm}}_{\gamma_1}$};
 \draw(6.5,-0.5)node{$\underbrace{\hspace*{12mm}}_{\gamma_2-e}$};
 \draw(10.5,-0.5)node{$\underbrace{\hspace*{12mm}}_{e}$};
 \draw (0,4)node[above]{$0$}--(0,0);
 \draw (1,4)node[above]{$1$}--(1,0);
 \draw[dots] (1.2,4)--(2.8,4);
 \draw[dots] (1.2,0)--(2.8,0);
 \draw (3,4)node[above]{$e-3$}--(3,0);
 \draw (4,4)node[above]{$e-2$}--(7,2)--(8,0);
 \draw (5,4)node[above]{$e-1$}--(4,2)--(5,0);
 \draw[dots] (5.2,4)--(6.8,4);
 \draw[dots] (5.2,0)--(6.8,0);
 \draw (7,4)node[above]{$e-3$}--(6,2)--(7,0);
 \draw (8,4)node[above]{$e-2$}--(7,2)--(4,0);
 \draw (9,4)node[above]{$e-1$}--(9,0);
 \draw[dots] (9.2,4)--(10.8,4);
 \draw[dots] (9.2,0)--(10.8,0);
 \draw (11,4)node[above]{$e-3$}--(11,0);
 \draw (12,4)node[above]{$e-2$}--(12,0);
 \node[greendot] at (12,2){};
\end{braid}\notag \\
&& +
\begin{braid}
 \draw(2,-0.5)node{$\underbrace{\hspace*{16mm}}_{\gamma_1}$};
 \draw(6.5,-0.5)node{$\underbrace{\hspace*{12mm}}_{\gamma_2-e}$};
 \draw(10.5,-0.5)node{$\underbrace{\hspace*{12mm}}_{e}$};
 \draw (0,4)node[above]{$0$}--(0,0);
 \draw (1,4)node[above]{$1$}--(1,0);
 \draw[dots] (1.2,4)--(2.8,4);
 \draw[dots] (1.2,0)--(2.8,0);
 \draw (3,4)node[above]{$e-3$}--(3,0);
 \draw (4,4)node[above]{$e-2$}--(9,2)--(12,0);
 \draw (5,4)node[above]{$e-1$}--(4,2)--(5,0);
 \draw[dots] (5.2,4)--(6.8,4);
 \draw[dots] (5.2,0)--(6.8,0);
 \draw (7,4)node[above]{$e-3$}--(6,2)--(7,0);
 \draw (8,4)node[above]{$e-2$}--(7,2)--(8,0);
 \draw (9,4)node[above]{$e-1$}--(10,2)--(9,0);
 \draw[dots] (9.2,4)--(10.8,4);
 \draw[dots] (9.2,0)--(10.8,0);
 \draw (11,4)node[above]{$e-3$}--(12,2)--(11,0);
 \draw (12,4)node[above]{$e-2$}--(9,2)--(4,0);
 \node[greendot] at (8,4){};
\end{braid}
+
\begin{braid}
 \draw(2,-0.5)node{$\underbrace{\hspace*{16mm}}_{\gamma_1}$};
 \draw(6.5,-0.5)node{$\underbrace{\hspace*{12mm}}_{\gamma_2-e}$};
 \draw(10.5,-0.5)node{$\underbrace{\hspace*{12mm}}_{e}$};
 \draw (0,4)node[above]{$0$}--(0,0);
 \draw (1,4)node[above]{$1$}--(1,0);
 \draw[dots] (1.2,4)--(2.8,4);
 \draw[dots] (1.2,0)--(2.8,0);
 \draw (3,4)node[above]{$e-3$}--(3,0);
 \draw (4,4)node[above]{$e-2$}--(7,2)--(8,0);
 \draw (5,4)node[above]{$e-1$}--(4,2)--(5,0);
 \draw[dots] (5.2,4)--(6.8,4);
 \draw[dots] (5.2,0)--(6.8,0);
 \draw (7,4)node[above]{$e-3$}--(6,2)--(7,0);
 \draw (8,4)node[above]{$e-2$}--(9,2)--(12,0);
 \draw (9,4)node[above]{$e-1$}--(10,2)--(9,0);
 \draw[dots] (9.2,4)--(10.8,4);
 \draw[dots] (9.2,0)--(10.8,0);
 \draw (11,4)node[above]{$e-3$}--(12,2)--(11,0);
 \draw (12,4)node[above]{$e-2$}--(9,2)--(4,0);
\end{braid} \notag\\
& = &
(y_n +
\begin{braid}
 \draw(2,-0.5)node{$\underbrace{\hspace*{16mm}}_{\gamma_1}$};
 \draw(6.5,-0.5)node{$\underbrace{\hspace*{12mm}}_{\gamma_2-e}$};
 \draw(10.5,-0.5)node{$\underbrace{\hspace*{12mm}}_{e}$};
 \draw (0,4)node[above]{$0$}--(0,0);
 \draw (1,4)node[above]{$1$}--(1,0);
 \draw[dots] (1.2,4)--(2.8,4);
 \draw[dots] (1.2,0)--(2.8,0);
 \draw (3,4)node[above]{$e-3$}--(3,0);
 \draw (4,4)node[above]{$e-2$}--(4,0);
 \draw (5,4)node[above]{$e-1$}--(5,0);
 \draw[dots] (5.2,4)--(6.8,4);
 \draw[dots] (5.2,0)--(6.8,0);
 \draw (7,4)node[above]{$e-3$}--(7,0);
 \draw (8,4)node[above]{$e-2$}--(9,2)--(12,0);
 \draw (9,4)node[above]{$e-1$}--(10,2)--(9,0);
 \draw[dots] (9.2,4)--(10.8,4);
 \draw[dots] (9.2,0)--(10.8,0);
 \draw (11,4)node[above]{$e-3$}--(12,2)--(11,0);
 \draw (12,4)node[above]{$e-2$}--(9,2)--(8,0);
\end{braid}){\cdot}
\begin{braid}
 \draw(2,-0.5)node{$\underbrace{\hspace*{16mm}}_{\gamma_1}$};
 \draw(6.5,-0.5)node{$\underbrace{\hspace*{12mm}}_{\gamma_2-e}$};
 \draw(10.5,-0.5)node{$\underbrace{\hspace*{12mm}}_{e}$};
 \draw (0,4)node[above]{$0$}--(0,0);
 \draw (1,4)node[above]{$1$}--(1,0);
 \draw[dots] (1.2,4)--(2.8,4);
 \draw[dots] (1.2,0)--(2.8,0);
 \draw (3,4)node[above]{$e-3$}--(3,0);
 \draw (4,4)node[above]{$e-2$}--(7,2)--(8,0);
 \draw (5,4)node[above]{$e-1$}--(4,2)--(5,0);
 \draw[dots] (5.2,4)--(6.8,4);
 \draw[dots] (5.2,0)--(6.8,0);
 \draw (7,4)node[above]{$e-3$}--(6,2)--(7,0);
 \draw (8,4)node[above]{$e-2$}--(7,2)--(4,0);
 \draw (9,4)node[above]{$e-1$}--(9,0);
 \draw[dots] (9.2,4)--(10.8,4);
 \draw[dots] (9.2,0)--(10.8,0);
 \draw (11,4)node[above]{$e-3$}--(11,0);
 \draw (12,4)node[above]{$e-2$}--(12,0);
\end{braid} \label{I-problem: induction2: example2: equation3}\\
&& +
\begin{braid}
 \draw(2,-0.5)node{$\underbrace{\hspace*{16mm}}_{\gamma_1}$};
 \draw(6.5,-0.5)node{$\underbrace{\hspace*{12mm}}_{\gamma_2-e}$};
 \draw(10.5,-0.5)node{$\underbrace{\hspace*{12mm}}_{e}$};
 \draw (0,4)node[above]{$0$}--(0,0);
 \draw (1,4)node[above]{$1$}--(1,0);
 \draw[dots] (1.2,4)--(2.8,4);
 \draw[dots] (1.2,0)--(2.8,0);
 \draw (3,4)node[above]{$e-3$}--(3,0);
 \draw (4,4)node[above]{$e-2$}--(9,2)--(12,0);
 \draw (5,4)node[above]{$e-1$}--(4,2)--(5,0);
 \draw[dots] (5.2,4)--(6.8,4);
 \draw[dots] (5.2,0)--(6.8,0);
 \draw (7,4)node[above]{$e-3$}--(6,2)--(7,0);
 \draw (8,4)node[above]{$e-2$}--(7,2)--(8,0);
 \draw (9,4)node[above]{$e-1$}--(10,2)--(9,0);
 \draw[dots] (9.2,4)--(10.8,4);
 \draw[dots] (9.2,0)--(10.8,0);
 \draw (11,4)node[above]{$e-3$}--(12,2)--(11,0);
 \draw (12,4)node[above]{$e-2$}--(9,2)--(4,0);
 \node[greendot] at (8,4){};
\end{braid} \notag,
\end{eqnarray}
where by \autoref{I-problem: short}
\begin{eqnarray}
&&y_n +
\begin{braid}
 \draw(2,-0.5)node{$\underbrace{\hspace*{16mm}}_{\gamma_1}$};
 \draw(6.5,-0.5)node{$\underbrace{\hspace*{12mm}}_{\gamma_2-e}$};
 \draw(10.5,-0.5)node{$\underbrace{\hspace*{12mm}}_{e}$};
 \draw (0,4)node[above]{$0$}--(0,0);
 \draw (1,4)node[above]{$1$}--(1,0);
 \draw[dots] (1.2,4)--(2.8,4);
 \draw[dots] (1.2,0)--(2.8,0);
 \draw (3,4)node[above]{$e-3$}--(3,0);
 \draw (4,4)node[above]{$e-2$}--(4,0);
 \draw (5,4)node[above]{$e-1$}--(5,0);
 \draw[dots] (5.2,4)--(6.8,4);
 \draw[dots] (5.2,0)--(6.8,0);
 \draw (7,4)node[above]{$e-3$}--(7,0);
 \draw (8,4)node[above]{$e-2$}--(9,2)--(12,0);
 \draw (9,4)node[above]{$e-1$}--(10,2)--(9,0);
 \draw[dots] (9.2,4)--(10.8,4);
 \draw[dots] (9.2,0)--(10.8,0);
 \draw (11,4)node[above]{$e-3$}--(12,2)--(11,0);
 \draw (12,4)node[above]{$e-2$}--(9,2)--(8,0);
\end{braid}\notag\\
& = &
\begin{braid}
 \draw(2,-0.5)node{$\underbrace{\hspace*{16mm}}_{\gamma_1}$};
 \draw(6.5,-0.5)node{$\underbrace{\hspace*{12mm}}_{\gamma_2-e}$};
 \draw(10.5,-0.5)node{$\underbrace{\hspace*{12mm}}_{e}$};
 \draw (0,4)node[above]{$0$}--(0,0);
 \draw (1,4)node[above]{$1$}--(1,0);
 \draw[dots] (1.2,4)--(2.8,4);
 \draw[dots] (1.2,0)--(2.8,0);
 \draw (3,4)node[above]{$e-3$}--(3,0);
 \draw (4,4)node[above]{$e-2$}--(4,0);
 \draw (5,4)node[above]{$e-1$}--(5,0);
 \draw[dots] (5.2,4)--(6.8,4);
 \draw[dots] (5.2,0)--(6.8,0);
 \draw (7,4)node[above]{$e-3$}--(7,0);
 \draw (8,4)node[above]{$e-2$}--(8,0);
 \draw (9,4)node[above]{$e-1$}--(9,0);
 \draw[dots] (9.2,4)--(10.8,4);
 \draw[dots] (9.2,0)--(10.8,0);
 \draw (11,4)node[above]{$e-3$}--(11,0);
 \draw (12,4)node[above]{$e-2$}--(12,0);
 \node[greendot] at (12,2){};
\end{braid}
+
\begin{braid}
 \draw(2,-0.5)node{$\underbrace{\hspace*{16mm}}_{\gamma_1}$};
 \draw(6.5,-0.5)node{$\underbrace{\hspace*{12mm}}_{\gamma_2-e}$};
 \draw(10.5,-0.5)node{$\underbrace{\hspace*{12mm}}_{e}$};
 \draw (0,4)node[above]{$0$}--(0,0);
 \draw (1,4)node[above]{$1$}--(1,0);
 \draw[dots] (1.2,4)--(2.8,4);
 \draw[dots] (1.2,0)--(2.8,0);
 \draw (3,4)node[above]{$e-3$}--(3,0);
 \draw (4,4)node[above]{$e-2$}--(4,0);
 \draw (5,4)node[above]{$e-1$}--(5,0);
 \draw[dots] (5.2,4)--(6.8,4);
 \draw[dots] (5.2,0)--(6.8,0);
 \draw (7,4)node[above]{$e-3$}--(7,0);
 \draw (8,4)node[above]{$e-2$}--(11,2)--(12,0);
 \draw (9,4)node[above]{$e-1$}--(8,2)--(9,0);
 \draw[dots] (9.2,4)--(10.8,4);
 \draw[dots] (9.2,0)--(10.8,0);
 \draw (11,4)node[above]{$e-3$}--(10,2)--(11,0);
 \draw (12,4)node[above]{$e-2$}--(11,2)--(8,0);
\end{braid}\notag
\end{eqnarray}

\begin{eqnarray}
&&
-
\begin{braid}
 \draw(2,-0.5)node{$\underbrace{\hspace*{16mm}}_{\gamma_1}$};
 \draw(6.5,-0.5)node{$\underbrace{\hspace*{12mm}}_{\gamma_2-e}$};
 \draw(10.5,-0.5)node{$\underbrace{\hspace*{12mm}}_{e}$};
 \draw (0,4)node[above]{$0$}--(0,0);
 \draw (1,4)node[above]{$1$}--(1,0);
 \draw[dots] (1.2,4)--(2.8,4);
 \draw[dots] (1.2,0)--(2.8,0);
 \draw (3,4)node[above]{$e-3$}--(3,0);
 \draw (4,4)node[above]{$e-2$}--(4,0);
 \draw (5,4)node[above]{$e-1$}--(5,0);
 \draw[dots] (5.2,4)--(6.8,4);
 \draw[dots] (5.2,0)--(6.8,0);
 \draw (7,4)node[above]{$e-3$}--(7,0);
 \draw (8,4)node[above]{$e-2$}--(8,0);
 \draw (9,4)node[above]{$e-1$}--(9,0);
 \draw[dots] (9.2,4)--(10.8,4);
 \draw[dots] (9.2,0)--(10.8,0);
 \draw (11,4)node[above]{$e-3$}--(11,0);
 \draw (12,4)node[above]{$e-2$}--(12,0);
 \node[greendot] at (12,2){};
\end{braid}
+
\begin{braid}
 \draw(2,-0.5)node{$\underbrace{\hspace*{16mm}}_{\gamma_1}$};
 \draw(6.5,-0.5)node{$\underbrace{\hspace*{12mm}}_{\gamma_2-e}$};
 \draw(10.5,-0.5)node{$\underbrace{\hspace*{12mm}}_{e}$};
 \draw (0,4)node[above]{$0$}--(0,0);
 \draw (1,4)node[above]{$1$}--(1,0);
 \draw[dots] (1.2,4)--(2.8,4);
 \draw[dots] (1.2,0)--(2.8,0);
 \draw (3,4)node[above]{$e-3$}--(3,0);
 \draw (4,4)node[above]{$e-2$}--(4,0);
 \draw (5,4)node[above]{$e-1$}--(5,0);
 \draw[dots] (5.2,4)--(6.8,4);
 \draw[dots] (5.2,0)--(6.8,0);
 \draw (7,4)node[above]{$e-3$}--(7,0);
 \draw (8,4)node[above]{$e-2$}--(8,0);
 \draw (9,4)node[above]{$e-1$}--(9,0);
 \draw[dots] (9.2,4)--(10.8,4);
 \draw[dots] (9.2,0)--(10.8,0);
 \draw (11,4)node[above]{$e-3$}--(11,0);
 \draw (12,4)node[above]{$e-2$}--(12,0);
 \node[greendot] at (11,2){};
\end{braid}\notag \\
&& +
\begin{braid}
 \draw(2,-0.5)node{$\underbrace{\hspace*{16mm}}_{\gamma_1}$};
 \draw(6.5,-0.5)node{$\underbrace{\hspace*{12mm}}_{\gamma_2-e}$};
 \draw(10.5,-0.5)node{$\underbrace{\hspace*{12mm}}_{e}$};
 \draw (0,4)node[above]{$0$}--(0,0);
 \draw (1,4)node[above]{$1$}--(1,0);
 \draw[dots] (1.2,4)--(2.8,4);
 \draw[dots] (1.2,0)--(2.8,0);
 \draw (3,4)node[above]{$e-3$}--(3,0);
 \draw (4,4)node[above]{$e-2$}--(4,0);
 \draw (5,4)node[above]{$e-1$}--(5,0);
 \draw[dots] (5.2,4)--(6.8,4);
 \draw[dots] (5.2,0)--(6.8,0);
 \draw (7,4)node[above]{$e-3$}--(7,0);
 \draw (8,4)node[above]{$e-2$}--(8,0);
 \draw (9,4)node[above]{$e-1$}--(9,0);
 \draw[dots] (9.2,4)--(10.8,4);
 \draw[dots] (9.2,0)--(10.8,0);
 \draw (11,4)node[above]{$e-3$}--(11,0);
 \draw (12,4)node[above]{$e-2$}--(12,0);
 \node[greendot] at (9,2){};
\end{braid}
-
\begin{braid}
 \draw(2,-0.5)node{$\underbrace{\hspace*{16mm}}_{\gamma_1}$};
 \draw(6.5,-0.5)node{$\underbrace{\hspace*{12mm}}_{\gamma_2-e}$};
 \draw(10.5,-0.5)node{$\underbrace{\hspace*{12mm}}_{e}$};
 \draw (0,4)node[above]{$0$}--(0,0);
 \draw (1,4)node[above]{$1$}--(1,0);
 \draw[dots] (1.2,4)--(2.8,4);
 \draw[dots] (1.2,0)--(2.8,0);
 \draw (3,4)node[above]{$e-3$}--(3,0);
 \draw (4,4)node[above]{$e-2$}--(4,0);
 \draw (5,4)node[above]{$e-1$}--(5,0);
 \draw[dots] (5.2,4)--(6.8,4);
 \draw[dots] (5.2,0)--(6.8,0);
 \draw (7,4)node[above]{$e-3$}--(7,0);
 \draw (8,4)node[above]{$e-2$}--(8,0);
 \draw (9,4)node[above]{$e-1$}--(9,0);
 \draw[dots] (9.2,4)--(10.8,4);
 \draw[dots] (9.2,0)--(10.8,0);
 \draw (11,4)node[above]{$e-3$}--(11,0);
 \draw (12,4)node[above]{$e-2$}--(12,0);
 \node[greendot] at (8,2){};
\end{braid}\notag\\
& = &
\begin{braid}
 \draw(2,-0.5)node{$\underbrace{\hspace*{16mm}}_{\gamma_1}$};
 \draw(6.5,-0.5)node{$\underbrace{\hspace*{12mm}}_{\gamma_2-e}$};
 \draw(10.5,-0.5)node{$\underbrace{\hspace*{12mm}}_{e}$};
 \draw (0,4)node[above]{$0$}--(0,0);
 \draw (1,4)node[above]{$1$}--(1,0);
 \draw[dots] (1.2,4)--(2.8,4);
 \draw[dots] (1.2,0)--(2.8,0);
 \draw (3,4)node[above]{$e-3$}--(3,0);
 \draw (4,4)node[above]{$e-2$}--(4,0);
 \draw (5,4)node[above]{$e-1$}--(5,0);
 \draw[dots] (5.2,4)--(6.8,4);
 \draw[dots] (5.2,0)--(6.8,0);
 \draw (7,4)node[above]{$e-3$}--(7,0);
 \draw (8,4)node[above]{$e-2$}--(11,2)--(12,0);
 \draw (9,4)node[above]{$e-1$}--(8,2)--(9,0);
 \draw[dots] (9.2,4)--(10.8,4);
 \draw[dots] (9.2,0)--(10.8,0);
 \draw (11,4)node[above]{$e-3$}--(10,2)--(11,0);
 \draw (12,4)node[above]{$e-2$}--(11,2)--(8,0);
\end{braid}
+
\begin{braid}
 \draw(2,-0.5)node{$\underbrace{\hspace*{16mm}}_{\gamma_1}$};
 \draw(6.5,-0.5)node{$\underbrace{\hspace*{12mm}}_{\gamma_2-e}$};
 \draw(10.5,-0.5)node{$\underbrace{\hspace*{12mm}}_{e}$};
 \draw (0,4)node[above]{$0$}--(0,0);
 \draw (1,4)node[above]{$1$}--(1,0);
 \draw[dots] (1.2,4)--(2.8,4);
 \draw[dots] (1.2,0)--(2.8,0);
 \draw (3,4)node[above]{$e-3$}--(3,0);
 \draw (4,4)node[above]{$e-2$}--(4,0);
 \draw (5,4)node[above]{$e-1$}--(5,0);
 \draw[dots] (5.2,4)--(6.8,4);
 \draw[dots] (5.2,0)--(6.8,0);
 \draw (7,4)node[above]{$e-3$}--(7,0);
 \draw (8,4)node[above]{$e-2$}--(8,0);
 \draw (9,4)node[above]{$e-1$}--(9,0);
 \draw[dots] (9.2,4)--(10.8,4);
 \draw[dots] (9.2,0)--(10.8,0);
 \draw (11,4)node[above]{$e-3$}--(11,0);
 \draw (12,4)node[above]{$e-2$}--(12,0);
 \node[greendot] at (11,2){};
\end{braid} \label{I-problem: induction2: example2: equation4} \\
&& +
\begin{braid}
 \draw(2,-0.5)node{$\underbrace{\hspace*{16mm}}_{\gamma_1}$};
 \draw(6.5,-0.5)node{$\underbrace{\hspace*{12mm}}_{\gamma_2-e}$};
 \draw(10.5,-0.5)node{$\underbrace{\hspace*{12mm}}_{e}$};
 \draw (0,4)node[above]{$0$}--(0,0);
 \draw (1,4)node[above]{$1$}--(1,0);
 \draw[dots] (1.2,4)--(2.8,4);
 \draw[dots] (1.2,0)--(2.8,0);
 \draw (3,4)node[above]{$e-3$}--(3,0);
 \draw (4,4)node[above]{$e-2$}--(4,0);
 \draw (5,4)node[above]{$e-1$}--(5,0);
 \draw[dots] (5.2,4)--(6.8,4);
 \draw[dots] (5.2,0)--(6.8,0);
 \draw (7,4)node[above]{$e-3$}--(7,0);
 \draw (8,4)node[above]{$e-2$}--(8,0);
 \draw (9,4)node[above]{$e-1$}--(9,0);
 \draw[dots] (9.2,4)--(10.8,4);
 \draw[dots] (9.2,0)--(10.8,0);
 \draw (11,4)node[above]{$e-3$}--(11,0);
 \draw (12,4)node[above]{$e-2$}--(12,0);
 \node[greendot] at (9,2){};
\end{braid}
-
\begin{braid}
 \draw(2,-0.5)node{$\underbrace{\hspace*{16mm}}_{\gamma_1}$};
 \draw(6.5,-0.5)node{$\underbrace{\hspace*{12mm}}_{\gamma_2-e}$};
 \draw(10.5,-0.5)node{$\underbrace{\hspace*{12mm}}_{e}$};
 \draw (0,4)node[above]{$0$}--(0,0);
 \draw (1,4)node[above]{$1$}--(1,0);
 \draw[dots] (1.2,4)--(2.8,4);
 \draw[dots] (1.2,0)--(2.8,0);
 \draw (3,4)node[above]{$e-3$}--(3,0);
 \draw (4,4)node[above]{$e-2$}--(4,0);
 \draw (5,4)node[above]{$e-1$}--(5,0);
 \draw[dots] (5.2,4)--(6.8,4);
 \draw[dots] (5.2,0)--(6.8,0);
 \draw (7,4)node[above]{$e-3$}--(7,0);
 \draw (8,4)node[above]{$e-2$}--(8,0);
 \draw (9,4)node[above]{$e-1$}--(9,0);
 \draw[dots] (9.2,4)--(10.8,4);
 \draw[dots] (9.2,0)--(10.8,0);
 \draw (11,4)node[above]{$e-3$}--(11,0);
 \draw (12,4)node[above]{$e-2$}--(12,0);
 \node[greendot] at (8,2){};
\end{braid}\notag.
\end{eqnarray}

Then by $\lambda\in\mathscr S_n^\Lambda$ and \autoref{send}, for the first term of (\ref{I-problem: induction2: example2: equation3}),
\begin{eqnarray}
&&
(y_n +
\begin{braid}
 \draw(2,-0.5)node{$\underbrace{\hspace*{16mm}}_{\gamma_1}$};
 \draw(6.5,-0.5)node{$\underbrace{\hspace*{12mm}}_{\gamma_2-e}$};
 \draw(10.5,-0.5)node{$\underbrace{\hspace*{12mm}}_{e}$};
 \draw (0,4)node[above]{$0$}--(0,0);
 \draw (1,4)node[above]{$1$}--(1,0);
 \draw[dots] (1.2,4)--(2.8,4);
 \draw[dots] (1.2,0)--(2.8,0);
 \draw (3,4)node[above]{$e-3$}--(3,0);
 \draw (4,4)node[above]{$e-2$}--(4,0);
 \draw (5,4)node[above]{$e-1$}--(5,0);
 \draw[dots] (5.2,4)--(6.8,4);
 \draw[dots] (5.2,0)--(6.8,0);
 \draw (7,4)node[above]{$e-3$}--(7,0);
 \draw (8,4)node[above]{$e-2$}--(9,2)--(12,0);
 \draw (9,4)node[above]{$e-1$}--(10,2)--(9,0);
 \draw[dots] (9.2,4)--(10.8,4);
 \draw[dots] (9.2,0)--(10.8,0);
 \draw (11,4)node[above]{$e-3$}--(12,2)--(11,0);
 \draw (12,4)node[above]{$e-2$}--(9,2)--(8,0);
\end{braid}){\cdot}
\begin{braid}
 \draw(2,-0.5)node{$\underbrace{\hspace*{16mm}}_{\gamma_1}$};
 \draw(6.5,-0.5)node{$\underbrace{\hspace*{12mm}}_{\gamma_2-e}$};
 \draw(10.5,-0.5)node{$\underbrace{\hspace*{12mm}}_{e}$};
 \draw (0,4)node[above]{$0$}--(0,0);
 \draw (1,4)node[above]{$1$}--(1,0);
 \draw[dots] (1.2,4)--(2.8,4);
 \draw[dots] (1.2,0)--(2.8,0);
 \draw (3,4)node[above]{$e-3$}--(3,0);
 \draw (4,4)node[above]{$e-2$}--(7,2)--(8,0);
 \draw (5,4)node[above]{$e-1$}--(4,2)--(5,0);
 \draw[dots] (5.2,4)--(6.8,4);
 \draw[dots] (5.2,0)--(6.8,0);
 \draw (7,4)node[above]{$e-3$}--(6,2)--(7,0);
 \draw (8,4)node[above]{$e-2$}--(7,2)--(4,0);
 \draw (9,4)node[above]{$e-1$}--(9,0);
 \draw[dots] (9.2,4)--(10.8,4);
 \draw[dots] (9.2,0)--(10.8,0);
 \draw (11,4)node[above]{$e-3$}--(11,0);
 \draw (12,4)node[above]{$e-2$}--(12,0);
\end{braid} \notag  \hspace*{4mm}\text{by induction}\\
& =_\gamma &
(y_n +
\begin{braid}
 \draw(2,-0.5)node{$\underbrace{\hspace*{16mm}}_{\gamma_1}$};
 \draw(6.5,-0.5)node{$\underbrace{\hspace*{12mm}}_{\gamma_2-e}$};
 \draw(10.5,-0.5)node{$\underbrace{\hspace*{12mm}}_{e}$};
 \draw (0,4)node[above]{$0$}--(0,0);
 \draw (1,4)node[above]{$1$}--(1,0);
 \draw[dots] (1.2,4)--(2.8,4);
 \draw[dots] (1.2,0)--(2.8,0);
 \draw (3,4)node[above]{$e-3$}--(3,0);
 \draw (4,4)node[above]{$e-2$}--(4,0);
 \draw (5,4)node[above]{$e-1$}--(5,0);
 \draw[dots] (5.2,4)--(6.8,4);
 \draw[dots] (5.2,0)--(6.8,0);
 \draw (7,4)node[above]{$e-3$}--(7,0);
 \draw (8,4)node[above]{$e-2$}--(9,2)--(12,0);
 \draw (9,4)node[above]{$e-1$}--(10,2)--(9,0);
 \draw[dots] (9.2,4)--(10.8,4);
 \draw[dots] (9.2,0)--(10.8,0);
 \draw (11,4)node[above]{$e-3$}--(12,2)--(11,0);
 \draw (12,4)node[above]{$e-2$}--(9,2)--(8,0);
\end{braid}){\cdot}
\begin{braid}
 \draw(2,-0.5)node{$\underbrace{\hspace*{16mm}}_{\gamma_1}$};
 \draw(6.5,-0.5)node{$\underbrace{\hspace*{12mm}}_{\gamma_2-e}$};
 \draw(10.5,-0.5)node{$\underbrace{\hspace*{12mm}}_{e}$};
 \draw (0,4)node[above]{$0$}--(0,0);
 \draw (1,4)node[above]{$1$}--(1,0);
 \draw[dots] (1.2,4)--(2.8,4);
 \draw[dots] (1.2,0)--(2.8,0);
 \draw (3,4)node[above]{$e-3$}--(3,0);
 \draw (4,4)node[above]{$e-2$}--(4,0);
 \draw (5,4)node[above]{$e-1$}--(5,0);
 \draw[dots] (5.2,4)--(6.8,4);
 \draw[dots] (5.2,0)--(6.8,0);
 \draw (7,4)node[above]{$e-3$}--(7,0);
 \draw (8,4)node[above]{$e-2$}--(8,0);
 \draw (9,4)node[above]{$e-1$}--(9,0);
 \draw[dots] (9.2,4)--(10.8,4);
 \draw[dots] (9.2,0)--(10.8,0);
 \draw (11,4)node[above]{$e-3$}--(11,0);
 \draw (12,4)node[above]{$e-2$}--(12,0);
 \node[greendot] at (8,2){};
\end{braid} \notag \hspace*{4mm} \text{by (\ref{I-problem: induction2: example2: equation4})}\\
& = &
\begin{braid}
 \draw(2,-0.5)node{$\underbrace{\hspace*{16mm}}_{\gamma_1}$};
 \draw(6.5,-0.5)node{$\underbrace{\hspace*{12mm}}_{\gamma_2-e}$};
 \draw(10.5,-0.5)node{$\underbrace{\hspace*{12mm}}_{e}$};
 \draw (0,4)node[above]{$0$}--(0,0);
 \draw (1,4)node[above]{$1$}--(1,0);
 \draw[dots] (1.2,4)--(2.8,4);
 \draw[dots] (1.2,0)--(2.8,0);
 \draw (3,4)node[above]{$e-3$}--(3,0);
 \draw (4,4)node[above]{$e-2$}--(4,0);
 \draw (5,4)node[above]{$e-1$}--(5,0);
 \draw[dots] (5.2,4)--(6.8,4);
 \draw[dots] (5.2,0)--(6.8,0);
 \draw (7,4)node[above]{$e-3$}--(7,0);
 \draw (8,4)node[above]{$e-2$}--(11,3)--(12,1)--(12,0);
 \draw (9,4)node[above]{$e-1$}--(8,3)--(8,1)--(9,0);
 \draw[dots] (9.2,4)--(10.8,4);
 \draw[dots] (9.2,0)--(10.8,0);
 \draw (11,4)node[above]{$e-3$}--(10,3)--(10,1)--(11,0);
 \draw (12,4)node[above]{$e-2$}--(12,3)--(11,1)--(8,0);
 \node[greendot] at (8,0){};
 \draw[loosely dashed,red] (-0.5,3) -- (8.5,3)--(8.5,1)--(-0.5,1)--(-0.5,3);
\end{braid}
+
\begin{braid}
 \draw(2,-0.5)node{$\underbrace{\hspace*{16mm}}_{\gamma_1}$};
 \draw(6.5,-0.5)node{$\underbrace{\hspace*{12mm}}_{\gamma_2-e}$};
 \draw(10.5,-0.5)node{$\underbrace{\hspace*{12mm}}_{e}$};
 \draw (0,4)node[above]{$0$}--(0,0);
 \draw (1,4)node[above]{$1$}--(1,0);
 \draw[dots] (1.2,4)--(2.8,4);
 \draw[dots] (1.2,0)--(2.8,0);
 \draw (3,4)node[above]{$e-3$}--(3,0);
 \draw (4,4)node[above]{$e-2$}--(4,0);
 \draw (5,4)node[above]{$e-1$}--(5,0);
 \draw[dots] (5.2,4)--(6.8,4);
 \draw[dots] (5.2,0)--(6.8,0);
 \draw (7,4)node[above]{$e-3$}--(7,0);
 \draw (8,4)node[above]{$e-2$}--(8,0);
 \draw (9,4)node[above]{$e-1$}--(9,0);
 \draw[dots] (9.2,4)--(10.8,4);
 \draw[dots] (9.2,0)--(10.8,0);
 \draw (11,4)node[above]{$e-3$}--(11,0);
 \draw (12,4)node[above]{$e-2$}--(12,0);
 \node[greendot] at (8,2){};
 \node[greendot] at (11,2){};
\end{braid} \notag \text{ by \autoref{send}} \\
&& +
\begin{braid}
 \draw(2,-0.5)node{$\underbrace{\hspace*{16mm}}_{\gamma_1}$};
 \draw(6.5,-0.5)node{$\underbrace{\hspace*{12mm}}_{\gamma_2-e}$};
 \draw(10.5,-0.5)node{$\underbrace{\hspace*{12mm}}_{e}$};
 \draw (0,4)node[above]{$0$}--(0,0);
 \draw (1,4)node[above]{$1$}--(1,0);
 \draw[dots] (1.2,4)--(2.8,4);
 \draw[dots] (1.2,0)--(2.8,0);
 \draw (3,4)node[above]{$e-3$}--(3,0);
 \draw (4,4)node[above]{$e-2$}--(4,0);
 \draw (5,4)node[above]{$e-1$}--(5,0);
 \draw[dots] (5.2,4)--(6.8,4);
 \draw[dots] (5.2,0)--(6.8,0);
 \draw (7,4)node[above]{$e-3$}--(7,0);
 \draw (8,4)node[above]{$e-2$}--(8,0);
 \draw (9,4)node[above]{$e-1$}--(9,0);
 \draw[dots] (9.2,4)--(10.8,4);
 \draw[dots] (9.2,0)--(10.8,0);
 \draw (11,4)node[above]{$e-3$}--(11,0);
 \draw (12,4)node[above]{$e-2$}--(12,0);
 \node[greendot] at (8,2){};
 \node[greendot] at (9,2){};
\end{braid}
-
\begin{braid}
 \draw(2,-0.5)node{$\underbrace{\hspace*{16mm}}_{\gamma_1}$};
 \draw(6.5,-0.5)node{$\underbrace{\hspace*{12mm}}_{\gamma_2-e}$};
 \draw(10.5,-0.5)node{$\underbrace{\hspace*{12mm}}_{e}$};
 \draw (0,4)node[above]{$0$}--(0,0);
 \draw (1,4)node[above]{$1$}--(1,0);
 \draw[dots] (1.2,4)--(2.8,4);
 \draw[dots] (1.2,0)--(2.8,0);
 \draw (3,4)node[above]{$e-3$}--(3,0);
 \draw (4,4)node[above]{$e-2$}--(4,0);
 \draw (5,4)node[above]{$e-1$}--(5,0);
 \draw[dots] (5.2,4)--(6.8,4);
 \draw[dots] (5.2,0)--(6.8,0);
 \draw (7,4)node[above]{$e-3$}--(7,0);
 \draw (8,4)node[above]{$e-2$}--(8,0);
 \draw (9,4)node[above]{$e-1$}--(9,0);
 \draw[dots] (9.2,4)--(10.8,4);
 \draw[dots] (9.2,0)--(10.8,0);
 \draw (11,4)node[above]{$e-3$}--(11,0);
 \draw (12,4)node[above]{$e-2$}--(12,0);
 \node[greendot] at (8,3){};
 \node[greendot] at (8,1){};
\end{braid} \in \Blam[\gamma] \label{I-problem: induction2: example2: equation5}.
\end{eqnarray}

For the second term of (\ref{I-problem: induction2: example2: equation3}), by induction, \autoref{I-problem: short}, $\lambda\in \mathscr S_n^\Lambda$ and \autoref{send},
\begin{eqnarray}
&&
\begin{braid}
 \draw(2,-0.5)node{$\underbrace{\hspace*{16mm}}_{\gamma_1}$};
 \draw(6.5,-0.5)node{$\underbrace{\hspace*{12mm}}_{\gamma_2-e}$};
 \draw(10.5,-0.5)node{$\underbrace{\hspace*{12mm}}_{e}$};
 \draw (0,4)node[above]{$0$}--(0,0);
 \draw (1,4)node[above]{$1$}--(1,0);
 \draw[dots] (1.2,4)--(2.8,4);
 \draw[dots] (1.2,0)--(2.8,0);
 \draw (3,4)node[above]{$e-3$}--(3,0);
 \draw (4,4)[densely dotted] node[above]{$e-2$}--(8,3)--(9,2)--(12,1)--(12,0);
 \draw (5,4)node[above]{$e-1$}--(4,3)--(4,1)--(5,0);
 \draw[dots] (5.2,4)--(6.8,4);
 \draw[dots] (5.2,0)--(6.8,0);
 \draw (7,4)node[above]{$e-3$}--(6,3)--(6,1)--(7,0);
 \draw (8,4)[densely dotted] node[above]{$e-2$}--(7,3)--(7,1)--(8,0);
 \draw (9,4)node[above]{$e-1$}--(9,3)--(10,2)--(9,1)--(9,0);
 \draw[dots] (9.2,4)--(10.8,4);
 \draw[dots] (9.2,0)--(10.8,0);
 \draw (11,4)node[above]{$e-3$}--(11,3)--(12,2)--(11,1)--(11,0);
 \draw (12,4)[densely dotted] node[above]{$e-2$}--(12,3)--(9,2)--(8,1)--(4,0);
 \node[greendot] at (8,4){};
\end{braid}
\overset{(\ref{dia:psipsi})}=
\begin{braid}
 \draw(2,-0.5)node{$\underbrace{\hspace*{16mm}}_{\gamma_1}$};
 \draw(6.5,-0.5)node{$\underbrace{\hspace*{12mm}}_{\gamma_2-e}$};
 \draw(10.5,-0.5)node{$\underbrace{\hspace*{12mm}}_{e}$};
 \draw (0,4)node[above]{$0$}--(0,0);
 \draw (1,4)node[above]{$1$}--(1,0);
 \draw[dots] (1.2,4)--(2.8,4);
 \draw[dots] (1.2,0)--(2.8,0);
 \draw (3,4)node[above]{$e-3$}--(3,0);
 \draw (4,4)node[above]{$e-2$}--(4,3)--(7,2)--(8,1)--(12,0);
 \draw (5,4)node[above]{$e-1$}--(5,3)--(4,2)--(5,1)--(5,0);
 \draw[dots] (5.2,4)--(6.8,4);
 \draw[dots] (5.2,0)--(6.8,0);
 \draw (7,4)node[above]{$e-3$}--(7,3)--(6,2)--(7,1)--(7,0);
 \draw (8,4)node[above]{$e-2$}--(9,3)--(9,1)--(8,0);
 \draw (9,4)node[above]{$e-1$}--(10,3)--(10,1)--(9,0);
 \draw[dots] (9.2,4)--(10.8,4);
 \draw[dots] (9.2,0)--(10.8,0);
 \draw (11,4)node[above]{$e-3$}--(12,3)--(12,1)--(11,0);
 \draw (12,4)node[above]{$e-2$}--(8,3)--(7,2)--(4,1)--(4,0);
 \node[greendot] at (8,4){};
 \draw[loosely dashed,red] (-0.5,3) -- (8.5,3)--(8.5,1)--(-0.5,1)--(-0.5,3);
\end{braid} \hspace*{4mm} \text{by induction} \notag \\
& =_\gamma &
\begin{braid}
 \draw(2,-0.5)node{$\underbrace{\hspace*{16mm}}_{\gamma_1}$};
 \draw(6.5,-0.5)node{$\underbrace{\hspace*{12mm}}_{\gamma_2-e}$};
 \draw(10.5,-0.5)node{$\underbrace{\hspace*{12mm}}_{e}$};
 \draw (0,4)node[above]{$0$}--(0,0);
 \draw (1,4)node[above]{$1$}--(1,0);
 \draw[dots] (1.2,4)--(2.8,4);
 \draw[dots] (1.2,0)--(2.8,0);
 \draw (3,4)node[above]{$e-3$}--(3,0);
 \draw (4,4)node[above]{$e-2$}--(4,0);
 \draw (5,4)node[above]{$e-1$}--(5,0);
 \draw[dots] (5.2,4)--(6.8,4);
 \draw[dots] (5.2,0)--(6.8,0);
 \draw (7,4)node[above]{$e-3$}--(7,0);
 \draw (8,4)[densely dotted] node[above]{$e-2$}--(9,2)--(8,0);
 \draw (9,4)node[above]{$e-1$}--(10,2)--(9,0);
 \draw[dots] (9.2,4)--(10.8,4);
 \draw[dots] (9.2,0)--(10.8,0);
 \draw (11,4)node[above]{$e-3$}--(12,2)--(11,0);
 \draw (12,4)[densely dotted] node[above]{$e-2$}--(8,2)--(12,0);
 \node[greendot] at (8,4){};
 \node[greendot] at (8,2){};
 \draw[->] (8,2) -- (9.5,2.75);
\end{braid}
\overset{\substack{(\ref{dia:y-psi com})\\ (\ref{dia:psipsi})}}=
-\begin{braid}
 \draw(2,-0.5)node{$\underbrace{\hspace*{16mm}}_{\gamma_1}$};
 \draw(6.5,-0.5)node{$\underbrace{\hspace*{12mm}}_{\gamma_2-e}$};
 \draw(10.5,-0.5)node{$\underbrace{\hspace*{12mm}}_{e}$};
 \draw (0,4)node[above]{$0$}--(0,0);
 \draw (1,4)node[above]{$1$}--(1,0);
 \draw[dots] (1.2,4)--(2.8,4);
 \draw[dots] (1.2,0)--(2.8,0);
 \draw (3,4)node[above]{$e-3$}--(3,0);
 \draw (4,4)node[above]{$e-2$}--(4,0);
 \draw (5,4)node[above]{$e-1$}--(5,0);
 \draw[dots] (5.2,4)--(6.8,4);
 \draw[dots] (5.2,0)--(6.8,0);
 \draw (7,4)node[above]{$e-3$}--(7,0);
 \draw (8,4)node[above]{$e-2$}--(8,3)--(9,2)--(12,1)--(12,0);
 \draw (9,4)node[above]{$e-1$}--(9,3)--(10,2)--(9,1)--(9,0);
 \draw[dots] (9.2,4)--(10.8,4);
 \draw[dots] (9.2,0)--(10.8,0);
 \draw (11,4)node[above]{$e-3$}--(11,3)--(12,2)--(11,1)--(11,0);
 \draw (12,4)node[above]{$e-2$}--(12,3)--(9,2)--(8,1)--(8,0);
 \node[greendot] at (8,4){};
 \draw[loosely dashed,red] (7.5,3) -- (12.5,3)--(12.5,1)--(7.5,1)--(7.5,3);
\end{braid} \hspace*{4mm} \text{by \autoref{I-problem: short}} \notag\\
& = &
\begin{braid}
 \draw(2,-0.5)node{$\underbrace{\hspace*{16mm}}_{\gamma_1}$};
 \draw(6.5,-0.5)node{$\underbrace{\hspace*{12mm}}_{\gamma_2-e}$};
 \draw(10.5,-0.5)node{$\underbrace{\hspace*{12mm}}_{e}$};
 \draw (0,4)node[above]{$0$}--(0,0);
 \draw (1,4)node[above]{$1$}--(1,0);
 \draw[dots] (1.2,4)--(2.8,4);
 \draw[dots] (1.2,0)--(2.8,0);
 \draw (3,4)node[above]{$e-3$}--(3,0);
 \draw (4,4)node[above]{$e-2$}--(4,0);
 \draw (5,4)node[above]{$e-1$}--(5,0);
 \draw[dots] (5.2,4)--(6.8,4);
 \draw[dots] (5.2,0)--(6.8,0);
 \draw (7,4)node[above]{$e-3$}--(7,0);
 \draw (8,4)node[above]{$e-2$}--(8,0);
 \draw (9,4)node[above]{$e-1$}--(9,0);
 \draw[dots] (9.2,4)--(10.8,4);
 \draw[dots] (9.2,0)--(10.8,0);
 \draw (11,4)node[above]{$e-3$}--(11,0);
 \draw (12,4)node[above]{$e-2$}--(12,0);
 \node[greendot] at (8,2){};
 \node[greendot] at (12,2){};
\end{braid}
-
\begin{braid}
 \draw(2,-0.5)node{$\underbrace{\hspace*{16mm}}_{\gamma_1}$};
 \draw(6.5,-0.5)node{$\underbrace{\hspace*{12mm}}_{\gamma_2-e}$};
 \draw(10.5,-0.5)node{$\underbrace{\hspace*{12mm}}_{e}$};
 \draw (0,4)node[above]{$0$}--(0,0);
 \draw (1,4)node[above]{$1$}--(1,0);
 \draw[dots] (1.2,4)--(2.8,4);
 \draw[dots] (1.2,0)--(2.8,0);
 \draw (3,4)node[above]{$e-3$}--(3,0);
 \draw (4,4)node[above]{$e-2$}--(4,0);
 \draw (5,4)node[above]{$e-1$}--(5,0);
 \draw[dots] (5.2,4)--(6.8,4);
 \draw[dots] (5.2,0)--(6.8,0);
 \draw (7,4)node[above]{$e-3$}--(7,0);
 \draw (8,4)node[above]{$e-2$}--(8,0);
 \draw (9,4)node[above]{$e-1$}--(9,0);
 \draw[dots] (9.2,4)--(10.8,4);
 \draw[dots] (9.2,0)--(10.8,0);
 \draw (11,4)node[above]{$e-3$}--(11,0);
 \draw (12,4)node[above]{$e-2$}--(12,0);
 \node[greendot] at (8,2){};
 \node[greendot] at (11,2){};
\end{braid} \notag \\
&&
-
\begin{braid}
 \draw(2,-0.5)node{$\underbrace{\hspace*{16mm}}_{\gamma_1}$};
 \draw(6.5,-0.5)node{$\underbrace{\hspace*{12mm}}_{\gamma_2-e}$};
 \draw(10.5,-0.5)node{$\underbrace{\hspace*{12mm}}_{e}$};
 \draw (0,4)node[above]{$0$}--(0,0);
 \draw (1,4)node[above]{$1$}--(1,0);
 \draw[dots] (1.2,4)--(2.8,4);
 \draw[dots] (1.2,0)--(2.8,0);
 \draw (3,4)node[above]{$e-3$}--(3,0);
 \draw (4,4)node[above]{$e-2$}--(4,0);
 \draw (5,4)node[above]{$e-1$}--(5,0);
 \draw[dots] (5.2,4)--(6.8,4);
 \draw[dots] (5.2,0)--(6.8,0);
 \draw (7,4)node[above]{$e-3$}--(7,0);
 \draw (8,4)node[above]{$e-2$}--(8,0);
 \draw (9,4)node[above]{$e-1$}--(9,0);
 \draw[dots] (9.2,4)--(10.8,4);
 \draw[dots] (9.2,0)--(10.8,0);
 \draw (11,4)node[above]{$e-3$}--(11,0);
 \draw (12,4)node[above]{$e-2$}--(12,0);
 \node[greendot] at (8,2){};
 \node[greendot] at (9,2){};
\end{braid}
+
\begin{braid}
 \draw(2,-0.5)node{$\underbrace{\hspace*{16mm}}_{\gamma_1}$};
 \draw(6.5,-0.5)node{$\underbrace{\hspace*{12mm}}_{\gamma_2-e}$};
 \draw(10.5,-0.5)node{$\underbrace{\hspace*{12mm}}_{e}$};
 \draw (0,4)node[above]{$0$}--(0,0);
 \draw (1,4)node[above]{$1$}--(1,0);
 \draw[dots] (1.2,4)--(2.8,4);
 \draw[dots] (1.2,0)--(2.8,0);
 \draw (3,4)node[above]{$e-3$}--(3,0);
 \draw (4,4)node[above]{$e-2$}--(4,0);
 \draw (5,4)node[above]{$e-1$}--(5,0);
 \draw[dots] (5.2,4)--(6.8,4);
 \draw[dots] (5.2,0)--(6.8,0);
 \draw (7,4)node[above]{$e-3$}--(7,0);
 \draw (8,4)node[above]{$e-2$}--(8,0);
 \draw (9,4)node[above]{$e-1$}--(9,0);
 \draw[dots] (9.2,4)--(10.8,4);
 \draw[dots] (9.2,0)--(10.8,0);
 \draw (11,4)node[above]{$e-3$}--(11,0);
 \draw (12,4)node[above]{$e-2$}--(12,0);
 \node[greendot] at (8,1){};
 \node[greendot] at (8,3){};
\end{braid} \notag \\
&&
-
\begin{braid}
 \draw(2,-0.5)node{$\underbrace{\hspace*{16mm}}_{\gamma_1}$};
 \draw(6.5,-0.5)node{$\underbrace{\hspace*{12mm}}_{\gamma_2-e}$};
 \draw(10.5,-0.5)node{$\underbrace{\hspace*{12mm}}_{e}$};
 \draw (0,4)node[above]{$0$}--(0,0);
 \draw (1,4)node[above]{$1$}--(1,0);
 \draw[dots] (1.2,4)--(2.8,4);
 \draw[dots] (1.2,0)--(2.8,0);
 \draw (3,4)node[above]{$e-3$}--(3,0);
 \draw (4,4)node[above]{$e-2$}--(4,0);
 \draw (5,4)node[above]{$e-1$}--(5,0);
 \draw[dots] (5.2,4)--(6.8,4);
 \draw[dots] (5.2,0)--(6.8,0);
 \draw (7,4)node[above]{$e-3$}--(7,0);
 \draw (8,4)node[above]{$e-2$}--(11,3)--(12,1)--(12,0);
 \draw (9,4)node[above]{$e-1$}--(8,3)--(8,1)--(9,0);
 \draw[dots] (9.2,4)--(10.8,4);
 \draw[dots] (9.2,0)--(10.8,0);
 \draw (11,4)node[above]{$e-3$}--(10,3)--(10,1)--(11,0);
 \draw (12,4)node[above]{$e-2$}--(12,3)--(11,1)--(8,0);
 \node[greendot] at (8,4){};
 \draw[loosely dashed,red] (-0.5,3) -- (8.5,3)--(8.5,1)--(-0.5,1)--(-0.5,3);
\end{braid} \notag \hspace*{4mm}\text{by \autoref{send}}\\
& =_\gamma &
\begin{braid}
 \draw(2,-0.5)node{$\underbrace{\hspace*{16mm}}_{\gamma_1}$};
 \draw(6.5,-0.5)node{$\underbrace{\hspace*{12mm}}_{\gamma_2-e}$};
 \draw(10.5,-0.5)node{$\underbrace{\hspace*{12mm}}_{e}$};
 \draw (0,4)node[above]{$0$}--(0,0);
 \draw (1,4)node[above]{$1$}--(1,0);
 \draw[dots] (1.2,4)--(2.8,4);
 \draw[dots] (1.2,0)--(2.8,0);
 \draw (3,4)node[above]{$e-3$}--(3,0);
 \draw (4,4)node[above]{$e-2$}--(4,0);
 \draw (5,4)node[above]{$e-1$}--(5,0);
 \draw[dots] (5.2,4)--(6.8,4);
 \draw[dots] (5.2,0)--(6.8,0);
 \draw (7,4)node[above]{$e-3$}--(7,0);
 \draw (8,4)node[above]{$e-2$}--(8,0);
 \draw (9,4)node[above]{$e-1$}--(9,0);
 \draw[dots] (9.2,4)--(10.8,4);
 \draw[dots] (9.2,0)--(10.8,0);
 \draw (11,4)node[above]{$e-3$}--(11,0);
 \draw (12,4)node[above]{$e-2$}--(12,0);
 \node[greendot] at (8,2){};
 \node[greendot] at (12,2){};
\end{braid}
= e_\gamma y_\gamma \label{I-problem: induction2: example2: equation6}.
\end{eqnarray}

Substitute the results of (\ref{I-problem: induction2: example2: equation5}) and (\ref{I-problem: induction2: example2: equation6}) to (\ref{I-problem: induction2: example2: equation3}), we have
$$
\begin{braid}
 \draw(2,-0.5)node{$\underbrace{\hspace*{16mm}}_{\gamma_1}$};
 \draw(6.5,-0.5)node{$\underbrace{\hspace*{12mm}}_{\gamma_2-e}$};
 \draw(10.5,-0.5)node{$\underbrace{\hspace*{12mm}}_{e}$};
 \draw (0,4)node[above]{$0$}--(0,0);
 \draw (1,4)node[above]{$1$}--(1,0);
 \draw[dots] (1.2,4)--(2.8,4);
 \draw[dots] (1.2,0)--(2.8,0);
 \draw (3,4)node[above]{$e-3$}--(3,0);
 \draw (4,4)node[above]{$e-2$}--(11,2)--(12,0);
 \draw (5,4)node[above]{$e-1$}--(4,2)--(5,0);
 \draw[dots] (5.2,4)--(6.8,4);
 \draw[dots] (5.2,0)--(6.8,0);
 \draw (7,4)node[above]{$e-3$}--(6,2)--(7,0);
 \draw (8,4)node[above]{$e-2$}--(7,2)--(8,0);
 \draw (9,4)node[above]{$e-1$}--(8,2)--(9,0);
 \draw[dots] (9.2,4)--(10.8,4);
 \draw[dots] (9.2,0)--(10.8,0);
 \draw (11,4)node[above]{$e-3$}--(10,2)--(11,0);
 \draw (12,4)node[above]{$e-2$}--(11,2)--(4,0);
 \node[greendot] at (7,2){};
\end{braid} =_\gamma e_\gamma y_\gamma.
$$

\textbf{Case \ref{I-problem: induction2}c:} $i = e-1$.
The method to prove this is the same as for Case \ref{I-problem: induction2}a so it is left as an exercise. Then by induction, this completes the proof.\endproof

Finally, we can use (\ref{I-problem: non-addable: help}) to prove our main result of this subsection.

\begin{Proposition} \label{I-problem: weak result2}
Suppose $m$ is a positive integer, $\lambda = (m,m,1) \in \mathscr S_n^\Lambda$ and $\gamma = (\gamma_1,\gamma_2) = (m,m+1)$. Recall $\lambda_- = (m,m)$. Write $\bi_{\lambda_-} = (i_1, i_2, \ldots, i_{n-1})$. If $k = i_{n-1}+1 \in I$, we have
$$
e_\gamma y_\gamma = e(\bi_{\lambda_-} \vee k) y_{\lambda_-} y_n^{b_k^{\lambda_-}} \in \Blam.
$$
\end{Proposition}

\proof Without loss of generally we assume $\Lambda = \Lambda_0$. When $m = 1$, then $\gamma = (1,2)$ and
$$
e_\gamma y_\gamma =
\begin{braid}
 \draw (0,4)node[above]{$0$}--(0,0);
 \draw (1,4)node[above]{$e-1$}--(1,0);
 \draw (2,4)node[above]{$0$}--(2,0);
\end{braid}
\overset{\substack{(\ref{dia:braid})\\(\ref{dia:ii3})}}=
-\begin{braid}
 \draw (0,4)node[above]{$0$}--(1,2)--(0,0);
 \draw (1,4)node[above]{$e-1$}--(2,2)--(1,0);
 \draw (2,4)node[above]{$0$}--(0,2)--(2,0);
 \node[greendot] at (0,2){};
\end{braid}
-
\begin{braid}
 \draw (0,4)node[above]{$0$}--(2,0);
 \draw (1,4)node[above]{$e-1$}--(0,2)--(1,0);
 \draw (2,4)node[above]{$0$}--(0,0);
\end{braid} \in N_3^{\Lambda_0} \subseteq \Blam.
$$

When $m > 1$, write $\bi_{(\gamma_1)} = (i_1, i_2, \ldots, i_m)$. Set $\sigma = (\gamma_1-1,\gamma_2) = (m-1,m)$. Then by \autoref{I-problem: induction2}, we have
$$
e_\gamma y_\gamma =
\begin{braid}
 \draw (0,4)node[above]{$0$}--(0,0);
 \draw (1,4)node[above]{$1$}--(1,0);
 \draw[dots] (1.2,4)--(2.8,4);
 \draw[dots] (1.2,0)--(2.8,0);
 \draw (3,4)node[above]{$i_m-1$}--(3,0);
 \draw (4,4)node[above]{$i_m$}--(4,0);
 \draw (5,4)node[above]{$e-1$}--(5,0);
 \draw (6,4)node[above]{$0$}--(6,0);
 \draw[dots] (6.2,4)--(7.8,4);
 \draw[dots] (6.2,0)--(7.8,0);
 \draw (8,4)node[above]{$i_m-1$}--(8,0);
 \draw (9,4)node[above]{$i_m$}--(9,0);
\end{braid}
=_\gamma
\begin{cases}
\begin{braid}
 \draw (0,4)node[above]{$0$}--(0,0);
 \draw (1,4)node[above]{$1$}--(1,0);
 \draw[dots] (1.2,4)--(2.8,4);
 \draw[dots] (1.2,0)--(2.8,0);
 \draw (3,4)node[above]{$i_m-1$}--(3,0);
 \draw (4,4)node[above]{$i_m$}--(8,3)--(9,1)--(9,0);
 \draw (5,4)node[above]{$e-1$}--(4,3)--(4,1)--(5,0);
 \draw (6,4)node[above]{$0$}--(5,3)--(5,1)--(6,0);
 \draw[dots] (6.2,4)--(7.8,4);
 \draw[dots] (6.2,0)--(7.8,0);
 \draw (8,4)node[above]{$i_m-1$}--(7,3)--(7,1)--(8,0);
 \draw (9,4)node[above]{$i_m$}--(9,3)--(8,1)--(4,0);
 \draw[loosely dashed,red] (-0.5,3) -- (7.5,3)--(7.5,1)--(-0.5,1)--(-0.5,3);
\end{braid}, & \text{if $i_m \neq e-1$},\\
&\\
\begin{braid}
 \draw (0,4)node[above]{$0$}--(0,0);
 \draw (1,4)node[above]{$1$}--(1,0);
 \draw[dots] (1.2,4)--(2.8,4);
 \draw[dots] (1.2,0)--(2.8,0);
 \draw (3,4)node[above]{$e-2$}--(3,0);
 \draw (4,4)node[above]{$e-1$}--(4,0);
 \draw (5,4)node[above]{$e-1$}--(8,3)--(9,1)--(9,0);
 \draw (6,4)node[above]{$0$}--(5,3)--(5,1)--(6,0);
 \draw[dots] (6.2,4)--(7.8,4);
 \draw[dots] (6.2,0)--(7.8,0);
 \draw (8,4)node[above]{$e-2$}--(7,3)--(7,1)--(8,0);
 \draw (9,4)node[above]{$e-1$}--(9,3)--(8,1)--(5,0);
 \draw[loosely dashed,red] (-0.5,3) -- (7.5,3)--(7.5,1)--(-0.5,1)--(-0.5,3);
\end{braid}, & \text{if $i_m = e-1$}.
\end{cases}
$$

In both cases, the parts bounded by square are both $e_\sigma y_\sigma$. As $|\sigma| = n-2$ and $\lambda \in \mathscr S_n^\Lambda$, by induction, $e_\sigma y_\sigma \in\Blam[\sigma]$. By the definition of $\sigma$ and \autoref{residue sequence has to be the same}, it forces that $e_\sigma y_\sigma \in \Blam[\lambda|_{n-2}]$. Then by \autoref{send}, we have $e_\gamma y_\gamma \in  \Blam $.    \endproof

\subsection{Final part of $y$-problem}

In the last subsection we have proved that if $\lambda = (m,m,1) \in \mathscr S_n^\Lambda$, then
$$
e(\bi_{\lambda_-} \vee k) y_{\lambda_-} y_n^{b_k^{\lambda_-}} \in \Blam
$$
with $k = i_n + 1$. In this subsection we will gradually remove the restrictions on $\lambda$ and $k$. First we are going to introduce a useful homomorphism and use it to prove a few more properties of $\mathscr R_n$ and $\R$. After that we are going to show that if $\lambda \in \mathscr P_I^\Lambda$, then we can extend $\lambda$ to a $\l+1$ multipartition by adding an $\emptyset$ at the end and thus the new multipartition is in $\mathscr P_I^{\Lambda+\Lambda_i}$ for any $i\in I$. Analogous results are also true for $\mathscr P_y^\Lambda$ and $\mathscr P_\psi^\Lambda$. These will allow us to extend the result to an arbitrary multipartition $\lambda$.

For any $\bj \in I^m$, we can define a linear map $\hat\theta_\bj \map{\mathscr{R}_n}{\mathscr{R}_{n+m}^\Lambda}$ sending $e(\bi)$ to $e(\bj\vee \bi)$, $y_r$ to $y_{r+m}$ and $\psi_r$ to $\psi_{r+m}$.\label{notation: hat theta} This map $\hat\theta_{\bj}$ works as embedding from $\mathscr R_n$ to $\mathscr R_{n+m}$ followed by the projection onto $\mathscr R_{n+m}^\Lambda$.

\begin{Lemma} \label{theta_j: homomorphism}
For $\bj\in I^m$, the map $\hat\theta_\bj$ is a homomorphism.
\end{Lemma}

\proof The map is defined to be linear. Hence we only have to check the relations. Since the relations of $\mathscr{R}_n$ and $\mathscr{R}_{n+m}^\Lambda$ from \autoref{notation: KLR} are independent of the value of $r$, we see that $\hat\theta_\bj$ is a homomorphism. \endproof

It will be necessary to cut a multicomposition $\lambda$ into one multicomposition $\mu$ and a composition $\gamma$ for our later work. Note that in our work we will mainly set $\mu$ to be a multipartition and $\gamma$ to be a partition, but generally we don't have such restriction.

\begin{Example}

Fix $e = 4$, $\Lambda = 2\Lambda_0 + \Lambda_1$, $\kappa_\Lambda = (0,1,0)$. Suppose $\lambda = (4,2|2^2,1|3^2,2)$. So
$$
[\lambda] = \Bigg(\hspace*{1mm}\ydiag(4,2)\hspace*{1mm}\Bigg|\hspace*{1mm}%
    \ydiag(2,2,1)\hspace*{1mm}\Bigg|\hspace*{1mm}\ydiag(3,3,2)\hspace*{1mm}\Bigg).
$$

We want to divide the last partition of $\lambda$ after the first row. This is called the cut row of $\lambda$. This gives us a multipartition $\mu$ with Young diagram
$$
[\mu] = \Bigg(\hspace*{1mm}\ydiag(4,2)\hspace*{1mm}\Bigg|\hspace*{1mm}%
    \ydiag(2,2,1)\hspace*{1mm}\Bigg|\hspace*{1mm}\ydiag(3)\hspace*{1mm}\Bigg),
$$
and a partition $\gamma$ with diagram
$$
[\gamma] = \ydiag(3,2).
$$

We call $\mu$ and $\gamma$ the cut part and the remaining part, respectively.

Moreover we want to preserve the following data. The value $|\mu|$ is called the cut of $\lambda$ which is 14 in this case.  The residue of the top left node of $\gamma$ as a subdiagram of $\lambda^{(3)}$ is called the cut residue, which in this case is 3.
\end{Example}

Now we give a formal definition.

\begin{Definition} \label{Definition: cut}

Suppose $\lambda= (\lambda^{(1)},\lambda^{(2)},\ldots,\lambda^{(\l)}) \in \mathscr C_n^\Lambda$ with $\lambda^{(\l)} = (\lambda^{(\l)}_1,\ldots,\lambda^{(\l)}_{k_\l})$ and $a$ is an integer such that $0\leq a < k_\l$. We call $m$ a \textbf{cut} of $\lambda$ and $a$ the \textbf{cut row} associated to $m$ where $m = \sum_{i = 1}^{\l-1}|\lambda^{(i)}| + \sum_{j = 1}^a \lambda_j^{(\l)}$. Define $\Lambda' = \Lambda_s$, where $s = \kappa_\l + 1 - (a+1) = \kappa_\l - a$, the residue of the node at position $(a+1,1,\l)$. We call $s$ to be \textbf{cut residue} associated to $m$ and $\Lambda'$ to be \textbf{cut weight} associated to $m$. We then define $\mu = \lambda|_m \in \mathscr C^\Lambda_m$ and $\gamma = (\lambda_{a+1}^{(\l)},\lambda_{a+2}^{(\l)},\ldots,\lambda_{k_\l}^{(\l)})\in \mathscr C^{\Lambda'}_{n-m}$ and call $\mu$ and $\gamma$ to be \textbf{cut part} and \textbf{remaining part} of $\lambda$ associated to $m$, respectively.

\end{Definition}

Note we can either remove a portion of the last tableau, or cut out the whole partition.

We will start to work with $\hat\theta_{\bi}$, which involving elements in both $\Rn$ and $\R$. Recall that $\hat e(\bi)$, $\hat y_r$, $\hat \psi_s$ and $\hat\psi_{\s\t}$ are elements from $\mathscr R_n$ and $e(\bi)$, $y_r$, $\psi_s$ and $\psi_{\s\t}$ are elements from $\R$.

\begin{Lemma} \label{theta_j: go up}
Suppose $\lambda \in \mathscr S_n^\Lambda$. Let $m$ be a cut of $\lambda$ with $m < n-1$, $\nu = \lambda|_m$ and $\Lambda'$ be the cut weight associated to $m$. Consider $N_{n-m}^{\Lambda'} \subseteq \mathscr{R}_{n-m}$. If $\hat\theta_{\bi_\nu}\map{\mathscr{R}_{n-m}}{\mathscr{R}_n^\Lambda}$, then $\hat\theta_{\bi_\nu}(N_{n-m}^{\Lambda'}) y_\nu \subseteq \Blam$.
\end{Lemma}
\proof Consider $r \in N_{n-m}^{\Lambda'}$. Then by (\ref{basis of KLR}),
$$
r = \sum_{\bj = (j_1,j_2,\ldots,j_{n-m}) \in I^{n-m}} c_\bj R'_\bj\hat e(\bj) \hat y_1^{(\Lambda',\alpha_{j_1})} R_\bj,
$$
where $R_\bj$ and $R'_\bj$ are some elements in $\mathscr{R}_{n-m}$ and $c_\bj \in \Z$. Therefore
\begin{eqnarray*}
\hat\theta_{\bi_\nu}(r)y_\nu & = & \sum_{\bj = (j_1,j_2,\ldots,j_{n-m}) \in I^{n-m}} c_\bj \hat\theta_{\bi_\nu}(R'_\bj)\hat\theta_{\bi_\nu}(\hat e(\bj) \hat y_1^{(\Lambda',\alpha_{j_1})}) \hat\theta_{\bi_\nu}(R_\bj) y_\nu\\
& = & \sum_{\bj = (j_1,j_2,\ldots,j_{n-m}) \in I^{n-m}} c_\bj  \hat\theta_{\bi_\nu}(R'_\bj)e(\bi_\nu\vee j_1 \vee j_2 j_3 \ldots j_{n-m}) y_\nu y_{m+1}^{(\Lambda',\alpha_{j_1})} \hat\theta_{\bi_\nu}(R_\bj)\\
& = & \sum_{\bj = (j_1,j_2,\ldots,j_{n-m}) \in I^{n-m}} c_\bj  \hat\theta_{\bi_\nu}(R'_\bj)\theta_{(j_2,j_3,\ldots,j_{n-m})}(e(\bi_\nu\vee j_1) y_\nu y_{m+1}^{(\Lambda',\alpha_{j_1})}) \hat\theta_{\bi_\nu}(R_\bj).
\end{eqnarray*}

Next we consider $e(\bi_\nu\vee j_1) y_\nu y_{m+1}^{(\Lambda',\alpha_{j_1})} \in \mathscr R_{m+1}^\Lambda$.

Recall that we can write $\nu = (\lambda^{(1)},\lambda^{(2)},\ldots,\lambda^{(\l-1)},\nu^{(\l)})$ and $\nu^{(\l)} = (\nu^{(\l)}_1,\ldots,\nu^{(\l)}_l)$. Let $\mu = \lambda|_{m+1}$. As $m < n-1$, $|\mu| = m+1 < n = |\lambda|$, and $\lambda\in \mathscr S_n^\Lambda$, we have $\mu \in \P_I$. Write $\bi_\nu = (i_1,i_2,\ldots,i_m)$. Notice that $(\Lambda',\alpha_{j_1}) = |\mathscr A_{\t^\nu}^{j_1}|$. We consider two cases.

Suppose ${j_1} = i_m +1$ and $\nu^+$ is a multipartition. By \autoref{notation: b k lambda} we have $|\mathscr A_{\t^\nu}^{j_1}| = b_{j_1}^\nu - 1$. Then by \autoref{I-problem: downstair cases} we have $e(\bi_\nu \vee {j_1})y_\nu y_{m+1}^{(\Lambda',\alpha_{j_1})} = e(\bi_\nu \vee {j_1})y_\nu y_{m+1}^{b_{j_1}^\nu - 1} = e_{\nu^+}y_{\nu^+}$. Because $m$ is a cut of $\lambda$ and $\nu = \lambda|_m$, $\mu = \lambda|_{m+1}$, we must have $\nu^+ > \mu$. Therefore $e_{\nu^+}y_{\nu^+} \in \Blam[\mu]$. So
$$
e(\bi_\nu \vee {j_1})y_\nu y_{m+1}^{(\Lambda',\alpha_{j_1})} \in \Blam[\mu].
$$

Otherwise, by \autoref{notation: b k lambda} we have $|\mathscr A_{\t^\nu}^{j_1}| = b_{j_1}^\nu$. Then by $\mu \in \P_I$ and the definition of $\P_I$, for any ${j_1}\in I$, we have $e(\bi_\nu \vee {j_1})y_\nu y_{m+1}^{b_{j_1}^\nu} \in \Blam[\mu]$ because $\nu = \mu_-$. Therefore
$$
e(\bi_\nu \vee {j_1})y_\nu y_{m+1}^{(\Lambda',\alpha_{j_1})} = e(\bi_\nu \vee {j_1})y_\nu y_{m+1}^{b_{j_1}^\nu} \in \Blam[\mu].
$$

Therefore for any ${j_1}\in I$ we have $e(\bi_\nu \vee {j_1})y_\nu y_{m+1}^{(\Lambda',\alpha_{j_1})} \in \Blam[\mu]$. Hence by \autoref{send} and \autoref{B_lambda is an ideal},
$$
 \hat\theta_{\bi_\nu}(R'_\bj)\theta_{(j_2,j_3,\ldots,j_{n-m})}(e(\bi_\nu\vee j_1) y_\nu y_{m+1}^{(\Lambda',\alpha_{j_1})}) \hat\theta_{\bi_\nu}(R_\bj) \subseteq \Blam.
$$

Therefore $\hat\theta_{\bi_\nu}(r)y_\nu = \sum_{\bj = (j_1,j_2,\ldots,j_{n-m}) \in I^{n-m}} c_\bj  \hat\theta_{\bi_\nu}(R'_\bj)\theta_{(j_2,j_3,\ldots,j_{n-m})}(e(\bi_\nu\vee j_1) y_\nu y_{m+1}^{(\Lambda',\alpha_{j_1})}) \hat\theta_{\bi_\nu}(R_\bj) \subseteq \Blam$.
\endproof

\begin{Definition} \label{notation: concatenate multipartition}
Suppose $\lambda$ is a multicomposition of $m$ and $\mu$ is a composition. If there exists a multicomposition $\gamma$ such that $\lambda$ and $\mu$ are cut part and remaining part of $\gamma$ associated $m$, we write $\gamma = \lambda \vee \mu$ and say $\gamma$ is the \textbf{concatenation} of $\lambda$ and $\mu$.
\end{Definition}

For example, suppose $\lambda = (2^2,1|3^3|2)$ and $\mu = (4,2)$, then $\gamma = \lambda\vee\mu = (2^2,1|3^3|2,4,2)$. Notice that in general $\gamma$ is not a multipartition.

The following Corollaries follows by the definition of $\lambda\vee\mu$.

\begin{Corollary} \label{concatenation: cor1}

Suppose $\lambda$ is a multipartition of $n$ and $\mu, \gamma$ are partitions of $m$. Then $\mu > \gamma$ if and only if $\lambda\vee\mu > \lambda\vee\gamma$.

\end{Corollary}

\begin{Corollary} \label{concatenation: cor2}

Suppose $\lambda$ is a multipartition of $n$ and $\mu$ is a partition of $m$. If $\gamma = \lambda\vee\mu$, $\hat\theta_{\bi_\lambda}(\hat e_\mu \hat y_\mu)y_\lambda = e_\gamma y_\gamma$.
\end{Corollary}

\begin{Corollary} \label{concatenation: cor3}
Suppose $\lambda$ and $\mu$ are multipartitions and $\gamma$ is a partition such that $\lambda = \mu \vee \gamma$. If $\u$ and $\v$ are standard $\gamma$-tableaux, there exist standard $\lambda$-tableaux $\dot\u$ and $\dot\v$ such that $\hat\theta_{\bi_\mu}(\hat\psi_{\u\v})y_\mu = \psi_{\dot\u\dot\v}$.
\end{Corollary}

\proof Suppose $\lambda\in\P_n$ and $\mu\in\P_m$. By \autoref{Definition: cut}, $\mu$ is the cut part of $\lambda$ associated to $m$. Let $a$ be the cut row associated to $m$. Define $\dot\u$ to be the standard $\lambda$-tableau such that $\dot\u|_m = \t^\mu$, and for any $k > m$, if $\dot\u^{-1}(k) = (r_1,c_1,\l_1)$ and $\u^{-1}(k-m) = (r_2,c_2,1)$, then
$$
\begin{array}{rl}
&c_1  =  c_2,\\
&r_1 = r_2 + a.
\end{array}
$$

Define $\dot\v$ similarly. It is trivial that $\hat\theta_{\bi_\mu} (\hat\psi_{d(\u)}) = \psi_{d(\dot\u)}$ and $\hat\theta_{\bi_\mu} (\hat\psi_{d(\v)}) = \psi_{d(\dot\v)}$. Therefore by \autoref{concatenation: cor2},
$$
\hat\theta_{\bi_\mu}(\hat\psi_{\u\v}) = \hat\theta_{\bi_\mu}(\hat\psi_{d(\u)}^*)\hat\theta_{\bi_\mu}(\hat e_\gamma \hat y_\gamma) \hat\theta_{\bi_\mu}(\hat \psi_{d(\v)}) = \psi_{d(\dot\u)}^* e_\lambda y_\lambda \psi_{d(\dot\v)} = \psi_{\dot\u\dot\v}.
$$ \endproof

\begin{Lemma} \label{concatenation: cor4}

Suppose $\lambda \in \mathscr S_n^\Lambda$ and $\mu \in \mathscr C_n^\Lambda$ with $\mu > \lambda$. If $\mu_- \neq \lambda_-$, then $e_\mu y_\mu \in \Blam$.

\end{Lemma}

\proof As $\mu > \lambda$ and $\mu_- \neq \lambda_-$, there exists $m < n$ such that $\mu|_m > \lambda|_m$ and $\mu|_{m-1} = \lambda|_{m-1}$. Set $\nu = \mu|_m$. If $\nu \in \P_m$, we have $e_\nu y_\nu = \psi_{\t^\nu \t^\nu} \in \Blam[\lambda|_m]$, so by \autoref{send} we have $e_\mu y_\mu \in \Blam$.

If $\nu \not\in \P_m$, because $\lambda\in \mathscr S_n^\Lambda$ and $|\nu| = m < n$, we have $\nu \in \P_I$. Notice that if we write $\nu = (\nu^{(1)},\ldots,\nu^{(l)},\emptyset,\ldots,\emptyset)$ with $\nu^{(l)} = (\nu^{(l)}_1,\ldots,\nu^{(l)}_{k-1},\nu^{(l)}_k)$, because $\nu|_{m-1} = \mu|_{m-1} = \lambda|_{m-1} \in \P_{m-1}$ and $\nu \not\in \P_m$, we must have $\nu_{k-1}^{(l)} + 1 = \nu_k^{(l)}$. Therefore if we write $\bi_\nu = (i_1,i_2,\ldots,i_m)$, we have
$$
e_\nu y_\nu = e(\bi_{\nu_-}\vee i_m) y_{\nu_-} y_m^{b_{i_m}^{\nu_-}} \in \Blam[\nu] \subseteq \Blam[\lambda|_m].
$$

Then by \autoref{send}, we have $e_\mu y_\mu \in \Blam$. This completes the proof. \endproof

Now we are ready to start proving the main result of this section. We start by proving two more specialized Propositions. After that we will introduce a Proposition which removes these restrictions and leads to the main Theorem of this section.

\begin{Proposition} \label{I-problem: no empty tail: part 1}

Suppose $\lambda\in\mathscr S_n^\Lambda$ and $\lambda_- = (\lambda^{(1)},\ldots,\lambda^{(\l)}_-)$ with $\lambda^{(\l)}_- = (\lambda^{(\l)}_1,\ldots,\lambda^{(\l)}_{l-1},\lambda^{(\l)}_l) \neq \emptyset$. Write $\bi_{\lambda_-} = (i_1,i_2,\ldots,i_{n-1})$. For $k\in I$, if $k \neq i_{n-1}+1$, or $k = i_{n-1}+1$ and $\lambda^{(\l)}_{l-1} > \lambda^{(\l)}_l$, we have
$$
e(\bi_{\lambda_-}\vee k)y_{\lambda_-} y_n^{b_k^{\lambda_-}} \in \Blam.
$$

\end{Proposition}

\proof For convenience set $m = \lambda^{(\l)}_l$ and $\mu = \lambda|_{n-m-1}$. Therefore $\mu = (\lambda^{(1)},\ldots,\lambda^{(\l-1)},\mu^{(\l)})$ where
$$
\mu^{(\l)} = \begin{cases}
(\lambda^{(\l)}_1,\ldots,\lambda^{(\l)}_{l-1}), & \text{ if $l > 1$,}\\
\emptyset, & \text{ if $l = 1$.}
\end{cases}
$$

Suppose $i$ is the residue of node $(l,1,\l)$ in $\lambda$ and $\Lambda' = \Lambda_i$. Define $\gamma = (m+1)\in \mathscr P^{\Lambda'}_{m+1}$. Notice that $\lambda_- = \mu\vee\gamma_-$. Because $k\neq i_{n-1}+1$ or $k = i_{n-1}+1$ and $\lambda_{k-1}^{(\l)} > \lambda_k^{(\l)}$, we have $b_k^{\gamma_-} = b_k^{\lambda_-}$. By \autoref{onerow: I}, in $\mathscr R_{m+1}^{\Lambda'}$ we have $e(\bi_{\gamma_-}\vee k)y_{\gamma_-}y_{m+1}^{b_k^{\gamma_-}} \in \Blam[\gamma]$. This implies that in $\mathscr R_{m+1}$, $\hat e(\bi_{\gamma_-}\vee k)\hat y_{\gamma_-}\hat y_{m+1}^{b_k^{\gamma_-}} \in N_{m+1}^{\Lambda'}$. Then let $\hat\theta_{\bi_\mu}\map{\mathscr R_{m+1}}{\R}$, by \autoref{theta_j: go up},
$$
e(\bi_{\lambda_-}\vee k)y_{\lambda_-} y_n^{b_k^{\lambda_-}} = e(\bi_\mu \vee \bi_{\gamma_-}\vee k)y_{\lambda_-} y_n^{b_k^{\gamma_-}} = \hat\theta_{\bi_\mu}(\hat e(\bi_{\gamma_-}\vee k) \hat y_{\gamma_-}\hat y_{m+1}^{b_k^{\gamma_-}})y_\mu \in \hat\theta_{\bi_\mu}(N_{m+1}^{\Lambda'}) y_\mu \subseteq\Blam,
$$
which completes the proof.\endproof

\begin{Proposition}\label{I-problem: no empty tail: part 2}

Suppose $\lambda = (\lambda^{(1)},\ldots,\lambda^{(\l)})\in\mathscr S_n^\Lambda$ with $\lambda^{(\l)} = (\lambda^{(\l)}_1,\ldots,\lambda^{(\l)}_{l-1},\lambda^{(\l)}_l,1)$ and $l \geq 2$, where $\lambda^{(\l)}_{l-1} = \lambda^{(\l)}_l$. Write $\bi_{\lambda_-} = (i_1,i_2,\ldots,i_{n-1})$. Suppose $k\in I$ and $k \equiv i_{n-1}+1\pmod{e}$. Then
$$
e(\bi_{\lambda_-}\vee k)y_{\lambda_-} y_n^{b_k^{\lambda_-}}\in \Blam.
$$

\end{Proposition}

\proof For convenience set $m = \lambda^{(\l)}_{l-1} = \lambda^{(\l)}_l$, and $\mu = \lambda|_{n-2m-1}$. Therefore $\mu = (\lambda^{(1)},\ldots,\lambda^{(\l-1)},\mu^{(\l)})$ where
$$
\mu^{(\l)} = \begin{cases}
(\lambda^{(\l)}_1,\ldots,\lambda^{(\l)}_{l-2}), & \text{ if $l > 2$,}\\
\emptyset, & \text{ if $l = 2$.}
\end{cases}
$$

Suppose $i$ is the residue of node $(l-1,1,\l)$ in $\lambda$ and $\Lambda' = \Lambda_i$. Define $\gamma = (m,m+1)\in \mathscr P^{\Lambda'}_{2m+1}$. Notice that $\lambda_- = \mu\vee \gamma_-$. Because $k \equiv i_{n-1}+1\pmod{e}$, we have $b_k^{\gamma_-} = b_k^{\lambda_-}$ and $e(\bi_{\gamma_-}\vee k)y_{\gamma_-} y_{2m+1}^{b_k^{\gamma_-}} = e_\gamma y_\gamma$. By \autoref{I-problem: weak result2}, we have $e_\gamma y_\gamma \in \Blam[\gamma]$. Therefore we can write $e_\gamma y_\gamma = \sum_{\substack{\u,\v\in \Std(\sigma)\\ \sigma > \gamma}} c_{\u\v}\psi_{\u\v}$ with $\sigma = (\sigma_1,\sigma_2)$ where $\sigma_2 \geq 0$ and $\sigma_1 > \gamma_1 = m$. Therefore in $\mathscr R_{2m+1}$, we have
$$
\hat e_\gamma \hat y_\gamma = \sum_{\substack{\u,\v\in \Std(\sigma)\\ \sigma > \gamma}} c_{\u\v}\hat\psi_{\u\v} + r,
$$
with $r \in N_{2m+1}^{\Lambda'}$ and $c_{\u\v}\in\Z$. Therefore
\begin{align}
e(\bi_{\lambda_-}\vee k)y_{\lambda_-}y_n^{b_k^{\lambda_-}} & = e(\bi_\mu\vee \bi_{\gamma_-}\vee k)y_{\lambda_-} y_n^{b_k^{\gamma_-}} = \hat\theta_{\bi_\mu} (\hat e(\bi_{\gamma_-}\vee k)\hat y_{\gamma_-}\hat y_{2m+1}^{b_k^{\gamma_-}})y_\mu = \hat\theta_{\bi_\mu}(\hat e_\gamma \hat y_\gamma) y_\mu \notag \\
& = \sum_{\substack{\u,\v\in\Std(\sigma)\\ \sigma > \gamma}} c_{\u\v}\hat\theta_{\bi_\mu}(\hat \psi_{\u\v})y_\mu + \hat\theta_{\bi_\mu}(r) \label{I-problem: no empty tail: part2: help}.
\end{align}

For the first term of (\ref{I-problem: no empty tail: part2: help}), define $\alpha = \mu\vee\sigma \in \mathscr C_n^\Lambda$. Because $\sigma > \gamma$, by \autoref{concatenation: cor1} we have $\alpha = \mu\vee\sigma > \mu\vee\gamma > \lambda$. Therefore
$$
\hat\theta_{\bi_\mu}(\hat\psi_{\u\v})y_\mu = \hat\theta_{\bi_\mu}(\hat\psi_{d(\u)}^*)\hat\theta_{\bi_\mu}(\hat e_\sigma \hat y_\sigma) y_\mu \hat\theta_{\bi_\mu}(\hat\psi_{d(\v)}) = \hat\theta_{\bi_\mu}(\hat\psi_{d(\u)}^*)e_{\mu\vee\sigma}y_{\mu\vee\sigma}\hat\theta_{\bi_\mu}(\hat\psi_{d(\v)}) = \hat\theta_{\bi_\mu}(\hat\psi_{d(\u)}^*)e_\alpha y_\alpha \hat\theta_{\bi_\mu}(\hat\psi_{d(\v)})
$$

Because $\sigma = (\sigma_1,\sigma_2) > \gamma = (m,m+1)$, we must have $\sigma_1 > m$. Therefore $\alpha_- = \mu\vee\sigma_- \neq \mu\vee\gamma_- = \lambda_-$. Then by \autoref{concatenation: cor4}, $e_\alpha y_\alpha \in \Blam$. By \autoref{B_lambda is an ideal}, we have $\hat\theta_{\bi_\mu}(\hat\psi_{\u\v})y_\mu = \hat\theta_{\bi_\mu}(\hat\psi_{d(\u)}^*)e_\alpha y_\alpha \hat\theta_{\bi_\mu}(\hat\psi_{d(\v)})\in \Blam$ which yields $\sum_{\substack{\u,\v\in\Std(\sigma)\\ \sigma > \gamma}} c_{\u\v}\hat\theta_{\bi_\mu}(\hat \psi_{\u\v})y_\mu\in \Blam$.

For the second term of (\ref{I-problem: no empty tail: part2: help}), by \autoref{theta_j: go up}, we have $\hat\theta_{\bi_\mu}(r) \in \Blam$. Therefore $$e(\bi_{\lambda_-}\vee k)y_{\lambda_-} y_n^{b_k^{\lambda_-}}\in \Blam.$$\endproof

Suppose $\lambda \in \mathscr S_n^\Lambda$. If $\lambda_- = (\lambda^{(1)},\lambda^{(2)},\ldots,\lambda^{(\l-1)},\lambda_-^{(\l)})$ with $\lambda_-^{(\l)} \neq \emptyset$ by \autoref{I-problem: no empty tail: part 1} and~\autoref{I-problem: no empty tail: part 2} we have $\lambda\in\P_I$. In the rest of the subsection we are going to prove the result is still true if $\lambda_-^{(\l)} = \emptyset$.

Suppose $\mu = (\mu^{(1)},\ldots,\mu^{(\l)}) \in \P_n$ and $\kappa_\Lambda = (\kappa_1,\ldots,\kappa_\l)$, if $\mu^{(\l)} = \emptyset$, we define $\bar\Lambda = \Lambda - \Lambda_{\kappa_\l}$, $\kappa_{\bar\Lambda} = (\kappa_1,\ldots,\kappa_{\l-1})$ and $\bar\mu = (\mu^{(1)},\ldots,\mu^{(\l-1)}) \in \mathscr P^{\bar\Lambda}_n$. Suppose $\u,\v$ are two standard $\mu$-tableaux, define $\bar \u$ and $\bar \v$ to be standard $\bar\mu$-tableaux obtained by removing the $\emptyset$ at the end of $\u$ and $\v$ respectively. Write $k = \kappa_\l$ for convenience. If $\bi = (i_1,\ldots, i_n) \in I^n$, define
\begin{equation*} \label{notation: y i k}
y_{\bi,k} = y_1^{\delta_{i_1,k}}y_2^{\delta_{i_2,k}}\ldots y_n^{\delta_{i_n,k}}.
\end{equation*}

\begin{Lemma} \label{remove empty tableaux: help}
Suppose the notations are defined as above and $\bi_\v$ is the residue sequence of $\v$, then
$$
\psi_{\u\v} = \psi_{\bar \u \bar \v} y_{\bi_\nu,k},
$$
where $\psi_{\u\v}$ is an element in $\mathscr{R}_n^\Lambda$ and $\psi_{\bar \u \bar \v}$ is an element in $\mathscr{R}_n^{\bar\Lambda}$.
\end{Lemma}

\proof Without loss of generality, assume $\u = \t^\mu$. By the definition of $\mu$ and $\bar\mu$, writing $\bi_\mu = (i_1, i_2, \ldots, i_n)$, we have $y_\mu = y_{\bar\mu} y_1^{\delta_{i_1,k}}y_2^{\delta_{i_2,k}}\ldots y_n^{\delta_{i_n,k}} = y_{\bar\mu} y_{\bi_\mu,k}$.

Now for any residue sequence $\bi = (i_1,i_2,\ldots, i_n) \in I^n$ and any $r$, If $i_r \neq i_{r+1}$
\begin{eqnarray*}
e(\bi)y_1^{\delta_{i_1,k}} y_2^{\delta_{i_2,k}}\ldots y_n^{\delta_{i_n,k}}\psi_r & = & (e(\bi)y_r^{\delta_{i_r,k}}y_{r+1}^{\delta_{i_{r+1},k}} \psi_r) y_1^{\delta_{i_1,k}} \ldots y_{r-1}^{\delta_{i_{r-1},k}}y_{r+2}^{\delta_{i_{r+1},k}}\ldots y_n^{\delta_{i_n,k}}\\
& = &e(\bi) \psi_r y_r^{\delta_{i_{r+1},k}} y_{r+1}^{\delta_{i_r,k}} y_1^{\delta_{i_1,k}} \ldots y_{r-1}^{\delta_{i_{r-1},k}}y_{r+2}^{\delta_{i_{r+1},k}}\ldots y_n^{\delta_{i_n,k}}\\
& = & e(\bi) \psi_r y_1^{\delta_{s_r(i_1),k}} \ldots y_{n}^{\delta_{s_r(i_{n}),k}} = e(\bi) \psi_r y_{\bi{\cdot}s_r,k}.
\end{eqnarray*}

If $i_r = i_{r+1}$, then by relation (\ref{alg:ypsi}), as $\delta_{i_r,k} = \delta_{i_{r+1},k}$, we have the same result.

Hence
\begin{eqnarray*}
\psi_{\t^\mu\v} & = & e_\mu y_\mu \psi_{d(\v)} = e(\bi_\mu)y_{\bar\mu} y_{\bi_\mu,k} \psi_{d(\v)}\\
& = & e(\bi_\mu)y_{\bar\mu}\psi_{d(\v)} y_{\bi_\mu{\cdot}d(\v),k} = e_\mu y_{\bar\mu} \psi_{d(v)} y_{\bi_\v,k}.
\end{eqnarray*}

As $e_\mu = e_{\bar\mu}$ and $\psi_{d(\v)} = \psi_{d(\bar \v)}$, this completes the proof. \endproof

\begin{Proposition} \label{remove empty tableaux}

Suppose $\mu$, $\kappa_\Lambda$, $\bar\mu$, and $\kappa_{\bar\Lambda}$ are defined as above. Then $\bar\mu \in \mathscr P_I^{\bar\Lambda}\cap \mathscr P_y^{\bar\Lambda}\cap  \mathscr P_\psi^{\bar\Lambda}$ implies $\mu \in \P_I \cap \P_y \cap \P_\psi$.

\end{Proposition}

\proof We are only going to prove that $\bar\mu \in \mathscr P_I^{\bar\Lambda}$ implies $\mu \in \P_I$. The other two cases are similar.

Suppose $\bar\mu \in \mathscr P^{\bar\Lambda}_I$. Then for any $s \in I$, by the definition of $\mathscr P^{\bar\Lambda}_I$,
$$
e(\bi_{\bar\mu_-} \vee s)y_{\bar\mu_-} y_n^{b_s^{\bar\mu_-}}  = \sum_{\bar\u,\bar\v \in \Std( > \bar\mu)} c_{\bar\u \bar\v} \psi_{\bar\u \bar \v},
$$
where $\bi_{\bar\v} = \bi_{\bar\mu_-}\vee s = \bi_{\mu_-}\vee s$ and $c_{\bar\u \bar\v}\in \Z$.

Also we have $e(\bi_{\bar\mu_-} \vee s)y_{\bar\mu_-} y_n^{b_s^{\bar\mu_-}} = \theta_s(\psi_{t^{\bar\mu_-} t^{\bar\mu_-}})  y_n^{b_s^{\bar\mu_-}}$. Therefore we have
$$
\theta_s(\psi_{t^{\bar\mu_-} t^{\bar\mu_-}})  y_n^{b_s^{\bar\mu_-}} = \sum_{\bar\u,\bar\v \in \Std( > \bar\mu)} c_{\bar\u \bar\v} \psi_{\bar\u \bar \v}.
$$

Notice that $\t^{\bar\mu_-} = \overline {\t^{\mu_-}}$. Recall $k = \kappa_\l$, the last term of the multicharge $\kappa_\Lambda$. We consider two cases, $s \neq k$ and $s = k$ in the rest of the proof.

If $s \neq k$, then $b_s^{\mu_-} = b_s^{\bar\mu_-}$. Hence by \autoref{remove empty tableaux: help}
\begin{eqnarray*}
e(\bi_{\mu_-} \vee s)y_{\mu_-} y_n^{b_s^{\mu_-}} & = & \theta_s(\psi_{t^{\mu_-} t^{\mu_-}})  y_n^{b_s^{\mu_-}}\\
& = & \theta_s(\psi_{\bar \t^{\mu_-} \bar \t^{\mu_-}} y_{\bi_{\mu_-},k}) y_n^{b_s^{\bar\mu_-}}\\
& = & \theta_s(\psi_{\t^{\bar\mu_-} \t^{\bar\mu_-}}) y_n^{b_s^{\bar\mu_-}}  y_{\bi_{\mu_-},k}\\
& = & \sum_{\bar\u,\bar\v \in \Std( > \bar\mu)} c_{\bar\u \bar\v} \psi_{\bar\u \bar\v} y_{\bi_{\mu_-},k},
\end{eqnarray*}
and as $s \neq k$, $\delta_{s,k} = 0$. Hence $y_{\bi_{\mu_-},k} = y_{\bi_{\mu_-}\vee s, k} = y_{\bi_{\bar\v},k} = y_{\bi_\v,k}$. By \autoref{remove empty tableaux: help},
$$
e(\bi_{\mu_-} \vee s)y_{\mu_-} y_n^{b_s^{\mu_-}} = \sum_{\bar\u,\bar\v \in \Std( > \bar\mu)} c_{\bar\u \bar\v} \psi_{\bar\u \bar\v} y_{\bi_\v,k} = \sum_{\bar\u,\bar\v \in \Std( > \bar\mu)} c_{\bar\u \bar\v} \psi_{\u\v} \in \Blam,
$$
because $\bar\u,\bar\v \in \Std( > \bar \mu)$ implies $\u,\v \Std( > \mu)$.\\

If $s = k$, then $b_s^{\mu_-} = b_s^{\bar\mu_-} + 1$. Hence by \autoref{remove empty tableaux: help}
\begin{eqnarray*}
e(\bi_{\mu_-} \vee s)y_{\mu_-} y_n^{b_s^{\mu_-}} & = & \theta_s(\psi_{t^{\mu_-} t^{\mu_-}})  y_n^{b_s^{\mu_-}}\\
& = & \theta_s(\psi_{\bar \t^{\mu_-} \bar \t^{\mu_-}} y_{\bi_{\mu_-},k}) y_n^{b_s^{\bar\mu_-}}y_n\\
& = & \theta_s(\psi_{\t^{\bar\mu_-} \t^{\bar\mu_-}}) y_n^{b_s^{\bar\mu_-}}  y_{\bi_{\mu_-},k}y_n\\
& = & \sum_{\bar\u,\bar\v \in \Std( > \bar\mu)} c_{\bar\u \bar\v} \psi_{\bar\u \bar\v} y_{\bi_{\mu_-},k}y_n,
\end{eqnarray*}
and as $s = k$, $\delta_{s,k} = 1$. Hence $y_{\bi_{\mu_-},k}y_n = y_{\bi_{\mu_-}\vee s, k} = y_{\bi_{\bar\v},k} = y_{\bi_\v,k}$. By \autoref{remove empty tableaux: help}
$$
e(\bi_{\mu_-} \vee s)y_{\mu_-} y_n^{b_s^{\mu_-}} = \sum_{\bar\u,\bar\v \in \Std( > \bar\mu)} c_{\bar\u \bar\v} \psi_{\bar\u \bar\v} y_{\bi_\v,k} = \sum_{\bar\u,\bar\v \in \Std( > \bar\mu)} c_{\bar\u \bar\v} \psi_{\u\v} \in \Blam.
$$

These implies that $\mu \in \mathscr P_I$. \endproof

Now we are ready to prove \autoref{I-problem: final}.\\

\hspace*{-6mm}\textbf{Proof of \autoref{I-problem: final}.} Write $\mu = \lambda_- = (\lambda^{(1)},\lambda^{(2)},\ldots,\lambda^{(\l-1)},\mu^{(\l)})$. If $\mu^{(\l)} \neq \emptyset$, by \autoref{I-problem: no empty tail: part 1} and \autoref{I-problem: no empty tail: part 2}, we have $\lambda\in\P_I$.

If $\mu^{(\l)} = \emptyset$, write $\lambda^{(\l-1)} = (\lambda^{(\l-1)}_1,\lambda^{(\l-1)}_2,\ldots,\lambda^{(\l-1)}_{k_{\l-1}})$ and define $\gamma = (\lambda^{(1)},\ldots,\lambda^{(\l-2)},\gamma^{(\l-1)},\emptyset) \in \mathscr P_n^\Lambda$ with $\gamma^{(\l-1)} = (\lambda^{(\l-1)}_1,\lambda^{(\l-1)}_2,\ldots,\lambda^{(\l-1)}_{k_{\l-1}},1)$ and $\bar\gamma = (\lambda^{(1)},\ldots,\lambda^{(\l-2)},\gamma^{(\l-1)})$. As $l(\bar\gamma) < l(\lambda) = \l$, by the definition of $\mathscr S_n^\Lambda$, $\bar\gamma \in \mathscr P_I^{\bar\Lambda}$. Then by \autoref{remove empty tableaux} we have $\gamma \in \P_I$. Since $\gamma_- = \mu = \lambda_-$ and $\gamma > \lambda$, for any $k\in I$,
$$
e(\bi_{\lambda_-}\vee k) y_{\lambda_-} y_n^{b_k^{\lambda_-}} = e(\bi_{\gamma_-}\vee k) y_{\gamma_-} y_n^{b_k^{\gamma_-}} \in \Blam[\gamma] \subseteq \Blam,
$$
which yields that $\lambda \in \mathscr P_I$. This completes the proof. \endproof

The following Corollary is directly implied by \autoref{I-problem: final}. It will contribute to proving the $\psi$-problem.

\begin{Corollary} \label{composition in Blam}

Suppose $\lambda \in \mathscr S_n^\Lambda$ and $\mu \in \mathscr C_n^\Lambda$ where $\mu > \lambda$. Then we have $e_\mu y_\mu \in \Blam$.

\end{Corollary}

\proof If $\mu_- \neq \lambda_-$, using \autoref{concatenation: cor4} we have $e_\mu y_\mu \in \Blam$. Suppose then that $\mu_- = \lambda_-$. If $\mu \in \P_n$, then $e_\mu y_\mu = \psi_{\t^\mu \t^\mu} \in \Blam$. Finally, suppose that $\mu \not\in \P_n$. If we write $\mu = (\mu^{(1)},\ldots,\mu^{(l)},\emptyset,\ldots,\emptyset)$ with $\mu^{(l)} = (\mu^{(l)}_1,\ldots,\mu^{(l)}_{k-1},\mu^{(l)}_k)$, we must have $\mu^{(l)}_{k-1} + 1 = \mu^{(l)}_k$. If we write $\bi_\mu = (i_1,i_2,\ldots,i_n)$, we have $e_\mu y_\mu = e(\bi_{\mu_-}\vee i_n) y_{\mu_-} y_n^{b_{i_n}^{\mu_-}}$. By \autoref{I-problem: final}, as $\lambda\in \mathscr S_n^\Lambda$, we have $\lambda\in \P_I$. Since $\lambda_- = \mu_-$,
$$
e_\mu y_\mu = e(\bi_{\mu_-}\vee i_n) y_{\mu_-} y_n^{b_{i_n}^{\mu_-}} = e(\bi_{\lambda_-}\vee i_n) y_{\lambda_-} y_n^{b_{i_n}^{\lambda_-}}  \in \Blam.
$$\endproof

\section{Integral Basis Theorem II}

In this section our main purpose is to prove that $\R = \BZ$ by proving that $\psi_{\s\t}y_r$ and $\psi_{\s\t}\psi_r$ are both in $\BZ$. We first define an integer $m_\lambda$ such that if $\t\in\Std(\lambda)$ and $l(d(\t)) < m_\lambda$, we have $\psi_{\s\t}y_r \in \BZ$ and $\psi_{\s\t}\psi_r \in \BZ$. Our first step is to show that $m_\lambda > 0$. Then we prove if $l(d(\t)) \leq m_\lambda$, we will have $\psi_{\s\t}y_r \in \BZ$ and $\psi_{\s\t}\psi_r \in \BZ$ as well. By induction we will show that for any $\t\in\Std(\lambda)$, $l(d(\t)) < m_\lambda$, which indicates that $\psi_{\s\t}y_r \in \BZ$ and $\psi_{\s\t}\psi_r \in \BZ$ for any $\s,\t\in\Std(\lambda)$. Finally combining the results from the last section, we can prove that $\R = \BZ$.

\subsection{Base case of induction}
In this section we fix $\lambda\in \mathscr S_n^\Lambda$. First we will give a proper definition of $m_\lambda$.

\begin{Definition}\label{notation: W lambda}
Define $m_\lambda$ to be the smallest nonnegative integer such that for any standard $\lambda$-tableau $\t$ with $l(d(\t)) < m_\lambda$ we have
\begin{eqnarray*}
\psi_{\s\t} y_r & = & \sum_{(\u,\v) \rhd (\s,\t)} c_{\u\v} \psi_{\u\v},\\
\psi_{\s\t} \psi_r & = & \begin{cases}
\psi_{\s \w} + \sum_{(\u,\v)\rhd (\s,\t)}c_{\u\v}\psi_{\u\v},& \text{ if $\w = \t{\cdot}s_r$ is standard and $d(\u){\cdot}s_r$ is reduced,}\\
\sum_{(\u,\v)\rhd (\s,\t)}c_{\u\v} \psi_{\u\v}, & \text{ if $\u{\cdot}s_r$ is not standard or $d(\u){\cdot}s_r$ is not reduced.}
\end{cases}
\end{eqnarray*}
for some $c_{\u\v} \in \Z$.
\end{Definition}

We will use induction to prove that for any $\t\in\Std(\lambda)$, $l(d(\t)) < m_\lambda$ in this section. In this subsection we will prove that $m_\lambda > 0$, which is the base case of the induction.

\begin{Lemma} \label{y-problem: induction}

Suppose $\lambda\in \mathscr S_n^\Lambda$. For any $1 \leq r \leq n$, $e_\lambda y_\lambda y_r = \psi_{\t^\lambda \t^\lambda}y_r \in \Blam$.

\end{Lemma}

\proof If $ r < n$, write $\mu = \lambda|_r$. As $\lambda\in \mathscr S_n^\Lambda$ we have $\mu \in \P_y$. Therefore $e_\mu y_\mu y_r \in \Blam[\mu]$. By \autoref{send}, we have $e_\lambda y_\lambda y_r \in \Blam$.

If $r = n$, write $\bi_\lambda = (i_1,i_2,\ldots,i_n)$. There exists a positive integer $b$ such that $e_\lambda y_\lambda = e(\bi_{\lambda_-}\vee i_n)y_{\lambda_-}y_n^b$. By the definition of $b_{i_n}^{\lambda_-}$ we have $b < b_{i_n}^{\lambda_-}$. If $b + 1 < b_{i_n}^{\lambda_-}$, by \autoref{I-problem: downstair cases} there exists $\nu\in \P_n$ such that
$$
e_\lambda y_\lambda y_n = e(\bi_{\lambda_-}\vee i_n) y_{\lambda_-} y_n^{b+1} = e_\nu y_\nu,
$$
and it is trivial that $\nu > \lambda$. Therefore $e_\lambda y_\lambda y_n \in \Blam$. If we have $b + 1 = b_{i_n}^{\lambda_-}$, by \autoref{I-problem: final} we have
$$
e_\lambda y_\lambda y_n = e(\bi_{\lambda_-}\vee i_n) y_{\lambda_-} y_n^{b+1} = e(\bi_{\lambda_-}\vee i_n) y_{\lambda_-} y_n^{b_{i_n}^{\lambda_-}} \in \Blam,
$$
which completes the proof. \endproof

\begin{Lemma} \label{psi-problem: induction}

Suppose $\lambda\in\mathscr S_n^\Lambda$. For any $1\leq r < n$, $e_\lambda y_\lambda \psi_r = \psi_{\t^\lambda \t^\lambda}\psi_r \in \Bgelam$.

\end{Lemma}

\proof Suppose $\t^\lambda{\cdot}s_r = \t$ is standard, we have $e_\lambda y_\lambda \psi_r = \psi_{\t^\lambda \t} \in \Bgelam$. So we consider the case that $\t^\lambda{\cdot}s_r = \t$ is not standard. If $r < n-1$, as $\lambda\in\mathscr S_n^\Lambda$, we have $\mu = \lambda_- \in \P_\psi$. Because $\t^\mu{\cdot}s_r = \t|_{n-1}$ which is not standard, $e_\mu y_\mu \psi_r \in \Blam[\mu]$. Then by \autoref{send}, we have $e_\lambda y_\lambda \psi_r \in \Blam$. If $r = n-1$ and write $\lambda = (\lambda^{(1)},\ldots,\lambda^{(l)},\emptyset,\ldots,\emptyset)$ with $\lambda^{(l)} = (\lambda^{(l)}_1,\ldots,\lambda^{(l)}_{k-1},\lambda^{(l)}_k)$, we must have either $\lambda^{(l)}_k \geq 2$ or $\lambda^{(l)}_{k-1} = \lambda^{(l)}_k = 1$. Then set $\mu = (\lambda^{(1)},\ldots,\mu^{(l)},\emptyset,\ldots,\emptyset)$ with
$$
\mu^{(l)} = \begin{cases}
(\lambda^{(l)}_1,\ldots,\lambda^{(l)}_{k-1}), & \text{ if $\lambda^{(l)}_k \geq 2$ and $k > 1$,}\\
(\lambda^{(l)}_1,\ldots,\lambda^{(l)}_{k-2}), & \text{ if $\lambda^{(l)}_{k-1} = \lambda^{(l)}_k = 1$ and $k > 2$,}\\
\emptyset, & \text{ otherwise.}
\end{cases}
$$

Suppose $i$ is the residue of node $(k,1,l)$ in $\lambda$ or the residue of node $(k-1,1,l)$ in $\lambda$, $\Lambda' = \Lambda_i$, $m = \lambda^{(l)}_k$ or $m = 2$ and $\gamma = (m)\in \mathscr P^{\Lambda'}_m$ or $\gamma = (1,1)\in \mathscr P^{\Lambda'}_m$ if $\lambda^{(l)}_k \geq 2$ or $\lambda^{(l)}_{k-1} = \lambda^{(l)}_k = 1$ respectively. Therefore $\lambda = \mu\vee \gamma$. Because $\lambda \in \mathscr S_n^\Lambda$, we have $\gamma \in \P_\psi$. Hence because $\t^\gamma {\cdot}s_{m-1}$ is not standard, we have $e_\gamma y_\gamma \psi_{m-1} \in \Blam[\gamma] = N_m^{\Lambda'}$. Then by \autoref{theta_j: go up},
$$
e_\lambda y_\lambda \psi_r = e(\bi_\mu \vee \bi_\gamma) y_\lambda \psi_r = \hat\theta_{\bi_\mu}(e_\gamma y_\gamma \psi_{m-1}) y_\mu \in \Blam,
$$
which completes the proof. \endproof

\begin{Corollary} \label{m_lambda is positive}
For $\lambda \in \mathscr S_n^\Lambda$, we have $m_\lambda > 0$.

\end{Corollary}

\proof Combining the above two Lemmas, \autoref{residue sequence has to be the same} and \autoref{basis to hecke}, the Corollary follows.\endproof

\subsection{Completion of the $y$-problem}

In this subsection we are going to prove that for any $\t\in\Std(\lambda)$, if $l(d(\t)) \leq m_\lambda$, then for any $1\leq r\leq n-1$ and any $\s\in\Std(\lambda)$, if $\t{\cdot}s_r$ is standard and $d(\t){\cdot}s_r$ is reduced, $\psi_{\s\t}\psi_r \in \Bgelam$ and for any $1\leq r\leq n$, $\psi_{\s\t}y_r \in \Bgelam$.

First we introduce the following Lemma.

\begin{Lemma} \label{length of psi}
Suppose $m$ is a positive integer such that $m \leq m_\lambda$, then
\begin{eqnarray*}
e_\lambda y_\lambda \psi_{r_1}\psi_{r_2}\ldots \psi_{r_m}
& =_\lambda &\sum_{\substack{\v \in \text{Std}(\lambda)\\ l(d(\v)) \leq m}} c_{\t^\lambda\v} \psi_{t^\lambda \v}.
\end{eqnarray*}
\end{Lemma}

\proof We apply induction on $m$. Suppose $m = 0$ then there is nothing to prove. Assume for any $m' < m$ the Lemma holds. Therefore $e_\lambda y_\lambda\psi_{r_1}\psi_{r_2}\ldots \psi_{r_{m-1}} =_\lambda \sum_{\substack{\u \in \text{Std}(\lambda)\\ l(d(\u)) \leq m-1}} c_{\t^\lambda \u} \psi_{t^\lambda \u}$ which yields
$$
e_\lambda y_\lambda\psi_{r_1}\psi_{r_2}\ldots \psi_{r_{m-1}}\psi_{r_m} =_\lambda\sum_{\substack{\u \in \text{Std}(\lambda)\\ l(d(\u)) \leq m-1}} c_{\t^\lambda \u} \psi_{t^\lambda \u}\psi_{r_m}.
$$

For $\u\in\Std(\lambda)$ and $l(d(\u)) \leq m-1 < m_\lambda$, if $\u{\cdot}s_r$ is standard and $s_{d(\v)}{\cdot}s_{r_m}$ is reduced, by the definition of $m_\lambda$,
$$
\psi_{\t^\lambda \u}\psi_{r_m} = \psi_{\t^\lambda, \u{\cdot}s_r} + \sum_{(\x,\y)\rhd(\t^\lambda,\u)} c_{\x \y} \psi_{\x \y} =_\lambda \psi_{\t^\lambda, \u{\cdot}s_r} + \sum_{\v \rhd \u} c_{\t^\lambda \v} \psi_{\t^\lambda \v},
$$
where $l(d(\u{\cdot}s_r)) = 1 + l(d(\u)) \leq m$ and $l(d(\v)) < l(d(\u)) < m$ as $\v \rhd \u$. Hence
$$
\psi_{\t^\lambda \u}\psi_{r_m} =_\lambda \sum_{\substack{\v \in \text{Std}(\lambda)\\ l(d(\v)) \leq m}} c_\v \psi_{t^\lambda \v}.
$$

If $\u{\cdot}s_r$ is not standard or $s_{d(\v)}{\cdot}s_{r_m}$ is not reduced, we have
$$
\psi_{\t^\lambda \u} \psi_{r_m} = \sum_{(\x,\y)\rhd(\t^\lambda,\u)} c_{\x \y} \psi_{\x \y} =_\lambda \sum_{\v \rhd \u} c_{\t^\lambda \v} \psi_{\t^\lambda \v},
$$
where $l(d(\v)) < l(d(\u)) \leq m-1 < m$ as $\v \rhd \u$. Hence
$$
\psi_{\t^\lambda \u}\psi_{r_m} =_\lambda \sum_{\substack{\v \in \text{Std}(\lambda)\\ l(d(\v)) \leq m}} c_{\t^\lambda\v} \psi_{t^\lambda \v},
$$
which completes the proof. \endproof

Now we can start to prove that $\psi_{\s\t}\psi_r\in \Bgelam$ when $\t{\cdot}s_r$ is standard and $d(\t){\cdot}s_r$ is reduced.

\begin{Lemma}\label{basis change}
Suppose $\t$ is a standard $\lambda$-tableau with $d(\t) = s_{r_1}s_{ r_2} \ldots s_{r_l}$ where $l \leq m_\lambda$, and $d'(t) = s_{t_1} s_{t_2} \ldots s_{t_l}$ is another reduced decomposition of $d(\t)$, then
$$
e_\lambda y_\lambda \psi_{d(\t)} - e_\lambda y_\lambda \psi_{d'(\t)} = \sum_{(\u,\v)\rhd (\t^\lambda \t)} c_{\u\v} \psi_{\u\v}.
$$
\end{Lemma}

\proof By \cite[Proposition 2.5]{BKW:GradedSpecht}, we have
$$
y_\lambda e_\lambda \psi_{d(\t)} - y_\lambda e_\lambda \psi_{d'(\t)} = \sum_{u < d(\t)}c_{u,f}y_\lambda e_\lambda f(y) \psi_u,
$$
where $f(y)$ is a polynomial in $y_r$'s and $u$ is a word in $\mathfrak S_n$. If $f(y)\neq 1$, by \autoref{y-problem: induction} we have $e_\lambda y_\lambda f(y)\in \Blam $. Hence $y_\lambda e_\lambda f(y) \psi_u \in \Blam $. If $f(y) = 1$, as $u <  d(\t)$ then $l(u) <l \leq m_\lambda$, by \autoref{length of psi} we have $e_\lambda y_\lambda \psi_u \in \Bgelam $. Henceforth
$$
y_\lambda e_\lambda \psi_w - y_\lambda e_\lambda \psi_{w'} \in \Bgelam.
$$

Then by \autoref{basis to hecke} and \cite[Lemma 5.7]{HuMathas:GradedCellular}, we complete the proof.\endproof

The following Corollary is straightforward by \autoref{basis change} which explains the action of $\psi_r$ to $\psi_{\s\t}$ from right when $\t{\cdot}s_r$ is standard and $d(\t){\cdot}s_r$ is reduced.

\begin{Corollary}\label{psi problem: part2}
Suppose $\t$ is a standard $\lambda$-tableau with $l(d(\t)) \leq m_\lambda$, if $\w = \t{\cdot}s_r$ is standard and $d(\t){\cdot}s_r$ is reduced,
$$
\psi_{\s\t} \psi_r = \psi_{\s\w} + \sum_{(\u,\v)\rhd (\s,\t)}c_{\u\v} \psi_{\u\v}.
$$
\end{Corollary}

Now we start to prove that $\psi_{\s\t}y_r \in \Bgelam$.

\begin{Lemma} \label{yup: help}
Suppose $\t$ is a standard $\lambda$-tableau with $l(d(\t)) < m_\lambda$. For any $1\leq k \leq n$, $1\leq r \leq n-1$ and any standard $\lambda$-tableau $\s$, we have
$$
\psi_{\s\t}y_k\psi_r \in \Bgelam.
$$

\proof As $l(d(\t)) < m_\lambda$, we have
$$
\psi_{\s \t}y_k = \sum_{(\u,\v)\rhd (\s,\t)}c_{\u\v} \psi_{\u\v} = \sum_{\v\rhd \t}c_{\s \v} \psi_{\s \v} + \sum_{\u,\v\in\Std(>\lambda)}c_{\u\v} \psi_{\u\v}.
$$

For $\Shape(\v) = \lambda$, since $\v \rhd \t$, $l(d(\v)) < l(d(\t)) < m_\lambda$. Then we have $\psi_{\s\v}\psi_{r} \in \Bgelam$.

For $\u,\v\in\Std(>\lambda)$, $\psi_{\u\v}\in \Blam $. As $\lambda \in \mathscr S_n^\Lambda$, $\Blam $ is an ideal by \autoref{B_lambda is an ideal}. Hence $\psi_{\u\v}\psi_{r_l} \in \Blam$ and this completes the proof. \endproof

\end{Lemma}

\begin{Proposition}\label{yup}
Suppose $\t$ is a standard $\lambda$-tableau with $l(d(\t)) \leq m_\lambda$, for any $1\leq r\leq n$ and any standard $\lambda$-tableau $\s$, we have
$$
\psi_{\s\t} y_r = \sum_{(\u,\v)\rhd (\s,\t)}c_{\u\v} \psi_{\u\v}.
$$
\end{Proposition}

\proof Write $d(\t) = s_{r_1}s_{r_2}\ldots s_{r_{l-1}}s_{r_l}$ and $\w = \t{\cdot}s_{r_l} = s_{r_1}s_{r_2}\ldots s_{r_{l-1}}$. We prove this Proposition by considering different values of $r$.

If $r\not\in\{r_l,r_l+1\}$, then $\psi_{r_l}$ and $y_r$ commute. Hence
$$
\psi_{\s \t} y_r =\psi_{d(\s)}^* e_\lambda y_\lambda \psi_{d(\t)}y_r = \psi_{d(\s)}^* e_\lambda y_\lambda \psi_{d(\w)}y_r\psi_{r_l} = \psi_{\s\w}y_r\psi_{r_l}.
$$

As $l(d(\w)) = l(d(\t)) -1 < m_\lambda$, by \autoref{yup: help} we have $\psi_{\s\t}y_r \in \Bgelam $.\\

If $r = r_l$, let $\bj$ be a sequence such that $e(\bi_\lambda)\psi_{d(\t)} = \psi_{d(\w)}e(\bj)\psi_{r_l}$. We separate this case further into $j_{r_l} \neq j_{r_l+1}$ and $j_{r_l} = j_{r_l+1}$. First suppose $j_{r_l} \neq j_{r_l+1}$, then
$$
\psi_{\s \t} y_r =\psi_{d(\s)}^* e_\lambda y_\lambda \psi_{d(\t)}y_r = \psi_{d(\s)}^* e_\lambda y_\lambda \psi_{d(\w)}y_{r+1}\psi_{r_l} = \psi_{\s\w}y_{r+1}\psi_{r_l}.
$$

Hence as $l(d(\w)) < m_\lambda$, by \autoref{yup: help} we have $\psi_{\s\t}y_r \in \Bgelam $ when $j_{r_l} \neq j_{r_l+1}$. Now suppose $j_{r_l} = j_{r_l+1}$, we have
$$
\psi_{\s\t}y_r = \psi_{d(\s)}^* e_\lambda y_\lambda \psi_{d(\t)}\y_r = \psi^*_{d(\s)}e_\lambda y_\lambda \psi_{d(\w)} + \psi_{d(\s)}^* e_\lambda y_\lambda \psi_{d(\w)}y_{r+1}\psi_{r_l} = \psi_{\s\w} + \psi_{\s\w}y_{r+1}\psi_{r_l}.
$$

As $l(d(\w)) < m_\lambda$, by \autoref{yup: help} we have $\psi_{\s\w}y_{r+1}\psi_{r_l} \in \Bgelam $. As $\psi_{\s\w}\in \Bgelam $ as well, we have $\psi_{\s\t}y_r \in \Bgelam $. So for $r = r_l$, we have $\psi_{\s\t}y_r \in \Bgelam$.\\

If $r = r_l + 1$, the method is the same as $r = r_l$.

Therefore in all the cases, we have $\psi_{\s\t}y_r \in \Bgelam $. So
$$
\psi_{\s\t}y_r = \sum_{\u,\v\in\Std(\geq\lambda)} c_{\u\v}\psi_{\u\v},
$$
and by \autoref{basis to hecke} we complete the proof.\endproof

\subsection{Properties of $m_\lambda$}

In the rest of this section we will prove that if $\t\in\Std(\lambda)$ and $l(d(\t)) \leq m_\lambda$, then for any $1\leq r\leq n-1$ and any $\s\in\Std(\lambda)$, we have $\psi_{\s\t} \psi_r\in\Bgelam$. In this subsection we will give some properties for $m_\lambda$ which will be used in proving the above argument.

\begin{Lemma} \label{does not change by psi}
Suppose $\lambda\in \mathscr{S}_n^\Lambda$. For any permutation $w \in \mathfrak{S}_n$ with reduced expression $w = s_{r_1}s_{r_2}\ldots s_{r_{m-1}}s_{r_m}$ and $r = \text{min}\{r_1,r_2,\ldots,r_m\}$, if we write
$$
e_\lambda y_\lambda \psi_w = e_\lambda y_\lambda \psi_{r_1}\psi_{r_2}\ldots \psi_{r_m} = \sum_{\u,\v\in\Std(\geq\lambda)} c_{\u \v} \psi_{\u \v},
$$
then $c_{\u \v} \neq 0$ implies $\v|_k \unrhd \t^{\lambda|_k}$ for any $k < r$.

\end{Lemma}

\proof We prove the Lemma by induction. If $m = 1$, then $r_1 = r$.
$$
e_\lambda y_\lambda \psi_r = \begin{cases}
\psi_{\t^\lambda \v}, & \text{ if $\v = \t^\lambda{\cdot}s_r$ is standard,}\\
\sum_{(\u,\v) \rhd (\t^\lambda,\t^\lambda)} c_{\u\v} \psi_{\u\v}, & \text{ otherwise.}
\end{cases}
$$

If $\v = \t^\lambda{\cdot}s_r$ is standard, then by the definition of $\v$, $\v|_k = \t^\lambda|_k = \t^{\lambda|_k}$ for $k < r$. If it is the other case, as $\v \rhd \t^\lambda$, then $\v|_k \unrhd \t^\lambda|_k = \t^{\lambda|_k}$.\\

Assume for any $m' < m$ the Corollary holds. Then
$$
e_\lambda y_\lambda \psi_{r_1}\psi_{r_2}\ldots \psi_{r_{m-1}} = \sum_{\u_1,\v_1 \in \Std(\geq\lambda)} c_{\u_1\v_1}\psi_{\u_1\v_1},
$$
where $\v_1|_k \unrhd \t^{\lambda|_k}$ for $k < r$. Then
$$
e_\lambda y_\lambda \psi_{r_1}\psi_{r_2}\ldots \psi_{r_{m-1}}\psi_{r_m} = \sum_{\u_1,\v_1 \in \Std(\geq\lambda)}c_{\u_1\v_1}\psi_{\u_1\v_1}\psi_{r_m}.
$$

Since
$$
\psi_{\u_1\v_1}\psi_{r_m} = \begin{cases}
\psi_{\u_1\v} + \sum_{(\u,\v_2) \rhd (\u_1,\v_1)} c_{\u\v_2}\psi_{\u\v_2}, & \text{ if $\v = \v_1{\cdot}s_{r_m}$ is standard and $d(\v_1){\cdot}s_{r_m}$ is reduced,}\\
\sum_{(\u,\v) \rhd (\u_1,\v_1)} c_{\u\v}\psi_{\u\v}, & \text{ otherwise.}
\end{cases}
$$

If $\v = \v_1{\cdot}s_{r_m}$ is standard and $d(\v_1){\cdot}s_{r_m}$ is reduced, recall $\v_1|_k \unrhd \t^{\lambda|_k}$ for $k < r$, as $\v = \v_1{\cdot}s_{r_m}$, $\v|_k = \v_1|_k \unrhd \t^{\lambda|_k}$ for $k < r \leq r_m$. For $\v_2 \rhd \v_1$, we have $\v_2|_k \rhd \v_1|_k \unrhd \t^{\lambda|_k}$ for $k < r$.

If it is of the other case, as $\v \rhd \v_1$, $\v|_k \rhd \v_1|_k \unrhd \t^{\lambda|_k}$. Therefore
$$
e_\lambda y_\lambda \psi_w = e_\lambda y_\lambda \psi_{r_1}\psi_{r_2}\ldots \psi_{r_m} = \sum_{\u,\v\in\Std(\geq\lambda)} c_{\u \v} \psi_{\u \v},
$$
and $c_{\u \v} \neq 0$ implies $\v|_k \unrhd \t^{\lambda|_k}$ for any $k < r$. This completes the proof. \endproof

\begin{Lemma} \label{already in P}

Suppose $\lambda \in \mathscr{S}_n^\Lambda \cap (\P_I \cap \P_y \cap \P_\psi)$. Then for any $1 \leq r_1,r_2,\ldots,r_m \leq n-1$
$$
e_\lambda y_\lambda \psi_{r_1}\psi_{r_2} \ldots \psi_{r_m} = \sum_{\v\in\Std(\lambda)} c_{\t^\lambda\v}\psi_{\t^\lambda\v} + \sum_{\substack{\u,\v\in\Std(>\lambda)\\\u \rhd \t^\lambda}} c_{\u\v}\psi_{\u\v}.
$$

\end{Lemma}

\proof When $m = 1$, we have
$$
e_\lambda y_\lambda \psi_{r_1} = \begin{cases}
\psi_{\t^\lambda \v}, & \text{ if $\v = \t^\lambda{\cdot}s_{r_1}$ is standard,}\\
\sum_{(\u,\v) \rhd (\t^\lambda, \t^\lambda)} c_{\u\v} \psi_{\u\v} = \sum_{\substack{\u,\v\in\Std(>\lambda)\\ \u \rhd \t^\lambda}} c_{\u\v}\psi_{\u\v}, & \text{ if $\v = \t^\lambda{\cdot}s_{r_1}$ is not standard.}
\end{cases}
$$
which follows the Lemma.

\hh
Suppose for $m' < m$ the Lemma holds. Then by induction
\begin{equation}\label{psi_n-1 is not involved: help1}
e_\lambda y_\lambda \psi_{r_1}\ldots\psi_{r_{m-1}} \psi_{r_m} = \sum_{\v_1 \in \Std(\lambda)} c_{\t^\lambda\v_1}\psi_{\t^\lambda\v_1}\psi_{r_m} + \sum_{\substack{\u_1,\v_1 \in \Std(>\lambda)\\ \u_1 \rhd \t^\lambda}} c_{\u_1\v_1}\psi_{\u_1\v_1}\psi_{r_m}.
\end{equation}

For $\v_1 \in \Std(\lambda)$, as $\lambda\in \P_\psi$,
$$
\psi_{\t^\lambda\v_1}\psi_{r_m} = \begin{cases}
\psi_{\t^\lambda \v_2} + \sum_{(\u_2,\v_2)\rhd (\t^\lambda,\v_1)}c_{\u_2\v_2}\psi_{\u_2\v_2}, & \text{ if $\v_2 = \v_1{\cdot}s_{r_m}$ is standard}\\
& \hspace*{10mm}\text{and $d(\v_2) = d(\v_1){\cdot}s_{r_m}$ is reduced},\\
\sum_{(\u_2,\v_2)\rhd (\t^\lambda,\v_1)}c_{\u_2\v_2}\psi_{\u_2\v_2}, & \text{ if $\v_2 = \v_1{\cdot}s_{r_m}$ is not standard}\\
& \hspace*{10mm}\text{or $d(\v_2) = d(\v_1){\cdot}s_{r_m}$ is not reduced}.
\end{cases}
$$
where in both cases, we can write
\begin{equation}\label{psi_n-1 is not involved: help2}
\psi_{\t^\lambda\v_1}\psi_{r_m} = \sum_{\v_2 \in \Std(\lambda)} c_{\t^\lambda\v_2}\psi_{\t^\lambda\v} + \sum_{\substack{\u_2,\v_2 \in \Std(>\lambda)\\ \u_2 \rhd \t^\lambda}} c_{\u_2\v_2}\psi_{\u_2\v_2}.
\end{equation}

\hh
For $\u_1,\u_2 \in \Std(>\lambda)$,
$$
\psi_{\u_1\v_1}\psi_{r_m} = \begin{cases}
\psi_{\u_1 \v_2} + \sum_{(\u_2,\v_2)\rhd (\u_1,\v_1)}c_{\u_2\v_2}\psi_{\u_2\v_2}, & \text{ if $\v_2 = \v_1{\cdot}s_{r_m}$ is standard}\\
& \hspace*{10mm}\text{and $d(\v_2) = d(\v_1){\cdot}s_{r_m}$ is reduced},\\
\sum_{(\u_2,\v_2)\rhd (\u_1,\v_1)}c_{\u_2\v_2}\psi_{\u_2\v_2}, & \text{ if $\v_2 = \v_1{\cdot}s_{r_m}$ is not standard}\\
&\hspace*{10mm}\text{or $d(\v_2) = d(\v_1){\cdot}s_{r_m}$ is not reduced}.
\end{cases}
$$
where since $\u_1 \rhd \t^\lambda$, we can always write
\begin{equation} \label{psi_n-1 is not involved: help3}
\psi_{\u_1\v_1}\psi_{r_m} = \sum_{\substack{\u_2,\v_2 \in \Std(>\lambda)\\\u_2 \rhd \t^\lambda}} c_{\u_2\v_2}\psi_{\u_2\v_2}.
\end{equation}

\hh
Therefore, substitute (\ref{psi_n-1 is not involved: help2}) and (\ref{psi_n-1 is not involved: help3}) back to (\ref{psi_n-1 is not involved: help1}), we have
$$
e_\lambda y_\lambda \psi_{r_1}\psi_{r_2} \ldots \psi_{r_m} = \sum_{\v\in\Std(\lambda)} c_{\t^\lambda\v}\psi_{\t^\lambda\v} + \sum_{\substack{\u,\v\in\Std(>\lambda)\\ \u \rhd \t^\lambda}} c_{\u\v}\psi_{\u\v},
$$
which completes the proof. \endproof

\begin{Lemma} \label{psi_n-1 is not involved}

Suppose $\lambda \in \mathscr S_n^\Lambda$ and $r_1,r_2,\ldots,r_m$ are positive integers such that $r_1,\ldots,r_m < n-1$. Then
$$
e_\lambda y_\lambda \psi_{r_1}\psi_{r_2} \ldots \psi_{r_m} \in \Bgelam.
$$

\end{Lemma}

\proof Define $\mu = \lambda|_{n-1}$. As $\lambda \in \mathscr S_n^\Lambda$, $\mu \in \mathscr S_{n-1}^\Lambda\cap(\P_I \cap \P_y \cap \P_\psi)$. Define $i \in I$ such that $\bi_\lambda = \bi_\mu\vee i$. As $r_1,r_2,\ldots,r_m < n-1$, we have
$$
e_\lambda y_\lambda \psi_{r_1}\psi_{r_2}\ldots \psi_{r_m} = \theta_i(e_\mu y_\mu \psi_{r_1}\psi_{r_2} \ldots \psi_{r_m}),
$$
where
$$
e_\mu y_\mu \psi_{r_1}\psi_{r_2} \ldots \psi_{r_m} = \sum_{\dot\v\in\Std(\mu) } c_{\t^\mu \dot\v} \psi_{\t^\mu \dot\v} + \sum_{\dot\u,\dot\v\in\Std(>\mu)} c_{\dot\u\dot\v} \psi_{\dot\u\dot\v}.
$$

As $\sum_{\dot\u,\dot\v\in\Std(>\mu)} c_{\dot\u\dot\v} \psi_{\dot\u\dot\v} \in\Blam[\mu] = \Blam[\lambda|_{n-1}]$, by \autoref{send}, $\theta_i(\sum_{\dot\u,\dot\v\in\Std(>\mu)} c_{\dot\u\dot\v} \psi_{\dot\u\dot\v}) \in \Blam $.

For $\dot\v \in \Std(\mu) = \Std(\lambda|_{n-1})$ and $\bi_\mu\vee i = \bi_\lambda$, define $\v$ to be the standard $\lambda$-tableau with $\v|_{n-1} = \dot\v$. Hence $\theta_i(\psi_{\t^\mu \dot\v}) = \psi_{\t^\lambda \v}$. Therefore
$$
\theta_i(\sum_{\dot\v\in\Std(\mu) } c_{\t^\mu \dot\v} \psi_{\t^\mu \dot\v}) = \sum_{\v\in\Std(\lambda) } c_{\t^\mu \dot\v} \psi_{\t^\lambda \v} \in \Bgelam.
$$

So
$$
e_\mu y_\mu \psi_{r_1}\psi_{r_2} \ldots \psi_{r_m} = \theta_i(\sum_{\dot\v\in\Std(\mu) } c_{\t^\mu \dot\v} \psi_{\t^\mu \dot\v}) + \theta_i(\sum_{\dot\u,\dot\v\in\Std(>\mu)} c_{\dot\u\dot\v} \psi_{\dot\u\dot\v})  \in \Bgelam.
$$
\endproof

\subsection{Garnir tableaux}

In the following subsections we will prove that $\psi_{\s\t}\psi_r \in \Bgelam$ for $l(d(\t)) \leq m_\lambda$. Generally, if $\t{\cdot}s_r$ is standard and $d(\t){\cdot}s_r$ is reduced or $l(d(\t)){\cdot}s_r$ is not reduced, it is comparatively easy to prove that $\psi_{\s\t}\psi_r\in\Bgelam$. Our main difficulty is to prove that when $\t{\cdot}s_r$ is not standard then $\psi_{\s\t}\psi_r\in\Bgelam$. In order to prove this we consider different types of $\t$. Among these cases the hardest part is that when $\t$ is a special kind of tableaux which is called the Garnir tableau and $\t{\cdot}s_r$ is not standard. In this subsection we will prove that in such case $\psi_{\s\t}\psi_r \in \Bgelam$.

The method of proving the argument in this subsection is first assuming that $\Shape(\t)$ is a partition of two rows, and using the similar argument we used in the last section to extend the result to general multipartitions. First we give a detailed definition of garnir tableaux.

We introduce a special kind of tableaux, the Garnir tableaux, which was first introduced by Murphy~\cite{Murphy:basis}. Let $(a,b,m)$ be a node of $\lambda$ such that $(a+1,b,m)$ is also a node of $\lambda$. The \textbf{$(a,b,m)$-Garnir belt} of $\lambda$ consists of the nodes $(a,c,m)$ for $b \leq c\leq \lambda_a^{(m)}$ and the nodes $(a+1,g,m)$ for $1\leq g\leq b$. For example here is a picture of the $(2,3,2)$-Garnir belt for $\lambda = (3,1|7,6,5,2)$.
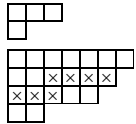
\begin{eqnarray*}
&&\tab(\ \ \ ,\ )\\
&&\tab(\ \ \ \ \ \ \ ,\ \ \times\times\times\times,\times\times\times\ \ ,\ \ )
\end{eqnarray*}

The \textbf{$(a,b,m)$-Garnir tableau} of shape $\lambda$ is the unique maximal standard $\lambda$-tableau with respect to the Bruhat order ($\rhd$) among the standard $\lambda$-tableaux which agree with $\t^\lambda$ outside the $(a,b,m)$-Garnir belt. For example the following is the $(2,3,2)$-Garnir tableau for $\lambda = (3,1|7,6,5,2)$.
$$
\Bigg(\hspace*{1mm}\tab(123,4)\hspace*{1mm}\Bigg|\hspace*{1mm}
          \tab(56789\ten\eleven,\twelve\thirteen\sixteen\eighteen\nineteen\twenty,\fourteen\fifteen\seventeen\twentyone\twentytwo,\twentythree,\twentyfour)\hspace*{1mm}\Bigg)
$$

Suppose $\lambda = (\lambda^{(1)},\ldots,\lambda^{(\l)})$ and $\lambda^{(\l)} = (\lambda_1^{(\l)},\ldots,\lambda_k^{(\l)})$. Let $(a,b,m) = (k-1,\lambda_k^{(\l)},\l)$ and $\t$ be the $(a,b,m)$-Garnir tableau. Let the entry in node $(a,b,m)$ of $\t$ be $r$. Then $\t{\cdot}s_r$ is not standard.

\begin{Definition}

Suppose $\lambda \in \P_n$ with $\lambda = (\lambda^{(1)},\ldots,\lambda^{(\l)})$ and $\lambda^{(\l)} = (\lambda^{(\l)}_1,\ldots,\lambda^{(\l)}_k)$. If $k \geq 2$, and $\t$ is the $(k-1,\lambda^{(\l)}_k,\l)$-Garnir tableau, then we call $\t$ the \textbf{last Garnir tableau} of shape $\lambda$, and $r = \t(k-1,\lambda^{(\l)}_k,\l)$ the \textbf{last Garnir entry} of $\t$.

\end{Definition}

For example
$$
\Bigg(\hspace*{1mm}\tab(123,4)\hspace*{1mm}\Bigg|\hspace*{1mm}
          \tab(56789\ten\eleven,\twelve\thirteen\sixteen\eighteen\nineteen\twenty,\fourteen\fifteen\seventeen)\hspace*{1mm}\Bigg)
$$
is the last $(2,3,2)$-Garnir tableau, and
$$
\Bigg(\hspace*{1mm}\tab(123,4)\hspace*{1mm}\Bigg|\hspace*{1mm}
          \tab(56789\ten\eleven,\twelve\thirteen\sixteen\eighteen\nineteen\twenty,\fourteen\fifteen\seventeen\twentyone\twentytwo,\twentythree\twentyfour)\hspace*{1mm}\Bigg)
$$
is not the last one. Notice that $\t{\cdot}s_r$ is not standard.\\

Because we are going to play around with $\psi_{d(\t)}$ a lot, we introduce more detailed notation for these elements in the next Lemma.

\begin{Lemma} \label{standard expression of garnir tableau}
Suppose $\lambda = (\lambda^{(1)},\ldots,\lambda^{(\l)})$. Let $\t$ be a $(a,b,m)$-Garnir tableau of shape $\lambda$ and $\lambda^{(m)} = (\lambda_1^{(m)},\ldots,\lambda_k^{(m)})$. Suppose
$$
\begin{cases}
\t^\lambda(a,b,m) & = l,\\
\t^\lambda(a,\lambda_a^{(m)},m) & = s,\\
\t^\lambda(a+1,b,m) & = t.
\end{cases}
$$

Then $l \leq s < t$. Write $t - s = c$,
$$
\psi_{d(\t)} = \psi_s\psi_{s+2} \ldots \psi_{t-1}{\cdot} \psi_{s-1}\psi_s\ldots \psi_{t-2}{\cdot}\ldots{\cdot}\psi_{l+1}\psi_{l+2}\ldots\psi_{l+c}{\cdot}\psi_l\psi_{l+1}\ldots\psi_{l+c-2}
$$
where
$$
l(\psi_s\psi_{s+1}\ldots \psi_{t-1}) = l(\psi_{s-1}\psi_{s}\ldots\psi_{t-2}) = \ldots = l(\psi_{l+1}\psi_{l+2}\ldots\psi_{l+c}) = c
$$
and
$$
l(\psi_l\psi_{l+1}\ldots\psi_{l+c-2}) = c-1.
$$
\end{Lemma}

\proof The Lemma follows by direct calculation. \endproof

\begin{Example}
Suppose $\lambda = (3,1|7,6,5,2)$ and $(a,b,m) = (2,3,2)$. Let $\t$ be the $(2,3,2)$-Garnir tableau of shape $\lambda$. Then
$$
\t =
\Bigg(\hspace*{1mm}\tab(123,4)\hspace*{1mm}\Bigg|\hspace*{1mm}
          \tab(56789\ten\eleven,\twelve\thirteen\sixteen\eighteen\nineteen\twenty,\fourteen\fifteen\seventeen\twentyone\twentytwo,\twentythree\twentyfour)\hspace*{1mm}\Bigg),
$$
and
$$
\t^\lambda =
\Bigg(\hspace*{1mm}\tab(123,4)\hspace*{1mm}\Bigg|\hspace*{1mm}
          \tab(56789\ten\eleven,\twelve\thirteen\fourteen\fifteen\sixteen\seventeen,\eighteen\nineteen\twenty\twentyone\twentytwo,\twentythree\twentyfour)\hspace*{1mm}\Bigg).
$$

Then
$$
\begin{array}{lll}
\t^\lambda(a,b,m) & = \t^\lambda(2,3,2) & = 14,\\
\t^\lambda(a,\lambda_a^{(m)},m) & = \t^\lambda(2,6,2) & = 17,\\
\t^\lambda(a+1,b,m) & = \t^\lambda(3,3,2) & = 20,
\end{array}
$$
and $c = t - s = 3$. Therefore
$$
\psi_{d(\t)} = \psi_{17}\psi_{18}\psi_{19}{\cdot}\psi_{16}\psi_{17}\psi_{18}{\cdot}\psi_{15}\psi_{16}\psi_{17}{\cdot}\psi_{14}\psi_{15}
$$
with
$$
l(\psi_{17}\psi_{18}\psi_{19}) = l(\psi_{16}\psi_{17}\psi_{18}) = l(\psi_{15}\psi_{16}\psi_{17}) = 3 = c
$$
and
$$
l(\psi_{14}\psi_{15}) = 2 = c-1.
$$
\end{Example}

\begin{Remark}
For $a \leq b-1$, we will write $\psi_{a,b} = \psi_a \psi_{a+1} \psi_{a+2} \ldots \psi_{b-2} \psi_{b-1}$ and $\psi_{b,a} = \psi^*_{a,b}$ in order to simplify our notations.
\end{Remark}

Our first step is to prove that when $\lambda$ is a partition with two rows and $\t$ is the last Garnir tableau of shape $\lambda$ with $r$ as its last Garnir entry, then $\psi_{\s\t}\psi_r \in \Bgelam$ for any $\s\in\Std(\lambda)$. We set $\lambda = (\lambda_1,\lambda_2)$ and without loss of generality, set $\Lambda = \Lambda_0$. Therefore $\lambda\in \P_n$ with $n = \lambda_1 + \lambda_2$. Also we set $\mu = (\lambda_1,\lambda_2-1,1)$, $\dot\lambda = (\lambda_1-1,\lambda_2)$ and $\dot\mu = (\lambda_1-1,\lambda_2-1,1)$. Furthermore, let $i = \res(\gamma_1)$, $j = \res(\gamma_2)$, where $\gamma_1 = (1,\lambda_1,1)$ and $\gamma_2 = (2,\lambda_2,1)$.

First we prove a few useful Lemmas.
\begin{Lemma}\label{swap}
Suppose $\lambda$, $\dot\lambda$, $\dot\mu$, $i$ and $j$ are defined as above, we have
$$
e_\lambda y_\lambda \psi_{\lambda_1}\psi_{\lambda_1+1}\ldots \psi_{n-2}\psi_{n-1} =_\lambda \begin{cases}
\psi_{\lambda_1}\psi_{\lambda_1+1}\ldots\psi_{n-2}\psi_{n-1} e(\bi_{\dot{\lambda}}\vee i)y_{\dot{\lambda}}y_n\\
\hspace*{14mm} - \psi_{\lambda_1+1}\ldots\psi_{n-2}\psi_{n-1} e(\bi_{\dot{\lambda}}\vee i)y_{\dot{\lambda}}, & \text{ if $i = e-1$, $j\neq e-1$,}\\
\psi_{\lambda_1}\psi_{\lambda_1+1}\ldots\psi_{n-2}\psi_{n-1} e(\bi_{\dot{\lambda}}\vee i)y_{\dot{\lambda}}y_n\\
\hspace*{14mm} - \psi_{\lambda_1+1}\ldots\psi_{n-2}\psi_{n-1} e(\bi_{\dot{\lambda}}\vee i)y_{\dot{\lambda}}\\
\hspace*{28mm} - \psi_{\lambda_1}\ldots\psi_{n-2} e(\bi_{\dot{\lambda}}\vee i)y_{\dot{\lambda}}, & \text{ if $i = j = e-1$,}\\
\psi_{\lambda_1}\psi_{\lambda_1+1}\ldots\psi_{n-2}\psi_{n-1} e(\bi_{\dot{\lambda}}\vee i)y_{\dot{\lambda}}\\
\hspace*{14mm} + \psi_{\lambda_1}\psi_{\lambda_1+1}\ldots\psi_{n-2} e(\bi_{\dot{\mu}}\vee i)y_{\dot{\mu}}, & \text{ if $i = j = e-2$,}\\
\psi_{\lambda_1}\psi_{\lambda_1+1}\ldots\psi_{n-2}\psi_{n-1} e(\bi_{\dot{\lambda}}\vee i)y_{\dot{\lambda}}, & \text{ otherwise.}
\end{cases}
$$
\end{Lemma}

\proof In order to make the notation and diagrams clearer we set $e = 4$. For the other choices of $e$ one can check that the argument is the same because the proof does not depend on the value of $e$.

By using the diagrammatic notation we have
$$
e_\lambda y_\lambda \psi_{\lambda_1}\psi_{\lambda_1+1}\ldots \psi_{n-2}\psi_{n-1}
=
\begin{braid}
 \draw(3,-0.5)node{$\underbrace{\hspace*{24mm}}_{\lambda_1-1}$};
 \draw(10.5,-0.5)node{$\underbrace{\hspace*{27mm}}_{\lambda_2-1}$};
 \draw(3.5,5.2)node{$\overbrace{\hspace*{27mm}}^{\lambda_1}$};
 \draw(12,5.2)node{$\overbrace{\hspace*{31mm}}^{\lambda_2}$};
 \draw (0,4)node[above]{$0$}--(0,0);
 \draw (1,4)node[above]{$1$}--(1,0);
 \draw (2,4)node[above]{$2$}--(2,0);
 \draw (3,4)node[above]{$3$}--(3,0);
 \draw (4,4)node[above]{$0$}--(4,0);
 \draw[dots] (4.2,4)--(5.8,4);
 \draw[dots] (4.2,0)--(5.8,0);
 \draw (6,4)node[above]{$i-1$}--(6,0);
 \draw (7,4)node[above]{$i$}--(16,0);
 \draw (8,4)node[above]{$3$}--(7,0);
 \draw (9,4)node[above]{$0$}--(8,0);
 \draw (10,4)node[above]{$1$}--(9,0);
 \draw (11,4)node[above]{$2$}--(10,0);
 \draw (12,4)node[above]{$3$}--(11,0);
 \draw (13,4)node[above]{$0$}--(12,0);
 \draw[dots] (13.2,4)--(14.8,4);
 \draw[dots] (12.2,0)--(13.8,0);
 \draw (15,4)node[above]{$j-1$}--(14,0);
 \draw (16,4)node[above]{$j$}--(15,0);
 \node[greendot] at (3,4){};
 \node[greendot] at (7,4){};
 \node[below,scale = 1] at (7,4){$\delta_{i,3}$};
 \node[greendot] at (11,4){};
 \node[greendot] at (16,4){};
 \node[below right,scale = 1] at (16,4){$\delta_{j,2}$};
\end{braid}.
$$

Therefore, the Lemma is equivalent to moving all of the dots from the top to the bottom of the diagram. In order to do this we have to consider several cases.

\textbf{Case ~\ref{swap}a:} $i \neq e-2,e-1$.

Because we set $e = 4$ so in this case we have $i \neq 2,3$. As $i \neq 3$, $\delta_{i,3} = 0$. Therefore there are no dots on the strand labelled by $i$. And as $i \neq 2$, by relation~\ref{dia:y-psi com}, we have
\begin{eqnarray*}
&&e_\lambda y_\lambda \psi_{\lambda_1}\psi_{\lambda_1+1}\ldots \psi_{n-2}\psi_{n-1}\\
& = &
\begin{braid}
 \draw(3.5,5.2)node{$\overbrace{\hspace*{27mm}}^{\lambda_1}$};
 \draw (0,4)node[above]{$0$}--(0,0);
 \draw (1,4)node[above]{$1$}--(1,0);
 \draw (2,4)node[above]{$2$}--(2,0);
 \draw (3,4)node[above]{$3$}--(3,0);
 \draw (4,4)node[above]{$0$}--(4,0);
 \draw[dots] (4.2,4)--(5.8,4);
 \draw[dots] (4.2,0)--(5.8,0);
 \draw (6,4)node[above]{$i-1$}--(6,0);
 \draw (7,4)node[above]{$i$}--(16,0);
 \draw (8,4)node[above]{$3$}--(7,0);
 \draw (9,4)node[above]{$0$}--(8,0);
 \draw (10,4)node[above]{$1$}--(9,0);
 \draw (11,4)[densely dotted] node[above]{$2$}--(10,0);
 \draw (12,4)node[above]{$3$}--(11,0);
 \draw (13,4)node[above]{$0$}--(12,0);
 \draw[dots] (13.2,4)--(14.8,4);
 \draw[dots] (12.2,0)--(13.8,0);
 \draw (15,4)node[above]{$j-1$}--(14,0);
 \draw (16,4)[densely dotted] node[above]{$j$}--(15,0);
 \node[greendot] at (3,4){};
 \node[greendot] at (11,4){};
 \node[greendot] at (16,4){};
 \node[below right] at (16,4){$\delta_{j,2}$};
 \draw[->] (11,4) -- (10.5,2);
 \draw[->] (16,4) -- (15.02,0.08);
\end{braid}
 = 
\begin{braid}
 \draw (0,4)node[above]{$0$}--(0,0);
 \draw (1,4)node[above]{$1$}--(1,0);
 \draw (2,4)node[above]{$2$}--(2,0);
 \draw (3,4)node[above]{$3$}--(3,0);
 \draw (4,4)node[above]{$0$}--(4,0);
 \draw[dots] (4.2,4)--(5.8,4);
 \draw[dots] (4.2,0)--(5.8,0);
 \draw (6,4)node[above]{$i-1$}--(6,0);
 \draw (7,4)node[above]{$i$}--(16,0);
 \draw (8,4)node[above]{$3$}--(7,0);
 \draw (9,4)node[above]{$0$}--(8,0);
 \draw (10,4)node[above]{$1$}--(9,0);
 \draw (11,4)node[above]{$2$}--(10,0);
 \draw (12,4)node[above]{$3$}--(11,0);
 \draw (13,4)node[above]{$0$}--(12,0);
 \draw[dots] (13.2,4)--(14.8,4);
 \draw[dots] (12.2,0)--(13.8,0);
 \draw (15,4)node[above]{$j-1$}--(14,0);
 \draw (16,4)node[above]{$j$}--(15,0);
 \node[greendot] at (3,0){};
 \node[greendot] at (10,0){};
 \node[greendot] at (15,0){};
 \node[below] at (15,0){$\delta_{j,2}$};
\end{braid}\\
& = &
\psi_{\lambda_1}\psi_{\lambda_1+1}\ldots\psi_{n-2}\psi_{n-1} e(\bi_{\dot{\lambda}}\vee i)y_{\dot{\lambda}}.
\end{eqnarray*}

\textbf{Case ~\ref{swap}b:} $i = e-1$ and $j \neq e-1$.

Because we set $e = 4$ so in this case we have $i = 3$ and $j \neq 3$. Then $\delta_{i,3} = 1$. Hence
\begin{eqnarray*}
&&e_\lambda y_\lambda \psi_{\lambda_1}\psi_{\lambda_1+1}\ldots \psi_{n-2}\psi_{n-1}
 =
\begin{braid}
 \draw(3.5,5.2)node{$\overbrace{\hspace*{27mm}}^{\lambda_1}$};
 \draw (0,4)node[above]{$0$}--(0,0);
 \draw (1,4)node[above]{$1$}--(1,0);
 \draw (2,4)node[above]{$2$}--(2,0);
 \draw (3,4)node[above]{$3$}--(3,0);
 \draw (4,4)node[above]{$0$}--(4,0);
 \draw[dots] (4.2,4)--(5.8,4);
 \draw[dots] (4.2,0)--(5.8,0);
 \draw (6,4)node[above]{$2$}--(6,0);
 \draw (7,4)[densely dotted] node[above]{$3$}--(20,0);
 \draw (8,4)node[above]{$3$}--(7,0);
 \draw (9,4)node[above]{$0$}--(8,0);
 \draw (10,4)node[above]{$1$}--(9,0);
 \draw (11,4)node[above]{$2$}--(10,0);
 \draw (12,4)node[above]{$3$}--(11,0);
 \draw (13,4)node[above]{$0$}--(12,0);
 \draw (14,4)node[above]{$1$}--(13,0);
 \draw (15,4)node[above]{$2$}--(14,0);
 \draw (16,4)node[above]{$3$}--(15,0);
 \draw (17,4)node[above]{$0$}--(16,0);
 \draw[dots] (17.2,4)--(18.8,4);
 \draw[dots] (16.2,0)--(17.8,0);
 \draw (19,4)node[above]{$j-1$}--(18,0);
 \draw (20,4)node[above]{$j$}--(19,0);
 \node[greendot] at (3,4){};
 \node[greendot] at (7,4){};
 \node[greendot] at (11,4){};
 \node[greendot] at (15,4){};
 \node[greendot] at (20,4){};
 \node[below right] at (20,4){$\delta_{j,2}$};
 \draw[->] (7,4) -- (19.5, 4-3.815);
\end{braid}\\
&\hspace*{-5mm} \overset{(\ref{dia:y-psi com})}=&
-\begin{braid}
 \draw(3.5,5.2)node{$\overbrace{\hspace*{27mm}}^{\lambda_1}$};
 \draw (0,4)node[above]{$0$}--(0,0);
 \draw (1,4)node[above]{$1$}--(1,0);
 \draw (2,4)node[above]{$2$}--(2,0);
 \draw (3,4)node[above]{$3$}--(3,0);
 \draw (4,4)node[above]{$0$}--(4,0);
 \draw[dots] (4.2,4)--(5.8,4);
 \draw[dots] (4.2,0)--(5.8,0);
 \draw (6,4)node[above]{$2$}--(6,0);
 \draw (8,4)node[above]{$3$}--(20,0);
 \draw (7,4)node[above]{$3$}--(7,0);
 \draw (9,4)node[above]{$0$}--(8,0);
 \draw (10,4)node[above]{$1$}--(9,0);
 \draw (11,4)node[above]{$2$}--(10,0);
 \draw (12,4)node[above]{$3$}--(11,0);
 \draw (13,4)node[above]{$0$}--(12,0);
 \draw (14,4)node[above]{$1$}--(13,0);
 \draw (15,4)node[above]{$2$}--(14,0);
 \draw (16,4)node[above]{$3$}--(15,0);
 \draw (17,4)node[above]{$0$}--(16,0);
 \draw[dots] (17.2,4)--(18.8,4);
 \draw[dots] (16.2,0)--(17.8,0);
 \draw (19,4)node[above]{$j-1$}--(18,0);
 \draw (20,4)node[above]{$j$}--(19,0);
 \node[greendot] at (3,4){};
 \node[greendot] at (11,4){};
 \node[greendot] at (15,4){};
 \node[greendot] at (20,4){};
 \node[below right] at (20,4){$\delta_{j,2}$};
\end{braid}
 -
\begin{braid}
 \draw(3.5,5.2)node{$\overbrace{\hspace*{27mm}}^{\lambda_1}$};
 \draw (0,4)node[above]{$0$}--(0,0);
 \draw (1,4)node[above]{$1$}--(1,0);
 \draw (2,4)node[above]{$2$}--(2,0);
 \draw (3,4)node[above]{$3$}--(3,0);
 \draw (4,4)node[above]{$0$}--(4,0);
 \draw[dots] (4.2,4)--(5.8,4);
 \draw[dots] (4.2,0)--(5.8,0);
 \draw (6,4)node[above]{$2$}--(6,0);
 \draw (12,4)node[above]{$3$}--(20,0);
 \draw (8,4)node[above]{$3$}--(7,0);
 \draw (9,4)node[above]{$0$}--(8,0);
 \draw (10,4)node[above]{$1$}--(9,0);
 \draw (11,4)node[above]{$2$}--(10,0);
 \draw (7,4)node[above]{$3$}--(11,0);
 \draw (13,4)node[above]{$0$}--(12,0);
 \draw (14,4)node[above]{$1$}--(13,0);
 \draw (15,4)node[above]{$2$}--(14,0);
 \draw (16,4)node[above]{$3$}--(15,0);
 \draw (17,4)node[above]{$0$}--(16,0);
 \draw[dots] (17.2,4)--(18.8,4);
 \draw[dots] (16.2,0)--(17.8,0);
 \draw (19,4)node[above]{$j-1$}--(18,0);
 \draw (20,4)node[above]{$j$}--(19,0);
 \node[greendot] at (3,4){};
 \node[greendot] at (11,4){};
 \node[greendot] at (15,4){};
 \node[greendot] at (20,4){};
 \node[below right] at (20,4){$\delta_{j,2}$};
\end{braid}\\
&\hspace*{-5mm} - &
\begin{braid}
 \draw(3.5,5.2)node{$\overbrace{\hspace*{27mm}}^{\lambda_1}$};
 \draw (0,4)node[above]{$0$}--(0,0);
 \draw (1,4)node[above]{$1$}--(1,0);
 \draw (2,4)node[above]{$2$}--(2,0);
 \draw (3,4)node[above]{$3$}--(3,0);
 \draw (4,4)node[above]{$0$}--(4,0);
 \draw[dots] (4.2,4)--(5.8,4);
 \draw[dots] (4.2,0)--(5.8,0);
 \draw (6,4)node[above]{$2$}--(6,0);
 \draw (16,4)node[above]{$3$}--(20,0);
 \draw (8,4)node[above]{$3$}--(7,0);
 \draw (9,4)node[above]{$0$}--(8,0);
 \draw (10,4)node[above]{$1$}--(9,0);
 \draw (11,4)node[above]{$2$}--(10,0);
 \draw (12,4)node[above]{$3$}--(11,0);
 \draw (13,4)node[above]{$0$}--(12,0);
 \draw (14,4)node[above]{$1$}--(13,0);
 \draw (15,4)node[above]{$2$}--(14,0);
 \draw (7,4)node[above]{$3$}--(15,0);
 \draw (17,4)node[above]{$0$}--(16,0);
 \draw[dots] (17.2,4)--(18.8,4);
 \draw[dots] (16.2,0)--(17.8,0);
 \draw (19,4)node[above]{$j-1$}--(18,0);
 \draw (20,4)node[above]{$j$}--(19,0);
 \node[greendot] at (3,4){};
 \node[greendot] at (11,4){};
 \node[greendot] at (15,4){};
 \node[greendot] at (20,4){};
 \node[below right] at (20,4){$\delta_{j,2}$};
\end{braid}
 -  \ldots  +
\begin{braid}
 \draw(3.5,5.2)node{$\overbrace{\hspace*{27mm}}^{\lambda_1}$};
 \draw (0,4)node[above]{$0$}--(0,0);
 \draw (1,4)node[above]{$1$}--(1,0);
 \draw (2,4)node[above]{$2$}--(2,0);
 \draw (3,4)node[above]{$3$}--(3,0);
 \draw (4,4)node[above]{$0$}--(4,0);
 \draw[dots] (4.2,4)--(5.8,4);
 \draw[dots] (4.2,0)--(5.8,0);
 \draw (6,4)node[above]{$2$}--(6,0);
 \draw (7,4)node[above]{$3$}--(20,0);
 \draw (8,4)node[above]{$3$}--(7,0);
 \draw (9,4)node[above]{$0$}--(8,0);
 \draw (10,4)node[above]{$1$}--(9,0);
 \draw (11,4)node[above]{$2$}--(10,0);
 \draw (12,4)node[above]{$3$}--(11,0);
 \draw (13,4)node[above]{$0$}--(12,0);
 \draw (14,4)node[above]{$1$}--(13,0);
 \draw (15,4)node[above]{$2$}--(14,0);
 \draw (16,4)node[above]{$3$}--(15,0);
 \draw (17,4)node[above]{$0$}--(16,0);
 \draw[dots] (17.2,4)--(18.8,4);
 \draw[dots] (16.2,0)--(17.8,0);
 \draw (19,4)node[above]{$j-1$}--(18,0);
 \draw (20,4)node[above]{$j$}--(19,0);
 \node[greendot] at (3,4){};
 \node[greendot] at (20,0){};
 \node[greendot] at (11,4){};
 \node[greendot] at (15,4){};
 \node[greendot] at (20,4){};
 \node[below right] at (20,4){$\delta_{j,2}$};
\end{braid}\\
& \hspace*{-5mm} \overset{(\ref{dia:ii3})}= &
-\begin{braid}
 \draw (0,4)node[above]{$0$}--(0,0);
 \draw (1,4)node[above]{$1$}--(1,0);
 \draw (2,4)node[above]{$2$}--(2,0);
 \draw (3,4)node[above]{$3$}--(3,0);
 \draw (4,4)node[above]{$0$}--(4,0);
 \draw[dots] (4.2,4)--(5.8,4);
 \draw[dots] (4.2,0)--(5.8,0);
 \draw (6,4)node[above]{$2$}--(6,0);
 \draw (8,4)node[above]{$3$}--(20,0);
 \draw (7,4)node[above]{$3$}--(7,0);
 \draw (9,4)node[above]{$0$}--(8,0);
 \draw (10,4)node[above]{$1$}--(9,0);
 \draw (11,4)node[above]{$2$}--(10,0);
 \draw (12,4)node[above]{$3$}--(11,0);
 \draw (13,4)node[above]{$0$}--(12,0);
 \draw (14,4)node[above]{$1$}--(13,0);
 \draw (15,4)node[above]{$2$}--(14,0);
 \draw (16,4)node[above]{$3$}--(15,0);
 \draw (17,4)node[above]{$0$}--(16,0);
 \draw[dots] (17.2,4)--(18.8,4);
 \draw[dots] (16.2,0)--(17.8,0);
 \draw (19,4)node[above]{$j-1$}--(18,0);
 \draw (20,4)node[above]{$j$}--(19,0);
 \node[greendot] at (3,4){};
 \node[greendot] at (11,4){};
 \node[greendot] at (15,4){};
 \node[greendot] at (20,4){};
 \node[below right] at (20,4){$\delta_{j,2}$};
\end{braid} +
\begin{braid}
 \draw (0,4)node[above]{$0$}--(0,0);
 \draw (1,4)node[above]{$1$}--(1,0);
 \draw (2,4)node[above]{$2$}--(2,0);
 \draw (3,4)node[above]{$3$}--(3,0);
 \draw (4,4)node[above]{$0$}--(4,0);
 \draw[dots] (4.2,4)--(5.8,4);
 \draw[dots] (4.2,0)--(5.8,0);
 \draw (6,4)node[above]{$2$}--(6,0);
 \draw (7,4)node[above]{$3$}--(8,3)--(8,1)--(7,0);
 \draw (8,4)node[above]{$3$}--(7,3)--(7,1)--(11,0);
 \draw (9,4)node[above]{$0$}--(9,1)--(8,0);
 \draw (10,4)node[above]{$1$}--(10,1)--(9,0);
 \draw (11,4)node[above]{$2$}--(11,1)--(10,0);
 \draw (12,4)node[above]{$3$}--(20,0);
 \draw (13,4)node[above]{$0$}--(12,3)--(12,0);
 \draw (14,4)node[above]{$1$}--(13,0);
 \draw (15,4)node[above]{$2$}--(14,0);
 \draw (16,4)node[above]{$3$}--(15,0);
 \draw (17,4)node[above]{$0$}--(16,0);
 \draw[dots] (17.2,4)--(18.8,4);
 \draw[dots] (16.2,0)--(17.8,0);
 \draw (19,4)node[above]{$j-1$}--(18,0);
 \draw (20,4)node[above]{$j$}--(19,0);
 \node[greendot] at (3,2){};
 \node[greendot] at (7,2){};
 \node[greendot] at (11,2){};
 \node[greendot] at (15,4){};
 \node[greendot] at (20,4){};
 \node[below right] at (20,4){$\delta_{j,2}$};
 \draw[dashed, color = red] (-0.5,3) -- (12.5,3) -- (12.5,1) -- (-0.5,1) -- (-0.5,3);
\end{braid}\\
&\hspace*{-5mm} + &
\begin{braid}
 \draw (0,4)node[above]{$0$}--(0,0);
 \draw (1,4)node[above]{$1$}--(1,0);
 \draw (2,4)node[above]{$2$}--(2,0);
 \draw (3,4)node[above]{$3$}--(3,0);
 \draw (4,4)node[above]{$0$}--(4,0);
 \draw[dots] (4.2,4)--(5.8,4);
 \draw[dots] (4.2,0)--(5.8,0);
 \draw (6,4)node[above]{$2$}--(6,0);
 \draw (7,4)node[above]{$3$}--(8,3)--(8,1)--(7,0);
 \draw (8,4)node[above]{$3$}--(7,3)--(7,1)--(15,0);
 \draw (9,4)node[above]{$0$}--(9,1)--(8,0);
 \draw (10,4)node[above]{$1$}--(10,1)--(9,0);
 \draw (11,4)node[above]{$2$}--(11,1)--(10,0);
 \draw (12,4)node[above]{$3$}--(12,1)--(11,0);
 \draw (13,4)node[above]{$0$}--(13,1)--(12,0);
 \draw (14,4)node[above]{$1$}--(14,1)--(13,0);
 \draw (15,4)node[above]{$2$}--(15,1)--(14,0);
 \draw (16,4)node[above]{$3$}--(20,0);
 \draw (17,4)node[above]{$0$}--(16,3)--(16,0);
 \draw[dots] (17.2,4)--(18.8,4);
 \draw[dots] (16.2,0)--(17.8,0);
 \draw (19,4)node[above]{$j-1$}--(18,0);
 \draw (20,4)node[above]{$j$}--(19,0);
 \node[greendot] at (3,2){};
 \node[greendot] at (7,2){};
 \node[greendot] at (11,2){};
 \node[greendot] at (15,2){};
 \node[greendot] at (20,4){};
 \node[below right] at (20,4){$\delta_{j,2}$};
 \draw[dashed, color = red] (-0.5,3) -- (16.5,3) -- (16.5,1) -- (-0.5,1) -- (-0.5,3);
\end{braid} +\ \ \ \   \ldots\ldots\ldots\ldots \ \ \ \text{by \autoref{send}}\\
&\hspace*{-5mm} + &
\begin{braid}
 \draw (0,4)node[above]{$0$}--(0,0);
 \draw (1,4)node[above]{$1$}--(1,0);
 \draw (2,4)node[above]{$2$}--(2,0);
 \draw (3,4)node[above]{$3$}--(3,0);
 \draw (4,4)node[above]{$0$}--(4,0);
 \draw[dots] (4.2,4)--(5.8,4);
 \draw[dots] (4.2,0)--(5.8,0);
 \draw (6,4)node[above]{$2$}--(6,0);
 \draw (7,4)node[above]{$3$}--(20,0);
 \draw (8,4)node[above]{$3$}--(7,0);
 \draw (9,4)node[above]{$0$}--(8,0);
 \draw (10,4)node[above]{$1$}--(9,0);
 \draw (11,4)node[above]{$2$}--(10,0);
 \draw (12,4)node[above]{$3$}--(11,0);
 \draw (13,4)node[above]{$0$}--(12,0);
 \draw (14,4)node[above]{$1$}--(13,0);
 \draw (15,4)node[above]{$2$}--(14,0);
 \draw (16,4)node[above]{$3$}--(15,0);
 \draw (17,4)node[above]{$0$}--(16,0);
 \draw[dots] (17.2,4)--(18.8,4);
 \draw[dots] (16.2,0)--(17.8,0);
 \draw (19,4)node[above]{$j-1$}--(18,0);
 \draw (20,4)node[above]{$j$}--(19,0);
 \node[greendot] at (3,4){};
 \node[greendot] at (20,0){};
 \node[greendot] at (11,4){};
 \node[greendot] at (15,4){};
 \node[greendot] at (20,4){};
 \node[below right] at (20,4){$\delta_{j,2}$};
\end{braid}
\end{eqnarray*}

\begin{eqnarray*}
&\hspace*{-5mm} =_\lambda &
-\begin{braid}
 \draw (0,4)node[above]{$0$}--(0,0);
 \draw (1,4)node[above]{$1$}--(1,0);
 \draw (2,4)node[above]{$2$}--(2,0);
 \draw (3,4)node[above]{$3$}--(3,0);
 \draw (4,4)node[above]{$0$}--(4,0);
 \draw[dots] (4.2,4)--(5.8,4);
 \draw[dots] (4.2,0)--(5.8,0);
 \draw (6,4)node[above]{$2$}--(6,0);
 \draw (8,4)node[above]{$3$}--(20,0);
 \draw (7,4)node[above]{$3$}--(7,0);
 \draw (9,4)node[above]{$0$}--(8,0);
 \draw (10,4)node[above]{$1$}--(9,0);
 \draw (11,4)[densely dotted] node[above]{$2$}--(10,0);
 \draw (12,4)node[above]{$3$}--(11,0);
 \draw (13,4)node[above]{$0$}--(12,0);
 \draw (14,4)node[above]{$1$}--(13,0);
 \draw (15,4)[densely dotted] node[above]{$2$}--(14,0);
 \draw (16,4)node[above]{$3$}--(15,0);
 \draw (17,4)node[above]{$0$}--(16,0);
 \draw[dots] (17.2,4)--(18.8,4);
 \draw[dots] (16.2,0)--(17.8,0);
 \draw (19,4)node[above]{$j-1$}--(18,0);
 \draw (20,4)[densely dotted] node[above]{$j$}--(19,0);
 \node[greendot] at (3,4){};
 \node[greendot] at (11,4){};
 \node[greendot] at (15,4){};
 \node[greendot] at (20,4){};
 \node[below right] at (20,4){$\delta_{j,2}$};
 \draw[->] (11,4) -- (10.5,2);
 \draw[->] (15,4) -- (14.25,1);
 \draw[->] (20,4) -- (19.0,0);
\end{braid} +
\begin{braid}
 \draw (0,4)node[above]{$0$}--(0,0);
 \draw (1,4)node[above]{$1$}--(1,0);
 \draw (2,4)node[above]{$2$}--(2,0);
 \draw (3,4)node[above]{$3$}--(3,0);
 \draw (4,4)node[above]{$0$}--(4,0);
 \draw[dots] (4.2,4)--(5.8,4);
 \draw[dots] (4.2,0)--(5.8,0);
 \draw (6,4)node[above]{$2$}--(6,0);
 \draw (7,4)node[above]{$3$}--(20,0);
 \draw (8,4)node[above]{$3$}--(7,0);
 \draw (9,4)node[above]{$0$}--(8,0);
 \draw (10,4)node[above]{$1$}--(9,0);
 \draw (11,4)[densely dotted] node[above]{$2$}--(10,0);
 \draw (12,4)node[above]{$3$}--(11,0);
 \draw (13,4)node[above]{$0$}--(12,0);
 \draw (14,4)node[above]{$1$}--(13,0);
 \draw (15,4)[densely dotted] node[above]{$2$}--(14,0);
 \draw (16,4)node[above]{$3$}--(15,0);
 \draw (17,4)node[above]{$0$}--(16,0);
 \draw[dots] (17.2,4)--(18.8,4);
 \draw[dots] (16.2,0)--(17.8,0);
 \draw (19,4)node[above]{$j-1$}--(18,0);
 \draw (20,4)[densely dotted] node[above]{$j$}--(19,0);
 \node[greendot] at (3,4){};
 \node[greendot] at (20,0){};
 \node[greendot] at (11,4){};
 \node[greendot] at (15,4){};
 \node[greendot] at (20,4){};
 \node[below right] at (20,4){$\delta_{j,2}$};
 \draw[->] (11,4) -- (10.5,2);
 \draw[->] (15,4) -- (14.25,1);
 \draw[->] (20,4) -- (19,0);
\end{braid}\\
&\hspace*{-5mm} \overset{(\ref{dia:y-psi com})}= &
-\begin{braid}
 \draw (0,4)node[above]{$0$}--(0,0);
 \draw (1,4)node[above]{$1$}--(1,0);
 \draw (2,4)node[above]{$2$}--(2,0);
 \draw (3,4)node[above]{$3$}--(3,0);
 \draw (4,4)node[above]{$0$}--(4,0);
 \draw[dots] (4.2,4)--(5.8,4);
 \draw[dots] (4.2,0)--(5.8,0);
 \draw (6,4)node[above]{$2$}--(6,0);
 \draw (8,4)node[above]{$3$}--(20,0);
 \draw (7,4)node[above]{$3$}--(7,0);
 \draw (9,4)node[above]{$0$}--(8,0);
 \draw (10,4)node[above]{$1$}--(9,0);
 \draw (11,4)node[above]{$2$}--(10,0);
 \draw (12,4)node[above]{$3$}--(11,0);
 \draw (13,4)node[above]{$0$}--(12,0);
 \draw (14,4)node[above]{$1$}--(13,0);
 \draw (15,4)node[above]{$2$}--(14,0);
 \draw (16,4)node[above]{$3$}--(15,0);
 \draw (17,4)node[above]{$0$}--(16,0);
 \draw[dots] (17.2,4)--(18.8,4);
 \draw[dots] (16.2,0)--(17.8,0);
 \draw (19,4)node[above]{$j-1$}--(18,0);
 \draw (20,4)node[above]{$j$}--(19,0);
 \node[greendot] at (3,0){};
 \node[greendot] at (10,0){};
 \node[greendot] at (14,0){};
 \node[greendot] at (19,0){};
 \node[below right] at (19,0){$\delta_{j,2}$};
\end{braid} +
\begin{braid}
 \draw (0,4)node[above]{$0$}--(0,0);
 \draw (1,4)node[above]{$1$}--(1,0);
 \draw (2,4)node[above]{$2$}--(2,0);
 \draw (3,4)node[above]{$3$}--(3,0);
 \draw (4,4)node[above]{$0$}--(4,0);
 \draw[dots] (4.2,4)--(5.8,4);
 \draw[dots] (4.2,0)--(5.8,0);
 \draw (6,4)node[above]{$2$}--(6,0);
 \draw (7,4)node[above]{$3$}--(20,0);
 \draw (8,4)node[above]{$3$}--(7,0);
 \draw (9,4)node[above]{$0$}--(8,0);
 \draw (10,4)node[above]{$1$}--(9,0);
 \draw (11,4)node[above]{$2$}--(10,0);
 \draw (12,4)node[above]{$3$}--(11,0);
 \draw (13,4)node[above]{$0$}--(12,0);
 \draw (14,4)node[above]{$1$}--(13,0);
 \draw (15,4)node[above]{$2$}--(14,0);
 \draw (16,4)node[above]{$3$}--(15,0);
 \draw (17,4)node[above]{$0$}--(16,0);
 \draw[dots] (17.2,4)--(18.8,4);
 \draw[dots] (16.2,0)--(17.8,0);
 \draw (19,4)node[above]{$j-1$}--(18,0);
 \draw (20,4)node[above]{$j$}--(19,0);
 \node[greendot] at (3,0){};
 \node[greendot] at (10,0){};
 \node[greendot] at (14,0){};
 \node[greendot] at (19,0){};
 \node[below right] at (19,0){$\delta_{j,2}$};
 \node[greendot] at (20,0){};
\end{braid}\\
&\hspace*{-5mm} = & -\psi_{\lambda_1+1}\ldots\psi_{n-2}\psi_{n-1} e(\bi_{\dot{\lambda}}\vee i)y_{\dot{\lambda}} + \psi_{\lambda_1}\psi_{\lambda_1+1}\ldots\psi_{n-2}\psi_{n-1} e(\bi_{\dot{\lambda}}\vee i)y_{\dot{\lambda}}y_n.
\end{eqnarray*}

\textbf{Case ~\ref{swap}c:} $i = j = e-1$.

Because we set $e = 4$ so in this case we have $i = j = 3$. Similarly as in Case \ref{swap}b, we have
\begin{eqnarray*}
&&e_\lambda y_\lambda \psi_{\lambda_1}\psi_{\lambda_1+1}\ldots \psi_{n-2}\psi_{n-1}\\
&\hspace*{-10mm} = &
-\begin{braid}
 \draw(3.5,5.2)node{$\overbrace{\hspace*{27mm}}^{\lambda_1}$};
 \draw (0,4)node[above]{$0$}--(0,0);
 \draw (1,4)node[above]{$1$}--(1,0);
 \draw (2,4)node[above]{$2$}--(2,0);
 \draw (3,4)node[above]{$3$}--(3,0);
 \draw (4,4)node[above]{$0$}--(4,0);
 \draw[dots] (4.2,4)--(5.8,4);
 \draw[dots] (4.2,0)--(5.8,0);
 \draw (6,4)node[above]{$2$}--(6,0);
 \draw (8,4)node[above]{$3$}--(20,0);
 \draw (7,4)node[above]{$3$}--(7,0);
 \draw (9,4)node[above]{$0$}--(8,0);
 \draw (10,4)node[above]{$1$}--(9,0);
 \draw (11,4)node[above]{$2$}--(10,0);
 \draw (12,4)node[above]{$3$}--(11,0);
 \draw (13,4)node[above]{$0$}--(12,0);
 \draw (14,4)node[above]{$1$}--(13,0);
 \draw (15,4)node[above]{$2$}--(14,0);
 \draw (16,4)node[above]{$3$}--(15,0);
 \draw (17,4)node[above]{$0$}--(16,0);
 \draw[dots] (17.2,4)--(18.8,4);
 \draw[dots] (16.2,0)--(17.8,0);
 \draw (19,4)node[above]{$2$}--(18,0);
 \draw (20,4)node[above]{$3$}--(19,0);
 \node[greendot] at (3,4){};
 \node[greendot] at (11,4){};
 \node[greendot] at (15,4){};
 \node[greendot] at (19,4){};
\end{braid} -
\begin{braid}
 \draw(3.5,5.2)node{$\overbrace{\hspace*{27mm}}^{\lambda_1}$};
 \draw (0,4)node[above]{$0$}--(0,0);
 \draw (1,4)node[above]{$1$}--(1,0);
 \draw (2,4)node[above]{$2$}--(2,0);
 \draw (3,4)node[above]{$3$}--(3,0);
 \draw (4,4)node[above]{$0$}--(4,0);
 \draw[dots] (4.2,4)--(5.8,4);
 \draw[dots] (4.2,0)--(5.8,0);
 \draw (6,4)node[above]{$2$}--(6,0);
 \draw (7,4)node[above]{$3$}--(8,3)--(8,1)--(7,0);
 \draw (8,4)node[above]{$3$}--(7,3)--(7,1)--(11,0);
 \draw (9,4)node[above]{$0$}--(9,1)--(8,0);
 \draw (10,4)node[above]{$1$}--(10,1)--(9,0);
 \draw (11,4)node[above]{$2$}--(11,1)--(10,0);
 \draw (12,4)node[above]{$3$}--(20,0);
 \draw (13,4)node[above]{$0$}--(12,3)--(12,0);
 \draw (14,4)node[above]{$1$}--(13,0);
 \draw (15,4)node[above]{$2$}--(14,0);
 \draw (16,4)node[above]{$3$}--(15,0);
 \draw (17,4)node[above]{$0$}--(16,0);
 \draw[dots] (17.2,4)--(18.8,4);
 \draw[dots] (16.2,0)--(17.8,0);
 \draw (19,4)node[above]{$2$}--(18,0);
 \draw (20,4)node[above]{$3$}--(19,0);
 \node[greendot] at (3,2){};
 \node[greendot] at (7,2){};
 \node[greendot] at (11,2){};
 \node[greendot] at (15,4){};
 \node[greendot] at (19,4){};
 \draw[dashed, color = red] (-0.5,3) -- (12.5,3) -- (12.5,1) -- (-0.5,1) -- (-0.5,3);
\end{braid}\\
&\hspace*{-10mm} - &
\begin{braid}
 \draw(3.5,5.2)node{$\overbrace{\hspace*{27mm}}^{\lambda_1}$};
 \draw (0,4)node[above]{$0$}--(0,0);
 \draw (1,4)node[above]{$1$}--(1,0);
 \draw (2,4)node[above]{$2$}--(2,0);
 \draw (3,4)node[above]{$3$}--(3,0);
 \draw (4,4)node[above]{$0$}--(4,0);
 \draw[dots] (4.2,4)--(5.8,4);
 \draw[dots] (4.2,0)--(5.8,0);
 \draw (6,4)node[above]{$2$}--(6,0);
 \draw (7,4)node[above]{$3$}--(8,3)--(8,1)--(7,0);
 \draw (8,4)node[above]{$3$}--(7,3)--(7,1)--(15,0);
 \draw (9,4)node[above]{$0$}--(9,1)--(8,0);
 \draw (10,4)node[above]{$1$}--(10,1)--(9,0);
 \draw (11,4)node[above]{$2$}--(11,1)--(10,0);
 \draw (12,4)node[above]{$3$}--(12,1)--(11,0);
 \draw (13,4)node[above]{$0$}--(13,1)--(12,0);
 \draw (14,4)node[above]{$1$}--(14,1)--(13,0);
 \draw (15,4)node[above]{$2$}--(15,1)--(14,0);
 \draw (16,4)node[above]{$3$}--(20,0);
 \draw (17,4)node[above]{$0$}--(16,3)--(16,0);
 \draw[dots] (17.2,4)--(18.8,4);
 \draw[dots] (16.2,0)--(17.8,0);
 \draw (19,4)node[above]{$2$}--(18,0);
 \draw (20,4)node[above]{$3$}--(19,0);
 \node[greendot] at (3,2){};
 \node[greendot] at (7,2){};
 \node[greendot] at (11,2){};
 \node[greendot] at (15,2){};
 \node[greendot] at (19,4){};
 \draw[dashed, color = red] (-0.5,3) -- (16.5,3) -- (16.5,1) -- (-0.5,1) -- (-0.5,3);
\end{braid} +\ \ \  \ldots\ldots\ldots\ldots\ \ \ \text{by \autoref{send}} \\
&\hspace*{-10mm} - &
\begin{braid}
 \draw(3.5,5.2)node{$\overbrace{\hspace*{27mm}}^{\lambda_1}$};
 \draw (0,4)node[above]{$0$}--(0,0);
 \draw (1,4)node[above]{$1$}--(1,0);
 \draw (2,4)node[above]{$2$}--(2,0);
 \draw (3,4)node[above]{$3$}--(3,0);
 \draw (4,4)node[above]{$0$}--(4,0);
 \draw[dots] (4.2,4)--(5.8,4);
 \draw[dots] (4.2,0)--(5.8,0);
 \draw (6,4)node[above]{$2$}--(6,0);
 \draw (7,4)node[above]{$3$}--(19,0);
 \draw (8,4)node[above]{$3$}--(7,0);
 \draw (9,4)node[above]{$0$}--(8,0);
 \draw (10,4)node[above]{$1$}--(9,0);
 \draw (11,4)node[above]{$2$}--(10,0);
 \draw (12,4)node[above]{$3$}--(11,0);
 \draw (13,4)node[above]{$0$}--(12,0);
 \draw (14,4)node[above]{$1$}--(13,0);
 \draw (15,4)node[above]{$2$}--(14,0);
 \draw (16,4)node[above]{$3$}--(15,0);
 \draw (17,4)node[above]{$0$}--(16,0);
 \draw[dots] (17.2,4)--(18.8,4);
 \draw[dots] (16.2,0)--(17.8,0);
 \draw (19,4)node[above]{$2$}--(18,0);
 \draw (20,4)node[above]{$3$}--(20,0);
 \node[greendot] at (3,4){};
 \node[greendot] at (11,4){};
 \node[greendot] at (15,4){};
 \node[greendot] at (19,4){};
\end{braid} +
\begin{braid}
 \draw(3.5,5.2)node{$\overbrace{\hspace*{27mm}}^{\lambda_1}$};
 \draw (0,4)node[above]{$0$}--(0,0);
 \draw (1,4)node[above]{$1$}--(1,0);
 \draw (2,4)node[above]{$2$}--(2,0);
 \draw (3,4)node[above]{$3$}--(3,0);
 \draw (4,4)node[above]{$0$}--(4,0);
 \draw[dots] (4.2,4)--(5.8,4);
 \draw[dots] (4.2,0)--(5.8,0);
 \draw (6,4)node[above]{$2$}--(6,0);
 \draw (7,4)node[above]{$3$}--(20,0);
 \draw (8,4)node[above]{$3$}--(7,0);
 \draw (9,4)node[above]{$0$}--(8,0);
 \draw (10,4)node[above]{$1$}--(9,0);
 \draw (11,4)node[above]{$2$}--(10,0);
 \draw (12,4)node[above]{$3$}--(11,0);
 \draw (13,4)node[above]{$0$}--(12,0);
 \draw (14,4)node[above]{$1$}--(13,0);
 \draw (15,4)node[above]{$2$}--(14,0);
 \draw (16,4)node[above]{$3$}--(15,0);
 \draw (17,4)node[above]{$0$}--(16,0);
 \draw[dots] (17.2,4)--(18.8,4);
 \draw[dots] (16.2,0)--(17.8,0);
 \draw (19,4)node[above]{$2$}--(18,0);
 \draw (20,4)node[above]{$3$}--(19,0);
 \node[greendot] at (3,4){};
 \node[greendot] at (20,0){};
 \node[greendot] at (11,4){};
 \node[greendot] at (15,4){};
 \node[greendot] at (19,4){};
\end{braid}\\
&\hspace*{-5mm} =_\lambda &
-\begin{braid}
 \draw (0,4)node[above]{$0$}--(0,0);
 \draw (1,4)node[above]{$1$}--(1,0);
 \draw (2,4)node[above]{$2$}--(2,0);
 \draw (3,4)node[above]{$3$}--(3,0);
 \draw (4,4)node[above]{$0$}--(4,0);
 \draw[dots] (4.2,4)--(5.8,4);
 \draw[dots] (4.2,0)--(5.8,0);
 \draw (6,4)node[above]{$2$}--(6,0);
 \draw (8,4)node[above]{$3$}--(20,0);
 \draw (7,4)node[above]{$3$}--(7,0);
 \draw (9,4)node[above]{$0$}--(8,0);
 \draw (10,4)node[above]{$1$}--(9,0);
 \draw (11,4)[densely dotted] node[above]{$2$}--(10,0);
 \draw (12,4)node[above]{$3$}--(11,0);
 \draw (13,4)node[above]{$0$}--(12,0);
 \draw (14,4)node[above]{$1$}--(13,0);
 \draw (15,4)[densely dotted] node[above]{$2$}--(14,0);
 \draw (16,4)node[above]{$3$}--(15,0);
 \draw (17,4)node[above]{$0$}--(16,0);
 \draw[dots] (17.2,4)--(18.8,4);
 \draw[dots] (16.2,0)--(17.8,0);
 \draw (19,4)[densely dotted] node[above]{$2$}--(18,0);
 \draw (20,4)node[above]{$3$}--(19,0);
 \node[greendot] at (3,4){};
 \node[greendot] at (11,4){};
 \node[greendot] at (15,4){};
 \node[greendot] at (19,4){};
 \draw[->] (11,4) -- (10.5,2);
 \draw[->] (15,4) -- (14.2,0.8);
 \draw[->] (19,4) -- (18.02,0.08);
\end{braid} -
\begin{braid}
 \draw (0,4)node[above]{$0$}--(0,0);
 \draw (1,4)node[above]{$1$}--(1,0);
 \draw (2,4)node[above]{$2$}--(2,0);
 \draw (3,4)node[above]{$3$}--(3,0);
 \draw (4,4)node[above]{$0$}--(4,0);
 \draw[dots] (4.2,4)--(5.8,4);
 \draw[dots] (4.2,0)--(5.8,0);
 \draw (6,4)node[above]{$2$}--(6,0);
 \draw (7,4)node[above]{$3$}--(19,0);
 \draw (8,4)node[above]{$3$}--(7,0);
 \draw (9,4)node[above]{$0$}--(8,0);
 \draw (10,4)node[above]{$1$}--(9,0);
 \draw (11,4)[densely dotted] node[above]{$2$}--(10,0);
 \draw (12,4)node[above]{$3$}--(11,0);
 \draw (13,4)node[above]{$0$}--(12,0);
 \draw (14,4)node[above]{$1$}--(13,0);
 \draw (15,4)[densely dotted] node[above]{$2$}--(14,0);
 \draw (16,4)node[above]{$3$}--(15,0);
 \draw (17,4)node[above]{$0$}--(16,0);
 \draw[dots] (17.2,4)--(18.8,4);
 \draw[dots] (16.2,0)--(17.8,0);
 \draw (19,4)[densely dotted] node[above]{$2$}--(18,0);
 \draw (20,4)node[above]{$3$}--(20,0);
 \node[greendot] at (3,4){};
 \node[greendot] at (11,4){};
 \node[greendot] at (15,4){};
 \node[greendot] at (19,4){};
 \draw[->] (11,4) -- (10.5,2);
 \draw[->] (15,4) -- (14.2,0.8);
 \draw[->] (19,4) -- (18.02,0.08);
\end{braid}\\
&\hspace*{-5mm} + &
\begin{braid}
 \draw (0,4)node[above]{$0$}--(0,0);
 \draw (1,4)node[above]{$1$}--(1,0);
 \draw (2,4)node[above]{$2$}--(2,0);
 \draw (3,4)node[above]{$3$}--(3,0);
 \draw (4,4)node[above]{$0$}--(4,0);
 \draw[dots] (4.2,4)--(5.8,4);
 \draw[dots] (4.2,0)--(5.8,0);
 \draw (6,4)node[above]{$2$}--(6,0);
 \draw (7,4)node[above]{$3$}--(20,0);
 \draw (8,4)node[above]{$3$}--(7,0);
 \draw (9,4)node[above]{$0$}--(8,0);
 \draw (10,4)node[above]{$1$}--(9,0);
 \draw (11,4)[densely dotted] node[above]{$2$}--(10,0);
 \draw (12,4)node[above]{$3$}--(11,0);
 \draw (13,4)node[above]{$0$}--(12,0);
 \draw (14,4)node[above]{$1$}--(13,0);
 \draw (15,4)[densely dotted] node[above]{$2$}--(14,0);
 \draw (16,4)node[above]{$3$}--(15,0);
 \draw (17,4)node[above]{$0$}--(16,0);
 \draw[dots] (17.2,4)--(18.8,4);
 \draw[dots] (16.2,0)--(17.8,0);
 \draw (19,4)[densely dotted] node[above]{$2$}--(18,0);
 \draw (20,4)node[above]{$3$}--(19,0);
 \node[greendot] at (3,4){};
 \node[greendot] at (20,0){};
 \node[greendot] at (11,4){};
 \node[greendot] at (15,4){};
 \node[greendot] at (19,4){};
 \draw[->] (11,4) -- (10.5,2);
 \draw[->] (15,4) -- (14.2,0.8);
 \draw[->] (19,4) -- (18.02,0.08);
\end{braid}\\
&\hspace*{-5mm} =&
-\begin{braid}
 \draw (0,4)node[above]{$0$}--(0,0);
 \draw (1,4)node[above]{$1$}--(1,0);
 \draw (2,4)node[above]{$2$}--(2,0);
 \draw (3,4)node[above]{$3$}--(3,0);
 \draw (4,4)node[above]{$0$}--(4,0);
 \draw[dots] (4.2,4)--(5.8,4);
 \draw[dots] (4.2,0)--(5.8,0);
 \draw (6,4)node[above]{$2$}--(6,0);
 \draw (8,4)node[above]{$3$}--(20,0);
 \draw (7,4)node[above]{$3$}--(7,0);
 \draw (9,4)node[above]{$0$}--(8,0);
 \draw (10,4)node[above]{$1$}--(9,0);
 \draw (11,4)node[above]{$2$}--(10,0);
 \draw (12,4)node[above]{$3$}--(11,0);
 \draw (13,4)node[above]{$0$}--(12,0);
 \draw (14,4)node[above]{$1$}--(13,0);
 \draw (15,4)node[above]{$2$}--(14,0);
 \draw (16,4)node[above]{$3$}--(15,0);
 \draw (17,4)node[above]{$0$}--(16,0);
 \draw[dots] (17.2,4)--(18.8,4);
 \draw[dots] (16.2,0)--(17.8,0);
 \draw (19,4)node[above]{$2$}--(18,0);
 \draw (20,4)node[above]{$3$}--(19,0);
 \node[greendot] at (3,0){};
 \node[greendot] at (10,0){};
 \node[greendot] at (14,0){};
 \node[greendot] at (18,0){};
\end{braid} -
\begin{braid}
 \draw (0,4)node[above]{$0$}--(0,0);
 \draw (1,4)node[above]{$1$}--(1,0);
 \draw (2,4)node[above]{$2$}--(2,0);
 \draw (3,4)node[above]{$3$}--(3,0);
 \draw (4,4)node[above]{$0$}--(4,0);
 \draw[dots] (4.2,4)--(5.8,4);
 \draw[dots] (4.2,0)--(5.8,0);
 \draw (6,4)node[above]{$2$}--(6,0);
 \draw (7,4)node[above]{$3$}--(19,0);
 \draw (8,4)node[above]{$3$}--(7,0);
 \draw (9,4)node[above]{$0$}--(8,0);
 \draw (10,4)node[above]{$1$}--(9,0);
 \draw (11,4)node[above]{$2$}--(10,0);
 \draw (12,4)node[above]{$3$}--(11,0);
 \draw (13,4)node[above]{$0$}--(12,0);
 \draw (14,4)node[above]{$1$}--(13,0);
 \draw (15,4)node[above]{$2$}--(14,0);
 \draw (16,4)node[above]{$3$}--(15,0);
 \draw (17,4)node[above]{$0$}--(16,0);
 \draw[dots] (17.2,4)--(18.8,4);
 \draw[dots] (16.2,0)--(17.8,0);
 \draw (19,4)node[above]{$2$}--(18,0);
 \draw (20,4)node[above]{$3$}--(20,0);
 \node[greendot] at (3,0){};
 \node[greendot] at (10,0){};
 \node[greendot] at (14,0){};
 \node[greendot] at (18,0){};
\end{braid}
\end{eqnarray*}

\begin{eqnarray*}
&\hspace*{-5mm} + &
\begin{braid}
 \draw (0,4)node[above]{$0$}--(0,0);
 \draw (1,4)node[above]{$1$}--(1,0);
 \draw (2,4)node[above]{$2$}--(2,0);
 \draw (3,4)node[above]{$3$}--(3,0);
 \draw (4,4)node[above]{$0$}--(4,0);
 \draw[dots] (4.2,4)--(5.8,4);
 \draw[dots] (4.2,0)--(5.8,0);
 \draw (6,4)node[above]{$2$}--(6,0);
 \draw (7,4)node[above]{$3$}--(20,0);
 \draw (8,4)node[above]{$3$}--(7,0);
 \draw (9,4)node[above]{$0$}--(8,0);
 \draw (10,4)node[above]{$1$}--(9,0);
 \draw (11,4)node[above]{$2$}--(10,0);
 \draw (12,4)node[above]{$3$}--(11,0);
 \draw (13,4)node[above]{$0$}--(12,0);
 \draw (14,4)node[above]{$1$}--(13,0);
 \draw (15,4)node[above]{$2$}--(14,0);
 \draw (16,4)node[above]{$3$}--(15,0);
 \draw (17,4)node[above]{$0$}--(16,0);
 \draw[dots] (17.2,4)--(18.8,4);
 \draw[dots] (16.2,0)--(17.8,0);
 \draw (19,4)node[above]{$2$}--(18,0);
 \draw (20,4)node[above]{$3$}--(19,0);
 \node[greendot] at (3,0){};
 \node[greendot] at (20,0){};
 \node[greendot] at (10,0){};
 \node[greendot] at (14,0){};
 \node[greendot] at (18,0){};
\end{braid}\\
&\hspace*{-5mm}= &
-\psi_{\lambda_1+1}\ldots\psi_{n-2}\psi_{n-1} e(\bi_{\dot{\lambda}}\vee i)y_{\dot{\lambda}} - \psi_{\lambda_1}\ldots\psi_{n-2} e(\bi_{\dot{\lambda}}\vee i)y_{\dot{\lambda}} + \psi_{\lambda_1}\psi_{\lambda_1+1}\ldots\psi_{n-2}\psi_{n-1} e(\bi_{\dot{\lambda}}\vee i)y_{\dot{\lambda}}y_n.
\end{eqnarray*}

\textbf{Case ~\ref{swap}d:} $i = e - 2$ and $j \neq e - 2$.

Because we set $e = 4$ so in this case we have $i = 2$ and $j \neq 2$. Similarly we set $j = 3$ in this case in order to make the diagrams easier to read. For the other cases with $j \neq 2$ the argument is similar. By \autoref{send},
\begin{eqnarray*}
&&e_\lambda y_\lambda \psi_{\lambda_1}\psi_{\lambda_1+1}\ldots \psi_{n-2}\psi_{n-1}\\
&\hspace*{-5mm} = &
\begin{braid}
 \draw(3,5.2)node{$\overbrace{\hspace*{24mm}}^{\lambda_1}$};
 \draw (0,4)node[above]{$0$}--(0,0);
 \draw (1,4)node[above]{$1$}--(1,0);
 \draw (2,4)node[above]{$2$}--(2,0);
 \draw (3,4)node[above]{$3$}--(3,0);
 \draw[dots] (3.2,4)--(4.8,4);
 \draw[dots] (3.2,0)--(4.8,0);
 \draw (5,4)node[above]{$1$}--(5,0);
 \draw (6,4)node[above]{$2$}--(20,0);
 \draw (7,4)node[above]{$3$}--(6,0);
 \draw (8,4)node[above]{$0$}--(7,0);
 \draw (9,4)node[above]{$1$}--(8,0);
 \draw (10,4)node[above]{$2$}--(9,0);
 \draw (11,4)node[above]{$3$}--(10,0);
 \draw (12,4)node[above]{$0$}--(11,0);
 \draw (13,4)node[above]{$1$}--(12,0);
 \draw (14,4)node[above]{$2$}--(13,0);
 \draw (15,4)node[above]{$3$}--(14,0);
 \draw (16,4)node[above]{$0$}--(15,0);
 \draw[dots] (16.2,4)--(17.8,4);
 \draw[dots] (15.2,0)--(16.8,0);
 \draw (18,4)node[above]{$1$}--(17,0);
 \draw (19,4)[densely dotted] node[above]{$2$}--(18,0);
 \draw (20,4)node[above]{$3$}--(19,0);
 \node[greendot] at (3,4){};
 \node[greendot] at (10,4){};
 \node[greendot] at (14,4){};
 \node[greendot] at (19,4){};
 \draw[->] (19,4) -- (18.02,0.04);
\end{braid}\\
& \hspace*{-5mm}\overset{(\ref{dia:y-psi com})}= &
\begin{braid}
 \draw (0,4)node[above]{$0$}--(0,0);
 \draw (1,4)node[above]{$1$}--(1,0);
 \draw (2,4)node[above]{$2$}--(2,0);
 \draw (3,4)node[above]{$3$}--(3,0);
 \draw[dots] (3.2,4)--(4.8,4);
 \draw[dots] (3.2,0)--(4.8,0);
 \draw (5,4)node[above]{$1$}--(5,0);
 \draw (6,4)node[above]{$2$}--(20,0);
 \draw (7,4)node[above]{$3$}--(6,0);
 \draw (8,4)node[above]{$0$}--(7,0);
 \draw (9,4)node[above]{$1$}--(8,0);
 \draw (10,4)node[above]{$2$}--(9,0);
 \draw (11,4)node[above]{$3$}--(10,0);
 \draw (12,4)node[above]{$0$}--(11,0);
 \draw (13,4)node[above]{$1$}--(12,0);
 \draw (14,4)node[above]{$2$}--(13,0);
 \draw (15,4)node[above]{$3$}--(14,0);
 \draw (16,4)node[above]{$0$}--(15,0);
 \draw[dots] (16.2,4)--(17.8,4);
 \draw[dots] (15.2,0)--(16.8,0);
 \draw (18,4)node[above]{$1$}--(17,0);
 \draw (19,4)node[above]{$2$}--(18,0);
 \draw (20,4)node[above]{$3$}--(19,0);
 \node[greendot] at (3,4){};
 \node[greendot] at (10,4){};
 \node[greendot] at (14,4){};
 \node[greendot] at (18,0){};
\end{braid}
+
\begin{braid}
 \draw (0,4)node[above]{$0$}--(0,0);
 \draw (1,4)node[above]{$1$}--(1,0);
 \draw (2,4)node[above]{$2$}--(2,0);
 \draw (3,4)node[above]{$3$}--(3,0);
 \draw[dots] (3.2,4)--(4.8,4);
 \draw[dots] (3.2,0)--(4.8,0);
 \draw (5,4)node[above]{$1$}--(5,0);
 \draw (6,4)node[above]{$2$}--(6,1)--(18,0);
 \draw (7,4)node[above]{$3$}--(7,1)--(6,0);
 \draw (8,4)node[above]{$0$}--(8,1)--(7,0);
 \draw (9,4)node[above]{$1$}--(9,1)--(8,0);
 \draw (10,4)node[above]{$2$}--(10,1)--(9,0);
 \draw (11,4)node[above]{$3$}--(11,1)--(10,0);
 \draw (12,4)node[above]{$0$}--(12,1)--(11,0);
 \draw (13,4)node[above]{$1$}--(13,1)--(12,0);
 \draw (14,4)node[above]{$2$}--(14,1)--(13,0);
 \draw (15,4)node[above]{$3$}--(15,1)--(14,0);
 \draw (16,4)node[above]{$0$}--(16,1)--(15,0);
 \draw[dots] (16.2,4)--(17.8,4);
 \draw[dots] (15.2,0)--(16.8,0);
 \draw (18,4)node[above]{$1$}--(18,1)--(17,0);
 \draw (19,4)node[above]{$2$}--(20,3)--(20,0);
 \draw (20,4)node[above]{$3$}--(19,3)--(19,0);
 \node[greendot] at (3,2){};
 \node[greendot] at (10,2){};
 \node[greendot] at (14,2){};
 \draw[dashed, color = red] (-0.5,3) -- (19.5,3) -- (19.5,1) -- (-0.5,1) -- (-0.5,3);
\end{braid}\\
&\hspace*{-5mm} =_\lambda &
\begin{braid}
 \draw (0,4)node[above]{$0$}--(0,0);
 \draw (1,4)node[above]{$1$}--(1,0);
 \draw (2,4)node[above]{$2$}--(2,0);
 \draw (3,4)node[above]{$3$}--(3,0);
 \draw[dots] (3.2,4)--(4.8,4);
 \draw[dots] (3.2,0)--(4.8,0);
 \draw (5,4)node[above]{$1$}--(5,0);
 \draw (6,4)node[above]{$2$}--(20,0);
 \draw (7,4)node[above]{$3$}--(6,0);
 \draw (8,4)node[above]{$0$}--(7,0);
 \draw (9,4)node[above]{$1$}--(8,0);
 \draw (10,4)node[above]{$2$}--(9,0);
 \draw (11,4)node[above]{$3$}--(10,0);
 \draw (12,4)node[above]{$0$}--(11,0);
 \draw (13,4)node[above]{$1$}--(12,0);
 \draw (14,4)[densely dotted] node[above]{$2$}--(13,0);
 \draw (15,4)node[above]{$3$}--(14,0);
 \draw (16,4)node[above]{$0$}--(15,0);
 \draw[dots] (16.2,4)--(17.8,4);
 \draw[dots] (15.2,0)--(16.8,0);
 \draw (18,4)node[above]{$1$}--(17,0);
 \draw (19,4)node[above]{$2$}--(18,0);
 \draw (20,4)node[above]{$3$}--(19,0);
 \node[greendot] at (3,4){};
 \node[greendot] at (10,4){};
 \node[greendot] at (14,4){};
 \node[greendot] at (18,0){};
 \draw[->] (14,4) -- (13.2,0.8);
\end{braid} = \ \ \ \ \ldots\ldots\ldots\ldots\\
&\hspace*{-5mm} =_\lambda &
\begin{braid}
 \draw (0,4)node[above]{$0$}--(0,0);
 \draw (1,4)node[above]{$1$}--(1,0);
 \draw (2,4)node[above]{$2$}--(2,0);
 \draw (3,4)node[above]{$3$}--(3,0);
 \draw[dots] (3.2,4)--(4.8,4);
 \draw[dots] (3.2,0)--(4.8,0);
 \draw (5,4)node[above]{$1$}--(5,0);
 \draw (6,4)node[above]{$2$}--(20,0);
 \draw (7,4)node[above]{$3$}--(6,0);
 \draw (8,4)node[above]{$0$}--(7,0);
 \draw (9,4)node[above]{$1$}--(8,0);
 \draw (10,4)[densely dotted] node[above]{$2$}--(9,0);
 \draw (11,4)node[above]{$3$}--(10,0);
 \draw (12,4)node[above]{$0$}--(11,0);
 \draw (13,4)node[above]{$1$}--(12,0);
 \draw (14,4)node[above]{$2$}--(13,0);
 \draw (15,4)node[above]{$3$}--(14,0);
 \draw (16,4)node[above]{$0$}--(15,0);
 \draw[dots] (16.2,4)--(17.8,4);
 \draw[dots] (15.2,0)--(16.8,0);
 \draw (18,4)node[above]{$1$}--(17,0);
 \draw (19,4)node[above]{$2$}--(18,0);
 \draw (20,4)node[above]{$3$}--(19,0);
 \node[greendot] at (3,4){};
 \node[greendot] at (10,4){};
 \node[greendot] at (13,0){};
 \node[greendot] at (18,0){};
 \draw[->] (10,4) -- (9.5,2);
\end{braid}\\
&\hspace*{-5mm} \overset{(\ref{dia:y-psi com})}= &
\begin{braid}
 \draw (0,4)node[above]{$0$}--(0,0);
 \draw (1,4)node[above]{$1$}--(1,0);
 \draw (2,4)node[above]{$2$}--(2,0);
 \draw (3,4)node[above]{$3$}--(3,0);
 \draw[dots] (3.2,4)--(4.8,4);
 \draw[dots] (3.2,0)--(4.8,0);
 \draw (5,4)node[above]{$1$}--(5,0);
 \draw (6,4)node[above]{$2$}--(20,0);
 \draw (7,4)node[above]{$3$}--(6,0);
 \draw (8,4)node[above]{$0$}--(7,0);
 \draw (9,4)node[above]{$1$}--(8,0);
 \draw (10,4)node[above]{$2$}--(9,0);
 \draw (11,4)node[above]{$3$}--(10,0);
 \draw (12,4)node[above]{$0$}--(11,0);
 \draw (13,4)node[above]{$1$}--(12,0);
 \draw (14,4)node[above]{$2$}--(13,0);
 \draw (15,4)node[above]{$3$}--(14,0);
 \draw (16,4)node[above]{$0$}--(15,0);
 \draw[dots] (16.2,4)--(17.8,4);
 \draw[dots] (15.2,0)--(16.8,0);
 \draw (18,4)node[above]{$1$}--(17,0);
 \draw (19,4)node[above]{$2$}--(18,0);
 \draw (20,4)node[above]{$3$}--(19,0);
 \node[greendot] at (3,4){};
 \node[greendot] at (9,0){};
 \node[greendot] at (13,0){};
 \node[greendot] at (18,0){};
\end{braid}
+
\begin{braid}
 \draw (0,4)node[above]{$0$}--(0,0);
 \draw (1,4)node[above]{$1$}--(1,0);
 \draw (2,4)node[above]{$2$}--(2,0);
 \draw (3,4)node[above]{$3$}--(3,0);
 \draw[dots] (3.2,4)--(4.8,4);
 \draw[dots] (3.2,0)--(4.8,0);
 \draw (5,4)node[above]{$1$}--(5,0);
 \draw (6,4)node[above]{$2$}--(6,1)--(9,0);
 \draw (7,4)node[above]{$3$}--(7,1)--(6,0);
 \draw (8,4)node[above]{$0$}--(8,1)--(7,0);
 \draw (9,4)node[above]{$1$}--(9,1)--(8,0);
 \draw (10,4)node[above]{$2$}--(20,3)--(20,0);
 \draw (11,4)node[above]{$3$}--(10,3)--(10,0);
 \draw (12,4)node[above]{$0$}--(11,3)--(11,0);
 \draw (13,4)node[above]{$1$}--(12,3)--(12,0);
 \draw (14,4)node[above]{$2$}--(13,3)--(13,0);
 \draw (15,4)node[above]{$3$}--(14,3)--(14,0);
 \draw (16,4)node[above]{$0$}--(15,3)--(15,0);
 \draw[dots] (16.2,4)--(17.8,4);
 \draw[dots] (15.2,0)--(16.8,0);
 \draw (18,4)node[above]{$1$}--(17,3)--(17,0);
 \draw (19,4)node[above]{$2$}--(18,3)--(18,0);
 \draw (20,4)node[above]{$3$}--(19,3)--(19,0);
 \node[greendot] at (3,2){};
 \node[greendot] at (13,0){};
 \node[greendot] at (18,0){};
 \draw[dashed, color = red] (-0.5,3) -- (10.5,3) -- (10.5,1) -- (-0.5,1) -- (-0.5,3);
\end{braid}\\
&\hspace*{-5mm} =_\lambda &
\begin{braid}
 \draw (0,4)node[above]{$0$}--(0,0);
 \draw (1,4)node[above]{$1$}--(1,0);
 \draw (2,4)node[above]{$2$}--(2,0);
 \draw (3,4)node[above]{$3$}--(3,0);
 \draw[dots] (3.2,4)--(4.8,4);
 \draw[dots] (3.2,0)--(4.8,0);
 \draw (5,4)node[above]{$1$}--(5,0);
 \draw (6,4)node[above]{$2$}--(20,0);
 \draw (7,4)node[above]{$3$}--(6,0);
 \draw (8,4)node[above]{$0$}--(7,0);
 \draw (9,4)node[above]{$1$}--(8,0);
 \draw (10,4)node[above]{$2$}--(9,0);
 \draw (11,4)node[above]{$3$}--(10,0);
 \draw (12,4)node[above]{$0$}--(11,0);
 \draw (13,4)node[above]{$1$}--(12,0);
 \draw (14,4)node[above]{$2$}--(13,0);
 \draw (15,4)node[above]{$3$}--(14,0);
 \draw (16,4)node[above]{$0$}--(15,0);
 \draw[dots] (16.2,4)--(17.8,4);
 \draw[dots] (15.2,0)--(16.8,0);
 \draw (18,4)node[above]{$1$}--(17,0);
 \draw (19,4)node[above]{$2$}--(18,0);
 \draw (20,4)node[above]{$3$}--(19,0);
 \node[greendot] at (3,0){};
 \node[greendot] at (9,0){};
 \node[greendot] at (13,0){};
 \node[greendot] at (18,0){};
\end{braid}
=
\psi_{\lambda_1}\psi_{\lambda_1+1}\ldots\psi_{n-2}\psi_{n-1} e(\bi_{\dot{\lambda}}\vee i)y_{\dot{\lambda}}.
\end{eqnarray*}

\textbf{Case ~\ref{swap}e:} $i = j =e - 2$.

Because we set $e = 4$ so in this case we have $i = j = 2$. Then by \autoref{send},
\begin{eqnarray*}
&&\hspace*{-5mm}e_\lambda y_\lambda \psi_{\lambda_1}\psi_{\lambda_1+1}\ldots \psi_{n-2}\psi_{n-1} 
=
\begin{braid}
 \draw(3,5.2)node{$\overbrace{\hspace*{24mm}}^{\lambda_1}$};
 \draw (0,4)node[above]{$0$}--(0,0);
 \draw (1,4)node[above]{$1$}--(1,0);
 \draw (2,4)node[above]{$2$}--(2,0);
 \draw (3,4)node[above]{$3$}--(3,0);
 \draw[dots] (3.2,4)--(4.8,4);
 \draw[dots] (3.2,0)--(4.8,0);
 \draw (5,4)node[above]{$1$}--(5,0);
 \draw (6,4)node[above]{$2$}--(21,0);
 \draw (7,4)node[above]{$3$}--(6,0);
 \draw (8,4)node[above]{$0$}--(7,0);
 \draw (9,4)node[above]{$1$}--(8,0);
 \draw (10,4)node[above]{$2$}--(9,0);
 \draw (11,4)node[above]{$3$}--(10,0);
 \draw (12,4)node[above]{$0$}--(11,0);
 \draw (13,4)node[above]{$1$}--(12,0);
 \draw (14,4)node[above]{$2$}--(13,0);
 \draw[dots] (14.2,4)--(15.8,4);
 \draw[dots] (13.2,0)--(14.8,0);
 \draw (16,4)node[above]{$1$}--(15,0);
 \draw (17,4)[densely dotted] node[above]{$2$}--(16,0);
 \draw (18,4)node[above]{$3$}--(17,0);
 \draw (19,4)node[above]{$0$}--(18,0);
 \draw (20,4)node[above]{$1$}--(19,0);
 \draw (21,4)node[above]{$2$}--(20,0);
 \node[greendot] at (3,4){};
 \node[greendot] at (10,4){};
 \node[greendot] at (14,4){};
 \node[greendot] at (17,4){};
 \node[greendot] at (21,4){};
 \draw[->] (17,4) -- (16.2,0.8);
\end{braid}
\end{eqnarray*}

\begin{eqnarray*}
&\hspace*{-5mm} \overset{(\ref{dia:y-psi com})}= &
\begin{braid}
 \draw (0,4)node[above]{$0$}--(0,0);
 \draw (1,4)node[above]{$1$}--(1,0);
 \draw (2,4)node[above]{$2$}--(2,0);
 \draw (3,4)node[above]{$3$}--(3,0);
 \draw[dots] (3.2,4)--(4.8,4);
 \draw[dots] (3.2,0)--(4.8,0);
 \draw (5,4)node[above]{$1$}--(5,0);
 \draw (6,4)node[above]{$2$}--(21,0);
 \draw (7,4)node[above]{$3$}--(6,0);
 \draw (8,4)node[above]{$0$}--(7,0);
 \draw (9,4)node[above]{$1$}--(8,0);
 \draw (10,4)node[above]{$2$}--(9,0);
 \draw (11,4)node[above]{$3$}--(10,0);
 \draw (12,4)node[above]{$0$}--(11,0);
 \draw (13,4)node[above]{$1$}--(12,0);
 \draw (14,4)node[above]{$2$}--(13,0);
 \draw[dots] (14.2,4)--(15.8,4);
 \draw[dots] (13.2,0)--(14.8,0);
 \draw (16,4)node[above]{$1$}--(15,0);
 \draw (17,4)node[above]{$2$}--(16,0);
 \draw (18,4)node[above]{$3$}--(17,0);
 \draw (19,4)node[above]{$0$}--(18,0);
 \draw (20,4)node[above]{$1$}--(19,0);
 \draw (21,4)node[above]{$2$}--(20,0);
 \node[greendot] at (3,4){};
 \node[greendot] at (10,4){};
 \node[greendot] at (14,4){};
 \node[greendot] at (16,0){};
 \node[greendot] at (21,4){};
\end{braid}
+
\begin{braid}
 \draw (0,4)node[above]{$0$}--(0,0);
 \draw (1,4)node[above]{$1$}--(1,0);
 \draw (2,4)node[above]{$2$}--(2,0);
 \draw (3,4)node[above]{$3$}--(3,0);
 \draw[dots] (3.2,4)--(4.8,4);
 \draw[dots] (3.2,0)--(4.8,0);
 \draw (5,4)node[above]{$1$}--(5,0);
 \draw (6,4)node[above]{$2$}--(6,1)--(16,0);
 \draw (7,4)node[above]{$3$}--(7,1)--(6,0);
 \draw (8,4)node[above]{$0$}--(8,1)--(7,0);
 \draw (9,4)node[above]{$1$}--(9,1)--(8,0);
 \draw (10,4)node[above]{$2$}--(10,1)--(9,0);
 \draw (11,4)node[above]{$3$}--(11,1)--(10,0);
 \draw (12,4)node[above]{$0$}--(12,1)--(11,0);
 \draw (13,4)node[above]{$1$}--(13,1)--(12,0);
 \draw (14,4)node[above]{$2$}--(14,1)--(13,0);
 \draw[dots] (14.2,4)--(15.8,4);
 \draw[dots] (13.2,0)--(14.8,0);
 \draw (16,4)node[above]{$1$}--(16,1)--(15,0);
 \draw (17,4)node[above]{$2$}--(21,3)--(21,0);
 \draw (18,4)node[above]{$3$}--(17,3)--(17,0);
 \draw (19,4)node[above]{$0$}--(18,3)--(18,0);
 \draw (20,4)node[above]{$1$}--(19,3)--(19,0);
 \draw (21,4)node[above]{$2$}--(20,3)--(20,0);
 \node[greendot] at (3,2){};
 \node[greendot] at (10,2){};
 \node[greendot] at (14,2){};
 \node[greendot] at (21,4){};
 \draw[dashed, color = red] (-0.5,3) -- (17.5,3) -- (17.5,1) -- (-0.5,1) -- (-0.5,3);
\end{braid}\\
&\hspace*{-5mm} =_\lambda &
\begin{braid}
 \draw (0,4)node[above]{$0$}--(0,0);
 \draw (1,4)node[above]{$1$}--(1,0);
 \draw (2,4)node[above]{$2$}--(2,0);
 \draw (3,4)node[above]{$3$}--(3,0);
 \draw[dots] (3.2,4)--(4.8,4);
 \draw[dots] (3.2,0)--(4.8,0);
 \draw (5,4)node[above]{$1$}--(5,0);
 \draw (6,4)node[above]{$2$}--(21,0);
 \draw (7,4)node[above]{$3$}--(6,0);
 \draw (8,4)node[above]{$0$}--(7,0);
 \draw (9,4)node[above]{$1$}--(8,0);
 \draw (10,4)node[above]{$2$}--(9,0);
 \draw (11,4)node[above]{$3$}--(10,0);
 \draw (12,4)node[above]{$0$}--(11,0);
 \draw (13,4)node[above]{$1$}--(12,0);
 \draw (14,4)[densely dotted] node[above]{$2$}--(13,0);
 \draw[dots] (14.2,4)--(15.8,4);
 \draw[dots] (13.2,0)--(14.8,0);
 \draw (16,4)node[above]{$1$}--(15,0);
 \draw (17,4)node[above]{$2$}--(16,0);
 \draw (18,4)node[above]{$3$}--(17,0);
 \draw (19,4)node[above]{$0$}--(18,0);
 \draw (20,4)node[above]{$1$}--(19,0);
 \draw (21,4)node[above]{$2$}--(20,0);
 \node[greendot] at (3,4){};
 \node[greendot] at (10,4){};
 \node[greendot] at (14,4){};
 \node[greendot] at (16,0){};
 \node[greendot] at (21,4){};
 \draw[->] (14,4) -- (13.4,1.6);
\end{braid} = \ \ \ \ \ldots\ldots\ldots\ldots\\
&\hspace*{-5mm} =_\lambda &
\begin{braid}
 \draw (0,4)node[above]{$0$}--(0,0);
 \draw (1,4)node[above]{$1$}--(1,0);
 \draw (2,4)node[above]{$2$}--(2,0);
 \draw (3,4)node[above]{$3$}--(3,0);
 \draw[dots] (3.2,4)--(4.8,4);
 \draw[dots] (3.2,0)--(4.8,0);
 \draw (5,4)node[above]{$1$}--(5,0);
 \draw (6,4)node[above]{$2$}--(21,0);
 \draw (7,4)node[above]{$3$}--(6,0);
 \draw (8,4)node[above]{$0$}--(7,0);
 \draw (9,4)node[above]{$1$}--(8,0);
 \draw (10,4)[densely dotted] node[above]{$2$}--(9,0);
 \draw (11,4)node[above]{$3$}--(10,0);
 \draw (12,4)node[above]{$0$}--(11,0);
 \draw (13,4)node[above]{$1$}--(12,0);
 \draw (14,4)node[above]{$2$}--(13,0);
 \draw[dots] (14.2,4)--(15.8,4);
 \draw[dots] (13.2,0)--(14.8,0);
 \draw (16,4)node[above]{$1$}--(15,0);
 \draw (17,4)node[above]{$2$}--(16,0);
 \draw (18,4)node[above]{$3$}--(17,0);
 \draw (19,4)node[above]{$0$}--(18,0);
 \draw (20,4)node[above]{$1$}--(19,0);
 \draw (21,4)node[above]{$2$}--(20,0);
 \node[greendot] at (3,4){};
 \node[greendot] at (10,4){};
 \node[greendot] at (13,0){};
 \node[greendot] at (16,0){};
 \node[greendot] at (21,4){};
 \draw[->] (10,4) -- (9.6,2.4);
\end{braid}\\
&\hspace*{-5mm} \overset{(\ref{dia:y-psi com})}{=_\lambda} &
\begin{braid}
 \draw (0,4)node[above]{$0$}--(0,0);
 \draw (1,4)node[above]{$1$}--(1,0);
 \draw (2,4)node[above]{$2$}--(2,0);
 \draw (3,4)node[above]{$3$}--(3,0);
 \draw[dots] (3.2,4)--(4.8,4);
 \draw[dots] (3.2,0)--(4.8,0);
 \draw (5,4)node[above]{$1$}--(5,0);
 \draw (6,4)node[above]{$2$}--(21,0);
 \draw (7,4)node[above]{$3$}--(6,0);
 \draw (8,4)node[above]{$0$}--(7,0);
 \draw (9,4)node[above]{$1$}--(8,0);
 \draw (10,4)node[above]{$2$}--(9,0);
 \draw (11,4)node[above]{$3$}--(10,0);
 \draw (12,4)node[above]{$0$}--(11,0);
 \draw (13,4)node[above]{$1$}--(12,0);
 \draw (14,4)node[above]{$2$}--(13,0);
 \draw[dots] (14.2,4)--(15.8,4);
 \draw[dots] (13.2,0)--(14.8,0);
 \draw (16,4)node[above]{$1$}--(15,0);
 \draw (17,4)node[above]{$2$}--(16,0);
 \draw (18,4)node[above]{$3$}--(17,0);
 \draw (19,4)node[above]{$0$}--(18,0);
 \draw (20,4)node[above]{$1$}--(19,0);
 \draw (21,4)node[above]{$2$}--(20,0);
 \node[greendot] at (3,4){};
 \node[greendot] at (9,0){};
 \node[greendot] at (13,0){};
 \node[greendot] at (16,0){};
 \node[greendot] at (21,4){};
\end{braid}
+
\begin{braid}
 \draw (0,4)node[above]{$0$}--(0,0);
 \draw (1,4)node[above]{$1$}--(1,0);
 \draw (2,4)node[above]{$2$}--(2,0);
 \draw (3,4)node[above]{$3$}--(3,0);
 \draw[dots] (3.2,4)--(4.8,4);
 \draw[dots] (3.2,0)--(4.8,0);
 \draw (5,4)node[above]{$1$}--(5,0);
 \draw (6,4)node[above]{$2$}--(6,1)--(9,0);
 \draw (7,4)node[above]{$3$}--(7,1)--(6,0);
 \draw (8,4)node[above]{$0$}--(8,1)--(7,0);
 \draw (9,4)node[above]{$1$}--(9,1)--(8,0);
 \draw (10,4)node[above]{$2$}--(21,3)--(21,0);
 \draw (11,4)node[above]{$3$}--(10,3)--(10,0);
 \draw (12,4)node[above]{$0$}--(11,3)--(11,0);
 \draw (13,4)node[above]{$1$}--(12,3)--(12,0);
 \draw (14,4)node[above]{$2$}--(13,3)--(13,0);
 \draw[dots] (14.2,4)--(15.8,4);
 \draw[dots] (13.2,0)--(14.8,0);
 \draw (16,4)node[above]{$1$}--(15,3)--(15,0);
 \draw (17,4)node[above]{$2$}--(16,3)--(16,0);
 \draw (18,4)node[above]{$3$}--(17,3)--(17,0);
 \draw (19,4)node[above]{$0$}--(18,3)--(18,0);
 \draw (20,4)node[above]{$1$}--(19,3)--(19,0);
 \draw (21,4)node[above]{$2$}--(20,3)--(20,0);
 \node[greendot] at (3,2){};
 \node[greendot] at (13,0){};
 \node[greendot] at (16,0){};
 \node[greendot] at (21,4){};
 \draw[dashed, color = red] (-0.5,3) -- (10.5,3) -- (10.5,1) -- (-0.5,1) -- (-0.5,3);
\end{braid}\\
&\hspace*{-5mm} =_\lambda &
\begin{braid}
 \draw (0,4)node[above]{$0$}--(0,0);
 \draw (1,4)node[above]{$1$}--(1,0);
 \draw (2,4)node[above]{$2$}--(2,0);
 \draw (3,4)node[above]{$3$}--(3,0);
 \draw[dots] (3.2,4)--(4.8,4);
 \draw[dots] (3.2,0)--(4.8,0);
 \draw (5,4)node[above]{$1$}--(5,0);
 \draw (6,4)node[above]{$2$}--(21,0);
 \draw (7,4)node[above]{$3$}--(6,0);
 \draw (8,4)node[above]{$0$}--(7,0);
 \draw (9,4)node[above]{$1$}--(8,0);
 \draw (10,4)node[above]{$2$}--(9,0);
 \draw (11,4)node[above]{$3$}--(10,0);
 \draw (12,4)node[above]{$0$}--(11,0);
 \draw (13,4)node[above]{$1$}--(12,0);
 \draw (14,4)node[above]{$2$}--(13,0);
 \draw[dots] (14.2,4)--(15.8,4);
 \draw[dots] (13.2,0)--(14.8,0);
 \draw (16,4)node[above]{$1$}--(15,0);
 \draw (17,4)node[above]{$2$}--(16,0);
 \draw (18,4)node[above]{$3$}--(17,0);
 \draw (19,4)node[above]{$0$}--(18,0);
 \draw (20,4)node[above]{$1$}--(19,0);
 \draw (21,4)[densely dotted] node[above]{$2$}--(20,0);
 \node[greendot] at (3,0){};
 \node[greendot] at (9,0){};
 \node[greendot] at (13,0){};
 \node[greendot] at (16,0){};
 \node[greendot] at (21,4){};
 \draw[->] (21,4) -- (20.5,2);
\end{braid}\\
&\hspace*{-5mm}  \overset{(\ref{dia:y-psi com})}= &
\begin{braid}
 \draw (0,4)node[above]{$0$}--(0,0);
 \draw (1,4)node[above]{$1$}--(1,0);
 \draw (2,4)node[above]{$2$}--(2,0);
 \draw (3,4)node[above]{$3$}--(3,0);
 \draw[dots] (3.2,4)--(4.8,4);
 \draw[dots] (3.2,0)--(4.8,0);
 \draw (5,4)node[above]{$1$}--(5,0);
 \draw (6,4)node[above]{$2$}--(21,0);
 \draw (7,4)node[above]{$3$}--(6,0);
 \draw (8,4)node[above]{$0$}--(7,0);
 \draw (9,4)node[above]{$1$}--(8,0);
 \draw (10,4)node[above]{$2$}--(9,0);
 \draw (11,4)node[above]{$3$}--(10,0);
 \draw (12,4)node[above]{$0$}--(11,0);
 \draw (13,4)node[above]{$1$}--(12,0);
 \draw (14,4)node[above]{$2$}--(13,0);
 \draw[dots] (14.2,4)--(15.8,4);
 \draw[dots] (13.2,0)--(14.8,0);
 \draw (16,4)node[above]{$1$}--(15,0);
 \draw (17,4)node[above]{$2$}--(16,0);
 \draw (18,4)node[above]{$3$}--(17,0);
 \draw (19,4)node[above]{$0$}--(18,0);
 \draw (20,4)node[above]{$1$}--(19,0);
 \draw (21,4)node[above]{$2$}--(20,0);
 \node[greendot] at (3,0){};
 \node[greendot] at (9,0){};
 \node[greendot] at (13,0){};
 \node[greendot] at (16,0){};
 \node[greendot] at (20,0){};
\end{braid}
+
\begin{braid}
 \draw (0,4)node[above]{$0$}--(0,0);
 \draw (1,4)node[above]{$1$}--(1,0);
 \draw (2,4)node[above]{$2$}--(2,0);
 \draw (3,4)node[above]{$3$}--(3,0);
 \draw[dots] (3.2,4)--(4.8,4);
 \draw[dots] (3.2,0)--(4.8,0);
 \draw (5,4)node[above]{$1$}--(5,0);
 \draw (6,4)node[above]{$2$}--(20,0);
 \draw (7,4)node[above]{$3$}--(6,0);
 \draw (8,4)node[above]{$0$}--(7,0);
 \draw (9,4)node[above]{$1$}--(8,0);
 \draw (10,4)node[above]{$2$}--(9,0);
 \draw (11,4)node[above]{$3$}--(10,0);
 \draw (12,4)node[above]{$0$}--(11,0);
 \draw (13,4)node[above]{$1$}--(12,0);
 \draw (14,4)node[above]{$2$}--(13,0);
 \draw[dots] (14.2,4)--(15.8,4);
 \draw[dots] (13.2,0)--(14.8,0);
 \draw (16,4)node[above]{$1$}--(15,0);
 \draw (17,4)node[above]{$2$}--(16,0);
 \draw (18,4)node[above]{$3$}--(17,0);
 \draw (19,4)node[above]{$0$}--(18,0);
 \draw (20,4)node[above]{$1$}--(19,0);
 \draw (21,4)node[above]{$2$}--(21,0);
 \node[greendot] at (3,0){};
 \node[greendot] at (9,0){};
 \node[greendot] at (13,0){};
 \node[greendot] at (16,0){};
\end{braid}\\
&\hspace*{-5mm} = &
\psi_{\lambda_1}\psi_{\lambda_1+1}\ldots\psi_{n-2}\psi_{n-1} e(\bi_{\dot{\lambda}}\vee i)y_{\dot{\lambda}}
+ \psi_{\lambda_1}\psi_{\lambda_1+1}\ldots\psi_{n-2} e(\bi_{\dot{\mu}}\vee i)y_{\dot{\mu}},
\end{eqnarray*}
which completes the proof. \endproof

\begin{Remark}
If $\lambda_1 > \lambda_2$ and $\t$ is the last Garnir tableau of $\lambda$, by \autoref{standard expression of garnir tableau} we have
$$
\psi_{d(\t)}\psi_r = \psi_{a_n,n}\psi_{a_{n-1},n-1}\ldots \psi_{a_{r+2},r+2}\psi_{a_{r+1},r+1}.
$$
Define $\w$ to be the last Garnir tableau of shape $\dot{\lambda}$, we see that $\psi_{d(\w)} = \psi_{a_{n-1},n-1}\ldots \psi_{a_{r+2},r+2}\psi_{a_{r+1},r+1}$. Hence $e(\bi_{\dot{\lambda}}\vee i)y_{\dot{\lambda}}\psi_{d(\w)}\psi_r = \theta_i(\psi_{\t^{\dot{\lambda}} \w}\psi_r)$.
\end{Remark}

\begin{Lemma} \label{common in garnir}
Suppose $\t$ and $\dot\t$ are the last Garnir tableau of shape $\lambda$ and $\dot\lambda$ respectively with last Garnir entry $r$. Set
$$
\psi = \begin{cases}
\psi_{\lambda_1}\psi_{\lambda_1+1}\ldots \psi_{n-2}\psi_{n-1}y_n - \psi_{\lambda_1+1}\ldots \psi_{n-2}\psi_{n-1}, & \text{ if $i = e-1$, $j \neq e-1$.}\\
\psi_{\lambda_1}\psi_{\lambda_1+1}\ldots \psi_{n-2}\psi_{n-1}y_n - \psi_{\lambda_1+1}\ldots \psi_{n-2}\psi_{n-1}-\psi_{\lambda_1}\psi_{\lambda_1+1}\ldots\lambda_{n-2}, & \text{ if $i = j = e-1$.}\\
\psi_{\lambda_1}\psi_{\lambda_2}\ldots\psi_{n-2}\psi_{n-1}, & \text{ otherwise.}
\end{cases}
$$

For any standard $\dot\lambda$-tableau $\dot\v$, if $d(\t) \leq m_\lambda$ and $\dot\v \rhd \dot\t$, then
$$
\begin{cases}
\psi{\cdot} \theta_i(\psi_{\t^{\dot\lambda} \dot\v}) \in\Bgelam[\mu], & \text{ if $i = j = e-2$,}\\
\psi{\cdot} \theta_i(\psi_{\t^{\dot\lambda} \dot\v}) \in \Bgelam,  & \text{ otherwise.}
\end{cases}
$$
\end{Lemma}

\proof If it is not the case that $i = j = e-2$. By \autoref{swap} we have
$$
\psi{\cdot} e(\bi_{\dot\lambda}\vee i) y_{\dot\lambda} =_\lambda e(\bi_\lambda) y_\lambda \psi_{\lambda_1}\psi_{\lambda_1+1}\ldots\psi_{n-2}\psi_{n-1}.
$$

Then we have
$$
\psi{\cdot} \theta_i(\psi_{\t^{\dot\lambda} \dot\v}) = \psi{\cdot} e(\bi_{\dot\lambda}\vee i)y_{\dot\lambda}\psi_{d(\dot\v)} =_\lambda e_\lambda y_\lambda \psi_{\lambda_1}\psi_{\lambda_1+1}\ldots\psi_{n-2}\psi_{n-1}\psi_{d(\dot\v)},
$$
where as $\dot\v \rhd \dot\t$, then $d(\dot\v) < d(\dot\t)$ and
$$
l(\psi_{\lambda_1}\psi_{\lambda_1+1}\ldots\psi_{n-2}\psi_{n-1}\psi_{d(\dot\v)}) < l(\psi_{\lambda_1}\psi_{\lambda_1+1}\ldots\psi_{n-2}\psi_{n-1}) + l(d(\dot\t)) = l(d(t)) \leq m_\lambda.
$$

Then by \autoref{length of psi} we have $\psi{\cdot} \theta_i(\psi_{\t^{\dot\lambda} \dot\v}) \in \Bgelam $.\\

For $i = j = e-2$, set $\dot\mu = (\lambda_1-1,\lambda_2-1,1)$,$\gamma = \lambda|_{n-1} = (\lambda_1,\lambda_2-1)$ and $\dot\gamma = (\lambda_1-1,\lambda_2-1)$. Because $y_{\dot\gamma} = y_{\dot\mu}$. By \autoref{swap},
\begin{eqnarray*}
\psi{\cdot} \theta_i(\psi_{\t^{\dot\lambda}\dot\v}) & = & \psi_{\lambda_1,n}e(\bi_{\dot\lambda}\vee i)y_{\dot\lambda}\psi_{d(\dot\v)} =_\lambda e_\lambda y_\lambda \psi_{\lambda_1,n}\psi_{d(\dot\v)} - \psi_{\lambda_1,n-1}e(\bi_{\dot\mu}\vee i)y_{\dot\mu}\psi_{d(\dot\v)}\\
& = & e_\lambda y_\lambda \psi_{\lambda_1,n}\psi_{d(\dot\v)} - \theta_i(\psi_{\lambda_1,n-1}e(\bi_{\dot\gamma}\vee i)y_{\dot\mu}\psi_{d(\dot\v)})\\
& = &e_\lambda y_\lambda \psi_{\lambda_1,n}\psi_{d(\dot\v)} - \theta_i(\psi_{\lambda_1,n-1}e(\bi_{\dot\gamma}\vee i)y_{\dot\gamma}\psi_{d(\dot\v)}).
\end{eqnarray*}

Again by \autoref{swap}, $\psi_{\lambda_1,n-1} e(\bi_{\dot\gamma}\vee i)y_{\dot\gamma} \psi_{d(\dot\v)} =_{\gamma} e_\gamma y_\gamma\psi_{\lambda_1,n-1}\psi_{d(\dot\v)}$. Since $\gamma = \lambda|_{n-1}$, by \autoref{send},
$$
\theta_i(\psi_{\lambda_1,n-1} e(\bi_{\dot\gamma}\vee i)y_{\dot\gamma} \psi_{d(\dot\v)}) =_\lambda \theta_i(e_\gamma y_\gamma\psi_{\lambda_1,n-1}\psi_{d(\dot\v)}).
$$

Therefore
$$
\psi{\cdot} \theta_i(\psi_{\t^{\dot\lambda}\dot\v}) =_\lambda e_\lambda y_\lambda \psi_{\lambda_1,n}\psi_{d(\dot\v)} - \theta_i(e_\gamma y_\gamma\psi_{\lambda_1,n-1}\psi_{d(\dot\v)}).
$$

As $\lambda\in \mathscr S_n^\Lambda$ and $|\gamma| = n-1 < |\lambda|$
\begin{equation} \label{common in garnir: help1}
e_\gamma y_\gamma\psi_{\lambda_1,n-1}\psi_{d(\dot\v)} = \sum_{\y\in\Shape(\gamma)} c_{\t^\gamma \y} \psi_{\t^{\gamma} y} + \sum_{\x,\y\in\Std(>\gamma)} c_{\x\y}\psi_{\x\y}.
\end{equation}

For the first term of the left hand side of (\ref{common in garnir: help1}), because $\gamma = \lambda|_{n-1}$ and $j = i = e-2$, we have $b_i^\gamma = 2$. By \autoref{I-problem: downstair cases} and the definition of $\gamma$, $\theta_i(\psi_{\t^\gamma \y}) \in \Bgelam[\mu]$. For the second term of the left hand side of (\ref{common in garnir: help1}), as $\x,\y\in\Std(>\gamma) = \Std(>\lambda|_{n-1})$, $\psi_{\x\y} \in\Blam[\lambda|_{n-1}]$. By \autoref{send}, $\theta_i(\psi_{\x\y})\in \Blam \subseteq \Bgelam[\mu]$. Therefore,
$$
\theta_i(e_\gamma y_\gamma\psi_{\lambda_1,n-1}\psi_{d(\dot\v)}) \in \Bgelam[\mu].
$$

Finally, as
$$
l(\psi_{\lambda_1}\psi_{\lambda_1+1}\ldots\psi_{n-2}\psi_{n-1}\psi_{d(\dot\v)}) < l(\psi_{\lambda_1}\psi_{\lambda_1+1}\ldots\psi_{n-2}\psi_{n-1}) + l(d(\dot\t)) = l(d(\t)) \leq m_\lambda,
$$
by \autoref{length of psi} we have $e_\lambda y_\lambda \psi_{\lambda_1,n}\psi_{d(\dot\v)}\in \Blam  \subseteq \Bgelam[\mu]$. Hence $\psi{\cdot}\theta_i(\psi_{\t^{\dot\lambda}\dot\v})\in \Bgelam[\mu]$. This completes the proof. \endproof

\begin{Lemma} \label{common in garnir 2}
Suppose $\lambda_1 - \lambda_2 \equiv e-1 \pmod{e}$, i.e. $i = j$, and $(\lambda_1-\lambda_2+1)\lambda_2 - 1 \leq m_\lambda$. Let $\dot\u$ and $\dot\v$ be standard $\lambda|_{n-1}$-tableaux with $\dot\u \rhd \t^{\dot\lambda}$. Assume $i = j \neq e-2$. Then set
$$
\psi = \begin{cases}
\psi_{\lambda_1}\psi_{\lambda_1+1}\ldots \psi_{n-2}\psi_{n-1}y_n - \psi_{\lambda_1+1}\ldots \psi_{n-2}\psi_{n-1}-\psi_{\lambda_1}\psi_{\lambda_1+1}\ldots\psi_{n-2}, & \text{ if $i = j = e-1$,}\\
\psi_{\lambda_1}\psi_{\lambda_2}\ldots\psi_{n-2}\psi_{n-1}, & \text{ if $i = j \neq e-1, e-2$.}
\end{cases}
$$
and we have
$$
\psi{\cdot}\theta_i(\psi_{\dot\u\dot\v}) \in \Bgelam.
$$
\end{Lemma}

\proof We assume that $i = j \neq e-2$. First we need to introduce some properties of $\dot\u$. Because $\dot\u$ is a standard $\lambda|_{n-1}$-tableau and $\dot\u \rhd \t^{\dot\lambda}$, the only possible choice of $\dot\u$ is that $\dot\u|_{n-2} = \t^{(\lambda_1-1,\lambda_2-1)}$. Define $\u$ and $\v$ to be the unique standard $\lambda$-tableau with $\u|_{n-1} = \dot\u$ and $\v|_{n-1} = \dot\v$, respectively. For example, when $\lambda = (7,4)$ and $e = 3$, then $\dot{\u} = \tab(123456\ten,789)$ and $\u = \tab(123456\ten,789\eleven)$. From the definitions of $\dot\u$, $\dot\v$ and $\u$, $\v$ we see that $d(\v) = d(\dot\v)$ and $l(d(\u)) = l(d(\dot\u)) = \lambda_2 - 1$. Notice that if $i = j \neq e-2$, then $\bi_{\lambda} = \bi_{\lambda|_{n-1}}\vee i$ and $y_{\lambda|_{n-1}} = y_{\lambda}$.

\hh
Now we consider different cases for $i,j$. Suppose $i = j \neq e-1,e-2$, then
\begin{eqnarray*}
\psi{\cdot}\theta_i(\psi_{\dot\u\dot\v}) & = & \psi_{\lambda_1}\psi_{\lambda_2}\ldots\psi_{n-2}\psi_{n-1} \psi_{d(\dot\u)} e(\bi_{\lambda|_{n-1}}\vee i) y_{\lambda|_{n-1}} \psi_{d(\dot\v)}\\
& = & \psi_{\lambda_1}\psi_{\lambda_2}\ldots\psi_{n-2}\psi_{n-1} \psi_{d(\dot\u)} e_\lambda y_\lambda \psi_{d(\dot\v)}.
\end{eqnarray*}

Recall that $e \geq 3$. As $\lambda_1 - \lambda_2 \equiv e-1 \pmod{e}$, we must have $\lambda_1 - \lambda_2 \geq e-1 \geq 2$. Also because of $\lambda_2 \geq 1$,
$$
m_\lambda\geq(\lambda_1-\lambda_2+1)\lambda_2-1 \geq 3\lambda_2-1 \geq 2\lambda_2 > l(\psi_{\lambda_1,n}\psi_{d(\dot\u)}) = 2\lambda_2 - 1.
$$
Hence by \autoref{length of psi}
\begin{eqnarray*}
\psi_{\lambda_1}\psi_{\lambda_2}\ldots\psi_{n-2}\psi_{n-1} \psi_{d(\dot\u)} e_\lambda y_\lambda \psi_{d(\dot\v)} & = &
\psi_{\lambda_1}\psi_{\lambda_2}\ldots\psi_{n-2}\psi_{n-1} \psi_{d(\dot\u)} e_\lambda y_\lambda \psi_{d(\v)}\\
& = &\psi_{\lambda_1}\psi_{\lambda_2}\ldots\psi_{n-2}\psi_{n-1} \psi_{d(\dot\u)} \psi_{\t^\lambda \v} \in \Bgelam.
\end{eqnarray*}

\hh
Suppose $i = j = e-1$, then
\begin{eqnarray*}
&&\psi{\cdot}\theta_i(\psi_{\dot\u\dot\v})\\
 & = & (\psi_{\lambda_1}\psi_{\lambda_1+1}\ldots \psi_{n-2}\psi_{n-1}y_n - \psi_{\lambda_1+1}\ldots \psi_{n-2}\psi_{n-1}-\psi_{\lambda_1}\psi_{\lambda_1+1}\ldots\lambda_{n-2}) \psi_{d(\dot\u)}e(\bi_{\lambda|_{n-1}}\vee i) y_{\lambda|_{n-1}} \psi_{d(\dot\v)}\\
& = & \psi_{\lambda_1,n}y_n\psi_{d(\dot\u)} \psi_{\t^\lambda \v} - \psi_{\lambda_1+1,n}\psi_{d(\dot\u)} \psi_{\t^\lambda \v} - \psi_{\lambda_1,n-1}\psi_{d(\dot\u)}\psi_{\t^\lambda \v}.
\end{eqnarray*}

As $\psi_{d(\dot\u)}$ doesn't involve $\psi_{n-1}$, by \autoref{yup} and \autoref{B_lambda is an ideal},
$$
\psi_{\lambda_1,n}y_n\psi_{d(\dot\u)} \psi_{\t^\lambda \v} = \psi_{\lambda_1,n}\psi_{d(\dot\u)} y_n \psi_{\t^\lambda \v} \in \Blam.
$$

As $l(\psi_{\lambda_1+1,n}\psi_{d(\dot\u)}) = l(\psi_{\lambda_1,n-1}\psi_{d(\dot\u)}) = \lambda_2-1+\lambda_2-1 = 2\lambda_2-2 < m_\lambda$, by \autoref{length of psi}, $\psi_{\lambda_1+1,n}\psi_{d(\dot\u)} \psi_{\t^\lambda \v}$ and $\psi_{\lambda_1,n-1}\psi_{d(\dot\u)}\psi_{\t^\lambda \v}$ are both in $\Bgelam $. Hence $\psi{\cdot}\theta_i(\psi_{\dot\u\dot\v}) \in \Bgelam $. \endproof

Now we are ready to prove that $\psi_{\s\t}\psi_r \in \Bgelam$ when $\Shape(\t)$ has only two rows.

\begin{Proposition}\label{garnir}
Suppose $\lambda\in \mathscr S_n^\Lambda$ and $\t$ is the last Garnir tableau of shape $\lambda$ with $r$ to be the last Garnir entry and $l(d(\t)) \leq m_\lambda$, we have
$$
\psi_{\t^\lambda \t}\psi_r = \sum_{(\u,\v)\rhd (\t^\lambda,\t)}c_{\u \v} \psi_{\u \v}.
$$
\end{Proposition}
\proof By \autoref{standard expression of garnir tableau}, as $\t$ is the last Garnir tableau of shape $\lambda$, we have
$$
\psi_{d(\t)}\psi_r = \psi_{\lambda_1,n}\psi_{\lambda_1-1,n-1}\ldots \psi_{\lambda_2,r+1},
$$
where $l(\psi_{\lambda_1,n}) = l(\psi_{\lambda_1-1,n-1}) = \ldots = l(\psi_{\lambda_2,r+1}) = \lambda_2$. We prove the Proposition by induction on $\lambda_1$. Recall that we write $\dot{\lambda} = (\lambda_1-1,\lambda_2)$, $\mu = (\lambda_1,\lambda_2-1,1)$ and $\dot{\mu} = (\lambda_1-1,\lambda_2-1,1)$.

\hh
When $\lambda_1 = 1$, by definition of Garnir tableau, $\lambda_1 = \lambda_2$. Without loss of generality, we set $\Lambda = \Lambda_0$. In this case $i = 0$ and $j = e-1$. Hence
$$
\psi_{\t^\lambda \t}\psi_r = \psi_{\t^\lambda \t^\lambda}\psi_r = e(0,e-1)\psi_1 = \psi_1 e(e-1,0) = 0 \in \Bgelam.
$$

So, when $\lambda_1 = 1$ the Proposition is true.\\

As that the Proposition holds for any partition of two rows with the length of its first row less than $\lambda_1$, by \autoref{swap} we have
\begin{eqnarray}
\psi_{\t^\lambda \t}\psi_r & = & e_\lambda y_\lambda \psi_{d(\t)} \psi_r = e_\lambda y_\lambda \psi_{\lambda_1,n}\psi_{d(\dot\t)}\psi_r \notag \\
& = & \begin{cases}
\psi_{\lambda_1,n} e(\bi_{\dot{\lambda}}\vee i)y_{\dot{\lambda}}\psi_{d(\dot\t)}\psi_r + \psi_{\lambda_1,n-1}e(\bi_{\dot{\mu}}\vee i)y_{\dot{\mu}}\psi_{d(\dot\t)}\psi_r\\
\hspace*{14mm}=\psi_{\lambda_1,n}\theta_i(\psi_{\t^{\dot\lambda}\dot\t} \psi_r) + \psi_{\lambda_1,n-1}\theta_i(e_{\dot\mu}y_{\dot\mu}\psi_{d(\dot\t)}\psi_r), & \text{ if $i = j = e-2$,}\\
\psi{\cdot}\theta_i(e_{\dot\lambda}y_{\dot\lambda}\psi_{\lambda_1-1,n-1} \ldots \psi_{\lambda_2+1,r+2}\psi_{\lambda_2,r+1}) = \psi {\cdot} \theta_i(\psi_{\t^{\dot\lambda} \dot \t}\psi_r), & \text{ otherwise.}
\end{cases} \label{garnir: general form}
\end{eqnarray}
where $\dot\t$ is the last Garnir tableau with shape $\dot\lambda = (\lambda_1-1,\lambda_2)$, and
$$
\psi = \begin{cases}
\psi_{\lambda_1}\psi_{\lambda_1+1}\ldots \psi_{n-2}\psi_{n-1}y_n - \psi_{\lambda_1+1}\ldots \psi_{n-2}\psi_{n-1}, & \text{ if $i = e-1$, $j \neq e-1$,}\\
\psi_{\lambda_1}\psi_{\lambda_1+1}\ldots \psi_{n-2}\psi_{n-1}y_n - \psi_{\lambda_1+1}\ldots \psi_{n-2}\psi_{n-1} - \psi_{\lambda_1}\psi_{\lambda_1+1}\ldots\psi_{n-2}, & \text{ if $i = j = e-1$,}\\
\psi_{\lambda_1}\psi_{\lambda_2}\ldots\psi_{n-2}\psi_{n-1}, & \text{ otherwise.}
\end{cases}
$$

Now we separate the question into different cases.

\textbf{Case \ref{garnir}a:} $i\neq j$. By (\ref{garnir: general form}) we have
$$
\psi_{\t^\lambda \t}\psi_r = \psi{\cdot} \theta_i(\psi_{\t^{\dot\lambda} \dot\t} \psi_r).
$$

By induction, $\psi_{\t^{\dot\lambda}\dot\t}\psi_r = \sum_{\substack{\dot\v \in\Std(\dot\lambda)\\\dot\v \rhd \dot\t}} c_{\t^{\dot\lambda}\dot\v} \psi_{\t^{\dot\lambda}\dot\v} + \sum_{\dot\u,\dot\v\in\Std(>\dot\lambda)} c_{\dot\u \dot\v} \psi_{\dot\u \dot\v}$. Therefore
$$
\psi_{\t^\lambda \t}\psi_r = \sum_{\substack{\dot\v \in\Std(\dot\lambda)\\\dot\v \rhd \dot\t}} c_{\t^{\dot\lambda}\dot\v} \psi{\cdot}\theta_i(\psi_{\t^{\dot\lambda}\dot\v}) + \sum_{\dot\u,\dot\v\in\Std(>\dot\lambda)} c_{\dot\u \dot\v} \psi{\cdot}\theta_i(\psi_{\dot\u \dot\v}).
$$

For $\dot\u,\dot\v\in\Std(>\dot\lambda)$, by \autoref{residue sequence has to be the same}, $\res(\dot\u) = \res(\t^{\dot\lambda})$. Because $i \neq j$, we always have $\Shape(\dot \u) > \lambda|_{n-1}$. Hence $\psi_{\dot\u\dot\v} \in \Blam[\lambda|_{n-1}]$. Therefore by \autoref{send} and \autoref{B_lambda is an ideal}, $\psi{\cdot}\theta_i(\psi_{\dot\u \dot\v}) \in \Blam$. So $\sum_{\dot\u,\dot\v\in\Std(>\dot\lambda)} c_{\dot\u \dot\v} \psi{\cdot}\theta_i(\psi_{\dot\u \dot\v})\in\Blam$.

For $\dot\v\in\Std(\dot\lambda)$ with $\dot\v\unrhd \dot\t$, by \autoref{common in garnir}, $\psi{\cdot}\theta_i(\psi_{\t^{\dot\lambda}\dot\v}) \in \Bgelam $. Therefore $\sum_{\substack{\dot\v \in\Std(\dot\lambda)\\\dot\v \rhd \dot\t}} c_{\t^{\dot\lambda}\dot\v} \psi{\cdot}\theta_i(\psi_{\t^{\dot\lambda}\dot\v}) \in \Bgelam$. These yield $\psi_{\t^{\lambda}\t}\psi_r \in \Bgelam $.

\hh
\textbf{Case \ref{garnir}b:} $i = j \neq e-2$. By (\ref{garnir: general form}) we have
$$
\psi_{\t^\lambda \t}\psi_r = \psi{\cdot} \theta_i(\psi_{\t^{\dot\lambda} \dot\t} \psi_r).
$$

By induction, $\psi_{\t^{\dot\lambda}\dot\t}\psi_r = \sum_{\substack{\dot\v\in\Std(\dot\lambda)\\\dot\v \rhd \dot\t}}  c_{\t^{\dot\lambda}\dot\v} \psi_{\t^{\dot\lambda}\dot\v} + \sum_{\substack{\dot\u,\dot\v\in\Std(\lambda|_{n-1})\\ (\dot\u,\dot\v) \rhd (\t^{\dot\lambda},\t)}} c_{\dot\u \dot\v} \psi_{\dot\u \dot\v} + \sum_{\dot\u,\dot\v \in \Shape(> \lambda|_{n-1})} c_{\dot\u \dot\v} \psi_{\dot\u \dot\v}$. Therefore
$$
\psi_{\t^\lambda \t}\psi_r = \sum_{\substack{\dot\v\in\Std(\dot\lambda)\\\dot\v \rhd \dot\t}}  c_{\t^{\dot\lambda}\dot\v} \psi{\cdot}\theta_i(\psi_{\t^{\dot\lambda}\dot\v}) + \sum_{\substack{\dot\u,\dot\v\in\Std(\lambda|_{n-1})\\ (\dot\u,\dot\v) \rhd (\t^{\dot\lambda},\t)}} c_{\dot\u \dot\v} \psi{\cdot}\theta_i(\psi_{\dot\u \dot\v}) + \sum_{\dot\u,\dot\v \in \Shape(> \lambda|_{n-1})} c_{\dot\u \dot\v} \psi{\cdot}\theta_i(\psi_{\dot\u \dot\v}).
$$

For $\dot\u,\dot\v\in\Std(>\lambda|_{n-1})$, $\psi_{\dot\u\dot\v} \in\Blam[\lambda|_{n-1}]$. As $\lambda\in \mathscr S_n^\Lambda$, by \autoref{send} we have $\psi{\cdot}\theta_i(\psi_{\dot\u\dot\v}) \in \Blam $. Hence $\sum_{\dot\u,\dot\v \in \Std(> \lambda|_{n-1})} c_{\dot\u \dot\v} \psi{\cdot}\theta_i(\psi_{\dot\u \dot\v}) \in \Blam $.

For $\dot\u,\dot\v\in\Std(\lambda|_{n-1})$ with $\dot\u \rhd \t^{\dot\lambda}$, because $m_\lambda \geq d(\t) = (\lambda_1 - \lambda_2+1)\lambda_2-1$, by \autoref{common in garnir 2} we have $\psi{\cdot}\theta_i(\psi_{\dot\u\dot\v})\in \Bgelam $. So $\sum_{\substack{\dot\u,\dot\v\in\Std(\lambda|_{n-1})\\ (\dot\u,\dot\v) \rhd (\t^{\dot\lambda},\t)}} c_{\dot\u \dot\v} \psi{\cdot}\theta_i(\psi_{\dot\u \dot\v})\in\Bgelam$.

For $\dot\v\in\Std(\dot\lambda)$ with $\dot\v \rhd \dot\t$, by \autoref{common in garnir}, $\psi{\cdot} \theta_i(\psi_{\t^{\dot\lambda} \dot\v}) \in \Bgelam $. So $\sum_{\dot\u,\dot\v \in \Shape(> \lambda|_{n-1})} c_{\dot\u \dot\v} \psi{\cdot}\theta_i(\psi_{\dot\u \dot\v})\in\Bgelam$.

Therefore we have $\psi_{\t^{\lambda}\t}\psi_r \in \Bgelam $.

\hh
\textbf{Case \ref{garnir}c:} $i = j = e-2$. By (\ref{garnir: general form}) we have
\begin{equation} \label{garnir: part c}
\psi_{\t^\lambda \t}\psi_r = \psi_{\lambda_1,n}\theta_i(\psi_{\t^{\dot\lambda}\dot\t} \psi_r) + \psi_{\lambda_1,n-1}\theta_i(e_{\dot\mu}y_{\dot\mu}\psi_{d(\dot\t)}\psi_r).
\end{equation}

For the first term of (\ref{garnir: part c}), by induction,
\begin{eqnarray} \label{garnir: part c: first term: help1}
\psi_{\lambda_1,n}\theta_i(\psi_{\t^{\dot\lambda}\dot\t} \psi_r) & = & \sum_{\substack{\dot\v\in\Std(\dot\lambda)\\\dot\v \rhd \dot\t}}  c_{\t^{\dot\lambda}\dot\v} \psi_{\lambda_1,n}\theta_i(\psi_{\t^{\dot\lambda}\dot\v}) + \sum_{\substack{\dot\u,\dot\v\in\Std(\lambda|_{n-1}) \\ (\dot\u,\dot\v) \rhd (\t^{\dot\lambda},\dot\t)}} c_{\dot\u \dot\v} \psi_{\lambda_1,n}\theta_i(\psi_{\dot\u \dot\v})\\
&& + \sum_{\dot\u,\dot\v\in\Std(> \lambda|_{n-1})} c_{\dot\u \dot\v} \psi_{\lambda_1,n}\theta_i(\psi_{\dot\u \dot\v}).\notag
\end{eqnarray}

For $\dot\v\in\Std(\dot\lambda)$ with $\dot\v\rhd\dot\t$, by \autoref{common in garnir}, we have $\psi_{\lambda_1,n}\theta_i(\psi_{\t^{\dot\lambda}\dot\v})\in \Bgelam[\mu]$. Therefore
\begin{equation} \label{garnir: part c: first term: part1}
\sum_{\substack{\dot\v\in\Std(\dot\lambda)\\\dot\v \rhd \dot\t}}  c_{\t^{\dot\lambda}\dot\v} \psi_{\lambda_1,n}\theta_i(\psi_{\t^{\dot\lambda}\dot\v})\in \Bgelam[\mu].
\end{equation}

\hh
For $\dot\u,\dot\v\in\Std(\lambda|_{n-1})$ with $(\dot\u,\dot\v)\rhd(\t^{\dot\lambda},\dot\t)$, by \autoref{residue sequence has to be the same}, we have $\res(\dot\u) = \bi_{\dot\lambda}$. So the choice of $\dot\u$ is unique, where $d(\dot\u) = \psi_{\lambda_1,n-1}$. Hence as $\bi_\mu = \bi_{\lambda|_{n-1}}\vee i$ and $y_{\lambda|_{n-1}} = y_\mu$
\begin{equation} \label{garnir: part c: first term: part2: help1}
\psi_{\lambda_1,n}\theta_i(\psi_{\dot\u\dot\v}) = \psi_{\lambda_1,n}\psi_{n-1,\lambda_1} e(\bi_{\lambda|_{n-1}}\vee i)y_{\lambda|_{n-1}}\psi_{d(\dot\v)} = \psi_{\lambda_1,n}\psi_{n-1,\lambda_1} e_\mu y_\mu \psi_{d(\dot\v)}.
\end{equation}

We work with $\psi_{\lambda_1,n}\psi_{n-1,\lambda_1} e_\mu y_\mu = \psi_{\lambda_1}\psi_{\lambda_1+1}\ldots\psi_{n-2}\psi_{n-1}\psi_{n-2}\ldots\psi_{\lambda_1+1}\psi_{\lambda_1}e_\mu y_\mu$ first. We define a partition $\sigma = (\lambda_1,\lambda_2-2,1)$. Then
\begin{eqnarray}
&&\psi_{\lambda_1}\ldots \psi_{n-3}\psi_{n-2}\psi_{n-1}\psi_{n-2} \psi_{n-3}\ldots \psi_{\lambda_1} e_\mu y_\mu \notag \\
& = & \psi_{\lambda_1}\ldots \psi_{n-3}\psi_{n-1}\psi_{n-2}\psi_{n-1} \psi_{n-3}\ldots \psi_{\lambda_1} e_\mu y_\mu - \psi_{\lambda_1}\ldots \psi_{n-3} \psi_{n-3}\ldots \psi_{\lambda_1} e_\mu y_\mu\notag \\
& = & \psi_{n-1}\theta_{i-1}(\psi_{\lambda_1}\ldots \psi_{n-3}\psi_{n-2}\psi_{n-3}\ldots \psi_{\lambda_1} e_\sigma y_\sigma) \psi_{n-1} - \psi_{\lambda_1}\ldots \psi_{n-3} \psi_{n-3}\ldots \psi_{\lambda_1} e_\mu y_\mu. \label{garnir: part c: first term: part2: help2}
\end{eqnarray}

Consider the lefthand term in (\ref{garnir: part c: first term: part2: help2}). As $\lambda\in \mathscr S_n^\Lambda$ and $|\sigma| = n-1 < |\lambda|$, we have
$$
\psi_{\lambda_1}\ldots \psi_{n-3}\psi_{n-2}\psi_{n-3}\ldots \psi_{\lambda_1} e_\sigma y_\sigma = \sum_{\u \in \text{Std}(\sigma)} c_{\u \t^\sigma} \psi_{\u \t^\sigma} + \sum_{\u,\v\in\Std(>\sigma)} c_{\u \v} \psi_{\u \v},
$$
where $\res(\u) = \bi_\sigma{\cdot} s_{\lambda_1}s_{\lambda_1+1} \ldots s_{n-3}s_{n-2}s_{n-3}\ldots s_{\lambda_1+1}s_{\lambda_1} = \bi_\sigma$ by \autoref{residue sequence has to be the same}, and $\res(\v) = \bi_\sigma$. Since $\text{min}\{\lambda_1,\ldots,n-2\} = \lambda_1$, by \autoref{does not change by psi}, $c_{\u \t^\sigma} \neq 0$ implies $\u|_{\lambda_1-1} \unrhd \t^{\sigma|_{\lambda_1-1}}$. Then the unique choice for $\u$ is $\u = \t^\sigma$. Hence
$$
\psi_{\lambda_1}\ldots \psi_{n-3}\psi_{n-2}\psi_{n-3}\ldots \psi_{\lambda_1} e_\sigma y_\sigma = c{\cdot} e_\sigma y_\sigma + \sum_{\u,\v\in\Std(>\sigma)} c_{\u \v} \psi_{\u \v}.
$$

Further more if $\u$ is a standard tableau with $\Shape(\u) > \sigma$ and $\res(\u) = \bi_\sigma$, we must have $\Shape(\u) > \lambda|_{n-1}$. Hence by \autoref{send},
$$
\psi_{n-1} \theta_{i-1}(\sum_{\u,\v\in\Std(>\sigma)} c_{\u \v} \psi_{\u \v})\psi_{n-1} \in \Blam.
$$

Therefore
\begin{eqnarray*}
\psi_{n-1}\theta_{i-1}(\psi_{\lambda_1}\ldots \psi_{n-3}\psi_{n-2}\psi_{n-3}\ldots \psi_{\lambda_1} e_\sigma y_\sigma)\psi_{n-1} & =_\lambda & c{\cdot}\psi_{n-1}\theta_{i-1}(e_\sigma y_\sigma) \psi_{n-1}\\
& = & c{\cdot}\psi_{n-1}^2 e_\mu y_\mu = c{\cdot}(e_\lambda y_\lambda - e_\mu y_\mu y_{n-1}).
\end{eqnarray*}

By \autoref{yup} we have $e_\mu y_\mu y_{n-1}\in \Blam $, we have
\begin{equation} \label{garnir: help2}
\psi_{n-1}\theta_{i-1}(\psi_{\lambda_1}\ldots \psi_{n-3}\psi_{n-2}\psi_{n-3}\ldots \psi_{\lambda_1} e_\sigma y_\sigma) \psi_{n-1} =_\lambda c{\cdot}e_\lambda y_\lambda.
\end{equation}

For the righthand term in (\ref{garnir: part c: first term: part2: help2}), as $\lambda\in \mathscr S_n^\Lambda$, $\lambda|_{n-1} \in \mathscr{S}_{n-1}^\Lambda\cap(\P_I \cap \P_y \cap \P_\psi)$. By \autoref{already in P},
$$
\psi_{\lambda_1}\ldots\psi_{n-3}\psi_{n-3}\ldots\psi_{\lambda_1} e_{\lambda|_{n-1}} y_{\lambda|_{n-1}} =_{\lambda|_{n-1}} \sum_{\dot\u\in\Std(\lambda|_{n-1})} c_{\dot\u\t^{\lambda|_{n-1}}} \psi_{\dot\u\t^{\lambda|_{n-1}}}.
$$

Then by \autoref{send},
\begin{eqnarray}
\psi_{\lambda_1}\ldots \psi_{n-3} \psi_{n-3}\ldots \psi_{\lambda_1} e_\mu y_\mu & = & \theta_{i}(\psi_{\lambda_1}\ldots \psi_{n-3} \psi_{n-3}\ldots \psi_{\lambda_1} e_{\lambda|_{n-1}} y_{\lambda|_{n-1}}) \notag \\
& =_{\lambda} & \sum_{\dot\u\in\Std(\lambda|_{n-1})} c_{\dot\u \t^{\lambda|_{n-1}}} \theta_i(\psi_{\dot\u \t^{\lambda|_{n-1}}}) \notag \\
& = & \sum_{\u\in\Std(\mu)} c_{\dot\u \t^{\lambda|_{n-1}}} \psi_{\u\t^\mu} \label{garnir: help3},
\end{eqnarray}
where $\u$ is the unique $\mu$-tableau such that $\u|_{n-1} = \dot\u$.

So substitute (\ref{garnir: help2}) and (\ref{garnir: help3}) to (\ref{garnir: part c: first term: part2: help2}), we have
$$
\psi_{\lambda_1,n}\psi_{n-1,\lambda_1} e_\mu y_\mu =_\lambda \sum_{\u\in\Std(\mu)}c_{\u\t^\mu} \psi_{\u\t^\mu} \pm c{\cdot} \psi_{\t^\lambda \t^\lambda}.
$$

As $\dot\v$ is a standard tableau of shape $\lambda|_{n-1} = \mu|_{n-1}$, we can define $\v_1$ and $\v_2$ to be a standard $\mu$-tableau and $\lambda$-tableau where $\v_1|_{n-1} = \v_2|_{n-1} = \dot\v$, respectively. Henceforth $d(\v_1) = d(\v_2) = d(\dot\v)$ and by (\ref{garnir: part c: first term: part2: help1}),
\begin{eqnarray*}
\psi_{\lambda_1,n}\theta_i(\psi_{\dot\u\dot\v}) & = & \psi_{\lambda_1,n}\psi_{n-1,\lambda_1} e_\mu y_\mu \psi_{d(\dot\v)}=_\lambda \sum_{\u\in\Std(\mu)}c_{\u\t^\mu} \psi_{\u\t^\mu}\psi_{d(\dot\v)} \pm c{\cdot} \psi_{\t^\lambda \t^\lambda}\psi_{d(\dot\v)}\\
& = & \sum_{\u\in\Std(\mu)}c_{\u\t^\mu} \psi_{\u\t^\mu}\psi_{d(\v_1)} \pm c {\cdot} \psi_{\t^\lambda\t^\lambda}\psi_{d(\v_2)}\\
& = & \sum_{\u\in\Std(\mu)}c_{\u\t^\mu} \psi_{\u\v_1} \pm c {\cdot} \psi_{\t^\lambda\v_2} \in \Bgelam[\mu].
\end{eqnarray*}

Therefore,
\begin{equation} \label{garnir: part c: first term: part2}
\sum_{\substack{\dot\u,\dot\v\in\Std(\lambda|_{n-1}) \\ (\dot\u,\dot\v) \rhd (\t^{\dot\lambda},\dot\t)}} c_{\dot\u \dot\v} \psi_{\lambda_1,n}\theta_i(\psi_{\dot\u \dot\v}) \in \Bgelam[\mu].
\end{equation}

Finally, suppose $\dot\u,\dot\v\in\Shape(>\lambda|_{n-1})$, by \autoref{send}, we have $\psi_{\lambda_1,n}\theta_i(\psi_{\dot\u\dot\v}) \in \Blam$. Therefore
\begin{equation} \label{garnir: part c: first term: part3}
\sum_{\dot\u,\dot\v\in\Std(> \lambda|_{n-1})} c_{\dot\u \dot\v} \psi_{\lambda_1,n}\theta_i(\psi_{\dot\u \dot\v}) \in \Blam.
\end{equation}

Substitute (\ref{garnir: part c: first term: part1}), (\ref{garnir: part c: first term: part2}) and (\ref{garnir: part c: first term: part3}) to (\ref{garnir: part c: first term: help1}), we have
\begin{equation} \label{garnir: part c: first term}
\psi_{\lambda_1,n}\theta_i(\psi_{\t^{\dot\lambda}\dot\t} \psi_r) \in \Bgelam[\mu].
\end{equation}

For the second term of (\ref{garnir: part c}), by \autoref{swap}
\begin{eqnarray*}
\psi_{\lambda_1,n-1}\theta_i(e_{\dot\mu}y_{\dot\mu}\psi_{d(\dot\t)}\psi_r) & = & \theta_i(\psi_{\lambda_1,n-1} e_{\dot\mu}y_{\dot\mu}\psi_{d(\dot\t)}\psi_r) = \theta_i(e_{\lambda|_{n-1}}y_{\lambda|_{n-1}} \psi_{\lambda_1,n-1} \psi_{d(\dot\t)} \psi_r)\\
& = & e_\mu y_\mu \psi_{\lambda_1,n-1} \psi_{d(\dot\t)} \psi_r,
\end{eqnarray*}
where by \autoref{psi_n-1 is not involved}, because $\psi_{\lambda_1,n_1}\psi_{d(\dot\t)}\psi_r$ doesn't involve $\psi_{n-1}$,
\begin{equation} \label{garnir: part c: second term}
e_\mu y_\mu \psi_{\lambda_1,n-1} \psi_{d(\dot\t)} \psi_r \in \Bgelam[\mu].
\end{equation}

Therefore substitute (\ref{garnir: part c: first term}) and (\ref{garnir: part c: second term}) to (\ref{garnir: part c}), we have
$$
\psi_{\t^\lambda \t}\psi_r \in \Bgelam[\mu].
$$

Then by \autoref{basis to hecke} the proof is completed. \endproof

\begin{Example} \label{example: garnir c}
We give an example of Case \ref{garnir}c. Suppose $\lambda = (7,4)$, $e = 4$ and $\Lambda = \Lambda_0$. Therefore $i = j = 3$ and
$$
\t = \tab(12379\ten\eleven,4568) \hspace*{11mm} \t^\lambda = \tab(1234567,89\ten\eleven),
$$
with $d(\t) = s_7 s_8 s_9 s_{10} s_6 s_7 s_8 s_9 s_5 s_6 s_7 s_8 s_4 s_5 s_6$ and $r = 7$.

By \autoref{swap} we have
\begin{eqnarray*}
e_\lambda y_\lambda \psi_7 \psi_8 \psi_9 \psi_{10} & = & e(01230123012)y_4 y_{11} \psi_7 \psi_8 \psi_9 \psi_{10}\\
& = & \psi_7 \psi_8 \psi_9 \psi_{10} e(01230130122)y_4y_{11} + \psi_7 \psi_8 \psi_9 e(01230130122)y_4\\
& = & \psi_7 \psi_8 \psi_9 \psi_{10} e(\bi_{\dot\lambda}\vee i) y_{\dot\lambda} + \psi_7 \psi_8 \psi_9 e(\bi_{\dot\mu} \vee i) y_{\dot\mu}\\
& = & \psi_7 \psi_8 \psi_9 \psi_{10} \theta_i(e(\bi_{\dot\lambda}) y_{\dot\lambda}) + \psi_7 \psi_8 \psi_9 \theta_i(e(\bi_{\dot\mu}) y_{\dot\mu}),
\end{eqnarray*}
where $\dot\lambda = (6,4)$ and $\dot\mu = (6,3,1)$. Therefore
$$
\dot\t = \tab(12379\ten\eleven,4568)
$$
and $d(\dot\t) = s_6 s_7 s_8 s_9 s_5 s_6 s_7 s_8 s_4 s_5 s_6$, which indicates
\begin{eqnarray}
&&\psi_{\t^\lambda \t}\psi_r = e_\lambda y_\lambda \psi_7 \psi_8 \psi_9 \psi_{10} \psi_6 \psi_7 \psi_8 \psi_9 \psi_5 \psi_6 \psi_7 \psi_8 \psi_4 \psi_5 \psi_6 \psi_7 \notag\\
& = & \psi_7 \psi_8 \psi_9 \psi_{10} \theta_i(e(\bi_{\dot\lambda}) y_{\dot\lambda})\psi_6 \psi_7 \psi_8 \psi_9 \psi_5 \psi_6 \psi_7 \psi_8 \psi_4 \psi_5 \psi_6 \psi_7\notag \\
&& + \psi_7 \psi_8 \psi_9 \theta_i(e(\bi_{\dot\mu}) y_{\dot\mu})\psi_6 \psi_7 \psi_8 \psi_9 \psi_5 \psi_6 \psi_7 \psi_8 \psi_4 \psi_5 \psi_6 \psi_7\notag \\
& = & \psi_7 \psi_8 \psi_9 \psi_{10} \theta_i(e(\bi_{\dot\lambda}) y_{\dot\lambda})\psi_{d(\dot\t)}\psi_7 + \psi_7 \psi_8 \psi_9 \theta_i(e(\bi_{\dot\mu}) y_{\dot\mu})\psi_{d(\dot\t)} \psi_7 \notag \\
& = & \psi_{7,11} \theta_i(\psi_{\t^{\dot\lambda} \dot\t}\psi_7) + \psi_{7,10} \theta_i(e(\bi_{\dot\mu}) y_{\dot\mu}\psi_{d(\dot\t)} \psi_7) \label{help: garnir c: help1}.
\end{eqnarray}

For the first term of (\ref{help: garnir c: help1}),
\begin{eqnarray} \label{help: garnir c: help1: first term}
\psi_{7,11}\theta_i(\psi_{\t^{\dot\lambda}\dot\t} \psi_r) & = & \sum_{\substack{\dot\v\in\Std(\dot\lambda)\\\dot\v \rhd \dot\t}}  c_{\t^{\dot\lambda}\dot\v} \psi_{7,11}\theta_i(\psi_{\t^{\dot\lambda}\dot\v}) + \sum_{\substack{\dot\u,\dot\v\in\Std(\lambda|_{n-1}) \\ (\dot\u,\dot\v) \rhd (\t^{\dot\lambda},\dot\t)}} c_{\dot\u \dot\v} \psi_{7,11}\theta_i(\psi_{\dot\u \dot\v})\\
&& + \sum_{\dot\u,\dot\v\in\Std(> \lambda|_{n-1})} c_{\dot\u \dot\v} \psi_{7,11}\theta_i(\psi_{\dot\u \dot\v})\notag.
\end{eqnarray}

For $\dot\v \in \Std(\dot\lambda)$ with $\dot\v\rhd\dot\t$, by \autoref{common in garnir} we have
\begin{equation} \label{help: garnir c: help1: first term: part1}
\sum_{\substack{\dot\v\in\Std(\dot\lambda)\\\dot\v \rhd \dot\t}}  c_{\t^{\dot\lambda}\dot\v} \psi_{7,11}\theta_i(\psi_{\t^{\dot\lambda}\dot\v})\in \Bgelam[\mu].
\end{equation}

For $\dot\u, \dot\v \in \Std(\lambda|_{n-1})$ with $(\dot\u,\dot\v) \rhd (\t^{\dot\lambda},\t)$, then $\res(\dot\u) = \bi_{\dot\lambda} = 0123013012$, and because $\Shape(\dot\u) = \lambda|_{n-1} = (7,3)$ with residues
$$
\tab(0123012,301)
$$
and
$$
\dot\u \rhd \t^{\dot\lambda} = \tab(123456,789\ten).
$$

The only possible choice of $\dot\u$ is
$$
\dot\u = \tab(123456\ten,789),
$$
with $d(\dot\u) = s_7 s_8 s_9 = \psi_{\lambda_1,n-1}$. Hence
\begin{equation} \label{help: garnir c: help1: first term: part2: help1}
\psi_{7,11}\theta_i(\psi_{\dot\u \dot\v}) = \psi_7 \psi_8 \psi_9 \psi_{10} \psi_9 \psi_8 \psi_7 e(0123012301) y_4 \psi_{d(\dot\v)}.
\end{equation}

Notice we have
\begin{eqnarray}
&&\psi_7 \psi_8 \psi_9 \psi_{10} \psi_9 \psi_8 \psi_7 e(01230123012) y_4 \notag \\
& = & \psi_7 \psi_8 \psi_9 \psi_{10} \psi_9e(01230130212) \psi_8 \psi_7 y_4 \notag \\
& = & \psi_7 \psi_8 \psi_{10} \psi_9 \psi_{10} e(01230130212) \psi_8 \psi_7 y_4 - \psi_7 \psi_8 e(01230130212) \psi_8 \psi_7 y_4 \notag\\
& = & \psi_7 \psi_8 \psi_{10} \psi_9 \psi_{10} \psi_8 \psi_7 e(01230123012) y_4 - \psi_7 \psi_8 \psi_8 \psi_7 e(01230123012) y_4 \notag\\
& = & \psi_{10} \psi_7 \psi_8 \psi_9 \psi_8 \psi_7 e(01230123021) y_4 \psi_{10} - \psi_7 \psi_8 \psi_8 \psi_7 e(01230123012) y_4 \notag\\
& = & \psi_{10} \theta_1(\psi_7 \psi_8 \psi_9 \psi_8 \psi_7 e(0123012302) y_4) \psi_{10} - \psi_7 \psi_8 \psi_8 \psi_7 e(01230123012) y_4 \notag\\
& = & \psi_{10} \theta_1(\psi_7 \psi_8 \psi_9 \psi_8 \psi_7 e_\sigma y_\sigma) \psi_{10} - \psi_7 \psi_8 \psi_8 \psi_7 e_\mu y_\mu \label{help: garnir c: help1: first term: part2: help2},
\end{eqnarray}
where $\sigma = (7,2,1)$. Consider the left term of (\ref{help: garnir c: help1: first term: part2: help2}), because $|\sigma| < |\lambda|$ and $\lambda\in\mathscr S_n^\Lambda$, we have
$$
\psi_7 \psi_8 \psi_9 \psi_8 \psi_7 e_\sigma y_\sigma = \sum_{\u\in \Std(\sigma)} c_{\u\t^\sigma} \psi_{\u\t^\sigma} + \sum_{\u,\v\in\Std(>\sigma)} c_{\u\v} \psi_{\u\v}.
$$

For $\u \in \Std(\sigma)$, by \autoref{does not change by psi} and $\psi_7 \psi_8 \psi_9 \psi_8 \psi_7$ doesn't involve $\psi_s$ with $s \leq 6$, we have $u|_6 \rhd \t^\sigma|_6$. Then because $\res(\u) = \bi_\sigma{\cdot}s_7 s_8 s_9 s_8 s_7 = 0123012302$, by the definition of $\sigma$
$$
[\sigma] = \ydiag(7,2,1) \hspace*{5mm}\text{with residues}\hspace*{5mm} \tab(0123012,30,2).
$$

Then the only possible choice of $\u$ is $\t^\sigma = \tab(1234567,89,\ten)$. Hence
$$
\psi_7 \psi_8 \psi_9 \psi_8 \psi_7 e_\sigma y_\sigma = c{\cdot} \psi_{\t^\sigma \t^\sigma} + \sum_{\u,\v\in\Std(>\sigma)} c_{\u\v} \psi_{\u\v}.
$$

For $\u,\v\in \Std(>\sigma)$, we have $\res(\u) = \bi_\sigma = 0123012302$. It is impossible that $\Shape(\u) = \lambda|_{n-1}$ because $\bi_{\lambda|_{n-1}} = 0123012301$. Hence $\sum_{\u,\v\in\Std(>\sigma)} c_{\u\v} \psi_{\u\v} = \sum_{\u,\v\in\Std(>\lambda|_{n-1})} c_{\u\v} \psi_{\u\v} \in \Blam[\lambda|_{n-1}]$. So $\psi_7 \psi_8 \psi_9 \psi_8 \psi_7 e_\sigma y_\sigma = c{\cdot} \psi_{\t^\sigma \t^\sigma} + \Blam[\lambda|_{n-1}]$ and hence by \autoref{send},
\begin{eqnarray}
\psi_{10} \theta_1(\psi_7 \psi_8 \psi_9 \psi_8 \psi_7 e_\sigma y_\sigma) \psi_{10} & = & c{\cdot}\psi_{10} \theta_1(\psi_{\t^\sigma \t^\sigma})\psi_{10} + \psi_{10}\theta_1(\Blam[\lambda|_{n-1}])\psi_{10} \notag \\
& = & c{\cdot}\psi_{10} e(01230123021)y_4\psi_{10} + \Blam \notag \\
& =_\lambda & c{\cdot} e(01230123012)y_4 \psi_{10}^2 \notag \\
& = & c{\cdot} e(01230123012) y_4 y_{10} - c{\cdot} e(01230123012) y_4 y_9 \notag \\
& =_\lambda & c{\cdot} e(01230123012)y_4 y_{10} = c{\cdot} e_\lambda y_\lambda  \label{help: garnir c: help1: first term: part2: help3}.
\end{eqnarray}

For the right term of (\ref{help: garnir c: help1: first term: part2: help2}), because $\lambda|_{n-1}\in\mathscr S_{n-1}^\Lambda\cap(\P_I\cap\P_y\cap\P_\psi)$, by \autoref{already in P} we have
$$
\psi_7 \psi_8 \psi_8 \psi_7 e_{\lambda|_{n-1}} y_{\lambda|_{n-1}} = \sum_{\dot\u\in\Std(\lambda|_{n-1})} c_{\dot\u\t^{\lambda|_{n-1}}} \psi_{\dot\u\t^{\lambda|_{n-1}}} + \Blam[\lambda|_{n-1}].
$$

Then by \autoref{send},
\begin{eqnarray}
\psi_7 \psi_8 \psi_8 \psi_7 e_\mu y_\mu & = & \theta_2(\psi_7 \psi_8 \psi_8 \psi_7 e_{\lambda|_{n-1}} y_{\lambda|_{n-1}}) \notag \\
& = & \sum_{\dot\u\in\Std(\lambda|_{n-1})} c_{\dot\u\t^{\lambda|_{n-1}}} \theta_2(\psi_{\dot\u\t^{\lambda|_{n-1}}}) + \theta_2(\Blam[\lambda|_{n-1}]) \notag \\
& = & \sum_{\u\in\Std(\mu)} c_{\dot\u\t^{\lambda|_{n-1}}} \psi_{\u\t^{\mu}} + \Blam \label{help: garnir c: help1: first term: part2: help4}.
\end{eqnarray}

Substitute (\ref{help: garnir c: help1: first term: part2: help3}) and (\ref{help: garnir c: help1: first term: part2: help4}) back to (\ref{help: garnir c: help1: first term: part2: help2}), we have
$$
\psi_7 \psi_8 \psi_9 \psi_{10} \psi_9 \psi_8 \psi_7 e(01230123012) y_4 = \sum_{\u \in\Std(\mu)} c_{\u\t^\mu}\psi_{\u\t^\mu} + c{\cdot}e_\lambda y_\lambda + \Blam.
$$

Recall $\dot\v$ is a standard tableau of shape $\lambda|_{n-1} = \mu|_{n-1}$, we can define $\v_1 \in \Std(\mu)$ and $\v_2\in\Std(\lambda)$ such that $d(\v_1) = d(\v_2) = d(\dot\v)$. Hence by (\ref{help: garnir c: help1: first term: part2: help1}),
\begin{eqnarray*}
\psi_{7,11}\theta_i(\psi_{\dot\u \dot\v}) & = & \psi_7 \psi_8 \psi_9 \psi_{10} \psi_9 \psi_8 \psi_7 e(0123012301) y_4 \psi_{d(\dot\v)}\\
& = & \sum_{\u \in\Std(\mu)} c_{\u\v_1}\psi_{\u\v_1} + c{\cdot}\psi_{\t^\lambda\v_2} + \Blam \in \Bgelam[\mu],
\end{eqnarray*}
which yields
\begin{equation} \label{help: garnir c: help1: first term: part2}
\sum_{\substack{\dot\u,\dot\v\in\Std(\lambda|_{n-1}) \\ (\dot\u,\dot\v) \rhd (\t^{\dot\lambda},\dot\t)}} c_{\dot\u \dot\v} \psi_{7,11}\theta_i(\psi_{\dot\u \dot\v}) \in \Bgelam[\mu].
\end{equation}

Finally, suppose $\dot\u,\dot\v\in\Shape(>\lambda|_{n-1})$, by \autoref{send} we have $\psi_{7,11}\theta_i(\psi_{\dot\u \dot\v})\in \Blam$. Therefore
\begin{equation} \label{help: garnir c: help1: first term: part3}
\sum_{\dot\u,\dot\v\in\Std(> \lambda|_{n-1})} c_{\dot\u \dot\v} \psi_{7,11}\theta_i(\psi_{\dot\u \dot\v}) \in \Blam.
\end{equation}

Substitute (\ref{help: garnir c: help1: first term: part1}), (\ref{help: garnir c: help1: first term: part2}) and (\ref{help: garnir c: help1: first term: part3}) to (\ref{help: garnir c: help1: first term}), we have
\begin{equation} \label{help: garnir c: help1: part1}
\psi_{7,11}\theta_i(\psi_{\t^{\dot\lambda}\dot\t} \psi_r) \in \Bgelam[\mu].
\end{equation}

For the second term of (\ref{help: garnir c: help1}), by \autoref{swap}
\begin{eqnarray*}
\psi_{7,10} \theta_i(e(\bi_{\dot\mu}) y_{\dot\mu}\psi_{d(\dot\t)} \psi_7) & = & \theta_2(\psi_7 \psi_8 \psi_9e(0123013012) y_4\psi_{d(\dot\t)} \psi_7)\\
& = & \theta_2(e(0123012301)y_4 \psi_7 \psi_8 \psi_9 \psi_{d(\dot\t)}\psi_7 )\\
& = & e(01230123012)y_4 \psi_7 \psi_8 \psi_9 \psi_{d(\dot\t)}\psi_7 = e_\mu y_\mu \psi_{7,10}\psi_{d(\dot\t)}\psi_r.
\end{eqnarray*}

Then by \autoref{psi_n-1 is not involved}, because $\psi_{7,10}\psi_{d(\dot\t)}\psi_r$ doesn't involve $\psi_{10}$, we have $e_\mu y_\mu \psi_{7,10}\psi_{d(\dot\t)}\psi_r \in \Bgelam[\mu]$. Therefore
\begin{equation} \label{help: garnir c: help1: part2}
\psi_{7,10} \theta_i(e(\bi_{\dot\mu}) y_{\dot\mu}\psi_{d(\dot\t)} \psi_7) \in \Bgelam[\mu].
\end{equation}

Substitute (\ref{help: garnir c: help1: part1}) and (\ref{help: garnir c: help1: part2}) to (\ref{help: garnir c: help1}), we have $\psi_{\t^\lambda \t}\psi_r \in \Bgelam[\mu]$. Finally by \autoref{basis to hecke}, we have
$$
\psi_{\t^\lambda \t}\psi_r = \sum_{(\u,\v)\rhd (\t^\lambda,\t)}c_{\u \v} \psi_{\u \v}.
$$
\end{Example}

Finally, we can extend the above Proposition to arbitrary multipartition using arguments similar to those we used in the last section.

\begin{Corollary}
Suppose $\lambda\in \mathscr S_n^\Lambda$ and $\t$ is the last Garnir tableau of shape $\lambda$ with $r$ the last Garnir entry and $l(d(\t)) \leq m_\lambda$. Therefore, for any standard $\lambda$-tableau $\s$, $\psi_{\s\t} \psi_r = \sum_{(\u,\v)\rhd (\s,\t)}c_{\u\v} \psi_{\u\v}$.
\end{Corollary}

\proof Write $\lambda = (\lambda^{(1)},\ldots,\lambda^{(\l)})$ and $\lambda^{(\l)} = (\lambda_1^{(\l)},\ldots,\lambda_k^{(\l)})$. If $\lambda^{(\l)} = \emptyset$, then define $\bar\lambda = (\lambda^{(1)},\ldots,\lambda^{(\l-1)})$. As $l(\bar\lambda) = \l-1 < l(\lambda)$, we have $\bar\lambda \in \mathscr P^{\bar\Lambda}_I \cap \mathscr P^{\bar\Lambda}_y \cap \mathscr P^{\bar\Lambda}_\psi$. By \autoref{remove empty tableaux}, $\lambda \in \P_I \cap \P_y \cap \P_\psi$, and the Corollary follows.

Now suppose $\lambda^{(\l)} \neq \emptyset$, First we assume $\s = \t^\lambda$. As $\t$ is the last Garnir tableau of shape $\lambda$, $k \geq 2$. Setting $m = \lambda_{k-1}^{(\l)} + \lambda_k^{(\l)}$. As $\t$ is the last Garnir tableau, by the definition we see that $\t|_{n - m} = \t^\lambda|_{n - m}$ and $k \geq 2$. Define $i$ to be the residue of the node $(k-1,1,\l)$, $\Lambda' = \Lambda_i$, and $\tilde\t$ to be the last Garnir tableau of shape $(\lambda_{k-1}^{(\l)},\lambda_k^{(\l)})$. If we write $\mu = (\lambda^{(1)},\ldots,\lambda^{(\l-1)},\mu^{(\l)})$ with $\mu^{(\l)} = (\lambda_1^{(\l)},\ldots,\lambda_{k-2}^{(\l)})$ and $\gamma = (\lambda_{k-1}^{(\l)},\lambda_k^{(\l)})$, then
$$
\psi_{\t^\lambda \t} \psi_r = \hat\theta_{\bi_\mu}(\hat\psi_{\t^\gamma \tilde\t} \hat\psi_{r-(n-m)})y_\mu.
$$

Recall that $\hat\psi_{\t^\gamma \tilde\t}$ and $\hat\psi_{r-(n-m)}$ are elements of $\mathscr{R}_m$ and $\psi_{\t^\gamma \tilde\t}$ and $\psi_{r-(n-m)}$ are elements of $\mathscr R_m^{\Lambda'}$. Then by \autoref{garnir}, we have $\psi_{\t^\gamma \tilde\t} \psi_{r-(n-m)}\in \Bgelam[\gamma]$. Therefore we can write $\psi_{\t^\gamma \tilde\t} \psi_{r-(n-m)} = \sum_{\u,\v\in\Std(\gamma)} c_{\u\v}\psi_{\u\v} + \sum_{\u,\v\in\Std(>\gamma)} c_{\u\v} \psi_{\u\v}$ and hence $\hat\psi_{\t^\gamma \tilde\t} \hat\psi_{r-(n-m)} = \sum_{\u,\v\in\Std(\gamma)} c_{\u\v}\hat\psi_{\u\v} + \sum_{\u,\v\in\Std(>\gamma)} c_{\u\v} \hat\psi_{\u\v} + r$ where $r \in N_m^{\Lambda'}$. Therefore
$$
\psi_{\t^\lambda \t} \psi_r  = \sum_{\u,\v\in\Std(\gamma)} c_{\u\v}\hat\theta_{\bi_\mu}(\hat\psi_{\u\v})y_\mu + \sum_{\u,\v\in\Std(>\gamma)} c_{\u\v} \hat\theta_{\bi_\mu}(\hat\psi_{\u\v})y_\mu + \hat\theta_{\bi_\mu}(r)y_\mu.
$$

For $\u,\v\in\Std(\gamma)$, by \autoref{concatenation: cor3} we have $\hat\theta_{\bi_\mu}(\psi_{\u\v}) \in \Bgelam[\lambda]$. Hence $\sum_{\u,\v\in\Std(\gamma)} c_{\u\v}\hat\theta_{\bi_\mu}(\hat\psi_{\u\v})y_\mu\in\Bgelam$.

For $\u,\v\in\Std(>\gamma)$, write $\Shape(\u) = \Shape(\v) = \sigma$ and $\nu = \mu\vee\sigma$. By \autoref{concatenation: cor1} we have $\nu > \lambda = \mu\vee\gamma$. Then by \autoref{composition in Blam} and \autoref{send}, $\hat\theta_{\bi_\mu}(\psi_{\u\v})\in \Blam$. Hence $\sum_{\u,\v\in\Std(>\gamma)} c_{\u\v} \hat\theta_{\bi_\mu}(\hat\psi_{\u\v})y_\mu\in\Blam$.

Finally by \autoref{theta_j: go up}, $\hat\theta_{\bi_\mu} (r)y_\mu \in \Blam$. These yield that
$$
\psi_{\t^\lambda \t}\psi_r = \hat\theta_{\bi_\mu}(\psi_{\t^\gamma \tilde\t} \psi_{r-(n-m)})y_\mu \in \Bgelam.
$$

Now choose any $\s\in \Std(\lambda)$. Because $\psi_{\t^\lambda \t}\psi_r \in \Bgelam $, we have
$$
\psi_{\t^\lambda \t} \psi_r=_\lambda \sum_{\substack{\v\in \Std(\lambda)}} c_{\t^\lambda \v} \psi_{\t^\lambda \v}.
$$

Hence
$$
\psi_{\s\t} \psi_r= \psi^*_{d(\s)}\psi_{\t^\lambda \t}\psi_r =_\lambda \sum_{\substack{\v\in \Std(\lambda)}} c_{\t^\lambda \v} \psi^*_{d(\s)}\psi_{\t^\lambda \v} = \sum_{\substack{\v\in \Std(\lambda)}} c_{\t^\lambda \v} \psi_{\s\v}.
$$

Therefore, $\psi_{\s\t}\psi_r \in \Bgelam $. By \autoref{basis to hecke} we completes the proof.\endproof

\begin{Remark}
Generally it is not easy to find $c_{\u\v}$. Kleshchev-Mathas-Ram~\cite{KMR:UniversalSpecht} explicitly describes how to compute $c_{\u\v}$ where $\Shape(\u) = \Shape(\v) = \lambda$. This paper also gives an implicit method to compute these coefficients.
\end{Remark}

\subsection{Completion of the $\psi$-problem}

In this subsection we are going to prove that $\psi_{\s\t}\psi_r \in \Bgelam$. \autoref{psi problem: part2} shows that if $\t{\cdot}s_r$ is standard and $d(\t){\cdot}s_r$ is reduced then our claim is true. It remains to consider the case when $\t{\cdot}s_r$ is not standard or $d(\t){\cdot}s_r$ is not reduced.

First we introduce two special conditions on $\t \in \Std(\lambda)$ and $r$ for $1 \leq r \leq n-1$.

\begin{Definition} \label{definition: unlock}
Suppose $\t$ is a standard $\lambda$-tableau and $1 \leq r \leq n-1$. 

(a). If there exists a reduced expression $s_{r_1}s_{r_2}\ldots s_{r_{l-1}} s_{r_l}$ of $d(\t)$ such that $|r - r_l| > 1$, then $\t$ is \textbf{unlocked by $s_r$ in type I}.

(b). If there exists a reduced expression $s_{r_1}s_{r_2}\ldots s_{r_{l-1}} s_{r_l}$ of $d(\t)$ such that $r_{l-1} = r$ and $r = r_l \pm 1$, then $\t$ is \textbf{unlocked by $s_r$ in type II}.
\end{Definition}

The following Lemmas show that if $l(\t) \leq m_\lambda$ and $\t$ is unlocked by $s_r$ in type I or type II, then we have $\psi_{\s\t} \psi_r \in R_n^{\geq \lambda}$.

\begin{Lemma}\label{psiproblem:swap1}
Suppose $\t$ and $\s$ are two standard $\lambda$-tableaux with $d(\t) = d(\s){\cdot}s_k$ for some $k$ and $l(d(\t)) = l(d(\s)) + 1$. If for some $r\not\in\{ k-1, k, k+1\}$, $\t{\cdot}s_r$ is not standard or $d(\t){\cdot}s_r$ is not reduced, then $\s{\cdot}s_r$ is not standard or $d(\s){\cdot}s_r$ is not standard, respectively.
\end{Lemma}
\proof When $\t{\cdot}s_r$ is not standard, then $r$ and $r+1$ in $\t$ are adjacent, either in the same row or in the same column. Since $d(\t) = d(\s){\cdot}s_k$, we have $\s = \t{\cdot}s_k$. As $r\not\in \{k-1,k,k+1\}$, $r$ and $r+1$ are in the same nodes in $\t$ as in $\s$. Hence $\s{\cdot}s_r$ is not standard as well.

When $d(\t){\cdot}s_r$ is not reduced, $d(\t)(r) > d(\t)(r+1)$. As $d(\t) = d(\s){\cdot}s_k$ we have $d(\t)(r) = d(\s){\cdot}s_k(r) = d(\s)(r)$ and $d(\t)(r+1) = d(\s){\cdot}s_k(r+1) = d(\s)(r+1)$ as $r\not\in\{k-1,k,k+1\}$. Hence $d(\s)(r) > d(\s)(r+1)$. Therefore, $d(\s){\cdot}s_r$ is not reduced. This completes the proof.\endproof

\begin{Example} \label{example: unlock}
Suppose $\t = \tab(1279,358,46\ten)$ and $\s = \tab(1278,349,56\ten)$. Then
\begin{eqnarray*}
d(\t) & = & s_4 s_5 s_6 s_7 s_8 s_6 s_7 s_3 s_4 s_5 s_6 s_4,\\
d(\s) & = & s_7 s_8 s_4 s_5 s_6 s_7 s_3 s_4 s_5 s_6 s_4.
\end{eqnarray*}

Set $k = 8$, we have $d(\t) = d(\s) {\cdot}s_8$ and $l(d(\t)) = l(d(\s)) + 1$. Let $r = 3$, i.e. $r \not\in\{k-1,k,k+1\}$ and $\t{\cdot}s_r$ is not standard. We see that $\s{\cdot}s_r$ is not standard either. Similarly, let $r = 6$, i.e. $r \not\in \{k-1,k,k+1\}$ and $d(\t){\cdot}s_r$ is not reduced. Hence $d(\s){\cdot}s_r$ is not reduced either.
\end{Example}

\begin{Lemma} \label{psiproblem: unlock}

Suppose $\lambda\in \mathscr S_n^\Lambda$ and $\t$ is a standard $\lambda$-tableau with $l(d(\t)) \leq m_\lambda$. If $\t$ is unlocked by $s_r$ in type I, then $\psi_{\s\t}\psi_r \in \Bgelam $ for any standard $\lambda$-tableau $\s$.

\end{Lemma}

\proof Suppose $\t{\cdot}s_r$ is standard and $d(\t){\cdot}s_r$ is reduced, by \autoref{psi problem: part2}, $\psi_{\s\t}\psi_r \in \Bgelam$.

Suppose $\t{\cdot}s_r$ is not standard or $d(\t){\cdot}s_r$ is not reduced. Since $\t$ is a standard $\lambda$-tableau unlocked by $s_r$ in type I, by \autoref{definition: unlock}, there exists a reduced expression $s_{r_1}s_{r_2}\ldots s_{r_l}$ of $d(\t)$ such that $|r - r_l| > 1$. Define $\w = \t^\lambda{\cdot}s_{r_1}s_{r_2}\ldots s_{r_{l-1}}$. By \autoref{standard: short}, $\w$ is a standard $\lambda$-tableau. It is easy to see that $d(\t) = d(\w){\cdot}s_{r_l}$ and $l(d(\t)) = l(d(\w)) + 1$. Hence by \autoref{psiproblem:swap1}, $\w{\cdot}s_r$ is not standard or $d(\w){\cdot}s_r$ is not reduced. So $\psi_{\s\w}\psi_r =_\lambda \sum_{\substack{\v\in\Std(\lambda)\\\v\rhd\w}}c_{\s\v}\psi_{\s\v}$ because $l(d(\w)) = l(d(\t)) - 1 < w_\lambda$.

We can write $d(\t) = d(\w){\cdot}s_{r_l}$ as a reduced expression. By \autoref{basis change}, we have
$$
\sum_{\substack{\v\in\Std(\lambda)\\ \v\rhd\t}}c_{\s\v}\psi_{\s\v}\psi_r =_\lambda \psi_{\s\t}\psi_r - \psi_{d(\s)}^* e_\lambda y_\lambda \psi_{d(\w)}\psi_{r_l}\psi_r = \psi_{\s\t}\psi_r - \psi_{\s\w}\psi_r\psi_{r_l}.
$$

Because $\v \rhd \t$ and $l(d(\v)) < l(d(\t)) \leq m_\lambda$, we have $\sum_{\substack{\v\in\Std(\lambda)\\ \v\rhd\t}}c_{\s\v}\psi_{\s\v}\psi_r \in \Bgelam$. We can write $\psi_{\s\w}\psi_r\psi_{r_l} =_\lambda \sum_{\substack{\v\in\Std(\lambda)\\\v\rhd\w}}c_{\s\v}\psi_{\s\v}\psi_{r_l}$. Because $\v \rhd \w$, we have $l(d(\v)) < l(d(\w)) < m_\lambda$, which yields that $\psi_{\s\w}\psi_r\psi_{r_l} \in \Bgelam$. Therefore we have $\psi_{\s\t}\psi_r \in \Bgelam$. \endproof

\begin{Lemma}\label{psiproblem:swap2}
Suppose $\t$ is a standard $\lambda$-tableau and that there exists a standard $\lambda$-tableau $\w$ such that $d(\t) = d(\w){\cdot}s_r s_{r+1}$ for some $r$ and $l(d(\t)) = l(d(\w)) + 2$. If $\t{\cdot}s_r$ is not standard or $d(\t){\cdot}s_r$ is not reduced, then $\w{\cdot}s_{r+1}$ is not standard or $d(\w){\cdot}s_{r+1}$ is not reduced, respectively.

Similarly suppose $d(\t) = d(\w){\cdot}s_r s_{r-1}$ for some $r$ and $l(d(\t)) = l(d(\w)) + 2$. If $\t{\cdot}s_r$ is not standard or $d(\t){\cdot}s_r$ is not reduced, then $\w{\cdot}s_{r-1}$ is not standard or $d(\w){\cdot}s_{r-1}$ is not reduced, respectively.
\end{Lemma}

\proof Suppose $d(\t) = d(\w){\cdot}s_r s_{r+1}$. If $\t{\cdot}s_r$ is not standard, $r$ and $r+1$ are adjacent in $\t$. But $r$ and $r+1$ occupy the same positions as $r+1$ and $r+2$, respectively in $\w$. So $\w{\cdot}s_{r+1}$ is not standard. If $d(\t){\cdot}s_r$ is not reduced, as $d(\w)^{-1}(r+1) = d(\t)^{-1}(r)$ and $d(\w)^{-1}(r+2) = d(\t)^{-1}(r+1)$, by \autoref{change length of w}, then $d(\t){\cdot}s_r$ is not reduced implies $d(\w){\cdot}s_{r+1}$ is not reduced. The other case is similar.

\begin{Remark}
In \autoref{psiproblem:swap1} and \autoref{psiproblem:swap2}, when we say $d(\t) = d(\s){\cdot}s_r$ or $d(\t) = d(\s){\cdot}s_r s_{r+1}$, it means $d(\t)$ and $d(\s){\cdot}s_r$ or $d(\t)$ and $d(\s){\cdot}s_r s_{r+1}$ are the same as permutations.
\end{Remark}

\begin{Example} \label{example: unlock on tails}
Let $\t = \tab(123\twelve,456\thirteen,78\eleven,9\ten\fourteen)$. Suppose $\s = \tab(123\twelve,456\thirteen,789,\ten\eleven\fourteen)$, we have
\begin{eqnarray*}
d(\t) & = & s_8 s_9 s_{10}s_{11}s_{12}s_4 s_5 s_6 s_7 s_8 s_9 s_{10}s_{11}s_9 s_{10},\\
d(\s) & = & s_8 s_9 s_{10}s_{11}s_{12}s_4 s_5 s_6 s_7 s_8 s_9 s_{10}s_{11}.
\end{eqnarray*}

So we have $d(\t) = d(\s) s_9 s_{10}$ and therefore $r = 9$. We see that $\t{\cdot}s_r$ and $\s{\cdot}s_{r+1}$ are both non-standard.

Suppose $\s = \tab(123\eleven,456\thirteen,78\ten,9\twelve\fourteen)$, we have
\begin{eqnarray*}
d(\t) & = & s_8 s_9 s_{10}s_{11}s_{12}s_4 s_5 s_6 s_7 s_8 s_9 s_{10}s_{11}s_9 s_{10},\\
d(\s) & = & s_8 s_9 s_{10}s_{11}s_{12}s_4 s_5 s_6 s_7 s_8 s_9 s_{10}s_{9}.
\end{eqnarray*}

So we have $d(\t) = d(\s) s_{11}s_{10}$ and therefore $r = 11$. We see that $d(\t){\cdot}s_r$ and $d(\s){\cdot}s_{r-1}$ are both non-reduced because in $\t$, $r$ is below $r+1$ and in $\s$, $r-1$ is below $r$.
\end{Example}

%
%
%
%

\begin{Lemma} \label{psiproblem: unlock on tails}
Suppose $\lambda\in \mathscr S_n^\Lambda$ and $\t$ is a standard $\lambda$-tableau with $l(d(\t)) \leq m_\lambda$. If $\t$ is unlocked by $s_r$ in type II, then $\psi_{\s\t}\psi_r \in \Bgelam $ for any standard $\lambda$-tableau $\s$.

\end{Lemma}

\proof Suppose $\t{\cdot}s_r$ is standard and $d(\t){\cdot}s_r$ is reduced. Then, by \autoref{psi problem: part2}, $\psi_{\s\t}\psi_r \in \Bgelam$.

Suppose $\t{\cdot}s_r$ is not standard or $d(\t){\cdot}s_r$ is not reduced. Then, by \autoref{definition: unlock}, there exists a reduced expression $s_{r_1}s_{r_2}\ldots s_{r_{l-1}} s_{r_l}$ for $d(\t)$ such that $r_{l-1} = r$ and $r = r_l \pm 1$. Without loss of generality, set $r = r_l - 1$. Define $\w = \t^{\lambda}s_{r_1}s_{r_2}\ldots s_{r_{l-2}}$. By \autoref{standard: short}, $\w$ is a standard $\lambda$-tableau. It is easy to see that $d(\t) = d(\w){\cdot}s_{r}s_{r+1}$ and $l(d(\t)) = l(d(\w)) + 2$. Hence, by \autoref{psiproblem:swap2}, $\w{\cdot}s_{r+1}$ is not standard or $d(\w){\cdot}s_{r+1}$ is not reduced. So $\psi_{\s\w}\psi_{r+1} =_\lambda \sum_{\substack{\v\in\Std(\lambda)\\\v\rhd\w}}c_{\s\v}\psi_{\s\v}$ because $l(d(\w)) = l(d(\t)) - 2 < w_\lambda$.

Because $d(\t) = d(\w){\cdot}s_{r}s_{r+1}$ as a reduced expression, by \autoref{basis change}, we have
\begin{equation}\label{psiproblem: unlock on tails: 1}
\sum_{\substack{\v\in\Std(\lambda)\\ \v\rhd\t}}c_{\s\v}\psi_{\s\v}\psi_r =_\lambda \psi_{\s\t}\psi_r - \psi_{d(\s)}^* e_\lambda y_\lambda \psi_{d(\w)}\psi_r \psi_{r+1}\psi_r = \psi_{\s\t}\psi_r - \psi_{\s\w}\psi_r\psi_{r+1}\psi_{r}.
\end{equation}

Now $\v \rhd \t$, so $l(d(\v)) < l(d(\t)) \leq m_\lambda$. Hence we have
\begin{equation}\label{psiproblem: unlock on tails: 2}
\sum_{\substack{\v\in\Std(\lambda)\\ \v\rhd\t}}c_{\s\v}\psi_{\s\v}\psi_r \in \Bgelam.
\end{equation}

Let $\res(\w) = i_1 i_2 \ldots i_n$ be the residue sequence of $\w$. Then
$$
\psi_{\s \w}\psi_r\psi_{r+1}\psi_r = \begin{cases}
\psi_{\s \w}\psi_{r+1}\psi_r\psi_{r+1} \pm \psi_{\s\w}, & \text{ if $i_r = i_{r+2} = i_{r+1}\pm 1$,}\\
\psi_{\s \w}\psi_{r+1}\psi_r\psi_{r+1}, & \text{ otherwise.}
\end{cases}
$$

Because $\psi_{\s\w}\psi_{r+1} =_\lambda \sum_{\substack{\v\in\Std(\lambda)\\\v\rhd\w}}c_{\s\v}\psi_{\s\v}$,
$$
\psi_{\s\w}\psi_{r+1}\psi_r\psi_{r+1} =_\lambda \sum_{\v \rhd \w} c_{\s\v} \psi_{\s\v}\psi_r\psi_{r+1}.
$$

Since $\v \rhd \w$, $l(d(\v)) < l(d(\w)) = l(d(\t)) - 2 \leq m_\lambda - 2$. Hence, $l(\psi_{d(\v)}\psi_r\psi_{r+1}) = l(d(\v)) + 2 < m_\lambda$. By \autoref{length of psi} we have $\psi_{\s\v}\psi_r\psi_{r+1} \in \Bgelam $ if $\v \rhd \w$. Therefore, we always have
\begin{equation}\label{psiproblem: unlock on tails: 3}
\psi_{\s\w}\psi_{r+1}\psi_r\psi_{r+1} \in \Bgelam
\end{equation}
in both cases. Substituting (\ref{psiproblem: unlock on tails: 2}) and (\ref{psiproblem: unlock on tails: 3}) into (\ref{psiproblem: unlock on tails: 1}) shows that $\psi_{\s\t}\psi_r \in \Bgelam$.\endproof

The following Lemmas are technical results which we will use later.

\begin{Lemma} \label{help: case 3.5}
Suppose $\t\in\Std(\lambda)$ with $d(\t) = s_{n-1}s_{n-2}\ldots s_{r+1}$, and $\t{\cdot}s_r$ is not standard. Then $\t$ is the last Garnir tableau of shape $\lambda$.

\end{Lemma}

\proof As $d(\t)$ is the standard expression, we have $w_n = s_{n-1}$, $w_{n-1} = s_{n-2}, \ldots, w_{r+2} = s_{r+1}$ and $\t = \t^{(1)} = \t^{(2)} = \ldots = \t^{(r+1)}$. Write $\lambda = (\lambda^{(1)},\ldots,\lambda^{(\l)})$ and $\lambda^{(\l)} = (\lambda^{(\l)}_1,\ldots,\lambda^{(\l)}_k)$, as $\t^{(n)} = \t^{\lambda}{\cdot}w_n = \t^{\lambda}{\cdot}s_{n-1}$ is standard, $n-1$ and $n$ are not adjacent in $\t^\lambda$. This forces $\lambda^{(\l)}_k = 1$.

By (\ref{remark: standard expression}) we have
$$
\t^{-1}(k) = \begin{cases}
(\t^\lambda)^{-1}(k-1), & \text{ if $r+2 \leq k \leq n$,}\\
(\t^\lambda)^{-1}(n), &\text{ if $k = r+1$,}\\
(\t^\lambda)^{-1}(k), & \text{ otherwise.}
\end{cases}
$$
and since $\t{\cdot}s_r$ is not standard, $r$ and $r+1$ are adjacent in $\t$. As $\t^{-1}(r+1) = (\t^\lambda)^{-1}(n) = (k,\lambda^{(\l)}_k,\l) = (k,1,\l)$, we must have $\t^{-1}(r) = (k-1,1,\l)$. This shows that $\t$ is the last Garnir tableau with shape $\lambda$.\endproof

\begin{Example} \label{example: case 3.5}
Suppose $\lambda = (4,4,1)$ and $\t = \tab(1234,5789,6)$. Therefore $d(\t) = s_8 s_7 s_6$ and $\t{\cdot}s_5$ is not standard. Notice that $\t$ is the last Garnir tableau of shape $\lambda$.
\end{Example}

\begin{Lemma} \label{help: case 4.3}

Suppose $\t\in \Std(\lambda)$ with $d(\t) = s_r s_{r+1} \ldots s_{n-2}$, and $\t{\cdot}s_{n-1}$ is not standard. Then $\t$ is the last Garnir tableau of shape $\lambda$.

\end{Lemma}

\proof As $d(\t)$ is the standard expression, we have $w_n = s_r s_{r+1} \ldots s_{n-2}$ and $w_{n-1} = \ldots w_1 = 1$. Write $\lambda = (\lambda^{(1)},\ldots,\lambda^{(\l)})$ and $\lambda^{(\l)} = (\lambda^{(\l)}_1,\ldots,\lambda^{(\l)}_k)$, as $\t = \t^\lambda{\cdot}w_n$, we have $\t^{-1}(n) = (k,\lambda^{(\l)}_k,\l)$. By (\ref{remark: standard expression}) we have
$$
\t^{-1}(k) =
\begin{cases}
(\t^\lambda)^{-1}(k+1), & \text{ if $r \leq k \leq n-2$,}\\
(\t^\lambda)^{-1}(r),  & \text{ if $k = n - 1$,}\\
(\t^\lambda)^{-1}(k), & \text{ otherwise.}
\end{cases}
$$

As $\t{\cdot}s_{n-1}$ is not standard, $n-1$ and $n$ are adjacent in $\t$. So in $r$ and $n$ are adjacent in $\t^\lambda$. But $r \leq n-2$. Hence $r$ has to be on the above of $n$ in $\t^\lambda$. i.e. $(\t^\lambda)^{-1}(r) = \t^{-1}(n-1) = (k-1,\lambda^{(\l)}_k,\l)$. This shows that $\lambda^{(\l)}_{k-1} = \lambda^{(\l)}_k$ and $\t$ is the last Garnir tableau of shape $\lambda$. \endproof

\begin{Example} \label{example: case 4.3}
Suppose $\lambda = (4,3,3)$ and $\t = \tab(1234,569,89\ten)$. Therefore $d(\t) = s_7 s_8$ and $\t{\cdot}s_9$ is not standard. Notice that $\t$ is the last Garnir tableau of shape $\lambda$.
\end{Example}

Recall that we have a standard expression for $d(\t)$ such that $d(\t) = w_n w_{n-1} \ldots w_1$, where $w_i = 1$ or $w_i = s_{a_i} s_{a_i+1} \ldots s_{i-1}$ for $1 \leq a_i \leq i-1$ and $1 \leq i \leq n$.

\begin{Lemma} \label{help: case 4.5}

Suppose $\t \in \Std(\lambda)$ and $d(\t) = w_n w_{n-1} \ldots w_{1}$ is the standard expression where $w_i \neq 1$ if $i \geq r+2$ or $i = r$ and $w_i = 1$ if $i < r$ or $i = r+1$, i.e. $d(\t) = w_n w_{n-1} \ldots w_{r+2}w_r$. If $\t{\cdot}s_r$ is not standard, then $l(w_i) \geq l(w_{r})+1$ for $i \geq r+2$.

\end{Lemma}

\proof Write $\lambda = (\lambda^{(1)},\ldots,\lambda^{(\l)})$ and $\lambda^{(\l)} = (\lambda^{(\l)}_1,\ldots,\lambda^{(\l)}_k)$. Because $\t{\cdot}s_r$ is not standard, $r$ and $r+1$ are adjacent in $\t$. By (\ref{remark: standard expression}), as $w_i \neq 1$ for $i \geq r+2$
$$
(k,\lambda^{(\l)}_k,\l) = (\t^{(n+1)})^{-1}(n) = (\t^{(n)})^{-1}(n-1) = \ldots = (\t^{(r+3)})^{-1}(r+2) = (\t^{(r+2)})^{-1}(r+1).
$$
Notice that $w_{r+1} = 1$ and $w_r$ doesn't involve $s_r$ or $s_{r+1}$, we have
$$
(k,\lambda^{(\l)}_k,\l) = (\t^{(r+2)})^{-1}(r+1) = (\t^{(r+1)})^{-1}(r+1) = (\t^{(r)})^{-1}(r+1).
$$

Since $w_r \neq 1$ and $w_{r+1} = 1$, recall $w_r = s_{a_r}s_{a_r+1}\ldots s_{r-2}s_{r-1}$, by (\ref{remark: standard expression}) we have
$$
(\t^{(r+2)})^{-1}(a_r) = (\t^{(r+1)})^{-1}(a_r) = (\t^{(r)})^{-1}(r).
$$

Since $w_i = 1$ for $i < r$, we have $\t^{(r)} = \t$. Then $\t^{-1}(r+1) = (k,\lambda^{(\l)}_k,\l)$. Because $a_r \leq r-1 < r+1$, by (\ref{remark: standard expression}), $a_r$ is not on the left of $r+1$ in $\t^{(r+2)}$ because $\t^{(r+2)}|_{r+1} = \t^\mu$ with $\mu = \Shape(\t^{(r+2)}|_{r+1})$. As $r$ and $r+1$ are adjacent in $\t$ and $(\t^{(r+2)})^{-1}(a_r) = (\t^{(r)})^{-1}(r) = \t^{-1}(r)$, we must have $\t^{-1}(r) = (k-1,\lambda^{(\l)}_k,\l)$. Therefore by the definition of the standard expression, we have $l(w_r) = \lambda^{(\l)}_k-1$.

Since $(k,\lambda^{(\l)}_k,\l) = (\t^{(n+1)})^{-1}(n) = (\t^{(n)})^{-1}(n-1) = \ldots = (\t^{(r+2)})^{-1}(r+1)$ and (\ref{remark: standard expression}), we have $l(w_i) \geq \lambda^{(\l)}_k = l(w_r)+1$ for all $i \geq r+2$. \endproof

\begin{Lemma} \label{help: case 4.6}

Suppose $\t\in \Std(\lambda)$ and $d(\t) = w_n w_{n-1}\ldots w_1$ with $w_i \neq 1$ if $i > r+2$ or $i = r$ and $w_i = 1$ if $i < r$ or $i = r+1$. If $l(w_i) = l(w_r) + 1$ for all $i \geq r+2$, i.e. $d(\t) = w_n w_{n-1} \ldots w_{r+2}w_r$, and $\t{\cdot}s_r$ is not standard, then $\t$ is the last Garnir tableau of shape $\lambda$.

\end{Lemma}

\proof Write $\lambda = (\lambda^{(1)},\ldots,\lambda^{(\l)})$ and $\lambda^{(\l)} = (\lambda^{(\l)}_1,\ldots,\lambda^{(\l)}_k)$. From the proof of \autoref{help: case 4.5} we have seen that $l(w_i) = \lambda^{(\l)}_k$ for $i \geq r+2$ and $l(w_r) = \lambda^{(\l)}_k - 1$. Therefore if $\t^\lambda(k-1,\lambda_{k-1}^{(\l)},\l) = t$,
$$
\begin{cases}
w_n = s_t s_{t+1} \ldots s_{n-1}, &\\
w_{n-1} = s_{t-1} s_t \ldots s_{n-2}, &\\
\ldots\ldots\ldots &\\
w_{r+2} = s_{t - n + r + 2}s_{t-n+r+3}\ldots s_{r+1}, &\\
w_r = s_{t-n+r+1}s_{t-n+r+2} \ldots s_{r-1},
\end{cases}
$$
and by direct calculation we see that such $d(\t)$ is the last Garnir tableau of shape $\lambda$.\endproof

\begin{Example} \label{example: case 4.6}
Suppose $\lambda = (7,5,3)$ and $\t = \tab(1234567,89\twelve\fourteen\fifteen,\ten\eleven\thirteen)$. Then
$$d(\t) = s_{12}s_{13}s_{14}{\cdot}s_{11}s_{12}s_{13}{\cdot}s_{10}s_{11}.$$

So we can write $d(\t) = w_{15} w_{14} w_{13} w_{12}$ where $w_{15} = s_{12}s_{13}s_{14}$, $w_{14} = s_{11}s_{12}s_{13}$, $w_{13} = 1$ and $w_{12} = s_{10}s_{11}$. Notice $l(w_{15}) = l(w_{14}) = l(w_{12}) + 1$ and $\t{\cdot}s_{12}$ is not standard, and furthermore, $\t$ is the last Garnir tableau of shape $\lambda$.
\end{Example}

Finally we are ready to prove the most important result of this subsection.

\begin{Proposition}\label{psi problem: part1}
Suppose $\lambda\in \mathscr S_n^\Lambda$ and that $\t$ is a standard $\lambda$-tableau such that $l(d(\t)) \leq m_\lambda$ and $d(\t){\cdot}s_r$ is not reduced or $\t{\cdot}s_r$ is not standard for some $r$. Then
$$
\psi_{\s\t} \psi_r = \sum_{(\u,\v)\rhd (\s,\t)}c_{\u\v} \psi_{\u\v}.
$$
for any standard $\lambda$-tableau $\s$.
\end{Proposition}
\proof First we set $\s = \t^\lambda$. By the definition of $m_\lambda$ the Proposition holds when $l(d(\t)) < m_\lambda$. Hence we only have to consider the case when $l(d(\t)) = m_\lambda$. By~\autoref{m_lambda is positive} we have $l(d(\t)) > 0$. Therefore $\psi_{d(\t)} \neq 1$. 

Recall that the standard expression of $\psi_{d(\t)}$ has the form $\psi_{w_n}\psi_{w_{n-1}}\ldots \psi_{w_2}$ where $\psi_{w_i} = \psi_{a_i}\psi_{a_i+1}\psi_{a_i+2}\ldots \psi_{i-1}$ for some $a_i \leq i-1$ or $\psi_{w_i} = 1$. Let $k$ be the smallest positive integer such that $\psi_{w_k}\neq 1$ and $\psi_{w_i} = 1$ for all $i < k$. Because $\psi_{d(\t)} \neq 1$, the integer $k$ is well-defined. So $\psi_{d(\t)} = \psi_{w_n}\psi_{w_{n-1}}\ldots \psi_{w_k}$.

Recall that by \autoref{psiproblem: unlock} and~\ref{psiproblem: unlock on tails}, if $d(\t)$ is unlocked by $s_r$ in type I or type II, we have $\psi_{\t^\lambda \t} \psi_r \in \Bgelam $.

We separate the problem into several cases:

\textbf{Case \ref{psi problem: part1}a}: $k-1\not\in \{r-1,r,r+1\}$. Then
$$
\psi_{d(\t)} = \psi_{w_n}\ldots \psi_{w_k} = \psi_{w_n}\ldots\psi_{w_{k+1}}\psi_{a_k}\psi_{a_k+1}\ldots \psi_{k-2}\psi_{k-1}.
$$

In this case $\t$ is unlocked by $s_r$ in type I. Therefore by \autoref{psiproblem: unlock}, $\psi_{\t^\lambda \t}\psi_r \in \Bgelam $. \\

\textbf{Case \ref{psi problem: part1}b}: $k-1 = r$. Define $\w = \t{\cdot}s_r$. Hence $d(\t) = d(\w){\cdot}s_r$. Write $\bi_\w = (i_1 i_2 \ldots i_n)$.
$$
e_\lambda y_\lambda \psi_{d(\t)}\psi_r = y_\lambda \psi_{d(\w)}e(\bi_\w)\psi_r^2 = \begin{cases}
0, & \text{ if $i_r = i_{r+1}$,}\\
y_\lambda \psi_{d(\w)}e(\bi_\w) = \psi_{\t^\lambda \w}, & \text{ if $|i_r - i_{r+1}|>1$,}\\
\pm\ y_\lambda \psi_{d(\w)}e(\bi_\w)(y_r - y_{r+1})\\
\hspace*{10mm} = \pm\ \psi_{\t^\lambda \w}(y_r - y_{r+1}), & \text{ if $i_r = i_{r+1}\pm 1$.}
\end{cases}
$$

By \autoref{yup} we have $\psi_{\t^\lambda \t} \psi_r \in \Bgelam $.\\

\textbf{Case \ref{psi problem: part1}c}: $k - 1 = r + 1$.

\textbf{\ref{psi problem: part1}c.1:} $\psi_{w_i} = 1$ for some $i > k$ and $i \neq n$. Then we have
\begin{eqnarray*}
\psi_{d(\t)}\psi_r & = & \psi_{w_n}\psi_{w_{n-1}}\ldots\psi_{w_{i+2}}\psi_{w_{i+1}}\psi_{w_{i-1}}\ldots \psi_{w_k}\psi_r \\
& = & \psi_{w_n}\psi_{w_{n-1}}\ldots\psi_{w_{i+2}}(\psi_{a_{i+1}} \psi_{a_{i+1}+1}\ldots \psi_{i-1}\psi_i )\psi_{w_{i-1}}\ldots \psi_{w_k}\psi_r.
\end{eqnarray*}

As $i > k = r + 2 > r+1$, we have
$$
\psi_{d(\t)}\psi_r = (\psi_{w_n}\psi_{w_{n-1}}\ldots\psi_{w_{i+2}}\psi_{a_{i+1}} \psi_{a_{i+1}+1}\ldots \psi_{i-1}\psi_{w_{i-1}}\ldots \psi_{w_k}\psi_i)\psi_r,
$$
which shows that $\t$ is unlocked by $s_r$ in type I. By \autoref{psiproblem: unlock} we have $\psi_{\t^\lambda \t}\psi_r \in \Bgelam $.\\

\textbf{\ref{psi problem: part1}c.2:} $\psi_{w_n} = 1$. In this case $\psi_{n-1}$ is not involved in $\psi_{d(\t)}\psi_r$. By \autoref{psi_n-1 is not involved} we have $\psi_{\t^\lambda \t}\psi_r \in \Bgelam $.\\

\textbf{\ref{psi problem: part1}c.3:} $\psi_{w_i}\neq 1$ for $i > k$ and $l(\psi_{w_k}) > 1$. Then we see that $\t$ is unlocked by $s_r$ in type II. By \autoref{psiproblem: unlock on tails} we have $\psi_{\t^\lambda \t}\psi_r \in \Bgelam $.\\

\textbf{\ref{psi problem: part1}c.4:} $\psi_{w_i}\neq 1$ for $i > k$, and there exists $k < j < n$ such that $l(\psi_{w_k}) = l(\psi_{w_{k+1}}) = \ldots = l(\psi_{w_{j-1}}) = 1$ and $l(\psi_{w_j}) > 1$. Then we have
$$
w_j{\cdot}w_{j-1} = s_{a_j} s_{a_j+1}\ldots s_{j-3} s_{j-2}s_{j-1}{\cdot}s_{j-2} = s_{a_j} s_{a_j+1}\ldots s_{j-3}{\cdot}s_{j-1}s_{j-2}s_{j-1}.
$$

Therefore
\begin{eqnarray*}
d(\t) & = & w_n w_{n-1}\ldots w_{j+1}{\cdot}s_{a_j} s_{a_j+1}\ldots s_{j-3}{\cdot} s_{j-1}s_{j-2}s_{j-1}{\cdot}w_{j-2}\ldots w_k\\
& = & w_n w_{n-1}\ldots w_{j+1}s_{a_j} s_{a_j+1}\ldots s_{j-3}{\cdot} s_{j-1}s_{j-2}w_{j-2}\ldots w_k{\cdot}s_{j-1},
\end{eqnarray*}
and $j - 1 \geq k = r + 2 > r + 1$, $s_{j-1}$ and $s_r$ commute, which shows that $\t$ is unlocked by $s_r$ in type I. By \autoref{psiproblem: unlock}, we have $\psi_{\t^\lambda \t}\psi_r \in \Bgelam $.\\

\textbf{\ref{psi problem: part1}c.5:} $l(\psi_{w_k}) = l(\psi_{w_{k+1}}) = \ldots = l(\psi_{w_n}) = 1$. Then by \autoref{help: case 3.5}, $\t$ is the last Garnir tableau of shape $\lambda$. Hence by \autoref{garnir}, $\psi_{\t^\lambda \t}\psi_r \in \Bgelam $.\\

\textbf{Case \ref{psi problem: part1}d}: $k - 1 = r - 1$.

\textbf{\ref{psi problem: part1}d.1:} $\psi_{w_{k+1}} \neq 1$. Then
\begin{eqnarray*}
d(\t) & = & w_n w_{n-1} \ldots w_{k+2}{\cdot}w_{k+1} w_k\\
& = & w_n w_{n-1} \ldots w_{k+2}{\cdot} s_{a_{k+1}}s_{a_{k+1}+1}\ldots s_{k-1}s_k {\cdot} s_{a_k} s_{a_k+1}\ldots s_{k-2}s_{k-1}\\
& = & w_n w_{n-1} \ldots w_{k+2} s_{a_{k+1}}s_{a_{k+1}+1}\ldots s_{k-1} s_{a_k} s_{a_k+1}\ldots s_{k-2}{\cdot}s_ks_{k-1},
\end{eqnarray*}
and as $r = k$, we see that $\t$ is unlocked by $s_r$ in type II. Therefore $\psi_{\t^\lambda \t}\psi_r \in \Bgelam $ by \autoref{psiproblem: unlock on tails}.\\

\textbf{\ref{psi problem: part1}d.2:} $k = n-1$ and $\psi_{w_{k+1}} = \psi_{w_n} = 1$. The $\psi_{d(\t)}\psi_r = \psi_{w_{n-1}}\psi_{n-1} = \psi_{a_{n-1}}\psi_{a_{n-1}+1}\ldots\psi_{n-2}\psi_{n-1}$. Then by \autoref{help: case 4.3}, $\t$ is the last Garnir tableau of shape $\lambda$. Hence by \autoref{garnir}, $\psi_{\t^\lambda \t}^\lambda\psi_r \in \Bgelam $.\\

\textbf{\ref{psi problem: part1}d.3:} $k < n-1$, $\psi_{w_{k+1}} = 1$ and $\psi_{w_n} = 1$. Then $n-1 > k = r$. So $\psi_{d(\t)}\psi_r$ doesn't involve $\psi_{n-1}$. By \autoref{psi_n-1 is not involved} we have $\psi_{\t^\lambda \t}\psi_r \in \Bgelam $.\\

\textbf{\ref{psi problem: part1}d.4:} $k < n-1$, $\psi_{w_{k+1}} = 1$ and there exists $k+1 < j < n$ such that $\psi_{w_j} = 1$ and $\psi_{w_{j+1}} \neq 1$. In this case we have
\begin{eqnarray*}
d(\t) & = & w_n w_{n-1} \ldots w_{j+2}w_{j+1}w_j w_{j-1}\ldots w_k\\
& = & w_n w_{n-1} \ldots w_{j+2}{\cdot}s_{a_{j+1}}s_{a_{j+1}+1}\ldots s_{j-1}s_j{\cdot}w_{j-1}\ldots w_k\\
& = & (w_n w_{n-1} \ldots w_{j+2}{\cdot}s_{a_{j+1}}s_{a_{j+1}+1}\ldots s_{j-1}{\cdot}w_{j-1}\ldots w_k){\cdot}s_j.
\end{eqnarray*}

As $j > k+1 = r+1$, $\psi_j$ and $\psi_r$ commute. Therefore $\t$ is unlocked by $s_r$ in type I. By \autoref{psiproblem: unlock}, we have $\psi_{\t^\lambda\t}\psi_r \in \Bgelam $.\\

\textbf{\ref{psi problem: part1}d.5:} $k < n-1$, $\psi_{w_{k+1}} = 1$ and for any $j > k+1$, $\psi_{w_j} \neq 1$. Then by \autoref{help: case 4.5}, we have $l(\psi_{w_j}) \geq l(\psi_{w_k})+1$ for all $j \geq k+2$.\\

\textbf{\ref{psi problem: part1}d.5.1: } Suppose $l(\psi_{w_{k+2}}) > l(\psi_{w_k}) + 1$. So we have $a_{k+2} \leq a_{k}$, and therefore
\begin{eqnarray*}
w_{k+2} w_{k} & = & s_{a_{k+2}} s_{a_{k+2} + 1} \ldots s_{k} s_{k+1} {\cdot} s_{a_{k}}s_{a_{k}+1} \ldots s_{k-2} s_{k-1}\\
& = & s_{a_{k}+1} \ldots s_{k-1} s_{k}{\cdot} s_{a_{k+2}} s_{a_{k+2} + 1} \ldots s_{k} s_{k+1}.
\end{eqnarray*}

Therefore
\begin{eqnarray*}
d(\t) & = & w_n w_{n-1}\ldots w_{k+3}{\cdot}w_{k+2} w_{k}\\
& = & w_n w_{n-1} w_{k+3}{\cdot}s_{a_{k}+1} \ldots s_{k-1} s_{k}{\cdot} s_{a_{k+2}} s_{a_{k+2} + 1} \ldots s_{k} s_{k+1}.
\end{eqnarray*}

Then because $k = r$, $\t$ is unlocked by $s_r$ in type II. Therefore, by \autoref{psiproblem: unlock on tails}, $\psi_{\t^\lambda \t} \psi_r \in \Bgelam $.\\

\textbf{\ref{psi problem: part1}d.5.2: } There exists $j > k+2$ such that $l(\psi_{w_{k+2}}) = l(\psi_{w_{k+3}}) = \ldots  = l(\psi_{w_{j-1}}) = l(\psi_{w_k}) + 1$ and $l(\psi_{w_j}) > l(\psi_{w_k}) + 1$. So we have $l(\psi_{w_j}) > l(\psi_{w_{j-1}})$ and $a_j \leq a_{j-1}$, and therefore
\begin{eqnarray*}
w_j w_{j-1} & = & s_{a_j} s_{a_j + 1} \ldots s_{j-2} s_{j-1} {\cdot} s_{a_{j-1}}s_{a_{j-1}+1} \ldots s_{j-3} s_{j-2}\\
& = & s_{a_{j-1}+1} \ldots s_{j-2} s_{j-1}{\cdot} s_{a_j} s_{a_j + 1} \ldots s_{j-2} s_{j-1}.
\end{eqnarray*}

Therefore
\begin{eqnarray*}
d(\t) & = & w_n w_{n-1} \ldots w_{j+1} w_j w_{j-1} w_{j-2} \ldots w_{k}\\
& = & w_n w_{n-1} \ldots w_{j+1} {\cdot}s_{a_{j-1}+1} \ldots s_{j-1}{\cdot} s_{a_j} \ldots s_{j-2} s_{j-1} {\cdot} w_{j-2} \ldots w_k\\
& = & (w_n w_{n-1} \ldots w_{j+1} {\cdot}s_{a_{j-1}+1} \ldots s_{j-1}{\cdot} s_{a_j} \ldots s_{j-2} {\cdot} w_{j-2} \ldots w_k){\cdot} s_{j-1}.
\end{eqnarray*}

Then because $j-1 > k+1 = r+1$, $s_{j-1}$ and $s_r$ commutes. Hence $\t$ is unlocked by $s_r$ in type I and therefore, by \autoref{psiproblem: unlock}, $\psi_{\t^\lambda \t} \psi_r \in \Bgelam $.\\

\textbf{\ref{psi problem: part1}d.5.3:} $l(\psi_{w_{k+2}}) = l(\psi_{w_{k+3}}) = \ldots = l(\psi_{w_{n-1}}) = l(\psi_{w_k}) + 1$. By \autoref{help: case 4.6}, $\t$ is the last Garnir tableau of shape $\lambda$. By \autoref{garnir}, we have $\psi_{\t^\lambda \t}\psi_r \in \Bgelam $.\\

By the above cases, $\psi_{\t^\lambda \t} \psi_r$ is always in $\Bgelam $. Therefore by \autoref{basis to hecke}, we have
$$
\psi_{\t^\lambda \t} \psi_r  = \sum_{(\u,\v)\rhd (\t^\lambda, \t)}c_{\u\v}\psi_{\u\v} = \sum_{\v\rhd t}c_{\t^\lambda \v} \psi_{\t^\lambda \v} + \sum_{\u,\v\in\Std(>\lambda)}c_{\u\v} \psi_{\u\v}.
$$

Giving any standard $\lambda$-tableau $\s$, we have $\psi_{\s\t} \psi_r = \psi^*_{d(\s)}\psi_{\t^\lambda \t}\psi_r$. Notice $\psi_{d(\s)}^* \psi_{\t^\lambda \v} = \psi_{\s\v}$. For any $\u,\v\in\Std(>\lambda)$, $\psi_{\u\v}\in \Blam $. As $\lambda \in \mathscr S_n^\Lambda$, by \autoref{B_lambda is an ideal}, $\Blam $ is an ideal. Therefore $\psi^*_{d(\s)}\psi_{\u\v} \in \Blam $. These arguments yield that $\psi_{\s\t}\psi_r\in \Bgelam $. By \autoref{basis to hecke} we completes the proof.\endproof

The following Corollary is straightforward by \autoref{psi problem: part2} and \autoref{psi problem: part1}.

\begin{Corollary}\label{psi problem}
Suppose $\lambda\in \mathscr S_n^\Lambda$, for any standard $\lambda$-tableau $\t$ with $l(d(\t)) \leq m_\lambda$, then
$$
\psi_{\s\t} \psi_r = \begin{cases}
\psi_{\t^\lambda \w} + \sum_{(\u,\v)\rhd (\s,\t)}c_{\u\v} \psi_{\u\v}, & \text{ if $\w = \u{\cdot}s_r$ is standard and $d(\u){\cdot}s_r$ is reduced,}\\
\sum_{(\u,\v)\rhd (\s,\t)}c_{\u\v} \psi_{\u\v}, & \text{ if $\u{\cdot}s_r$ is not standard or $d(\u){\cdot}s_r$ is not reduced.}
\end{cases}
$$
for any standard $\lambda$-tableau $\s$.
\end{Corollary}

\subsection{Integral basis Theorem}

In this subsection we will complete the main Theorem of this paper.

\begin{Theorem}\label{main}
Suppose $\lambda\in \mathscr S_n^\Lambda$, we have $\lambda\in \P_I \cap \P_y \cap \P_\psi$.
\end{Theorem}

\proof By \autoref{I-problem: final} we have when $\lambda\in \mathscr S_n^\Lambda$ then $\lambda\in \P_I$. By \autoref{m_lambda is positive}, we have $0 < m_\lambda$, i.e. $1 \leq m_\lambda$. Assume $l = l(d(\u)$ for some $\u\in\Std(\lambda)$, by \autoref{yup} and \autoref{psi problem}, for any $\t\in\Std(\lambda)$ with $l(d(\t)) = l$, we have
\begin{eqnarray*}
\psi_{\s\t} y_r & = & \sum_{(\u,\v) \rhd (\s,\t)} c_{\u\v} \psi_{\u\v},\\
\psi_{\s\t} \psi_r & = & \begin{cases}
\psi_{\s \w} + \sum_{(\u,\v)\rhd (\s,\t)}c_{\u\v}\psi_{\u\v},& \text{ if $\w = \t{\cdot}s_r$ is standard and $d(\u){\cdot}s_r$ is reduced,}\\
\sum_{(\u,\v)\rhd (\s,\t)}c_{\u\v} \psi_{\u\v}, & \text{ if $\u{\cdot}s_r$ is not standard or $d(\u){\cdot}s_r$ is not reduced.}
\end{cases}
\end{eqnarray*}

This implies that $l < m_\lambda$, i.e. $l+1 \leq m_\lambda$. So by induction, for any $\t\in \Std(\lambda)$, we have $l(d(\t)) < m_\lambda$. Therefore $\lambda\in \P_y \cap \P_\psi$. This completes the proof.\endproof

\begin{Theorem}\label{integral basis theorem}
The set $\set{\psi_{st}^{\Z}|\s,\t \in \Std(\lambda) \text{ for $\lambda\in\P_n$}}$ is a graded cellular basis of $\mathscr{R}_n^\Lambda(\Z)$.
\end{Theorem}

\proof It's trivial that when $n = 1$ the Theorem holds. Assume for any $n' < n$ the Theorem follows. Suppose we can write all multipartitions of $n$ as $\lambda_{[1]}, \lambda_{[2]},\ldots,\lambda_{[k]}$ where $\lambda_{[1]} > \lambda_{[2]} > \ldots > \lambda_{[k]}$. As $\lambda_{[1]} = ((n),\emptyset,\ldots,\emptyset)$, by \autoref{onerow: I}, \autoref{onerow: psi} and~\ref{onerow: y}, we have $\lambda_{[1]}\in \P_I \cap \P_y \cap \P_\psi$. Hence $\lambda_{[2]}\in \mathscr S_n^\Lambda$. Now assume $\lambda_{[i]}\in \mathscr S_n^\Lambda$, by \autoref{main}, $\lambda_{[i]}\in \P_I \cap \P_y \cap \P_\psi$. Hence $\lambda_{[i+1]}\in \mathscr S_n^\Lambda$. Therefore for any $i$, $\lambda_{[i]}\in \mathscr S_n^\Lambda$. Hence, for any $\lambda\in \P_n$, $\lambda\in \P_I \cap \P_y \cap \P_\psi$. Recall that
$$
\BZ  =  \{r\in \mathscr{R}_n^\Lambda(\Z)\ |\ r = \sum_{\substack{\s,\t\in \Std(\mu)\\ \mu\in\P_n}} c_{\s\t}\psi_{\s\t},c_{\s\t}\in \Z\}.
$$
So $\BZ$ is an ideal.

Now, for any $\bi = (i_1,i_2,\ldots,i_n) \in I^n$, set $\bj = (i_1,i_2,\ldots,i_{n-1}) \in I^{n-1}$. Because $e(\bj) \in \mathscr R_{n-1}^\Lambda$, by assumption we have $e(\bj) = \sum_{\substack{\mu\in\P_{n-1}\\ \u,\v\in\Std(\mu)}}c_{\u\v}\psi_{\u\v} \in R_{n-1}^\Lambda$ and hence that $e(\bi) = \theta_{i_n}(e(\bj)) = \sum_{\substack{\mu\in\P_{n-1}\\ \u,\v\in\Std(\mu)}}c_{\u\v}\theta_{i_n}(\psi_{\u\v})$.

For any $\mu \in \mathscr P_{n-1}^\Lambda$ and $\u,\v\in\Std(\mu)$, we have
$$
\theta_{i_n} (\psi_{\u\v}) = \psi^*_{d(\u)} e(\bi_{\mu}\vee i_n) y_\mu \psi_{d(\v)}.
$$

By \autoref{I-problem: downstair cases} and \autoref{I-problem: final}, we have $e(\bi_\mu\vee i_n)y_\mu y_n^0 \in R_n^\Lambda$. Then because $\BZ$ is an ideal,
$$
e(\bi) = \sum_{\substack{\mu\in\P_{n-1}\\ \u,\v\in\Std(\mu)}}c_{\u\v}\theta_{i_n}(\psi_{\u\v}) \in \BZ.
$$

Then we have $\BZ = \mathscr R_n^\Lambda(\Z)$. By \autoref{linearly independent}, the set $\{\psi_{st}^{\Z}\ |\ \s,\t \in \text{Std}(\lambda) \text{ for $\lambda\in\P_n$}\}$ is linearly independent. Hence, $\{\psi_{st}^{\Z}\ |\ \s,\t \in \text{Std}(\lambda) \text{ for $\lambda\in\P_n$}\}$ is a basis of $\mathscr R_n^\Lambda(\Z)$.

By definition, the elements of the set are homogeneous. That this basis is cellular is trivial by \autoref{basis with field} and \autoref{basis to hecke}. This completes the proof.\endproof

The next Corollary is a straightforward application of~\autoref{integral basis theorem}.

\begin{Corollary} \label{bi is residue sequence}
For any $\bi\in I^n$, $e(\bi) \neq 0$ if and only if $\bi$ is the residue sequence of a standard tableau $\t$.
\end{Corollary}

\proof Suppose $\bi$ is the residue sequence of a standard tableau $\t$. By \autoref{integral basis theorem} we have $\psi^{\Z}_{\t\t}\neq 0$. Because $\psi^{\Z}_{\t\t} = \psi^{\Z}_{\t\t}e(\bi)$, we must have $e(\bi) \neq 0$.

Suppose $\bi$ is not the residue sequence of any standard tableau. By \autoref{integral basis theorem} we can write
$$
1 = \sum_{\s,\t} c_{\s\t}\psi^{\Z}_{\s\t},
$$
and hence $e(\bi) = 1{\cdot}e(\bi) = \sum_{\s,\t} c_{\s\t}\psi^{\Z}_{\s\t}e(\bi) = 0$, which completes the proof. \endproof

\section{A new basis of the Affine KLR Algebras}

In \autoref{integral basis theorem} we have shown that $\R(\Z)$ is a $\Z$-free algebra with basis $\{\psi_{\s\t}\ |\ \s,\t\in\Std(\lambda),\lambda\in\P_n\}$. In this section we extend this result and construct a graded cellular basis of $\mathscr R_n(\Z)$. Moreover, for any weight $\Lambda$ we obtain a homogeneous basis of the ideal $N_n^\Lambda$ as a subset of the graded cellular basis of $R_n(\Z)$. As an application we obtain a new classification of the complete set of simple $\Rn$-modules by using the graded cellular basis of $\Rn$. Unlike in previous sections, this section laces no constraints upon $e$ as we now allow the case $e = 2$.

\subsection{Infinite sequence of weights and basis of $\Ri$}

This subsection introduces a special kind of sequences of weights $\inftyweight$. We use these infinite sequences of weights to define a graded cellular basis for $\Ri$ by combining the graded cellular basis of the algebras $\Rni$.

We fix an integer $e\geq 0$ and $e \neq 1$, and consider the KLR algebras $\mathscr R_n = \mathscr R_n(\Z)$. Suppose $\Lambda = \sum_{i\in I}a_i \Lambda_i$ and $\Lambda' = \sum_{i\in I}a'_i \Lambda_i$ are two weights in $P_+$. Define a partial ordering on $P^+$ and write $\Lambda \leq \Lambda'$ if $a_i \leq a_i'$ for all $i\in I$. Write $\Lambda < \Lambda'$ if $\Lambda \leq \Lambda'$ and $\Lambda \neq \Lambda'$.

\begin{Definition} \label{standard sequence}
Suppose $\Lambda^\infty = (\Lambda^{(k)})_{k\geq 1}$ is a sequence of weights in $P^+$ where $\Lambda^{(k)} = \sum_{i \in I} a_i^{(k)} \Lambda_i$. It is an \textbf{increasing sequence} if $\Lambda^{(k)} <\Lambda^{(k+1)}$ for all $k \geq 1$. The sequence $\inftyweight$ is called \textbf{standard} if $\lim_{k\rightarrow\infty} a_i^{(k)} = \infty$, for all $i \in I$.
\end{Definition}

\begin{Example} \label{Example: standard sequence: e > 0}
Suppose $e > 0$. Let $\inftyweight$ be any increasing sequence of weights such that $\Lambda^{(1)} = \Lambda_1$ and $\Lambda^{(k)} = \Lambda^{(k-1)} + \Lambda_i$ whenever $k \equiv i\pmod{e}$. For example, when $e = 3$, we have
\begin{eqnarray*}
\Lambda^{(1)} & = & \Lambda_1,\\
\Lambda^{(2)} & = & \Lambda_1 + \Lambda_2,\\
\Lambda^{(3)} & = & \Lambda_1 + \Lambda_2 + \Lambda_0,\\
\Lambda^{(4)} & = & 2\Lambda_1 + \Lambda_2 + \Lambda_0,\\
\Lambda^{(5)} & = & 2\Lambda_1 + 2\Lambda_2 + \Lambda_0,\\
\Lambda^{(6)} & = & 2\Lambda_1 + 2\Lambda_2 + 2\Lambda_0,\\
\Lambda^{(7)} & = & 3\Lambda_1 + 2\Lambda_2 + 2\Lambda_0,\\
\ldots\ldots & \ldots& \ldots\ldots\ldots
\end{eqnarray*}
So in this case we have $\lim_{k\rightarrow\infty}a_i^{(k)} = \infty$ for any $i\in I$ and $\inftyweight$ is a standard sequence.
\end{Example}

\begin{Example}\label{Example: standard sequence: e = 0}
Suppose $e = 0$. Define $\inftyweight$ where $\Lambda^{(1)} = \Lambda_0$ and $\Lambda^{(k)} = \Lambda^{(k-1)} + \Lambda_i$ with $i = (k-1) - (n-1)^2 - (n-1) = k-n^2+n-1$ if $(n-1)^2 < k \leq n^2$. In more details,
\begin{eqnarray*}
\Lambda^{(1)} & = & \Lambda_0,\\
\Lambda^{(2)} & = & \Lambda_{-1} + \Lambda_0,\\
\Lambda^{(3)} & = & \Lambda_{-1} + 2\Lambda_0,\\
\Lambda^{(4)} & = & \Lambda_{-1} + 2\Lambda_0 + \Lambda_1,\\
\Lambda^{(5)} & = & \Lambda_{-2} + \Lambda_{-1} + 2\Lambda_0 + \Lambda_1,\\
\Lambda^{(6)} & = & \Lambda_{-2} + 2\Lambda_{-1} + 2\Lambda_0 + \Lambda_1,\\
\Lambda^{(7)} & = & \Lambda_{-2} + 2\Lambda_{-1} + 3\Lambda_0 + \Lambda_1,\\
\Lambda^{(8)} & = & \Lambda_{-2} + 2\Lambda_{-1} + 3\Lambda_0 + 2\Lambda_1,\\
\Lambda^{(9)} & = & \Lambda_{-2} + 2\Lambda_{-1} + 3\Lambda_0 + 2\Lambda_1 + \Lambda_2,\\
\Lambda^{(10)} & = & \Lambda_{-3} + \Lambda_{-2} + 2\Lambda_{-1} + 3\Lambda_0 + 2\Lambda_1 + \Lambda_2,\\
\Lambda^{(11)} & = & \Lambda_{-3} + 2\Lambda_{-2} + 2\Lambda_{-1} + 3\Lambda_0 + 2\Lambda_1 + \Lambda_2,\\
\ldots\ldots & \ldots & \ldots\ldots\ldots
\end{eqnarray*}
So in this case we have $\lim_{k\rightarrow\infty}a_i^{(k)} = \infty$ for any $i\in I$ and $\inftyweight$ is a standard sequence.
\end{Example}

Recall that for any weight $\Lambda = \sum_{i\in I} a_i \Lambda_i$, we can define the two-sided ideal $N_n^\Lambda$ of $\mathscr R_n$. The ideal $N_n^\Lambda$ is generated by elements $e(\bi)y_1^{a_{i_1}}$ for all $\bi = (i_1,i_2,\ldots,i_n) \in I^n$. By definition, $\R \cong \mathscr R_n/N_n^\Lambda$. There is a natural injective homomorphism sending $\R$ to $\mathscr R_n$ by sending $e(\bi)$, $y_r$ and $\psi_r$ to $\hat e(\bi)$, $\hat y_r$ and $\hat \psi_r$, respectively. Hence we can consider $\R$ as a $\Z$-submodule of $\mathscr R_n$ and write $\mathscr R_n \cong \R \oplus N_n^\Lambda$ as $\Z$-modules.

Recall $Q_+ = \sum_{i \in I}\N \alpha_i$ is defined in subsection 1.1. For $\alpha = \sum_{i\in I} a_i \alpha_i \in Q_+$, define $|\alpha| = \sum_{\i \in I}a_i $. Then for any $\alpha\in Q_+$ with $|\alpha| = n$, define $I^\alpha$ to be the set of all sequences $\bi = (i_1,i_2,\ldots,i_n)\in I^n$ such that $a_i = |\set{1 \leq r \leq n|i_r = i}|$. By definition, if $\bi,\bj\in I^\alpha$ then there exists $v \in \mathfrak S_n$ such that $\bi = \bj{\cdot}v$. Define $\hat e_\alpha = \sum_{\bi \in I^\alpha} \hat e(\bi) \in \Rn$ and $e_\alpha = \sum_{\bi \in I^\alpha}e(\bi) \in \R$.

The following result is trivial by the relations of $\Rn$.

\begin{Lemma} \label{alpha R beta}
Suppose $\alpha,\beta\in Q_+$. Then $\Rn \hat e_\alpha \neq 0$ and $\hat e_\beta \Rn \hat e_\alpha = \delta_{\alpha\beta}\Rn e_\alpha = \delta_{\alpha\beta}\hat e_\beta \Rn$.
\end{Lemma}

We then define $\Rni = \Rn \hat e_\alpha$, $\Ri = \R e_\alpha$ and $N_\alpha^\Lambda = N_n^\Lambda \hat e_\alpha$. We see that $\Rni \hat e(\bj) = 0$ if $\bj \not\in I^\alpha$. Finally, because
$$
\Rn = \bigoplus_{\alpha\in Q_+} \Rni \hspace*{5mm}\text{and}\hspace*{5mm} \R = \bigoplus_{\alpha\in Q_+} \Ri,
$$
and by the relations $\Rni$ and $\Ri$'s are subalgebras of $\Rn$ and $\R$, respectively. Hence we will mainly work in $\Rni$, $\Ri$ and $N_\alpha^\Lambda$ and extend the basis of $\Ri$ to $\Rni$ and hence generate a graded cellular basis of $\Rn$.

By \autoref{integral basis theorem} and the orthogonality of $e(\bi)$'s we can give a basis for $\Ri$.

\begin{Proposition} \label{basis of Ri}
Suppose $\bi \in I^n$ and $\Lambda \in P_+$. The set
$$
\set{\psi_{\s\t}|\lambda\in\P_n, \s,\t\in\Std(\lambda), \res(\t) \in I^\alpha}
$$
is a graded cellular basis of $\Ri$.
\end{Proposition}

\subsection{Minimum degree of $N_\alpha^\Lambda$}

Fix $\alpha \in Q_+$. In the last subsection we introduced a standard sequence of $\inftyweight$. For each $k$ and $\bi\in I^n$, we define the \textbf{minimum degree} of $N_\alpha^{\Lambda^{(k)}}$ to be the integer 
$$
m_\alpha^{\Lambda^{(k)}} = \text{min}\{\deg(r)\ |\ \text{$r$ is a homogeneous element in $N_\alpha^\Lambda$}\}. 
$$

We will prove that $m_\alpha^{\Lambda^{(k)}} \rightarrow \infty$ with $k\rightarrow\infty$. This result is quite important in the next subsection because it will allow us to extend the basis of $\Ri$ to $\Rni$.

First we need to find a general description of the homogeneous elements of $N_\alpha^\Lambda$.

\begin{Lemma} \label{span set of nilpotents}
For $\Lambda = \sum_{i\in I} a_i \Lambda_i \in P_+$ and $\alpha\in Q_+$, the ideal $N_\alpha^\Lambda$ is spanned by
$$
\set{\psi_u e(\bi) y_1^{a_{i_1}} f(y) \psi_v|u,v\in\mathfrak S_n, f(y) \in \Z[y_1,y_2,\ldots,y_n], \bi \in I^\alpha}.
$$
\end{Lemma}

\proof By the definition of $N_\alpha^\Lambda$, any element of $N_\alpha^\Lambda$ can be written as linear combination of elements of the form
\begin{equation} \label{span set of nilpotents: help1}
\psi_{u_k}f_k(y)\psi_{u_{k-1}}\ldots\psi_{u_2}f_2(y)\psi_{u_1}f_1(y)e(\bi)y_1^{a_{i_1}}g_1(y)\psi_{v_1}g_2(y)\psi_{v_2}\ldots\psi_{v_{l-1}}g_l(y)\psi_{v_l},
\end{equation}
where $u_i,v_i\in\mathfrak S_n$, $\bi \in I^\alpha$ and $f_i(y), g_i(y) \in \Z[y_1,\ldots,y_n]$. In the view of \autoref{basis change} and \cite[Proposition 2.5]{BKW:GradedSpecht}, every element in the form of (\ref{span set of nilpotents: help1}) can be written as linear combination of terms of the form $\psi_u e(\bi)y_1^{a_{i_1}}f(y)\psi_v$'s. Hence $N_\alpha^\Lambda$ is spanned by the elements given in the statement of the Lemma. \endproof

The next result is directly implied by the above Lemma.

\begin{Proposition} \label{infinite degree}
Suppose that $\inftyweight$  is a standard sequence and $\alpha \in Q_+$. Then $\lim_{k\rightarrow \infty}m_\alpha^{\Lambda^{(k)}} = \infty$.
\end{Proposition}

\proof By \autoref{span set of nilpotents}, for any $k \geq 1$ we have
$$
m_\alpha^{\Lambda^{(k)}} = \min\set{\deg(\psi_u e(\bi) y_1^{a_{i_1}^{(k)}} f(y) \psi_v)|u,v\in\mathfrak S_n, f(y) \in \Z[y_1,y_2,\ldots,y_n], \bi \in I^\alpha}.
$$

By definition, $\deg(\psi_u e(\bi)y_1^{a_{i_1}^{(k)}}f(y)\psi_v) = \deg(\psi_ue(\bi)) + \deg(y_1^{a_{i_1}^{(k)}}) + \deg(f(y)) + \deg(\psi_v e(\bi{\cdot}v))$. As $u$ and $v$ are reduced expressions of permutations in $\mathfrak S_n$, $l(u) \leq \frac{(n-1)n}{2}$, and $\deg(\psi_r e(\bi)) \geq -2$ for any $\bi$. Hence, $\deg(\psi_u e(\bi)) \geq -(n-1)n$. By the same reasoning, $\deg(\psi_v e(\bi{\cdot}v)) \geq -(n-1)n$. Then as $\deg(f(y)) \geq 0$, we have $\deg(\psi_u e(\bi)y_1^{a_{i_1}^{(k)}}f(y)\psi_v) \geq -2(n-1)n + 2a_{i_1}^{(k)}$.

Define $a_\alpha^{(k)} = \min_{\bi\in I^\alpha} a_{i_1}^{(k)}$. We have
$$
\deg(\psi_u e(\bi)y_1^{a_{i_1}^{(k)}}f(y)\psi_v) \geq -2(n-1)n + 2a_\alpha^{(k)},
$$
for any $u$, $v$ and $f$. Therefore $m_\alpha^{\Lambda^{(k)}} \geq 2a_\alpha^{(k)} - 2(n-1)n$.

Choose $\bj \in I^\alpha$. By definition, $I^\alpha = \set{\bi \in I^n|\bi = \bj{\cdot}v\text{ with }v\in \mathfrak S_n}$. Then $|I^\alpha| \leq |\mathfrak S_n| < \infty$. Then $a_{i_1}^{(k)} \rightarrow \infty$ as $k\rightarrow \infty$ for any $\bi \in I^\alpha$ implies $a_\alpha^{(k)} \rightarrow \infty$ as $k \rightarrow \infty$ because $I^\alpha$ is finite. Therefore $m^{\Lambda^{(k)}}_\alpha \rightarrow \infty$. \endproof

\begin{Remark}

That the set $I^\alpha$ is finite is important in the proof of \autoref{infinite degree}. If it were possible for $I^\alpha$ to be infinite then knowing that $a^{(k)}_{i_1} \rightarrow \infty$ for all $i \in I^\alpha$ is not strong enough to imply that $a_\alpha^{(k)} \rightarrow \infty$ as $k \rightarrow \infty$.

\end{Remark}

\subsection{A graded cellular basis of $\mathscr R_n$}

In this subsection we will prove the main result of this section. First we introduce a special kind of multicharge $\kappa$ corresponding to a standard sequence $\inftyweight$. Then for any $\alpha \in Q_+$, we will find a graded cellular basis $\Bi$ of $\Rni$ that corresponds to $\kappa$.

\begin{Definition}\label{inverse multicharge}
Suppose $\inftyweight$ is a standard sequence. An \textbf{inverse multicharge sequence} for $\inftyweight$ is a infinite sequence $\kappa = (\ldots,\kappa_3,\kappa_2,\kappa_1)$ such that for any $k\geq 1$, if $\l_k = l(\Lambda^{(k)})$, then $\kappa_{\Lambda^{(k)}} = (\kappa_{\l_k},\kappa_{\l_k-1},\ldots,\kappa_2,\kappa_1)$ is a multicharge corresponding to $\Lambda^{(k)}$.
\end{Definition}

\begin{Example}\label{Example: inverse multicharge: e > 0}
Suppose $e = 3$. Using the standard sequence $\inftyweight$ introduced in Example~\ref{Example: standard sequence: e > 0}, we define a multicharge $\kappa = (\ldots,\kappa_3,\kappa_2,\kappa_1)$ where $\kappa_k \equiv k \pmod{e}$ for $k \geq 1$. We can write $\kappa = (\ldots,0,2,1,0,2,1,0,2,1,0,2,1)$, and we have
\begin{eqnarray*}
\kappa_{\Lambda^{(1)}} & = & (1),\\
\kappa_{\Lambda^{(2)}} & = & (2,1),\\
\kappa_{\Lambda^{(3)}} & = & (0,2,1),\\
\kappa_{\Lambda^{(4)}} & = & (1,0,2,1),\\
\kappa_{\Lambda^{(5)}} & = & (2,1,0,2,1),\\
\ldots\ldots & \ldots & \ldots\ldots\ldots
\end{eqnarray*}

All of these multicharges correspond to $\Lambda^{(k)}$.
\end{Example}

Fix a standard sequence $\inftyweight$ and an inverse multicharge sequence $\kappa$ for $\inftyweight$. An \textbf{affine multipartition} of $n$ is an ordered sequence $\hat\lambda = (\ldots,\lambda^{(2)},\lambda^{(1)})$ of partitions such that $\sum_{i = 1}^\infty |\lambda^{(i)}| = n$. Let $\affineP$ be the set of all affine multipartitions of $n$. We define \textbf{young diagram} $[\hat\lambda]$ and \textbf{standard affine tableau} $\S$ for affine multipartitions in the same way as for multipartitions. Let $\Std(\hat\lambda)$ be the set of all standard affine tableaux of shape $\hat\lambda$.

Fix an inverse multicharge sequence $\kappa = (\ldots, \kappa_2, \kappa_1)$. For any $\l > 0$, let $(\kappa_\l, \ldots, \kappa_1)$ be a multicharge and $\Lambda$ be the unique weight corresponds to $(\kappa_\l, \ldots, \kappa_1)$. We define a map $p_\l\map{\affineP}{\P_n}$, sending $\hat\lambda = (\ldots,\lambda^{(2)},\lambda^{(1)}) \in \affineP$ to $\lambda = (\lambda^{(\l)},\lambda^{(\l-1)},\ldots,\lambda^{(2)},\lambda^{(1)}) \in \P_n$.

Define the \textbf{level} of $\hat\lambda \in \affineP$ to be $l(\hat\lambda) = \l$ if $\lambda^{(\l)} \neq \emptyset$ and $\lambda^{(i)} = \emptyset$ for $i > \l$. Suppose $\hat\lambda \in \affineP$ with level $\l$ and $\lambda = p_\l(\hat\lambda)$, define $q_\l \map{\Std(\hat\lambda)}{\Std(\lambda)}$ sending $\hat\t = (\ldots, \t^{(2)}, \t^{(1)}) \in \Std(\hat\lambda)$ to $\t = (\t^{(\l)}, \t^{(\l-1)}, \ldots, \t^{(1)}) \in \Std(\lambda)$.

In order to simplify the notations, we write $\lambda = p_\l(\hat \lambda)$ if $l(\hat\lambda) = \l$. Similarly, if $l(\hat\lambda) = \l$, for $\T \in \Std(\hat\lambda)$, we write $\t = q_\l(\T) \in \Std(\lambda)$. Define the degree of each standard affine tableau to be $\deg(\S) = \deg(\s)$ and the residue sequence of the affine tableau $\res(\S) = \res(\s)$.

Extend dominance ordering $\unrhd$ and lexicographic ordering $\geq$ to $\affineP$ by defining $\hat\lambda \unrhd \hat\mu$ if $l(\hat\lambda) > l(\hat\mu)$ or $l(\hat\lambda) = l(\hat\mu)$ and $\lambda \unrhd \mu$ and $\hat\lambda \rhd \hat\mu$ if $\hat\lambda \unrhd \hat\mu$ and $\hat\lambda \neq \hat\mu$ for $\hat\lambda, \hat\mu\in\affineP$. We define the lexicographic orderings on $\affineP$ similarly.

\begin{Example}
Suppose $\hat\lambda = (\ldots|0|0|0|4,3,1|2,1|3,3)$. Then $\lambda = (4,3,1|2,1|3,3)$ and
$$
\S = \Bigg(\ldots\hspace*{1mm}\Bigg|\hspace*{1mm}\emptyset\hspace*{1mm}\Bigg|\hspace*{1mm}\emptyset\hspace*{1mm}\Bigg|\hspace*{1mm}\tab(18\thirteen\sixteen,7\twelve\fifteen,\ten)\hspace*{1mm}\Bigg|\hspace*{1mm}%
    \tab(26,3)\hspace*{1mm}\Bigg|\hspace*{1mm}\tab(45\eleven,9\fourteen\seventeen)\hspace*{1mm}\Bigg) \in \Std(\hat\lambda),
$$
and
$$
\s = \Bigg(\hspace*{1mm}\tab(18\thirteen\sixteen,7\twelve\fifteen,\ten)\hspace*{1mm}\Bigg|\hspace*{1mm}%
    \tab(26,3)\hspace*{1mm}\Bigg|\hspace*{1mm}\tab(45\eleven,9\fourteen\seventeen)\hspace*{1mm}\Bigg) \in \Std(\lambda).
$$
\end{Example}

Suppose $\Lambda \in P_+$ and $\lambda = (\lambda^{(\l)},\ldots,\lambda^{(1)}) \in \P_n$. Then for any $\s,\t\in\Std(\lambda)$, in \autoref{notation: psi basis} we have defined $\hat\psi_{\s\t}$ and $\psi_{\s\t} = \hat\psi_{\s\t} + N_n^\Lambda \in \R$. For any standard affine tableau $\S$, $\T$ we define $\psi_{\S\T} = \hat\psi_{\s\t}$. Also we can define $\psi_{\S\T}^* = \psi_{\T\S}$.

\begin{Example}
Suppose $\kappa = (\ldots,0,2,1,0,2,1,0,2,1)$ as in Example~\ref{Example: inverse multicharge: e > 0}. For
$$
\S =\Bigg(\ldots\hspace*{1mm}\Bigg|\hspace*{1mm}\emptyset\hspace*{1mm}\Bigg|\hspace*{1mm}\emptyset\hspace*{1mm}\Bigg|\hspace*{1mm}\tab(123)\hspace*{1mm}\Bigg|\hspace*{1mm}%
    \tab(4,5)\hspace*{1mm}\Bigg|\hspace*{1mm}\tab(6)\hspace*{1mm}\Bigg) \hspace*{10mm}
\T = \Bigg(\ldots\hspace*{1mm}\Bigg|\hspace*{1mm}\emptyset\hspace*{1mm}\Bigg|\hspace*{1mm}\emptyset\hspace*{1mm}\Bigg|\hspace*{1mm}\tab(124)\hspace*{1mm}\Bigg|\hspace*{1mm}%
    \tab(3,6)\hspace*{1mm}\Bigg|\hspace*{1mm}\tab(5)\hspace*{1mm}\Bigg) ,
$$
with
$$
\s = \Bigg(\hspace*{1mm}\tab(123)\hspace*{1mm}\Bigg|\hspace*{1mm}%
    \tab(4,5)\hspace*{1mm}\Bigg|\hspace*{1mm}\tab(6)\hspace*{1mm}\Bigg)\hspace*{10mm}
\t = \Bigg(\hspace*{1mm}\tab(124)\hspace*{1mm}\Bigg|\hspace*{1mm}%
    \tab(3,6)\hspace*{1mm}\Bigg|\hspace*{1mm}\tab(5)\hspace*{1mm}\Bigg).
$$

Then $\psi_{\S\T} = \hat\psi_{\s\t} = e(012211)y_2 y_3^2 y_5 \psi_5\psi_3 \in \mathscr R_n$.
\end{Example}

The next Lemma is a straightforward application of the definition of $\psi_{\S\T}$ and $\deg(\S)$.

\begin{Lemma} \label{degree part}
Suppose $\hat\lambda \in \affineP$ and $\S,\T\in \Std(\hat\lambda)$. Then $\psi_{\S\T}$ are homogeneous elements of $\Rn$ and $\deg(\psi_{\S\T}) = \deg(\S) + \deg(\T)$.
\end{Lemma}

Fix $\alpha\in Q_+$, a standard sequence $\inftyweight$ and an inverse multicharge sequence $\kappa$ corresponding to $\inftyweight$. Define a set of homogeneous elements of $\Rni$ by
$$
\Bi = \set{\psi_{\S\T}|\hat\lambda\in \affineP, \S,\T\in \Std(\hat\lambda), \res(\T) \in I^\alpha}.
$$

Note that by definition $\Bi$ depends on the choice of $\kappa$ and hence on $\inftyweight$. Remarkably, the main results of this section are true for any inverse multicharge sequence corresponding to $\inftyweight$.

\begin{Proposition} \label{B is a basis}
The set $\Bi$ is a homogeneous basis of $\Rni$.
\end{Proposition}

\proof By \autoref{degree part}, all of the elements of $\Bi$ are homogeneous. So we only have to prove that $\Bi$ is a basis of $\Rni$. First of all we show that $\Bi$ spans $\Rni$.

If $r\in \Rni$, then we can write $r$ as a linear combination of homogeneous elements, i.e. $r = \sum_{i\in\N} c_i r_i$, where $c_i\in \Z$, $\deg(r_i) = i$ and there are only finite many $i\in \N$ with $c_i \neq 0$. So it is enough to prove that any homogeneous element $r\in \Rni$ is a linear combination of $\Bi$.

For any $\Lambda < \Lambda'$, it is obvious that $N_\alpha^{\Lambda'}\subseteq N_\alpha^\Lambda$. Moreover, $N_\alpha^{\Lambda'}$ is a $\Rni$-ideal of $N_\alpha^\Lambda$. Hence we can define an infinite filtration
$$
\Rni > N_\alpha^{\Lambda^{(1)}} > N_\alpha^{\Lambda^{(2)}} > N_\alpha^{\Lambda^{(3)}} > \ldots.
$$

By \autoref{infinite degree}, $\lim_{k\rightarrow\infty} m^{\Lambda^{(k)}}_\alpha = \infty$, so if $r \in \Rni$ is homogeneous then there exists an integer $k(r)$ such that $m_\alpha^{\Lambda^{(k)}} > \deg(r)$ whenever $k > k(r)$. Fix $k > k(r)$ and hence $r \not\in N_\alpha^{\Lambda^{(k)}}$.

By \autoref{basis of Ri}, choosing a multicharge $\kappa$ corresponding to $\Lambda$, $\Ri \cong \Rni/ N_\alpha^\Lambda$ has a homogeneous basis $\set{\psi_{\s\t}|\lambda\in\P_n, \s,\t\in\Std(\lambda),\res(\t)\in I^\alpha}$. Fix a multicharge $(\kappa_{\l_k},\kappa_{\l_{k-1}},\ldots,\kappa_2,\kappa_1)$ corresponding $\Lambda^{(k)}$. For any homogeneous element $r \in \Rni$, there exists $c_{\s\t}\in \Z$ with $\res(\t) \in I^\alpha$ such that
\begin{eqnarray*}
&&r + N_\alpha^{\Lambda^{(k)}} = \sum_{\s,\t}c_{\s\t}\psi_{\s\t} = \sum_{\s,\t} c_{\s\t}\hat\psi_{\s\t} + N_\alpha^{\Lambda^{(k)}} = \sum_{\S,\T} c_{\s\t}\psi_{\S\T} + N_\alpha^{\Lambda^{(k)}} \\
&\Rightarrow & r - \sum_{\S,\T} c_{\s\t}\psi_{\S\T} \in N_\alpha^{\Lambda^{(k)}}.
\end{eqnarray*}

But as $r$ is a homogeneous element which is not in $N_\alpha^{\Lambda^{(k)}}$, we must have $ r - \sum_{\S,\T} c_{\s\t}\psi_{\S\T} = 0$, i.e. $ r = \sum_{\S,\T} c_{\s\t}\psi_{\S\T}$ with $\res(\T) = \res(\t) \in I^\alpha$. This shows that $r$ belongs to the span of $\Bi$. Hence $\Rni$ is spanned by $\Bi$.

Next we will prove that $\Bi$ is linearly independent. Suppose $S_\alpha$ is a finite subset of $\Bi$. Write $m_{S_\alpha} = \max \set{\deg(\psi_{\S\T})|\psi_{\S\T}\in S_\alpha}$. By \autoref{infinite degree} there exists $k$ such that $m^{\Lambda^{(k)}}_\alpha > m_{S_\alpha}$. Hence $\psi_{\S\T} \not\in N_\alpha^{\Lambda^{(k)}}$ for any $\psi_{\S\T}\in S_\alpha$. This means that for any $\psi_{\S\T}\in S_\alpha$, $\psi_{\s\t} \in \mathscr R_\alpha^{\Lambda^{(k)}}$ is nonzero. As by the definition, $\set{\psi_{\s\t}|\psi_{\S\T}\in S_\alpha}$ is a subset of the basis of $\mathscr R_\alpha^{\Lambda^{(k)}}$. We have $\sum_{\psi_{\S\T}\in S_\alpha} c_{\S\T} \psi_{\S\T} \in N_\alpha^{\Lambda^{(k)}}$ if and only if $\sum_{\psi_{\S\T}\in S_\alpha} c_{\S\T} \psi_{\s\t} = 0$ if and only if all $c_{\S\T} = 0$. But $\psi_{\S\T} \not\in N_\alpha^{\Lambda^{(k)}}$ for any $\psi_{\S\T}\in S_\alpha$, the above result yields that $\sum_{\psi_{\S\T}\in S_\alpha} c_{\S\T} \psi_{\S\T} = 0$ if and only if $c_{\S\T} = 0$. This shows that $\Bi$ is linearly independent. Hence $\Bi$ is a basis of $\Rni$. \endproof

Notice that in the definition of $\Bi$, it is well-defined for any inverse multicharge sequence $\kappa$ corresponds to $\inftyweight$. Hence for any weight $\Lambda$ with $\l = l(\Lambda)$, by the definition of the standard sequence, we can set $\Lambda^{(1)} = \Lambda$. Therefore, we obtain a subset of $\Bi$:
$$
\mathscr B^{\inftyweight}_\Lambda = \set{\psi_{\S\T}|\text{$\hat\lambda\in \affineP$ with $l(\hat\lambda)\leq \l$}, \S,\T\in \Std(\hat\lambda),\res(\T) \in I^\alpha}.
$$

\begin{Corollary} \label{basis of Ni}
Suppose $\Lambda$ is a weight with level $\l$ and $\inftyweight$ is a standard sequence with $\Lambda^{(1)} = \Lambda$. Then
$$
\Bi\backslash \mathscr B^{\inftyweight}_\Lambda = \set{\psi_{\S\T}|\text{$\hat\lambda\in \affineP$ with $l(\hat\lambda)> \l = l(\Lambda)$}, \S,\T\in \Std(\hat\lambda),\res(\T)\in I^\alpha}
$$
is a basis of $N_\alpha^\Lambda$.
\end{Corollary}

\proof By \autoref{basis of Ri}, $\Ri$ has a basis $\set{\psi_{\s\t}|\lambda\in \P_n, \s,\t\in \Std(\lambda),\res(\t)\in I^\alpha}$. It is easy to see that when $\Lambda^{(1)} = \Lambda$,
$$
\set{\psi_{\s\t}|\lambda\in \P_n, \s,\t\in \Std(\lambda),\res(\t)\in I^\alpha} = \set{\psi_{\s\t}=\psi_{\S\T} + N_n^\Lambda|\psi_{\S\T} \in \mathscr B_\Lambda^{\inftyweight}}.
$$
So for $\psi_{\S\T} \in \mathscr B^{\inftyweight}_\Lambda$, we must have $\psi_{\S\T}\not\in N_\alpha^\Lambda$.

Now suppose $\psi_{\S\T} \in \Bi\backslash \mathscr B^{\inftyweight}_\Lambda$. Then $\S,\T\in \Std(\hat\lambda)$ with $l(\hat\lambda) > \l$. By the definition it is obvious that $\psi_{\S\T} \in N_\alpha^\Lambda$ when $\Lambda^{(1)} = \Lambda$. Then $N_\alpha^\Lambda$ is spanned by $\Bi\backslash \mathscr B^{\inftyweight}_\Lambda$. $\Bi$ is a basis implies the linearly independence of $\Bi\backslash \mathscr B^{\inftyweight}_\Lambda$. So $\Bi\backslash \mathscr B^{\inftyweight}_\Lambda$ is a basis of $N_\alpha^\Lambda$. \endproof

Recall for any $\S,\T\in \Std(\hat\lambda)$ with $\hat\lambda \in \affineP$, we define $\psi^*_{\S\T} = \psi_{\T\S}$. By \autoref{B is a basis}, $*$ can be defined to be a linear bijection from $\Rni$ to $\Rni$. The next Corollary is straightforward by \autoref{basis of Ni}.

\begin{Corollary} \label{restriction of *}
Suppose $*\map{\Rni}{\Rni}$ is defined as above. Then it can be restricted to a linear bijection $*\map{N_\alpha^\Lambda}{N_\alpha^\Lambda}$.
\end{Corollary}

Now we can prove the main result of this section.

\begin{Proposition} \label{graded cellular basis of Rni}
Suppose $\inftyweight$ is a standard sequence and $\alpha \in Q_+$. The set $\Bi$ is a graded cellular basis of $\Rni$.
\end{Proposition}

\proof Recall \autoref{definition: graded cellular basis} gives the definition of graded cellular basis. \autoref{B is a basis} shows that $\Bi$ is a homogeneous basis of $\Rni$. To prove the Theorem we need to establish properties \ref{definition: graded cellular basis}(b) and \ref{definition: graded cellular basis}(c) of $\Bi$.

Suppose $a$ is an element of $\Rni$ and $\psi_{\S\T}\in \Bi$ with $\S,\T\in \Std(\hat\lambda)$. We can write $a = \sum_{i\in \N} c_i a_i$ where $c_i \in \Z$ and $a_i$ are homogeneous elements in $\Rni$ with $\deg(a_i) = i$. Define $d_1 = \deg(\psi_{\S\T})$ and $d_2 = \max\{i\ |\ c_i\neq 0\}$. By \autoref{infinite degree} there exists $k$ such that $m_\alpha^{\Lambda^{(k)}} > \max\{d_1,d_2,d_1+d_2\}$. This means that $\psi_{\S\T}$, $a$ and $\psi_{\S\T}a$ are not elements of $N_\alpha^{\Lambda^{(k)}}$. This means that $\psi_{\s\t} = \psi_{\S\T} + N_\alpha^{\Lambda^{(k)}}$, $a + N_\alpha^{\Lambda^{(k)}}$ and $\psi_{\S\T}a + N_\alpha^{\Lambda^{(k)}}$ are nonzero elements of $\mathscr R_\alpha^{\Lambda^{(k)}}$. By \autoref{basis of Ri} and because $t$ is a bijection,
\begin{eqnarray*}
&&\psi_{\s\t}(a + N_\alpha^{\Lambda^{(k)}}) = (\psi_{\S\T} + N_\alpha^{\Lambda^{(k)}})(a + N_\alpha^{\Lambda^{(k)}}) = \sum_{\v \in \Std(\lambda)} c_{\s\v}\psi_{\s\v} + \sum_{\substack{\u,\v\in\Std(\mu)\\ \mu\rhd\lambda}}c_{\u\v} \psi_{\u\v}\\
&\Rightarrow& \psi_{\S\T}a + N_\alpha^{\Lambda^{(k)}} = \sum_{\V\in\Std(\hat\lambda)} c_{\S\V}\psi_{\S\V} + \sum_{\substack{\U,\V\in\Std(\hat\mu)\\ \hat\mu\rhd\hat\lambda}}c_{\U\V} \psi_{\U\V} + N_\alpha^{\Lambda^{(k)}}\\
&\Rightarrow& \psi_{\S\T}a - (\sum_{\V\in\Std(\hat\lambda)} c_{\S\V}\psi_{\S\V} + \sum_{\substack{\U,\V\in\Std(\hat\mu)\\ \hat\mu\rhd\hat\lambda}}c_{\U\V} \psi_{\U\V}) \in N_\alpha^{\Lambda^{(k)}}.
\end{eqnarray*}

Since the left hand side of the above equation is homogeneous to $d_1 + d_2$ and $m_\alpha^{\Lambda^{(k)}} > d_1 + d_2$, we see that
$$
\psi_{\S\T}a = \sum_{\V\in\Std(\hat\lambda)} c_{\S\V}\psi_{\S\V} + \sum_{\substack{\U,\V\in\Std(\hat\mu)\\ \hat\mu\rhd\hat\lambda}}c_{\U\V} \psi_{\U\V}.
$$
which shows that $\Bi$ satisfies \ref{definition: graded cellular basis}(b).

For \ref{definition: graded cellular basis}(c), choose arbitrary $\psi_{\S\T},\psi_{\U\V} \in \Bi$. Suppose $\deg(\psi_{\S\T}) = k_1$ and $\deg(\psi_{\U\V}) = d_2$. By \autoref{infinite degree} we can choose $k$ so that $m_\alpha^{\Lambda^{(k)}} > \max\{k_1,k_2,k_1 + k_2\}$. Then by \autoref{restriction of *},
\begin{eqnarray*}
(\psi_{\s\t}\psi_{\u\v})^* & = & ((\psi_{\S\T} + N_\alpha^{\Lambda^{(k)}})(\psi_{\U\V} + N_\alpha^{\Lambda^{(k)}}))^* = (\psi_{\S\T}\psi_{\U\V} + N_\alpha^{\Lambda^{(k)}})^* = (\psi_{\S\T} \psi_{\U\V})^* + N_\alpha^{\Lambda^{(k)}},\\
\psi_{\v\u}\psi_{\t\s} & = & (\psi_{\V\U} + N_\alpha^{\Lambda^{(k)}})(\psi_{\T\S} + N_\alpha^{\Lambda^{(k)}}) = \psi_{\V\U} \psi_{\T\S} + N_\alpha^{\Lambda^{(k)}},
\end{eqnarray*}
which implies that $(\psi_{\S\T}\psi_{\U\V})^* -\psi_{\V\U} \psi_{\T\S} = N_\alpha^{\Lambda^{(k)}}$. Then because $m_\alpha^{\Lambda^{(k)}} > k_1+k_2$, $(\psi_{\S\T}\psi_{\U\V})^* - \psi_{\V\U} \psi_{\T\S} = 0$, i.e. $(\psi_{\S\T}\psi_{\U\V})^* = \psi_{\U\V}^* \psi_{\T\S}^*$. Because $*$ is a linear bijection and $\Bi$ is a basis of $\Rni$, this shows that $*\map{\Rni}{\Rni}$ is an anti-isomorphism. Hence $*$ satisfies \ref{definition: graded cellular basis}(c). This completes the proof.\endproof

Combining the above two Propositions and \autoref{basis of Ni} we can get the following results.

\begin{Theorem} \label{graded cellular basis of R_n}
For any standard sequence $\inftyweight$, the set
$$
\B = \set{\psi_{\S\T}|\hat\lambda \in \affineP, \S,\T\in\Std(\hat\lambda)}
$$
is a graded cellular basis of $\Rn$.
\end{Theorem}

\proof By definition we have $\B = \bigoplus_{\alpha\in Q_+} \Bi$ and $\Rn = \bigoplus_{\alpha\in Q_+} \Rni$. By the relations of $\Rn$ we see that $\Rni$ are subalgebras. The Theorem follows by \autoref{graded cellular basis of Rni} straightforward. \endproof

\begin{Corollary} \label{basis of N}
Suppose $\Lambda$ is a weight with level $\l$ and $\inftyweight$ is a standard sequence with $\Lambda^{(1)} = \Lambda$. Then
$$
\set{\psi_{\S\T}|\text{$\hat\lambda\in \affineP$ with $l(\hat\lambda)> \l = l(\Lambda)$}, \S,\T\in \Std(\hat\lambda)}
$$
is a basis of $N_n^\Lambda$.
\end{Corollary}

\subsection{Graded simple $\mathscr R_n$-modules}

\autoref{graded cellular basis of R_n} gives a graded cellular basis of $\mathscr R_n$. Graham and Lehrer~\cite{GL} described a complete set of irreducible representations of any finite dimensional cellular algebra, however, their results do not apply to $\Rn$ because it is an infinite dimensional algebra. In this subsection we use Graham and Lehrer's results to construct a complete set of graded simple $\Rn$-modules. The graded simple $\Rn$-modules have been desccribed by Brundan and Kleshchev~\cite[Theorem 5.19]{BK:GradedDecomp}. See~\autoref{remark: Rn modules} for more detials.

First we need to state some properties of the simple $\Rn$-modules. We start by showing that the graded dimension of an simple $\Rn$-module is bounded below.

\begin{Lemma} \label{minimal degree of Rn}
Suppose $r \in \Rn$ is a homogeneous element. Then $\deg(r) \geq -n(n-1)$.
\end{Lemma}

\proof By (\ref{basis of KLR}) we have the following basis of $\Rn$:
$$
\set{\hat e(\bi)\hat y_1^{\l_1} \hat y_2^{\l_2} \ldots \hat y_n^{\l_n} \hat \psi_w|\bi\in I^n, w\in\mathfrak{S}_n, \l_1,\l_2,\ldots,\l_n \geq 0}.
$$

If $\bi\in I^n$, $w\in\mathfrak S_n$ and $\l_1,\l_2,\ldots,\l_n \geq 0$, then
\begin{eqnarray*}
\deg(\hat e(\bi)\hat y_1^{\l_1} \hat y_2^{\l_2} \ldots \hat y_n^{\l_n} \hat \psi_w) & \geq & \deg(\hat e(\bi)\hat y_1^{\l_1} \hat y_2^{\l_2} \ldots \hat y_n^{\l_n}) + \deg(\hat e(\bi) \hat \psi_w)\\
& = & 2(\l_1 + \l_2 + \ldots + \l_n) + \deg(\hat e(\bi) \hat \psi_w)\\
& \geq & \deg(\hat e(\bi) \hat \psi_w).
\end{eqnarray*}

As $w \in \mathfrak S_n$, we have $l(w) \leq \frac{n(n-1)}{2}$ and by definition, $\hat e(\bi)\hat \psi_r \geq -2$ for any $r$ and $\bi\in I^n$. Therefore
$$
\deg(\hat e(\bi)\hat \psi_w) \geq -2 \times \frac{n(n-1)}{2} = -n(n-1).
$$

Hence $\deg(\hat e(\bi)\hat y_1^{\l_1} \hat y_2^{\l_2} \ldots \hat y_n^{\l_n} \hat \psi_w) \geq -n(n-1)$. This completes the proof. \endproof

Recall that for $|\alpha| = \sum_{i \in I} a_i$ and $\hat e_\alpha = \sum_{\bi \in I^\alpha} \hat e(\bi) \in \Rn$.

\begin{Lemma} \label{S is S_i}
Suppose $S$ is a simple $\Rn$-module. Then there exists $\alpha \in Q_+$ with $|\alpha| = n$ such that for any $\beta \in Q_+$ with $|\beta| = n$, $\hat e_\beta S = \delta_{\alpha\beta} S$.
\end{Lemma}

\proof Suppose $S$ is a simple $\Rn$-module. Because $1 = \sum_{\bj \in I^n} \hat e(\bj) = \sum_{\substack{\beta\in Q_+\\ |\beta| = n}}\hat e_\beta$, we can write $S = \bigoplus_{\substack{\beta\in Q_+\\ |\beta| = n}} \hat e_\beta S$. Suppose $\hat e_\alpha S \neq 0$ for some $\alpha \in Q_+$. Choose any nonzero element $s \in S$ and $\beta \in Q_+$ with $\beta \neq \alpha$. By \autoref{alpha R beta}, $\hat e_\alpha \Rn \hat e_\beta = 0$. So we must have $\hat e_\beta{\cdot}s = 0$. Hence $\hat e_\beta S = 0$. Therefore $S = \bigoplus_{\substack{\beta\in Q_+\\ |\beta| = n}} \hat e_\beta S = \hat e_\alpha S$. This completes the proof. \endproof

It is well-known that the irreducible representations of the affine Hecke algebra are finite dimensional as, by Bernstein, the affine Hecke algebra is finite dimensional over its centre. See for example, Proposition 4.1 and Corollary 4.2 of Grojnowski~\cite{Groj:control}, or Proposition 2.12 of Khovanov-Lauda~\cite{KhovLaud:diagI}. The next Proposition gives a different approach.

\begin{Proposition} \label{Rn in R}
Suppose $S$ is a graded simple $\Rn$-module and $\alpha \in Q_+$ is such that $\hat e_\beta S = \delta_{\alpha \beta}S$ for $\beta\in Q_+$. If $\Lambda \in P_+$ with $m_\alpha^\Lambda > n(n-1)$, then $S$ is isomorphic to a graded simple $\R$-module.
\end{Proposition}

\proof By \autoref{S is S_i} there exists $\alpha \in Q_+$ such that $\hat e_\beta S = \delta_{\alpha \beta}S$ for $\beta \in Q_+$. Then we choose an arbitrary nonzero homogeneous element $s \in S$ and suppose $\deg(s) = d$. Now for any nonzero homogeneous element $t \in S$, because $S$ is simple, there exists a homogeneous element $a \in \Rni$ such that $t = a{\cdot}s$. Therefore
$$
\deg(t) = \deg(a{\cdot}s) = \deg(a) + \deg(s) \geq d - n(n-1)
$$
where by \autoref{minimal degree of Rn} we have $\deg(a) \geq -n(n-1)$. So for any homogeneous nonzero element $t \in S$, we have
\begin{equation} \label{Rn in R: help1}
\deg(t) \geq d - n(n-1).
\end{equation}

Similarly, since for any nonzero homogeneous element $t \in \Rn$ there exists a homogeneous element $a \in \Rni$ such that $s = a{\cdot}t$, we have
\begin{equation} \label{Rn in R: help2}
\deg(t) \leq d + n(n-1).
\end{equation}

Combining (\ref{Rn in R: help1}) and (\ref{Rn in R: help2}), we have $|\deg(s) - \deg(t)| \leq n(n-1)$ for any nonzero homogeneous element $t\in S$. Because $s$ is chosen arbitrarily, we have
\begin{equation} \label{Rn in R: help3}
|\deg(s) - \deg(t)| \leq n(n-1)
\end{equation}
for any nonzero homogeneous elements $s,t \in S$.

Suppose $\Lambda \in P_+$ with $m_\alpha^\Lambda > n(n-1)$. For any homogeneous element $a \in N_\alpha^\Lambda$ and $t \in S$, we have $a{\cdot}t = 0$ because $\deg(a{\cdot}t) - \deg(t) = \deg(a) > n(n-1)$ and (\ref{Rn in R: help3}).

For any $s \in S$, we can define a map $f\map{\Rn}{S}$ by sending $a$ to $a{\cdot}s$. It is a homomorphism and it is obvious that $N_\alpha^\Lambda \subseteq \text{ker}f$. If $\beta\in Q_+$ and $\beta \neq \alpha$, then by \autoref{S is S_i} we have $\hat e_\beta{\cdot}s = 0$. Therefore $N_\beta^\Lambda \subseteq \text{ker}f$. Hence $N_n^\Lambda \subseteq \text{ker}f$. Therefore we can consider $S$ as a simple $\Rn/N_n^\Lambda$-module, i.e. $\R$-module. This completes the proof. \endproof

\begin{Corollary} \label{simple S is finite dimensional}
Suppose $S$ is a simple graded $\Rn$-module. Then $S$ is finite-dimensional.
\end{Corollary}

Building on Ariki's~\cite{Ariki:class} work in the ungraded case, Hu and Mathas~\cite{HuMathas:GradedCellular} constructed all graded simple $\R$-modules in the sense of Graham-Lehrer~\cite{GL}. They proved that, up to shift, graded simple $\R$-modules are labeled by the \textbf{Kleshchev multipartitions} of $n$, which were introduced by Ariki and Mathas~\cite{AM:simples}. Readers may also refer to Brundan and Kleshchev~\cite[(3.27)]{BK:GradedDecomp} (where they are called \textbf{restricted multipartitions}).

Suppose $\lambda = (\lambda^{(1)},\lambda^{(2)},\ldots,\lambda^{(\l)})\in \P_n$ and we consider the Young diagram $[\lambda]$. Let $\gamma = (r,c,l)$ be a node in the Young diagram with residue $i$, i.e. $i \equiv r-c+\kappa_l\pmod{e}$. Then $\gamma$ is an addable $i$-node if $\gamma \not\in [\lambda]$ and $[\lambda]\cup\{\gamma\}$ is the Young diagram of a multipartition, and $\gamma$ is a removable $i$-node if $\gamma \in [\lambda]$ and $[\lambda]\backslash \{\gamma\}$ is the Young diagram of a multipartition.

For each $\lambda \in \P_n$, we read all addable and removable $i$-nodes in the following order: we start with the first row of $\lambda^{(1)}$, and then read rows in $\lambda^{(1)}$ downward. We then read the first row of $\lambda^{(2)}$, and repeat the same procedure, until we finish reading all rows of $\lambda$. We write $A$ for an addable $i$-node, and $R$ for a removable $i$-node. Hence we get a sequence of $A$ and $R$. We then delete $RA$ as many as possible. For example, if we have a sequence $RARARRAAARRAR$, the resulting sequence will be $--------AR--R$. The node corresponding to the leftmost R is the \textbf{good $i$-node}.

The Kleshchev multipartition can then be defined recursively as follows.

\begin{Definition} \cite[Definition 2.3]{Ariki:class}
We declare that $\emptyset$ is Kelshchev. Assume that we have already defined the set of Kleshchev multipartitions up to size $n-1$. Let $\lambda$ be a multipartition of $n$. We say that $\lambda$ is a Kleshchev multipartition if there is a good node $\gamma$ in $[\lambda]$ such that if $[\mu] = [\lambda]\backslash\{\gamma\}$ then $\mu$ is a Kleshchev multipartition.
\end{Definition}

Let $\P_0$ be the set of Kleshchev multipartitions in $\P_n$. Let $S^\lambda$ be the cell module of $\R$(it is called the \textbf{Specht module} in $\R$), which was introduced in subsection 1.2, and $D^\lambda = S^\lambda/\text{rad }S^\lambda$.  Hu-Mathas~\cite[Corollary 5.11]{HuMathas:GradedCellular} gave a set of complete non-isomorphic graded simple $\R$-modules. We note that the graded simple $\R$-modules were first constructed by Brundan and Kleshchev~\cite[Theorem 4.11]{BK:GradedDecomp} giving the same classification but without using cellular algebra techniques.

\begin{Theorem} \label{simple in R}
The set $\set{D^\lambda\langle k\rangle|\lambda\in\P_0,  k\in\Z}$ is a complete set of pairwise non-isomorphic graded simple $\R$-modules.
\end{Theorem}

We consider $S^\lambda$ and $D^\lambda\langle k\rangle$ as $\Rn$-modules. The actions of $\hat e(\bi)$, $\hat y_r$ and $\hat \psi_s$ on $S^\lambda$ and $D^\lambda\langle k\rangle$ are the same as the actions of $e(\bi)$, $y_r$ and $\psi_s$. Therefore $D^\lambda\langle k\rangle$ is a simple $\Rn$-module. Hence the modules in \autoref{simple in R} are a set of simple $\Rn$-modules.

The next Lemma is straightforward by the definition of $D^\lambda$.

\begin{Lemma} \label{Rn if and only if R}
Suppose $\lambda,\mu \in \P_n$. Then $D^\lambda \cong D^\mu$ as $\R$-modules if and only if $D^\lambda \cong D^\mu$ as $\Rn$-modules.
\end{Lemma}

Now we can classify all graded simple $\Rn$-modules. Following the definitions in Section 1.2, for each $\hat\lambda \in \affineP$ we define the cell module $S^{\hat\lambda}$ of $\Rn$, which we also call a \textbf{Specht module}, associated with a bilinear form $\langle\cdot ,\cdot \rangle$. Then we define $\text{rad }S^{\hat\lambda}$ and hence the graded simple module $D^{\hat\lambda} = S^{\hat\lambda}/\text{rad }S^{\hat\lambda}$.

\begin{Lemma} \label{remove empty in the front}
Suppose $\hat\lambda \in \affineP$ and $\mu = p_k(\hat\lambda)$ for some $k \geq l(\hat\lambda)$. Then $S^\mu \cong S^{\hat\lambda}$ as $\Rn$-modules.
\end{Lemma}

\proof This is trivial given the definition of Specht modules in $\Rn$ and $\R$. \endproof

The next Corollary is straightforward by ~\autoref{Rn if and only if R} and~\autoref{remove empty in the front}.

\begin{Corollary} \label{same D lambda}
Suppose $\hat\lambda \in \affineP$ and $\mu = p_k(\hat\lambda)$ for some $k \geq l(\hat\lambda)$. Then $D^\mu \cong D^{\hat\lambda}$ as $\Rn$-modules.
\end{Corollary}

Hence we can prove the following Lemma.

\begin{Lemma}\label{help for all simple}
Suppose $\hat\lambda,\hat\mu \in \affineP$. Then $D^{\hat\lambda} \cong D^{\hat\mu}$ if and only if $\hat\lambda = \hat\mu$.
\end{Lemma}

\proof The if part is trivial. Now suppose $D^{\hat\lambda} \cong D^{\hat\mu}$. Choose $k > \max\{l(\hat\lambda),l(\hat\mu)\}$ and set $\nu = p_k(\hat\lambda)$ and $\sigma = p_k(\hat\mu)$. Then by \autoref{same D lambda} we have $D^\nu \cong D^\sigma$ as $\Rn$-modules. Then \autoref{simple in R} and \autoref{Rn if and only if R} implies $\nu = \sigma$. Therefore by the definition of $k$ we have $\hat\lambda = \hat\mu$. This completes the proof. \endproof

Now we extend the Kleshchev multipartitions to affine multipartitions. Define $\hat\lambda\in\affineP$ to be an \textbf{affine Kleshchev multipartition} if $\lambda$ is a Kleshchev multipartition. Let $\mathscr P_0^\kappa$ be the set of all affine Kleshchev multipartitions in $\affineP$. We can now give a complete set of pairwise non-isomorphic graded simple $\Rn$-modules.

\begin{Theorem} \label{all simple}
The set $\set{D^{\hat\lambda}\langle k\rangle|\hat\lambda\in\mathscr P_0^\kappa, k\in\Z}$ is a complete set of pairwise non-isomorphic graded simple $\Rn$-modules.
\end{Theorem}

\proof By the definition of (affine) Kleshchev multipartitions, \cite[Corollary 5.11]{HuMathas:GradedCellular} and \autoref{same D lambda}, $D^{\hat\lambda}\langle k\rangle \cong D^{\lambda}\langle k\rangle \neq 0$ if and only if $\hat\lambda \in \mathscr P_0^\kappa$.

Suppose $S$ is a graded simple $\Rn$-module. By \autoref{S is S_i} there exists $\alpha \in Q_+$ such that $\hat e_\beta S = \delta_{\alpha\beta}S$ for $\beta \in Q_+$. Then by \autoref{infinite degree} we can choose $i$ such that $m_\alpha^{\Lambda^{(i)}} > n(n-1)$ and hence by \autoref{Rn in R}, $S$ is isomorphic to a graded simple $\Rn^{\Lambda^{(i)}}$-module. Therefore, by \autoref{simple in R} there exist some $\mu \in \mathscr P_n^{\Lambda^{(i)}}$ and $k\in \Z$ such that $S\cong D^\mu \langle k\rangle$ as $\Rn^{\Lambda^{(i)}}$-modules, and hence as $\Rn$-modules.  Suppose $l(\mu) = \l$. We can choose $\hat\lambda \in \affineP$ such that $p_\l(\hat\lambda) = \mu$ with $l(\hat\lambda) \leq \l$. By \autoref{same D lambda} we have $D^{\hat\lambda} \cong D^\mu$ as $\Rn$-modules. Therefore, $S\cong D^{\hat\lambda}\langle k\rangle$. So
$$
\set{D^{\hat\lambda}\langle k\rangle|\hat\lambda\in\mathscr P_0^\kappa, k\in\Z}
$$
is a complete set of graded simple $\Rn$-modules.

By \autoref{help for all simple}, the set $\set{D^{\hat\lambda}\langle k\rangle|\hat\lambda\in\mathscr P_0^\kappa, k\in\Z}$ is a set of pairwise non-isomorphic graded $\Rn$-modules. This completes the proof.   \endproof

\begin{Remark} \label{remark: Rn modules}
Ariki-Mathas~\cite{AM:simples} showed that the simple $H_n$-modules are indexed by aperiodic multisegments. Khovanov and Lauda~\cite{KhovLaud:diagI,KhovLaud:diagII} also give a classification of the graded simple $\Rn$-modules for KLR algebras of arbitrary type. Interested readers may also refer to~\cite{KLeshRam:PureIrredRepsKLR},~\cite{BK:GradedDecomp},~\cite{KleshchevRam:IrredKLR},~\cite{LV:Crystals},~\cite{Webster:KnotI},~\cite{KK:HighestWeight} and~\cite{McNamara:FiniteDimenKLRI}. As far as we are aware the construction and classification in \autoref{all simple} is new.
\end{Remark}

\bibliography{papers}
\end{document}